%% file: Main.tex
% plain TeX wrapper for Dynamics in One Complex Variable: Intro. Lectures
%
\let\ReallyTheEnd=\end
\def\end{\vfill\eject}
%
\input header
\input imsmark
\def\IMSpubltext{A significantly revised version published}
\SBIMSMark{1990/5}{May 1990}{revised September 1991}
%
\input cd0-2.tex

\input cd3-5.tex

\input cd6-8.tex

\input cd9-11.tex

\input cd12-14.tex

\input cd15-16.tex

\input cd17-18.tex

\input cda-c.tex

\input cdd-g.tex

\input cdref.tex

%% file: header.tex
%   Standard Plaintex Header File
% Make sure we haven't already been run

\ifx\HeaderIsLoaded\undefined
\let\HeaderIsLoaded=\relax		% define HeaderIsLoaded

%\magnification=1200
\input psfig
\def\QP{\narrower\smallskip\noindent}
\def\ref{\hangindent=1pc \hangafter=1 \noindent}
\def\QED{  \rlap{$\sqcup$}$\sqcap$ \smallskip}
\def\mod{{\rm mod\;}}
\def\Poincare{Poincar\'e }
\def\[{$\,}
\def\]{\,$}
\def\C{{\bf C}}
\def\R{{\bf R}}
\def\Z{{\bf Z}}
\def\Q{{\bf Q}}
\def\={\;=\; }
\def\lln{\smallskip\centerline{------------------------------------------------------}\smallskip}
\font\bit=cmssi12 at 12truept
\font\tenmsy=msym10

\textfont8=\tenmsy
\mathchardef\ssm="7872
\mathsurround = 1pt
\abovedisplayskip=6pt
\belowdisplayskip=6pt
\parskip=2pt

\fi		% end of \ifx\HeaderIsLoaded

%% file: imsmark.tex
\def\IMSmarkvadjust{0 pt}
\def\IMSmarkhadjust{0 pt}
\def\IMSmarkhpadding{0 pt}
\def\IMSpubltext{Published in modified form:}
\def\SBIMSMark#1#2#3{
 \font\SBF=cmss10 at 10 true pt
 \font\SBI=cmssi10 at 10 true pt
 \setbox0=\hbox{\SBF \hbox to \IMSmarkhpadding{\relax}
                Stony Brook IMS Preprint \##1}
 \setbox2=\hbox to \wd0{\hfil \SBI #2}
 \setbox4=\hbox to \wd0{\hfil \SBI #3}
 \setbox6=\hbox to \wd0{\hss
             \vbox{\hsize=\wd0 \parskip=0pt \baselineskip=10 true pt
                   \copy0 \break%
                   \copy2 \break% 
                   \copy4 \break}}
 \dimen0=\ht6   \advance\dimen0 by \vsize \advance\dimen0 by 8 true pt
                \advance\dimen0 by -\pagetotal
	        \advance\dimen0 by \IMSmarkvadjust
 \dimen2=\hsize \advance\dimen2 by .25 true in
	        \advance\dimen2 by \IMSmarkhadjust

%
%   Check for publication info
%
%  \newread\jref
  \openin2=publishd.tex
  \ifeof2\setbox0=\hbox to 0pt{}
  \else 
     \setbox0=\hbox to 3.1 true in{
                \vbox to \ht6{\hsize=3 true in \parskip=0pt  \noindent  
                {\SBI \IMSpubltext}\hfil\break
                \input publishd.tex 
                \vfill}}
  \fi
  \closein2
  \ht0=0pt \dp0=0pt
 \ht6=0pt \dp6=0pt
 \setbox8=\vbox to \dimen0{\vfill \hbox to \dimen2{\copy0 \hss \copy6}}
 \ht8=0pt \dp8=0pt \wd8=0pt
 \copy8
 \message{*** Stony Brook IMS Preprint #1, #2. #3 ***}
}

%% file: publishd.tex
{\SBI in August 1999 by F. Vieweg \& Sohn} \hfil\break
{\it ISBN 3-528-03130-1}

%% file: cd0-2.tex
% written in plaintex
\input header
\def\PSL{{\rm PSL}}
\def\PD{\rho}
\def\SD{\sigma}
\footline={\hss\tenrm i-\folio\hss}\pageno=1

%table of content format
\def\leaderfill{{\leaders\hbox to 1em{\hss.\hss}\hfill}\hfil}
\def\eleaderfill{\hfill\hbox to 1em{\hss\hss}\hfill}

\vbox{}\vskip -2\baselineskip

\centerline{\bf DYNAMICS IN ONE COMPLEX VARIABLE}\smallskip
\centerline{\bf Introductory Lectures}\smallskip
\centerline{(Partially revised version of 9-5-91)}
\smallskip

\centerline{\bf John Milnor}\smallskip
% (SUNY Stony Brook and IAS Princeton)}\smallskip % \medskip

\centerline{\bit Institute for Mathematical Sciences, SUNY, Stony Brook NY}
\medskip

\centerline{\bf Table of Contents.}\smallskip
\line{Preface\leaderfill }
\line{Chronological Table\leaderfill }
\smallskip
\line{{\bf Riemann Surfaces}.\eleaderfill}
\line{\kern15pt 1.  Simply Connected Surfaces\leaderfill (9 pages)}
\line{\kern15pt 2. The Universal Covering, Montel's Theorem\leaderfill (5)}
\smallskip
\line{{\bf The Julia Set}.\eleaderfill}
\line{\kern15pt 3. Fatou and Julia: Dynamics
on the Riemann Sphere\leaderfill (9)}
\line{\kern15pt 4. Dynamics on Other Riemann Surfaces\leaderfill (6)}
\line{\kern15pt 5. Smooth Julia Sets\leaderfill (4)}
\smallskip
\line{{\bf Local Fixed Point Theory}.\eleaderfill}
\line{\kern15pt 6. Attracting and Repelling Fixed Points\leaderfill (7)}
\line{\kern15pt 7. Parabolic Fixed Points: the Leau-Fatou Flower\leaderfill (8)}
\line{\kern15pt 8. Cremer Points and Siegel Disks \leaderfill (11)}
\smallskip
\line{{\bf Global Fixed Point Theory}.\eleaderfill}
\line{\kern15pt 9. The Holomorphic Fixed Point Formula\leaderfill (3)}
\line{\kern10pt 10. Most Periodic Orbits Repel\leaderfill (3)}
\line{\kern10pt 11. Repelling Cycles are Dense in \[J\]\leaderfill (5)}
\smallskip
\line{{\bf Structure of the Fatou Set}.\eleaderfill}
\line{\kern10pt 12. Herman Rings\leaderfill (4)}
\line{\kern10pt 13. The
Sullivan Classification of Fatou Components\leaderfill (4)}
\line{\kern10pt 14. Sub-hyperbolic and hyperbolic Maps\leaderfill (7)}
\smallskip
\line{{\bf Carath\'eodory Theory}.\eleaderfill}
\line{\kern10pt 15. Prime Ends\leaderfill (5)}
\line{\kern10pt 16. Local Connectivity\leaderfill (4)}
\smallskip
\line{{\bf Polynomial Maps}.\eleaderfill}
\line{\kern10pt 17. The Filled Julia Set \[K\,$\leaderfill (4)}
\line{\kern10pt 18. External Rays and Periodic Points\leaderfill (8)}
\smallskip
\line{Appendix A. Theorems from Classical Analysis\leaderfill (4)}
\line{Appendix B. Length-Area-Modulus Inequalities\leaderfill (6)}
\line{Appendix C. Continued Fractions\leaderfill (5)}
\line{Appendix D. Remarks on Two Complex Variables\leaderfill (2)}
\line{Appendix E. Branched Coverings and Orbifolds\leaderfill (4)}
\line{Appendix F. Parameter Space\leaderfill (3)}
\line{Appendix G. Remarks on Computer Graphics\leaderfill (2)}
\smallskip
\line{References\leaderfill (7)}
\line{Index\leaderfill (3)}

\eject

\centerline{\bf PREFACE.}\bigskip

These notes will study the dynamics of iterated holomorphic
mappings from a\break Riemann surface to itself, concentrating on
the classical case of rational maps of the\break Riemann sphere.
They are based on introductory lectures
given at Stony Brook during the Fall Term of 1989-90. I am grateful
to the audience
for a great deal of constructive criticism, and
to Branner, Douady, Hubbard, and Shishikura who
taught me most of what I know in this field. The surveys by
Blanchard, Devaney, Douady, Keen, and Lyubich have been extremely useful, and
are highly recommended
to the reader. (Compare the list of references at the end.) Also, I want to
thank A. Poirier for his criticisms of my first draft.\smallskip

This subject is large and rapidly growing. These lectures are intended to
introduce the reader to some key ideas in
the field, and to form a basis for further study.
The reader is assumed to be familiar with the rudiments of complex variable
theory and of two-dimensional differential geometry. % Compare the list
\bigskip
\rightline{John Milnor, Stony Brook, February 1990}\medskip
\vfil\eject  % \bigskip

\centerline{\bf CHRONOLOGICAL TABLE}\bigskip

Following is a list of some of the founders of the field of complex dynamics.
\smallskip

$$\vbox{\settabs 2 \columns
\+Hermann Amandus Schwarz & 1843--1921\cr
\+Henri Poincar\'e & 1854--1912\cr
\+Gabriel Ko$\!$enigs	& 1858--1931\cr
\+L\'eopold Leau &	1868--1940$\pm$\cr
% \+L. \`E. B\"otk$\!$her  & ??\cr %(= L. \`E. B\"otk$\!$her = &\cr
\+Lucyan Emil B\"ottcher  & 1872-- ~?\cr
\+Constantin Carath\'eodory & 1873--1950\cr
\+Samuel Latt\`es & 1875$\pm$-1918\cr
\+Paul Montel   & 1876--1975\cr
\+Pierre Fatou &	1878--1929\cr
\+Paul Koebe & 1882--1945\cr
\+Gaston Julia &        1893--1978 \cr
\+Carl Ludwig Siegel &  1896--1981\cr
\+Hubert Cremer &	1897--1983\cr
\+Charles Morrey & 1907--1984\cr}$$
\smallskip

\noindent Among the many present day workers in the field, let me mention
a few whose work is emphasized in these notes:
Lars Ahlfors (1907),
Lipman Bers (1914),
Adrien Douady (1935),
Vladimir I. Arnold (1937),
Dennis P. Sullivan (1941),
Michael R. Herman (1942),
Bodil Branner (1943),
John Hamal Hubbard    (1945),
William P. Thurston (1946),
Mary Rees (1953),
Jean-Christophe Yoccoz (1955),
Mikhail Y. Lyubich (1959),
and Mitsuhiro Shishikura (1960).
% \bigskip
\vfil
%\noindent ------------

%\noindent * B\"ottcher was born in
%Warsaw, took his doctorate in Leipzig (1898), and was active in Lvov
%and Kazan. The Russian form of his name is ~{\cyr L. \`E. B\"etkher\cdprime}.

\eject
\footline={\hss\tenrm 1-\folio\hss}\pageno=1

\centerline{\bf RIEMANN SURFACES}\medskip

\centerline{\bf \S1. Simply Connected Surfaces.}\medskip

The first two sections will present an overview of well known material.
\smallskip

By a {\bit Riemann surface} we mean a connected complex analytic manifold
of complex dimension one. Two such surfaces \[S\] and \[S'\] are
{\bit conformally isomorphic} if there is a homeomorphism
from \[S\] onto \[S'\] which is holomorphic, with holomorphic inverse.
(Thus our conformal maps must always preserve orientation.) According
to \Poincare and Koebe, there are only three simply connected Riemann surfaces,
up to isomorphism.

{\QP{\bf 1.1. Uniformization Theorem.} \it Any simply connected Riemann
surface is conformally isomorphic either

\noindent $(a)$ to the plane \[\C\]
consisting of all complex numbers \[z=x+iy\],

\noindent $(b)$ to the open unit disk
\[D\subset\C\] consisting of all \[z\] with \[|z|^2=x^2+y^2<1\], or

\noindent $(c)$ to the
Riemann sphere \[\hat\C\] consisting of \[\C\] together with a point at
infinity.\medskip}

The proof, which is quite difficult, may be found in Springer, or Farkas \&
Kra, or Ahlfors [A2], or in Beardon.\QED

For the rest of this section, we will discuss these three surfaces in more
detail. We begin with a study of the unit disk \[D\].

{\QP{\bf 1.2. Schwarz Lemma.} {\it If \[f:D\to D\] is a holomorphic map with
\[f(0)=0\], then the derivative at the origin satisfies \[|f'(0)|\le 1\].
If equality holds, then \[f\] is a rotation about the origin, that is
\[f(z)=\lambda z\] for some constant \[\lambda
=f'(0)\] on the unit circle. In particular,
it follows that \[f\] is a conformal automorphism of \[D\]. On the other
hand, if \[|f'(0)|<1\], then \[|f(z)|<|z|\] for all \[z\ne 0\], and
\[f\] is not a conformal automorphism.}\medskip}

{\bf Proof.} We use the {\bit maximum modulus principle}, which
asserts that a non-constant
holomorphic function cannot attain its maximum absolute value at 
any interior point of its region of definition. First note that the
quotient function \[g(z)=f(z)/z\] is well defined and holomorphic
throughout the disk \[D\]. Since \[|g(z)|<1/r\] when \[|z|=r<1\], it follows
by the maximum modulus principle that \[|g(z)|< 1/r\] for all \[z\] in the
disk \[|z|\le r\]. Taking the limit as \[r\to 1\], we see that \[|g(z)|\le
1\] for all \[z\in D\]. Again by the maximum modulus principle, we see that
the case \[|g(z)|=1\], with \[z\] in the open disk, can occur only if the
function \[g(z)\] is constant. If we exclude this case \[f(z)/z=c\],
then it follows that \[|g(z)|=|f(z)/z|<1\] for all \[z\ne 0\], and similarly
that \[|g(0)|=|f'(0)|<1\]. Evidently the composition of two such maps
must also satisfy \[|f_1(f_2(z))|<|z|\], and hence cannot be the identity
map of \[D\]. \QED

{\bf 1.3. Remarks.} The Schwarz Lemma was first proved in this generality
by Carath\'eo\-dory. If we map the disk \[D_r\] of radius \[r\] into the disk
\[D_s\] of radius \[s\],
with \[f(0)=0\], then a similar argument shows that \[|f'(0)|
\le s/r\]. Even if we drop the condition that \[f(0)=0\], we certainly get the
inequality
$$      |f'(0)|\le 2s/r\qquad{\rm whenever}\qquad f(D_r)\subset D_s     $$
since the difference \[f(z)-f(0)\] takes values in \[D_{2s}\].
(In fact the extra factor of 2 is unnecessary. Compare Problem 1-2.)
One easy corollary
is the {\bit Theorem of Liouville}, which says that a bounded function
which is defined and holomorphic everywhere on \[\C\] must be constant.
Another closely related statement is the following.

{\QP{\bf 1.4. Theorem of Weierstrass.} \it If a sequence of holomorphic
functions converges uniformly, then their derivatives also converge uniformly,
and the limit function is itself holomorphic.\smallskip}

\noindent In fact the convergence of first derivatives follows easily from the
discussion above. For the proof of the full theorem, see for example [A1].
\smallskip

It follows from this discussion that our three model surfaces
really are distinct. For there are natural inclusion maps \[D\to\C\to\hat\C\];
yet it follows from the maximum modulus principle and Liouville's Theorem
that every holomorphic map \[\hat\C\to\C\] or %\break
\[\C\to D\] must be constant.\medskip

It is often more convenient to work with the {\bit upper half-plane}
\[H\], consisting of all complex numbers \[w=u+iv\] with \[v>0\].

{\QP{\bf 1.5. Lemma.} {\it The half-plane \[H\] is conformally isomorphic to
the disk \[D\] under the holomorphic mapping \[z=(i-w)/(i+w)\], with
inverse \[w=i(1-z)/(1+z)\].}\medskip}

{\bf Proof.} We have \[|z|^2<1\] if and only if \[|i-w|^2=
u^2+(1-2v+v^2)\] is less than
\[|i+w|^2=u^2+(1+2v+v^2)\], or in other words if and only if \[0<v\]. \QED

{\QP{\bf 1.6. Lemma.} {\it Given any point \[z_0\] of \[D\], there exists
a conformal automorphism \[f\] of \[D\] which maps \[z_0\] to the origin.
Furthermore, \[f\] is uniquely determined up to composition with a rotation
which fixes the origin.}\medskip}

{\bf Proof.} Given any two points \[w_1\] and \[w_2\] of the upper
half-plane \[H\], it is easy to check that there
exists a automorphism of the form \[w\mapsto a+bw\] which carries \[w_1\]
to \[w_2\]. Here the coefficients \[a\] and \[b\]
are to be real, with \[b>0\]. Since \[H\]  is conformally isomorphic to
\[D\], it follows that there exists a conformal automorphism of \[D\]
carrying any given point to zero. As noted above, the Schwarz Lemma
implies that
any automorphism of \[D\] which fixes the origin is a rotation. \QED

{\QP{\bf 1.7. Theorem.} \it The group \[G(H)\] consisting of all
conformal automorphisms of the upper half-plane \[H\] can be identified
with the group of all fractional linear transformations \[w\mapsto
(aw+b)/(cw+d)\] with real coefficients and with determinant \[ad-bc>0\].
\medskip}

If we normalize so that \[ad-bc=1\], then these coefficients are
well defined up to a simultaneous change of sign. {\it
Thus \[G(H)\] is isomorphic to the group \[\PSL(2,\R)\], consisting of all
\[2\times 2\] real matrices with determinant \[+1\] modulo the subgroup
\[\{\pm I\}\].}  To prove 1.7, it is only necessary to
note that this group acts transitively on \[H\], and that it contains the
group of ``rotations''
$$	g(w)= (w\cos\theta+\sin\theta)/(-w\sin\theta+\cos\theta)\,,
 \eqno (1)$$
which fix the point \[g(i)=i\] with \[g'(i)=e^{2i\theta
}\]. By 1.6, there can be no further
automorphisms. (Compare Problem 1-2.)\QED\smallskip

Next we introduce the {\bit Poincar\'e metric} on the half-plane \[H\].

{\QP{\bf 1.8. Theorem.} {\it There exists one and, up to multiplication by
a constant, only one Riemannian metric on the half-plane \[H\] which is invariant
under every conformal automorphism of \[H\].}\medskip}

It follows immediately that the same statement is true for the disk
\[D\], or for any other Riemann surface which is conformally isomorphic
to \[H\].\smallskip

{\bf Proof of 1.8.} A Riemannian metric \[ds^2=g_{11}dx^2
+2g_{12}dxdy+g_{22}dy^2\] is said to be {\bit conformal}$\,$ if
\[g_{11}=g_{22}\] and \[g_{12}=0\], so that the matrix
\[[g_{ij}]\], evaluated at any point \[z\],
is some multiple of the identity matrix. Such a metric can be
written as\break \[ds^2=\gamma(x+iy)^2(dx^2+dy^2)\], or briefly as \[ds=\gamma(z)
|dz|\], where the function \[\gamma(z)\] is smooth and strictly positive.
By definition, such a metric is {\bit invariant} under a conformal
automorphism \[z'=f(z)\] if and only if it satisfies the identity
\[\gamma(z')|dz'|=\gamma(z)|dz|\], or in other words.
$$	\gamma(f(z))=\gamma(z)/|f'(z)|\,.	\eqno (2)	$$
Equivalently, we may say that \[f\] is an {\bit isometry} with respect to
the metric.

As an example, suppose that a conformal metric \[\gamma(w)
|dw|\] on the upper half-plane
is invariant under every linear automorphism \[f(w)=a+bw\]. Then we must
have\break \[\gamma(a+bi)=\gamma(i)/b\,\]. If we normalize by setting \[\,
\gamma(i)
=1\,\], then we are led to the\break formula \[\gamma(u+iv)=1/v\], or
in other words
$$	ds=|dw|/v\qquad {\rm for}\quad w=u+iv\in H \,. \eqno(3) $$
{\it In fact, the metric defined by this formula is invariant under every
conformal automorphism \[g\] of \[H\].} For, if we select some arbitrary point
\[w_1\in H\] and set \[g(w_1)=w_2\], then \[g\]
can be expressed as the composition of a linear automorphism of the form
\[g_1(w)= a+bw\] which maps \[w_1\] to \[w_2\]
and an automorphism \[g_2\] which fixes \[w_2\].
We have constructed the metric (3) so that \[g_1\] is an isometry,
and it follows from
1.2 and 1.5 that \[|g'_2(w_2)|=1\], so that \[g_2\] is an isometry at \[w_2\].
Thus our metric is invariant at an arbitrarily chosen point under an arbitrary
automorphism.

To complete the proof of 1.8, we must show that a metric which is invariant
under all automorphisms
of \[H\] is necessarily conformal. For this purpose, it is
convenient to pass to the equivalent problem on \[D\]. In fact
a brief computation shows that any metric on \[D\] which is 
preserved by all rotations about the origin must be conformal
at the origin. Further details will be left to the reader. \QED

{\bf Definition.} This metric \[ds=|dw|/v\] is called
the {\bit Poincar\'e metric} on the upper half-plane
\[H\]. It can be shown that this is
the unique conformal metric on \[H\] which is complete, with
constant Gaussian curvature equal to \[-1\].
(Compare Problems 1-9 and 2-2.)
The corresponding expression on \[D\] is
$$ ds = 2|dz|/(1-|z|^2)\qquad {\rm for}\quad z=x+iy\in D\,, \eqno (4)$$
as can be verified using 1.5 and (2).\smallskip

{\bf Caution:} Some authors call \[|dz|/(1-|z|^2)\] the Poincar\'e
metric on \[D\], and correspondingly call
\[{1\over 2}|dw|/v\] the Poincar\'e metric on \[H\].
These modified metrics have constant Gaussian curvature equal to \[-4\].
\smallskip

Thus there is a preferred Riemannian metric \[ds\] on \[D\] or on \[H\].
More generally, if \[S\] is any Riemann surface which is conformally
isomorphic to \[D\], then there is a corresponding Poincar\'e metric
\[ds\] on \[S\].
Hence we can define the {\bit Poincar\'e distance} \[\PD(z_1\,,\, z_2)
=\PD_S(z_1\,,\, z_2)\] between two points of \[S\]
as the minimum, over all paths \[P\] joining
\[z_1\] to \[z_2\], of the integral \[\,\int_P\,ds\]. 

{\QP{\bf 1.9. Lemma.} \it The disk \[D\]
with this \Poincare metric is {\bit complete}. That is:

\noindent $(a)$ every Cauchy sequence with respect to the metric \[\rho_D\] converges,

\noindent $(b)$ every closed neighborhood \[\;N_r(z_0\,,\,\rho_D)\;=\;\{z\in D : \rho_D
(z,z_0)\le r\}\;\] is a compact subset of \[D\], and

\noindent$(c)$ any path leading from \[z_0\] to a point of the boundary
\[\partial D=\bar D\ssm D\subset \hat\C\] has infinite \Poincare length.

\noindent Furthermore, any two points of \[D\] are joined by a unique
minimal geodesic.\medskip}

(Compare Willmore.) Equivalently, we
have exactly the same statements for the half-plane \[H\].\smallskip

{\bf Proof.} Given any two points of \[D\] we can first choose a conformal
automorphism which moves the first point to the origin and the second
to some point \[\xi\] on the positive real axis. For any path \[P\]
between \[0\] and \[\xi\] within \[D\] we have
$$	\int_P ds\;=\;\int_P\;{2|dz|\over 1-|z|^2}\;\ge\;\int_P\,{2|dx|\over
 1-x^2}	\;\ge\;\int_{[0,\xi]}\;{2dx\over 1-x^2}\;=\;\log{1+\xi\over 1-\xi}\;,
	$$
with equality if and only if \[P\] is the straight line segment \[[0,\xi]\].
For any \[z\in D\], it follows easily  that the \Poincare distance \[\rho
=\rho_D
(0,z)\] is equal to \[\log\big((1+|z|)/(1-|z|)\big)\], and that the straight line
segment from \[0\] to \[z\] is the unique minimal \Poincare geodesic. Clearly
\[\rho\to\infty\] as \[|z|\to 1\], which proves (c). The proof of (b), and
hence of (a) is now straightforward.\QED\smallskip

The \Poincare metric
has the marvelous property of never increasing under holomorphic maps.

{\QP {\bf 1.10. Theorem of Pick.} {\it If \[f:S\to T\] is a holomorphic
map between
Riemann surfaces, both of which are conformally isomorphic to \[D\], then
$$	  \PD_{T}(f(z_1), \,f(z_2))\;\le\;\PD_{S}(z_1,\,z_2) \,.	$$
Furthermore, if equality holds for some
\[z_1\ne z_2\] in \[S\], then \[f\] must be
a conformal isomorphism from \[S\] onto \[T\].}\medskip}

{\bf Proof.} Join \[z_1\] to \[z_2\] by a geodesic of length equal to
the distance \[\PD_{S}(z_1\,,
\,z_2)\]. Let \[ds\] be the element of \Poincare length along this
geodesic, and let \[ds'\] be the element of length along the image curve in
\[T\]. It follows from the Schwarz Lemma and the definition of the
\Poincare metric that \[|ds'/ds|\le 1\], with equality if and only if \[f\]
is a conformal isomorphism; and the conclusion follows.\QED

Note that this distance depends sharply on the choice of \[S\].
As an example, suppose that \[U\subset D\] is a
simply connected open subset with \[U\ne D\]. Then \[U\] is conformally
isomorphic to \[D\] by the Riemann Mapping Theorem. Applying 1.10 to the
inclusion\break \[U\to D\], we see that
$$	\PD_D(z_1\,,\,z_2)\;<\;	\PD_U(z_1\,,\,z_2)	$$
for any two distinct points \[z_1\] and \[z_2\] in \[U\].\smallskip

{\bf 1.11. Remark.} In the special case of a holomorphic map from \[S\]
to itself, if the sharper inequality
$$	\PD(f(z_1),f(z_2))\le c\PD(z_1\,,\,z_2)	$$
is satisfied for all \[z_1\] and \[z_2\], where \[0<c<1\] is constant,
then \[f\] necessarily has a unique fixed point. However,
the example \[z\mapsto z^2\] on the unit disk shows that a map with a
unique fixed point need not satisfy this sharper inequality; and
the example \[w\mapsto w+i\] on the upper half-plane shows that map which
simply reduces \Poincare distance need not have any fixed point. (See also
Problem 1-1.)\medskip

Next let us consider the {\bit Riemann sphere} \[\hat\C\],
that is the compact Riemann surface consisting of the complex
numbers \[\C\] together with a single point at infinity. The conformal
structure on this complex manifold can be specified by using
the usual coordinate \[z\] as uniformizing parameter in the finite plane
\[\C\], and using \[\zeta=1/z\] as uniformizing parameter in \[\hat\C\ssm
\{0\}\].

When studying \[\hat\C\], we need some {\bit spherical metric} \[\SD\]
which is adapted to the topology, so that the point at infinity has
finite distance from other points of \[\hat\C\]. The precise choice of such
a metric does not matter for our purposes. However, one good choice
would be the distance function \[\SD(z_1\,,\,z_2)\] which is associated with
the Riemannian metric
$$	ds=2|dz|/(1+|z|^2)\,.	\eqno (6)	$$
This metric is smooth and well behaved, even in a neighborhood
of the point at infinity, and has constant Gaussian curvature \[+1\].
It corresponds to the standard metric on the unit sphere in \[\R^3\]
under stereographic projection.
Note that the map \[z\mapsto 1/z\] is an isometry for this metric.
(However, it is certainly not true that every conformal self-map of \[\hat\C\]
is an isometry.)\smallskip

The group \[G(\hat\C)\] consisting of all conformal automorphisms of \[\hat
\C\] can be described as follows. By definition,
an automorphism \[g\in G(\hat\C)\] is called an {\bit involution}
if the composition \[g\circ g\] is the identity map of \[\hat\C\].

{\QP{\bf 1.12. Theorem.} \it 
Every conformal automorphism \[g\] of \[\hat\C\] can be
expressed as a {\bit fractional linear transformation} or {\bit M\"obius
transformation}
$$	g(z)\;=\; (az+b)/(cz+d)\,,	$$
where the coefficients are complex numbers with determinant \[ad-bc\ne 0\].
Every non-identity automorphism of \[\hat\C\] either has two distinct fixed
points or one double fixed point in \[\hat\C\]. In general,
two non-identity elements of \[G(\hat\C)\] commute if and only if they have
precisely the same fixed points: the only exceptions to this statement
are provided by pairs of commuting involutions.\medskip}

Here by a ``double'' fixed point we mean one at which the derivative \[g'(z)\]
is equal to \[+1\].
If we normalize so that \[ad-bc=1\], then the coefficients are
well defined up to a simultaneous change of sign.
{\it Thus the group \[G(\hat\C)\] of conformal
automorphisms can be identified with the group \[\PSL(2,\C)\] consisting of
all complex matrices with determinant \[+1\] modulo the subgroup \[\{\pm
I\}\].}\smallskip

{\bf 1.13. Remark.} The group \[G(\C)\] of all conformal automorphisms
of the complex plane can be identified with the subgroup of \[G(\hat\C)\]
consisting of automorphisms \[g\] which fix the point \[\infty\], since
every conformal automorphism of \[\C\] extends uniquely
to a conformal automorphism of \[\hat\C\]. (Compare Ahlfors [A1, p.124] .)
It follows easily that \[G(\C)\] consists
of all affine transformations
$$	g(z)\;=\;\lambda z+c	 $$
with complex coefficients \[\lambda\ne 0\] and \[c\].\smallskip

{\bf Proof of 1.12.} Clearly \[G(\hat\C)\] contains this group of fractional
linear transformations as a subgroup. After composing the given \[g\in G(\hat
\C)\] with a suitable element of this subgroup, we may assume that \[g(0)=0\]
and \[g(\infty)=\infty\]. But then the quotient \[g(z)/z\] is a bounded
holomorphic function, hence constant by Liouville's Theorem (\S1.3). Thus
\[g\] is linear, and hence itself is an element of \[\PSL(2,\C)\].

The fixed points of a fractional linear transformation can be
determined by solving a quadratic equation, so it is easy to check that there
must be at least one and
at most two distinct solutions in the extended plane \[\hat\C\].
In particular, if an automorphism of \[\hat\C\] 
fixes three distinct points, then it must be the identity map.

In general, if two automorphisms \[g\] and \[h\] commute, we
can say that \[g\] maps every fixed point of \[h\] to a fixed point of \[h\],
and that \[h\] maps every fixed point of \[g\] to a fixed point of \[g\].
However, this leaves open the possibility that \[g\] interchanges the two fixed
of \[h\] and that \[h\] interchanges the two fixed points of \[g\].
If this is indeed the case, then both \[g\circ g\]
and \[h\circ h\] have at least four fixed points, hence both must equal
the identity map. (An example of this phenomenon is provided by the two
commuting involutions \[g(z)= -z\] with fixed points \[\{0,\infty\}\], and
\[h(z)= 1/z\] with fixed points \[\{\pm 1\}\].) If we exclude this possibility,
then elements which commute must have exactly the same fixed points.

Conversely, let us consider the subgroup consisting of all
elements of \[G(\hat\C)\] which fix two specified points \[z_0\ne z_1\].
It is convenient to conjugate by an automorphism
which carries \[z_0\] to zero and \[z_1\] to infinity. The argument above
shows that an automorphism \[g\] fixes both zero and infinity if and only
if it has the form \[g(z)=\lambda z\] for some \[\lambda\ne 0\]. Thus the
subgroup consisting of all such elements is isomorphic to the multiplicative
group \[\C\ssm\{0\}\], and hence is commutative as required.

Finally, consider automorphisms \[g\] which fix only the point at infinity.
A similar argument shows that \[g(z)\] must be equal to \[z+c\]
for some constant \[c\ne 0\]. (Compare 1.13.)
These automorphisms, together with the identity map, form a subgroup of
\[G(\hat\C)\] which is isomorphic to the additive group of \[\C\]; which again
is commutative as required. Further details will be left to the reader.\QED

The corresponding discussion for the group \[G(H)\] (or \[G(D)\])
will be important in \S3. To fix our ideas, let
us concentrate on the case of the half-plane. It follows from 1.7 that every
automorphism of \[H\] extends uniquely to an automorphism of \[\hat
\C\]. Hence 1.12 applies immediately. However, we must consider not only
fixed points inside of \[H\] but also fixed points which lie on the
boundary \[\partial H=\R\cup\infty\]. A priori, we should also consider
fixed points which lie in the lower half-plane, completely outside of the
closure \[\bar H\]. However, it is easy to check that \[w\] is a fixed point
of a fractional linear transformation with real coefficients
if and only if the complex
conjugate \[\bar w\] is also a fixed point. Thus each fixed point in the lower
half-plane, outside of \[\bar H\], is paired with a fixed point
which is inside the upper half-plane \[H\].\smallskip

{\bf 1.14. Definition.} The non-identity automorphisms of \[H\] fall into
three classes, as follows:

An automorphism \[g\in G(H)\] is said to be {\bit elliptic} if it has a
fixed point \[w_0\in H\]. We may also describe \[g\] as a {\bit rotation}
around \[w_0\]. (Compare (1) above.)

The automorphism \[g\in G(H)\] is {\bit hyperbolic} if it has two distinct
fixed points on the boundary \[\partial H\]. As an example,
\[g\] fixes the two points \[0\]
and \[\infty\] if and only if it has the form \[g(w)=\lambda w\] with
\[\lambda>0\].

The automorphism \[g\in G(H)\]
is {\bit parabolic} if it has just one double fixed
point, which necessarily belongs to the boundary \[\partial H=\R\cup
\infty\]. For example, if this double fixed point is the point at infinity,
then \[g\] is necessarily a translation: \[g(w)=w+c\;\] for some constant
\[c\ne 0\] in \[\R\].

{\QP{\bf 1.15. Lemma.} \it Two non-identity elements of \[G(H)\] commute if
and only if they have exactly
the same fixed points in \[\bar H=H\cup\R\cup\infty\]. The set
of all group elements with some specified
fixed point set, together with the identity
element, forms a commutative subgroup, which is isomorphic to a circle in the
elliptic case and to a real line in the parabolic or hyperbolic cases.
\medskip}

The proof is easily supplied.\QED

Again, we could equally well work with \[D\] in place of \[H\].
One defect of this exposition is that it requires going outside of the
half-plane \[H\] in order to distinguish between parabolic and hyperbolic
automorphisms. For a more intrinsic description of the difference between
these cases, see Problem 1-5.

\lln\medskip

We conclude this section with a number of problems for the reader.\medskip

{\bf Problem 1-1.} If a holomorphic map \[f:D\to D\] fixes the
origin, and is not a rotation, prove that the successive images \[f^
{\circ n}(z)\] converge to zero for all \[z\] in the open disk \[D\],
and prove that this convergence is uniform on compact subsets of \[D\].
(Here \[f^{\circ n}\] stands for the \[n$-fold iterate \[f\circ\cdots\circ
f\]. The example \[f(z)=z^2\] shows that convergence need not be uniform
on all of \[D\].) \smallskip

{\bf Problem 1-2.} Show that the group \[G(D)\] of conformal
automorphisms of the unit disk \[D\] consists of all maps
$$	g(z)\;=\; e^{i\theta} (z-a)/(1-\bar a z)	$$
with \[|e^{i\theta}|=1\], where \[a\in D\] is the unique point which maps
to zero. (Compare 1.7.) Check that \[|g'(0)|\le 1\], and conclude using 1.2
that any holomorphic map \[f:D\to D\] must satisfy \[|f'(0)|\le 1\].
\smallskip

{\bf Problem 1-3.} Show that the action of the group
\[G(\hat\C)\] on \[\hat\C\]
is {\bit simply 3-transitive}. That is, there is one and only one automorphism
which carries three distinct specified points of \[\hat\C\]
into three other specified
points. Similarly, show that the action of of \[G(\C)\] on \[\C\] is simply
2-transitive. Show that the action of \[G(D)\] on the boundary circle
\[\partial D\] carries three specified points into three other specified
points if and only if they have the same cyclic order; and show that the
action of \[G(D)\] carries two points of \[D\] into two other specified
points if and only if they have the same \Poincare distance.\smallskip

%{\bf Problem 1-4.} Using polar coordinates \[z=re^{i\theta}\] on the unit
%disk, show that the \Poincare metric is given by
%$$	ds=2|dz|/(1-r^2) =2\sqrt{(dr^2+r^2d\theta^2)}/(1-r^2)\ge 2|dr|
%/(1-r^2)\,,	$$
%where equality holds along any radial line \[\theta={\rm constant}\].
%Conclude that the associated distance function satisfies
%\[\PD(0,\, re^{i\theta})=\log((1+r)/(1-r))\]. In particular,
%show that this expression tends to infinity as \[r\to 1\], and hence check
%directly that the
%boundary of \[D\] has infinite \Poincare distance from any interior point, and
%that the closed neighborhood \[N_r(0,\rho_D)\] is a compact subset of \[D\].
%(\S1.8.)\smallskip

{\bf Problem 1-4.} By an {\bit anti-holomorphic} automorphism of \[\hat\C\]
we mean an orientation reversing self-homeomorphism of the form \[z\mapsto
\eta(\bar z)\], where \[\eta\] is holomorphic.
If \[L\] is a straight line or circle in
\[\hat\C\], show that there is one and only one anti-holomorphic involution
of \[\hat\C\] having \[L\] as fixed point set, and show that no other
non-vacuous fixed point sets can occur. Show that the automorphism group
\[G(\hat\C)\] acts transitively on the set of straight lines and circles
in \[\hat\C\]. Similarly, show that any
anti-holomorphic involution of \[D\] is the reflection in some \Poincare
geodesic; and show that \[G(D)\] acts transitively on the set of such
geodesics.\smallskip

{\bf Problem 1-5.}  Let \[g\] be an automorphism of \[D\] with \[g\circ g\]
not the identity map. Show that \[g\] is hyperbolic
if and only if there exists an automorphism \[h\] which satisfies
$$	h\circ g\circ h^{-1}\;=\;g^{-1}\,,	$$
or if and only if there exists some necessarily unique geodesic \[L\] with
respect to the \Poincare metric which is mapped onto itself by \[g\],
or if and only if \[g\] commutes with some anti-holomorphic involution.
(The possible choices for \[h\] are just the \[180
^\circ\] rotations about the points of \[L\].)\smallskip

{\bf Problem 1-6.} Define the {\bit infinite band} \[B\subset\C\] of width
\[\pi\] to be the set
of all \[z=x+iy\] with \[|y|<\pi/2\]. Show that the exponential map carries
\[B\] isomorphically onto the right half-plane \[\{u+iv:u>0\}\]. Show that
the \Poincare metric on \[B\] takes the form
$$	ds=|dz|/cos\,y\,. \eqno (7)	$$
Show that the
real axis is a geodesic whose \Poincare arclength coincides with its usual
Euclidean arclength, and show that each real translation \[z\mapsto z+c\] is
a hyperbolic automorphism of \[B\] having the real axis as its unique
invariant geodesic.\smallskip
\eject

{\bf Problem 1-7.} Given four distinct points \[z_j\] in
\[\hat\C\] show that the {\bit cross ratio}
$$	\chi(z_1\,,\,z_2\,;\,z_3\,,\,z_4)
	\;\;=\;\;{(z_1-z_3)(z_2-z_4)\over(z_1-z_4)(z_2-z_3)}\;\;\in
	\;\;\C\ssm\{0,1\} $$
is invariant under fractional linear transformations, and show that \[\chi\]
is real if and
only if the four points lie on a straight line or circle. \smallskip

{\bf Problem 1-8.} Show that each \Poincare
neighborhood \[N_r(w_0\,,\,\rho_H)\] in the upper half-plane is bounded
by a Euclidean circle, but that \[w_0\] is not its Euclidean center. Show that
each \Poincare geodesic in the upper half-plane is a straight line or
semi-circle which meets the real axis orthogonally. If the geodesic 
through \[w_1\] and \[w_2\] meets \[\partial H=\R\cup\infty\] at the points
\[\alpha\] and \[\beta\], show that the \Poincare distance between \[w_1\]
and \[w_2\] is equal to the logarithm of the cross ratio
\[\chi(w_1\,,\,w_2\,;\,\alpha\,,\,\beta)\]. Prove corresponding
statements for the unit disk.\smallskip

{\bf Problem 1-9.} The Gaussian curvature of a conformal metric \[ds=
\gamma(w)|dw|\] with \[w=u+iv\] is given by the formula
$$	 K\;\;=\;\; {\gamma_u^2+\gamma_v^2 -\gamma(\gamma_{uu}+\gamma_{vv})
	\over \gamma^4}	$$
where the subscripts stand for partial derivatives. (Compare Willmore, p. 79.)
Check that the \Poincare metrics (3), (4) and (7) have curvature \[K\equiv-1\],
and more generally that the metric
\[ds=c|dw|/v\] has curvature \[K\equiv-1/c^2\].
Check that the spherical metric (6) has curvature \[K\equiv +1\].\smallskip

{\bf Problem 1-10.} Classify conjugacy classes in the groups \[G(H)\cong
\PSL(2,\R)\] as follows.
Show that every automorphism of \[H\] without fixed point is conjugate to
a unique transformation of the form \[w\mapsto w+1\] or \[w\mapsto w-1\] or
\[w\mapsto\lambda w\] with \[\lambda>1\]; and show that the conjugacy class of
an automorphism \[g\] with fixed point \[w_0\in H\] is uniquely determined
by the derivative \[\lambda=g'(w_0)\], where \[|\lambda|=1\].
Show also that each non-identity element of \[\PSL_2(\R)\] belongs
to one and only one ``one-parameter subgroup'', and that each
one-parameter subgroup is conjugate to either
$$      t\;\;\mapsto\;\;
        \left[\matrix{1&t\cr 0&1\cr}\right]\quad{\rm or}\quad
        \left[\matrix{e^t&0\cr 0&e^{-t}\cr}\right]\quad{\rm or}\quad
        \left[\matrix{\cos t&\sin t\cr -\sin t&\cos t\cr}\right]      $$
according as its elements are parabolic or hyperbolic or elliptic.
Here \[t\] ranges over the additive group of real numbers.\smallskip

{\bf Problem 1-11.} For a non-identity automorphism \[g\in G(\hat\C)\], show
that the derivatives \[g'(z)\] at the two
fixed points are reciprocals, say \[\lambda\] and \[\lambda^{-1}\],
and show that the sum \[\lambda+\lambda^{-1}\] is a complete conjugacy
class invariant. (Here \[\lambda=1\] if and only if the two fixed points
coincide. In the
special case of a fixed point at infinity, one must evaluate the
derivative using the local coordinate \[\zeta=1/z\].)\smallskip

{\bf Problem 1-12.} Show that the conjugacy class of
a non-identity automorphism \[\;g(z)=\lambda z+c\;\]
in the group \[G(\C)\] is uniquely
determined by its image under the homomorphism \[g\mapsto \lambda\in\C\ssm
\{0\}\].

\vfil
%\rightline{jwm, version of 8-28-91}
\eject
% \bigskip\bigskip
\footline={\hss\tenrm 2-\folio\hss}\pageno=1

\centerline{\bf \S2. The  Universal Covering, Montel's Theorem.}\bigskip

If \[S\] is a completely arbitrary Riemann surface, then the universal
covering \[\tilde S\] is a well defined simply connected Riemann surface,
with a canonical projection map \[p:\tilde S\to S\].
(Compare Munkres, and also Appendix E.) According to the Uniformization
Theorem, this universal covering \[\tilde S\] must be
conformally isomorphic to one of the three model surfaces (\S1.1).
Thus we have the following.

{\QP{\bf 2.1. Lemma.} \it Every Riemann surface \[S\] is conformally
isomorphic to a quotient of the form \[\tilde S/\Gamma\], where
\[\tilde S\] is a simply connected surface which is isomorphic to
either \[D\], \[\C\], or \[\hat\C\]. Here \[\Gamma\] %subset G(\tilde S)\]
is a discrete group of conformal automorphisms of \[\tilde S\],
such that every non-identity element of \[\Gamma\] acts without fixed points
on \[\tilde S\].\medskip}

This discrete group \[\Gamma\subset G(\tilde S)\] can be identified with
the fundamental group \[\pi_1(S)\]. The elements of \[\Gamma\]
are called {\bit deck transformations}. They can be characterized as
maps \[\gamma:\tilde S\to\tilde S\] which satisfy the identity
\[\gamma=p\circ\gamma\], that is maps for which the diagram
$$      \matrix{\tilde S&\buildrel \gamma\over\longrightarrow& \tilde S\cr
       \;\;\;\downarrow p &&\;\;\;\downarrow p\cr
        S&\buildrel = \over\longrightarrow& S}  $$
is commutative.
Conversely, if we are given a group \[\Gamma\] of conformal automorphisms
of a simply connected surface \[\tilde S\] which is discrete as a subgroup of
\[G(\tilde S)\], and such that every non-identity element of \[\Gamma\] acts
without fixed points,
then it is not difficult to check that \[\Gamma\] is the group
of deck transformations for a universal covering map \[\tilde S\to\tilde
S/\Gamma\]. (Compare Problem 2-1.)\smallskip

We can analyze the three possibilities for \[\tilde S\] as follows.\smallskip
% (Note that (1) implies the sharper condition that \[\Gamma\] acts
% properly discontinuously --- that is, for any compact subset \[K\subset\tilde
% S\] there are only finitely many \[\gamma\] with \[\gamma(K)\cap K\ne
% \emptyset\].)

{\bf Spherical Case.} According to 1.12, every conformal automorphism of the
Riemann sphere \[\hat\C\] has
at least one fixed point. Therefore, if \[S\cong\hat\C/\Gamma\]
is a Riemann surface with universal covering
\[\tilde S\cong \hat\C\], then the group \[\Gamma\subset G(\hat\C)\] must be
trivial, hence \[S\] itself must be isomorphic to \[\hat\C\].\smallskip

{\bf Euclidean Case.} By 1.13, the group \[G(\C)\] of conformal
automorphisms of the complex plane
consists of all affine transformations \[z\mapsto \lambda z+c\] with
\[\lambda\ne 0\]. Every such transformation with \[\lambda
\ne 1\] has a fixed point.
Hence, if \[S\cong\C/\Gamma\] is a surface with universal covering
\[\tilde S\cong \C\], then \[\Gamma\] must be a discrete 
group of translations \[z\mapsto z+c\] of the complex plane \[\C\].
There are three subcases:

If \[\Gamma\] is trivial, then \[S\] itself is isomorphic to \[\C\].

If \[\Gamma\] has just one generator, then \[S\] is isomorphic to
the cylinder \[\C/\Z\], which in turn is isomorphic under the exponential
map to the {\bit punctured plane} \[\C\ssm\{0\}\].

If \[\Gamma\] has two generators, then it can be described as a two-dimensional
{\bit lattice}\break
\[\Lambda\subset\C\], that is an additive group generated by
two complex numbers which are linearly independent over \[\R\]. The
quotient \[{\bf T}=\C/ \Lambda\] is called a {\bit torus}.\smallskip

In all three of these cases, note that our surface inherits a locally
Euclidean geometry from the Euclidean metric \[|dz|\] on its universal
covering surface. For example the punctured plane \[\C\ssm\{0\}\], consisting
of points \[\exp(z)=w\], has a complete locally Euclidean metric \[|dz|
=|dw/w|\]. It is easy to check that such a metric is unique up to a
multiplicative constant. (Compare Problem 2-2.) We will refer to these
Riemann surfaces as [complete locally] {\bit
Euclidean surfaces}. The term ``{\bit parabolic surface\/}'' is also commonly used
in the literature.\smallskip

{\bf Hyperbolic Case.} In all other cases, the universal covering \[\tilde
S\] must be conformally isomorphic to the unit disk. Such Riemann surfaces
are said to be {\bit Hyperbolic}. {\it As an example, any Riemann surface
with non-abelian fundamental group, and in particular any
surface of higher genus, is necessarily Hyperbolic.}

{\bf Remark.} Here the word ``Hyperbolic'' is a reference to Hyperbolic
Geometry, that is the geometry of Lobachevsky.
Unfortunately the term ``hyperbolic'' has at least three
different well established meanings in conformal dynamics. (Compare \S1.14
and \S14.) In a crude attempt to avoid confusion, I will alway capitalize the
word when it is used with this geometric meaning.\smallskip

Every Hyperbolic surface \[S\] possesses a unique {\bit \Poincare metric},
which is complete, with Gaussian curvature identically equal to \[-1\].
To construct this metric, we note that the \Poincare metric on \[\tilde S\]
is invariant under the action of \[\Gamma\]. Hence there is one and only one
metric on \[S\] so that the projection \[\tilde S\cong D\to S\]
is a local isometry. Hence, just as in \S1.9, there is an
associated {\bit \Poincare distance function} \[\rho(z,z')=\rho_S(z,z')\] which
is equal to the length of a shortest geodesic from \[z\]to \[z'\]. As in 1.10
we have:

{\QP{\bf 2.2. Lemma.} \it If \[f:S\to T\] is a holomorphic map
between Hyperbolic surfaces, then
$$	 \PD_T(f(z),\,f(z'))\;\le\;\PD_S(z,z')\,.	$$
Furthermore, if equality holds for some \[z\ne z'\], then it follows
that \[f\] is a {\bit local isometry}. That is, \[f\] preserves the
infinitesimal distance element \[ds\], and hence preserves the distances
between nearby points.\medskip}

{\bf Proof.} This follows, just as in 1.10, since a minimal geodesic
joining \[z\] to \[z'\] must be mapped isometrically.\QED

{\bf Caution:} It no longer follows that \[f\] must be a conformal
isomorphism from \[S\] to \[T\]. However,
if \[f\] is a local isometry, then we can at least assert that \[f\] induces
an isometry \[\tilde S\buildrel\cong\over\longrightarrow\tilde T\] between
the universal covering surfaces, and hence that \[f\] is a covering map
from \[S\] onto \[T\].\smallskip

Here is an example. The map \[f(z)= z^2\] on the
{\bit punctured disk} \[D\ssm\{0\}\] is certainly not an automorphism, since
it is
two-to-one. However, the universal covering of \[D\ssm\{0\}\] can be identified
with the left half-plane, mapped to \[D\ssm\{0\}\] by the exponential map,
and \[f\] lifts to the automorphism \[w\mapsto 2w\] of this left half-plane.
Therefore, \[f\] is a covering map, and
preserves the \Poincare metric locally.\smallskip

Note that the punctured disk has abelian fundamental group \[\pi_1(D\ssm\{0\})
\cong\Z\]. Here is another such example.
For any \[ r>1\], the {\bit annulus}
$$	{\cal A}_r\;=\;\{\,z\; : \; 1<|z|<r\,\}		$$
is a Hyperbolic surface, since it admits a holomorphic map to the unit disk.
The fundamental group \[\pi_1({\cal A}_r)\] is evidently also free cyclic.
In fact annuli and the punctured disk
are the only Hyperbolic surfaces with abelian fundamental
group, other than the disk itself. (Compare Problem 2-3.)\smallskip

{\bf Maximal Hyperbolic Example.} If \[a_1\,,\,a_2\,,
\,a_3\] are three distinct points of \[\hat\C\], then the complement
\[\;\Sigma_3\,=\,\hat\C\ssm\{a_1\,,\,a_2\,,\,a_3\}\,\cong\,\C\ssm\{0,1\}\;\]
is called a {\bit thrice punctured sphere}. This is evidently a Hyperbolic
surface, for example since its fundamental group is not abelian.
One immediate corollary of this observation is the following.

{\QP {\bf 2.3. Picard's Theorem.}
{\it Any holomorphic map from \[\C\] to \[\hat\C\] which omits three
different values must be constant. More generally, if the Riemann surface
\[S\] admits some non-constant
holomorphic map to the thrice punctured
sphere \[\Sigma_3\], then \[S\] must be Hyperbolic.}\medskip}

For this map \[f:S\to \Sigma_3\] can be
lifted to a holomorphic map from the universal covering \[\tilde S\]
to the universal covering \[\tilde \Sigma_3\cong D\]. By Liouville's
Theorem (\S1.3) it follows that \[\tilde S\cong D\].\QED

Let \[U\subset \hat\C\] be any connected open set which omits at least
three points, and hence is Hyperbolic.
It is often useful to compare the \Poincare distance between two points
of \[U\] with the spherical distance between the same two points within
\[\hat\C\]. (See \S1(6).) Here is a crude estimate.
Again let \[N_r(z,\rho_U)\] be the closed neighborhood of
some fixed radius \[r\] with respect to the \Poincare metric about
the point \[z\] of \[U\].

{\QP {\bf 2.4. Lemma.} {\it 
As \[z\] converges towards the boundary \[\partial U\] in the spherical
metric,  the spherical diameter
of the neighborhood \[N_r(z,\rho_U)\] tends to zero.}\medskip}

\noindent (For a sharper statement in the simply connected case, see
A.8 in the Appendix.) \smallskip

{\bf Proof of 2.4.}
First consider the special case \[U=\Sigma_3=\hat\C\ssm\{a_1\,,\,
a_2\,,\,a_3\}\]. Choose some fixed base point \[z_0\] in \[U=\Sigma_3\].
For each fixed \[r>0\], the neighborhood \[N_r(z_0,\,\rho_U)\] is the image
under projection of a corresponding neighborhood in the universal covering
surface, and hence is compact and connected. (Compare 1.9.)
Now, as \[z\] tends to one of the three boundary points \[a_1\] of \[U\],
it must eventually
leave any compact subset of \[U\], hence the distance \[\PD_U(z,\,z_0)\] must
tend to infinity. For fixed \[r\], it follows that the neighborhood
\[N_r(z,\,\rho_U)\] will eventually be
disjoint from any given compact subset of \[\Sigma_3\]. In fact, since the set
\[N_r(z,\,\rho_U)\] is connected, this entire set must tend to just one
boundary point \[a_1\] with respect to the spherical metric.

Now consider an arbitrary Hyperbolic open set \[U\subset\hat\C\]. Given a
sequence of points of \[U\] tending to the boundary, we can choose a
subsequence \[\{z_j\}\] which converges to a single boundary point \[a_1\].
Choose two other boundary points \[a_2\]
and \[a_3\], and consider the inclusion map
from \[U\] to \[\Sigma_3=\hat\C\ssm\{a_1\,,\,a_2\,,\,a_3\}\]. Applying the
Pick inequality 2.2, we have \[N_r(z_j\,,\,\rho_U)\subset N_r(z_j\,,\,\rho_
{\Sigma_3})\], and it follows from the discussion above that
this entire neighborhood converges to the boundary point \[a_1\]
as \[j\to\infty\].\QED\smallskip

Using the \Poincare metric, we will develop another important tool.
A sequence of maps \[f_n:S\to\hat\C\] on a Riemann surface
\[S\] is said to {\bit converge locally uniformly} (or\break {\bit uniformly
on compact sets}) to the limit \[g:S\to\hat\C\] if for
every compact subset \[K\subset S\] the\eject
\noindent sequence
\[\{f_n|K\}\] of maps \[f_n\] restricted to \[K\] converges
uniformly to \[g|K\]. Here it is to be understood that we use the
spherical metric \[\SD(z_1\,,\,z_2)\] on the target space \[\hat\C\].
\smallskip

{\bf Definition.} Let \[S\] and \[S'\] be Riemann surfaces, with \[S'\]
compact. A collection \[{\cal F}\] of holomorphic maps
\[f_\alpha:S\to S'\] is {\bit normal}\[\] if every infinite
sequence of maps from \[{\cal F}\] contains a subsequence which converges
locally uniformly to a limit.\smallskip

(For the case of a non-compact target space \[S'\], see 2.6 below.)
Note that the limit function \[g\] must itself be holomorphic, by the Theorem
of Weierstrass. However, this limit \[g\] need not belong to the given family.
Roughly speaking, a family of maps is normal if and only if its closure,
in the space of all holomorphic maps from \[S\] to \[S'\], is a compact
set. (Ahlfors [A1, p.213].)

{\QP {\bf 2.5. Montel's Theorem.} {\it Let \[S\] be any Riemann surface,
which we may as well suppose to be Hyperbolic. If a collection \[{\cal F}\]
of holomorphic maps from \[S\] to the Riemann sphere \[\hat\C\] takes values
in some Hyperbolic open subset \[U\subset\hat\C\], or equivalently
if there are three distinct points of \[\hat\C\] which never occur as
values, then this collection \[{\cal F}\] is normal.} \medskip}

More explicitly, any sequence of holomorphic maps \[f_n:S\to U\]
contains a subsequence which converges, uniformly on compact subsets of \[S\],
to some holomorphic map \[g:S\to\bar U\]. Here it is essential that \[g\] be
allowed to take values in the closure \[\bar U\]. However, the proof will
show that the image \[g(S)\] is either contained in \[ U\], or is a single
point belonging to the boundary of \[ U\].\smallskip

{\bf Proof of 2.5.} First note by
2.3 that the surface \[S\] must be Hyperbolic, unless all of our maps
are constant. Hence \[S\] also has a \Poincare metric.
Choose a countable dense subset of points \[z_j\in S\,,\;j\ge 1\].
(It follows easily from \S1.1 that every Riemann surface possesses a countable
dense subset. This statement was first proved by Rad\'o. Compare
Ahlfors \& Sario.) Given any sequence
of holomorphic maps \[f_n:S\to U\subset\hat\C\], we can first choose an
infinite subsequence \[\{f_{n(p)}\}\] of the \[f_n\] so that the
images \[f_{n(p)}(z_1)\] converge to a limit within the closure of
\[U\]. Then choose a sub-sub-sequence
\[f_{n(p(q))}\] so that the \[f_{n(p(q))}(z_2)\] also converge to a limit,
and continue inductively. By a diagonal procedure,
taking the first element of the first subsequence, the second element of the
second subsequence, and so on, we construct a new infinite sequence of maps
\[g_m=f_{n_m}\] so that \[\lim_{m\to\infty}g_m(z_j)\in \bar U\] exists
for every choice of \[z_j\].

{\bf Case 1.} Suppose that every one of these limit points in \[\bar U\]
actually belongs to the set \[U\] itself. Given any compact set \[K\subset S\]
and any \[\epsilon>0\], we can choose finitely many \[z_j\] so that every
point \[z\in K\] has Poincar\'e distance \[\rho_S(z,z_j)<\epsilon\] from one of
these \[z_j\]. Further, we can choose \[m_0\] so that 
\[\rho_U(g_m(z_j)\,,\,g_n(z_j))<\epsilon\] for each of these finitely many
\[z_j\], whenever \[m,\,n>m_0\,\]. For any \[z\in K\] it then follows using
2.2 that \[\PD_U(g_m(z)\,,\,g_n(z))<3\epsilon\] whenever \[m,\,n>m_0\].
Thus the \[g_m(z)\] form a Cauchy sequence.
It follows that the sequence of functions \[\{g_m\}\] converges
to a limit, and that this convergence is uniform
on compact subsets of \[S\].

{\bf Case 2.} Suppose on the other hand that \[\lim_{m\to\infty}
g_m(z_j)\] is actually a boundary point \[a_0\in\partial U\] for some \[z_j\].
Then it follows
from 2.4 that \[g_m(z)\] converges to \[a_0\] for every \[z\in S\],
and that this convergence is uniform on compact subsets of \[S\].\QED
 \smallskip

{\bf 2.6. Remark.} If we consider maps \[S\to S'\] where the target space
\[S'\] is not compact, then the definition should be modified as follows.
We continue to assume that \[S\] is connected.
A collection \[\cal F\] of maps \[f_\alpha:S\to S'\] is {\bit normal~}
if every infinite sequence of maps from \[\cal F\] contains a subsequence
which either

\noindent $(1)$  converges locally uniformly
to a holomorphic map from \[S\] to \[S'\], or

\noindent $(2)$ diverges locally uniformly to infinity, in the sense that
the successive images of any compact subset of \[S\] eventually miss any given
compact subset of \[S'\].\smallskip

\lln

{\bf Concluding Problems.}\smallskip

{\bf Problem 2-1.} Let \[S\] be a simply connected Riemann surface, and let
\[\Gamma\subset G(S)\] be a discrete subgroup; that is, suppose that the
identity element
is an isolated point of \[\Gamma\]. If every non-identity element of
\[\Gamma\] acts on \[S\] without fixed points, show that the action of
\[\Gamma\] is
{\bit properly discontinuous}. That is, for every compact \[K\subset S\]
show that
only finitely many group elements satisfy \[K\cap\gamma(K)\ne
\emptyset\]. Show that each \[z\in S\] has a neighborhood \[U\] whose
translates \[\gamma(U)\] are pairwise disjoint. Conclude that \[S/\Gamma\]
is a well defined Riemann surface with \[S\] as its universal covering.
\smallskip

{\bf Problem 2-2.} Show that any Riemann surface which
is not conformally isomorphic to \[\hat\C\],
has one and up to a multiplicative constant only
one conformal Riemann metric which is complete, with constant Gaussian
curvature. (Make use of Hopf's Theorem which asserts that, for each
real number \[K\], there is one and only
one complete simply connected surface of constant curvature \[K\] up to
isometry. See Willmore p. 162.) On the other hand, show that \[\hat\C\]
has a 3-dimensional family of conformal metrics with curvature \[+1\].
\smallskip

{\bf Problem 2-3.} Using Problems 1-10 and 1-6, show that every Hyperbolic
surface
with abelian fundamental group is conformally isomorphic either to the disk
\[D\], or the punctured disk \[D\ssm\{0\}\], or to the annulus
$$	{\cal A}_r\;=\;\{z\in\C:1<|z|<r\}	$$
for some uniquely defined \[ r>1\]. Define the {\bit modulus} of
this annulus to be the number \[\mod\!({\cal A}_r)
=\log r/2\pi>0\].
Show that this annulus has a unique simple closed geodesic, which has
length \[\ell=\pi/\mod\!({\cal A}_r)=2\pi^2/\log r\,\]. On the other hand,
show that the punctured disk \[D\ssm\{0\}\] has no closed geodesic.
This punctured disk, which is conformally isomorphic to the set
\[{\cal A}_\infty\;=\;\{z\in\C:1<|z|<\infty\}\], is sometimes described
as an annulus of infinite modulus. (However this designation is ambiguous,
since the Euclidean surface \[\C-\{0\}\]
might also be described as an annulus of infinite modulus.)
\smallskip %By definition, this punctured disk has infinite modulus.\smallskip

{\bf Problem 2-4.} If \[S\] and \[S'\] are Hyperbolic Riemann surfaces
(not necessarily compact), show that {\bit every~} family of maps from
\[S\] to \[S'\] is normal.\smallskip

{\bf Problem 2-5.}
Show that normality is a local property. More precisely,
let \[S\] and \[S'\] be any Riemann surfaces,
and let \[\{f_\alpha\}\] be a family of holomorphic maps from \[S\] to \[S'\].
If every point of \[S\] has a neighborhood \[U\] such
that the collection \[\{f_\alpha|U\}\] of restricted maps
is normal, show by a diagonal argument as in the
proof of 2.5 that \[\{f_\alpha\}\] is normal.
\bigskip\vfill

%\rightline{jwm, version of 8-28-91}
\end

%% file: cd3-5.tex
\input header
\input ras

\def\GO{{G\!O}}
\def\TF{{\cal E}}
\def\T{{\bf T}}

\footline={\hss\tenrm 3-\folio\hss}\pageno=1

\centerline{\bf THE JULIA SET}\medskip

\centerline{\bf \S3. Fatou and Julia: Dynamics on the Riemann Sphere.}\bigskip

The local study of iterated holomorphic mappings, in a neighborhood of
a fixed point, was quite well developed in the late $19^{th}$ century.
(Compare \S\S6,7.) However, except for one very simple case studied
by Cayley, nothing was known about the global behavior of iterated
holomorphic maps
until 1906, when Pierre Fatou described the following startling example.
For the map \[z\mapsto z^2/(z^2+2)\], he showed that
almost every orbit under iteration
converges to zero, even though there is a Cantor set of exceptional points
for which the orbit remains bounded away from zero. (Problem 3-6.)
This aroused great interest. After a hiatus dug
the first world war, the subject was taken up in depth by Fatou, and also
by Gaston Julia and others such as S. Latt\'es and J. F. Ritt. The most
fundamental and incisive contributions were those of Fatou himself,
although Julia developed much closely related material at more or less the
same time. Julia, who had been badly wounded
during the war, was awarded the ``Grand Prix des Sciences math\'ematiques''
by the Paris Academy of Sciences in 1918 for his work.\medskip

{\bf Definition.} Let \[S\] be a Riemann surface, let \[f:S\to S\]
be a non-constant holomorphic mapping, and let \[f^{\circ n}:S\to S\] be its
\[n$-fold iterate.
Fixing some point \[z_0\in S\], we have the following basic dichotomy: {\it
If there exists some neighborhood \[U\] of \[z_0\] so that the sequence of
iterates \[\{f^{\circ n}\}\] restricted to \[U\] forms a normal family, then
we say that \[z_0\] is a {\bit regular\/} or {\bit normal\/} point, or that
\[z_0\] belongs to the {\bit Fatou set\/} of \[f\]. Otherwise, if no such
neighborhood exists, we say that \[z_0\] belongs
to the {\bit Julia set\/} \[J=J(f)\].} (For sharper formulations of this
dichotomy in the rational case, see \S11.8 and Problem 13-1.)\smallskip

Thus, by its very definition, the Julia set \[J\] is a closed subset of
\[S\], while the Fatou set \[S\ssm J\] is the complementary open
subset. (The choice
as to which of these two sets should be named after Julia and which after
Fatou is rather arbitrary. The term ``Julia set'' is firmly established,
but the Fatou set is often called by other names, such as ``stable set''
or ``normal set''.) {\it
Roughly speaking,
\[z_0\] belongs to the Fatou set if dynamics in some neighborhood of \[z_0\]
is relatively tame, and belongs to the Julia set, if dynamics in every
neighborhood of \[z_0\] is more wild.}\smallskip

The classical example, and the one which we will emphasize, is the
case where \[S\] is the Riemann sphere \[\hat\C=\C\cup\infty\].
Any holomorphic map \[f:\hat\C\to \hat\C\] on the Riemann sphere
can be expressed as a rational
function, that is as the quotient \[f(z)=p(z)/q(z)\] of two polynomials.
Here we may assume that \[p(z)\] and \[q(z)\] have no common roots.
The {\bit degree} \[d\] of \[f=p/q\] is then equal to the maximum of the
degrees of \[p\] and \[q\]. In particular, for almost every choice of constant
\[c\in\hat\C\] the equation \[f(z)=c\] has exactly \[d\] distinct solutions
in \[\hat\C\]. (For every choice of \[c\] it has at
least one solution, since we assume that \[d>0\].)\smallskip

As a simple example, consider the squaring map \[s:z\mapsto z^2\] on
\[\hat\C\]. The entire open disk \[D\] is contained in the Fatou set of
\[s\], since successive iterates on any compact subset converge uniformly
to zero. Similarly the exterior \[\hat\C\ssm \bar D\] is contained in the
Fatou set, since the iterates of \[s\] converge to the constant function
\[z\mapsto\infty\] outside of \[\bar D\]. On the other hand, if \[z_0\]
belongs to the unit circle, then in any neighborhood of \[z_0\]
any limit of iterates \[s^{\circ n}\] would necessarily
have a jump discontinuity as we cross the unit circle.
This shows that the Julia set \[J(s)\] is precisely equal to the unit circle.
\smallskip

Such smooth Julia sets are rather exceptional. (Compare \S5.)
See Figure 1 for some much more typical examples of Julia sets for polynomial
mappings. Figure 1a
shows a rather wild Jordan curve,
Figure 1b a rather thick Cantor set, Figure 1c a ``dendrite'', 
and Figure 1d a more complicated example, the ``airplane'', with a
superattracting period 3 orbit. (Further examples of Julia sets are shown
in Figures 2-5, 8-10, 12, and 17.)

\pageinsert % \eject

\vfil
\insertRaster 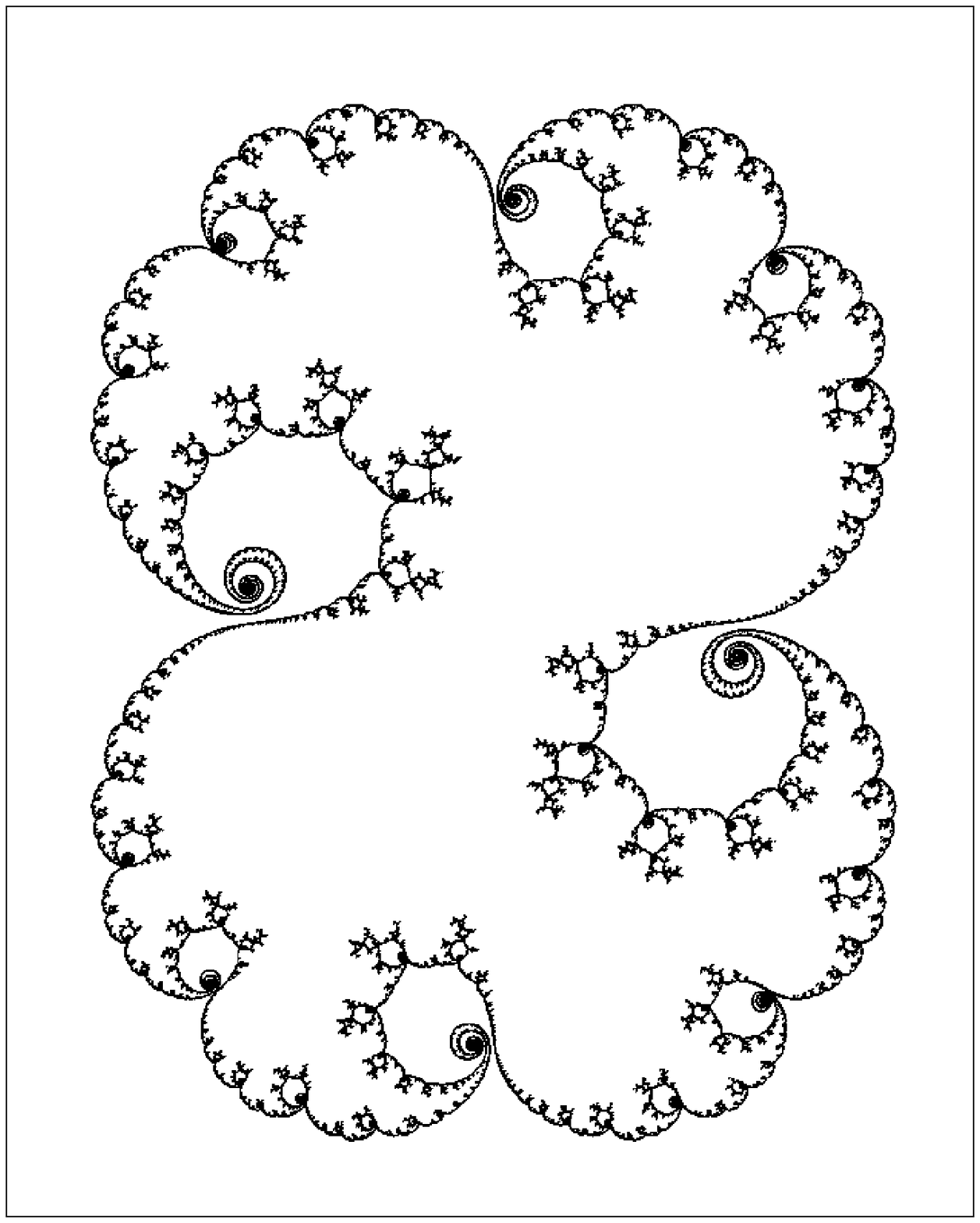 pixels 720 by 900 scaled 320
\smallskip

\centerline{Figure 1a. A simple closed curve,}
\centerline{Julia set for \[z\mapsto z^2+(.99+.14i)z\].}
\vfil
\insertRaster 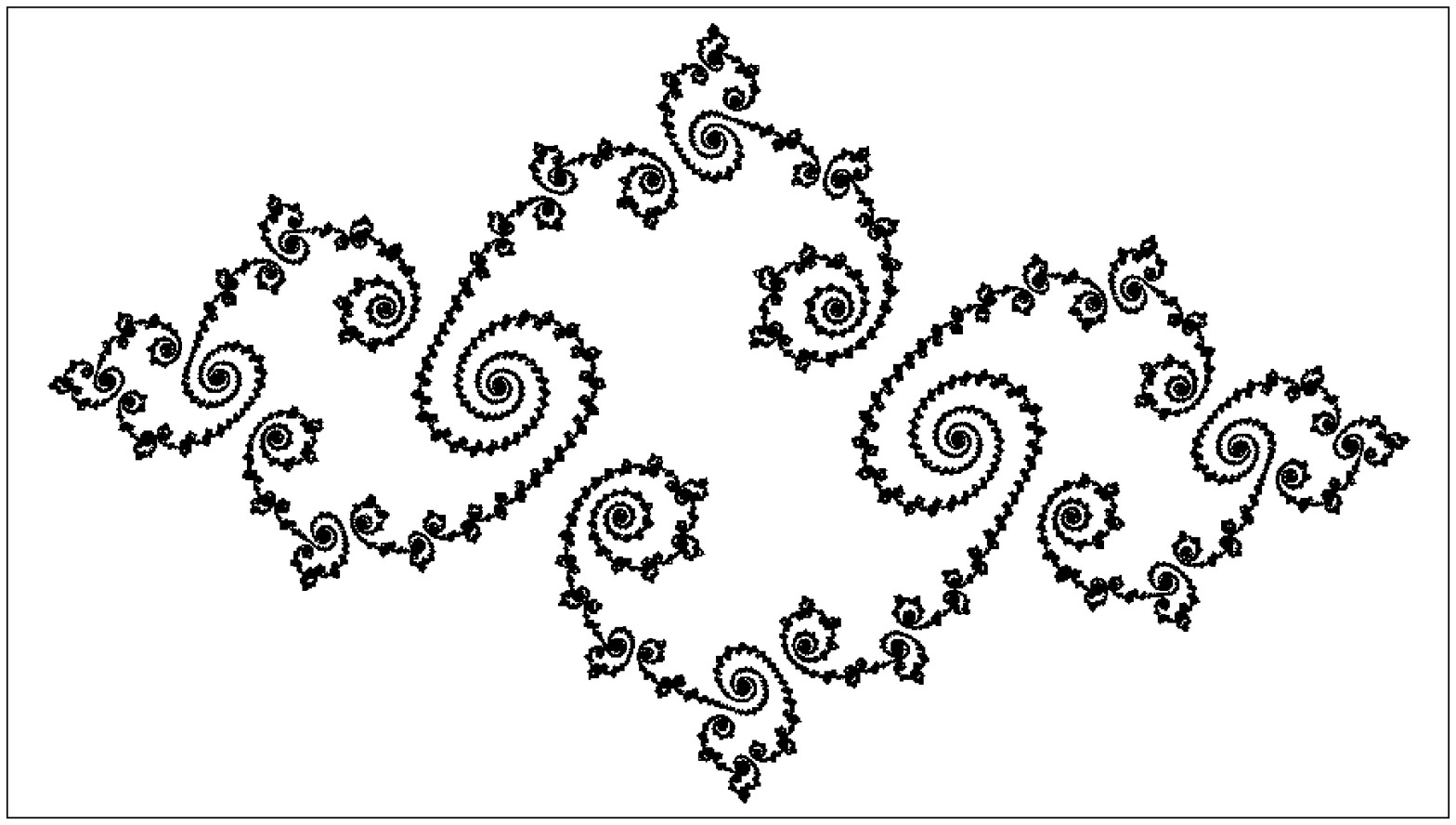 pixels 960 by 540 scaled 320
\smallskip
\centerline{Figure 1b. A totally disconnected Julia set,}
\centerline{\[z\mapsto z^2+(-.765+.12i)\].}
\vfil
\endinsert % \eject

\pageinsert
\vfil
\insertRaster 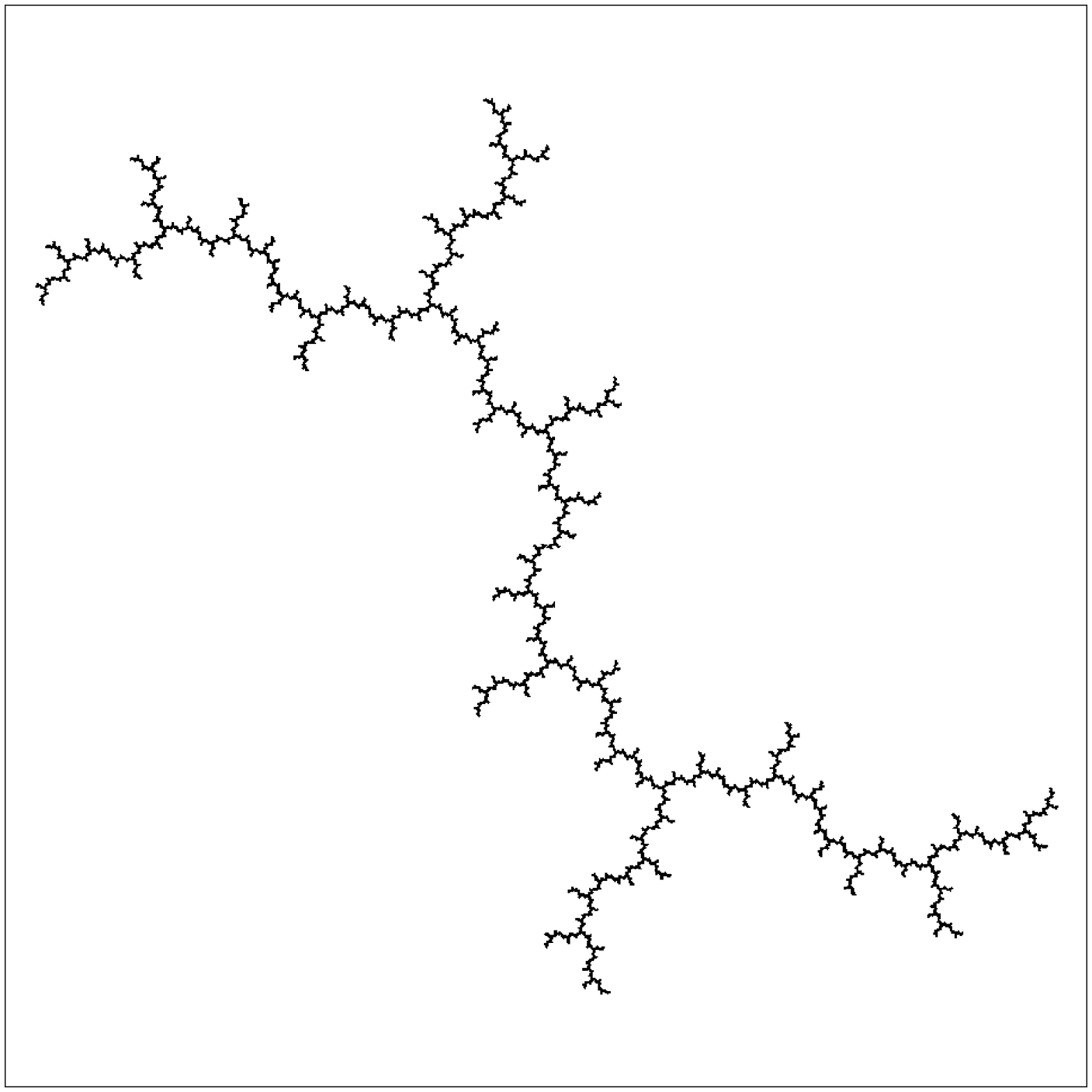 pixels 960 by 960 scaled 400
\smallskip

\centerline{Figure 1c. A ``dendrite'', Julia set for \[z\mapsto z^2+i\].}
\vfil%\endinsert % \eject

\insertRaster 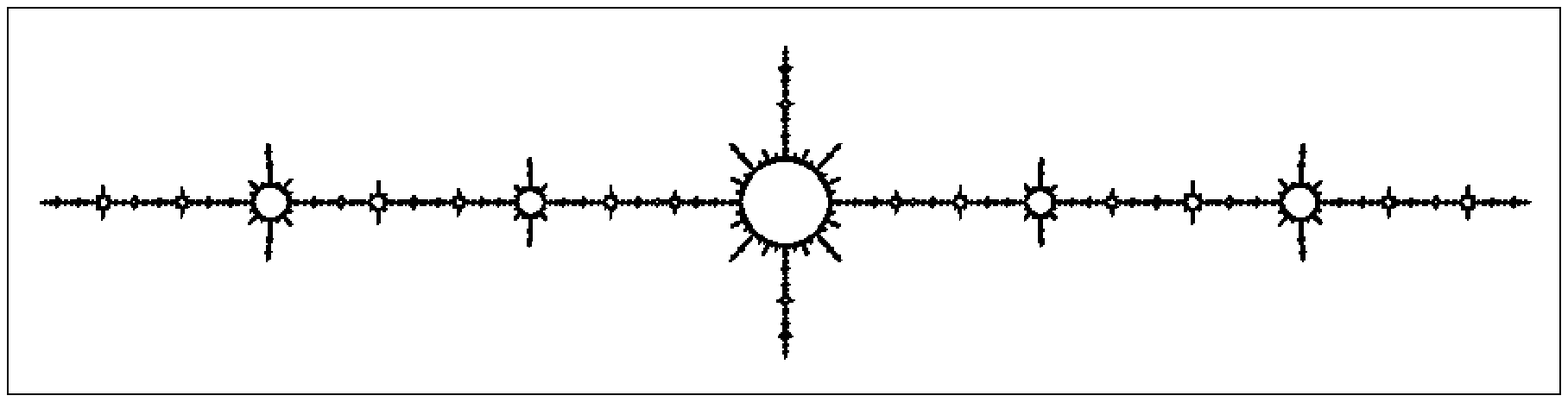 pixels 960 by 240 scaled 400
\smallskip
\centerline{Figure 1d. Julia set for \[z\mapsto z^2-1.75488\ldots\],
the ``airplane''.}
\vfil\endinsert % \eject

\smallskip

We will also need the following concept.\smallskip

{\bf Definition.} By the {\bit grand orbit} of a point
\[z\] under \[f:S\to S\] we mean the
set \[\GO(z,f)\] consisting of all points \[z'\in S\] whose orbits
eventually intersect the orbit of \[z\]. Thus \[z\] and \[z'\] have the same
grand orbit if and only if \[f^{\circ m}(z)=f^{\circ n}(z')\] for some
choice of \[m\ge 0\] and \[n\ge 0\].\smallskip

Here are some basic properties of the Julia set.

{\QP{\bf 3.1. Lemma.} {\it The Julia set \[J(f)\] of a
holomorphic map \[f:S\to S\] is {\bit fully invariant}
under \[f\]. That is,
if \[z\] belongs to \[J(f)\], then the entire grand orbit \[\GO(z,f)\]
is contained in \[J(f)\].}\medskip}

Evidently it suffices to prove that \[z\in J(f)\] if and only
if \[f(z)\in J(f)\].
A completely equivalent statement is that the Fatou set is fully invariant.
The proof, making use of the fact that a non-constant holomorphic map
takes open sets to open sets, is completely straightforward
and will be left to the reader.\QED

It follows that the Julia set possesses a great deal of {\bit
self-similarity\/}: {\it Whenever\break
\[f(z_1)=z_2\] in \[J(f)\], with derivative \[f'(z_1)\ne 0\], there
is an induced conformal\break isomorphism from a neighborhood \[N_1\] of \[z_1\]
to a neighborhood \[N_2\] of \[z_2\] which takes \[N_1\cap J(f)\] precisely
onto \[N_2\cap J(f)\].} (Compare Problem 3-7.)

{\QP {\bf 3.2. Lemma.} {\it For any \[n>0\], the Julia set \[J(f^{\circ n})
\] of the \[n$-fold iterate coincides with the Julia set \[J(f)\].}\medskip}

\noindent Again, the proof will be left to the reader.\QED

Now consider a {\bit periodic orbit\/} or ``{\bit cycle\/}''
$$	f\;:\;z_0\mapsto z_1\mapsto\;\cdots\;\mapsto z_{n-1}
	\mapsto z_n=z_0\,. $$
If the points \[z_1\,,\,\ldots\,,\,z_n\]
are all distinct, then the integer \[n\ge 1\] is called the {\bit period\/}.
If the Riemann surface \[S\] is \[\C\] (or an open subset
of \[\C\]), then the derivative
$$	\lambda\= (f^{\circ n})'(z_i)\=f'(z_1)\cdot f'(z_2)\cdots
	 f'(z_n) $$
is a well defined complex number called the {\bit multiplier~} or the
{\bit eigenvalue~}  of this periodic orbit. More generally, for
self-maps of an arbitrary Riemann surface the multiplier of a periodic
orbit can be defined using a local coordinate chart around any point of the
orbit. By definition, a
periodic orbit is either {\bit attracting~} or {\bit repelling~}
or {\bit indifferent\/} (\[=\]{\bit neutral~}) according as its multiplier
satisfies \[|\lambda|<1\] or \[|\lambda|>1\] or \[|\lambda|=1\]. The orbit
is called {\bit superattracting~} if \[\lambda=0\].

{\bf Caution:} In the special case where the point at infinity is periodic
under a rational map,
\[f^{\circ n}(\infty)=\infty\], this definition may be confusing.
The multiplier \[\lambda\] is {\it
not\/} equal to the limit as \[z\to\infty\] of the derivative of \[f^{\circ n}
(z)\], but is rather equal to
the reciprocal of this number. In fact if we introduce
the local uniformizing parameter \[w=1/z\] for \[z\] near \[\infty\],
then it is easy to check that the derivative of \[w\mapsto 1/f(1/w)\]
at \[w=0\] is
equal to \[\lim_{z\to\infty} 1/(f^{\circ n})'(z)\]. As an
example, if \[f(z)=2z\] then
\[\infty\] is an attracting fixed point with multiplier \[\lambda=1/2\].

In the
case of an attracting periodic orbit, we can define the {\bit basin of
attraction~} to be the open set \[\Omega\subset S\] consisting of all
points \[z\in S\] for which the successive iterates \[f^{\circ n}(z)\,,\,
f^{\circ 2n}(z)\,,\,\ldots\] converge towards some point of the periodic
orbit. In particular, this basin of attraction is defined
in the superattracting case. \smallskip

{\QP{\bf 3.3. Theorem.} {\it  Every attracting periodic orbit is contained
in the Fatou set. In fact the entire basin of
attraction \[\Omega\] for an attracting periodic orbit is contained in the
Fatou set. However the {\bit boundary} \[\partial \Omega\] is contained in the
Julia set, and every repelling periodic orbit is contained in the Julia set.}
\medskip}

{\bf Proof.} In view of 3.2, we need only consider the case of a fixed
point \[f(z_0)=z_0\]. If \[z_0\] is attracting, then it follows
from Taylor's Theorem that the successive iterates of
\[f\], restricted to a small neighborhood of \[z_0\], converge uniformly
to the constant function \[z\mapsto z_0\]. The corresponding statement for
any compact subset of the basin \[\Omega\] then follows easily.
On the other hand, around a
boundary point of this basin, that is a point which belongs to
the closure \[\bar \Omega\] but not to \[\Omega\] itself, it is clear that no
sequence of iterates can converge to a continuous limit. (See Problem
3-2 for a sharper statement.) If \[z_0\] is repelling, then
no sequence of iterates can converge uniformly near \[z_0\], since the
derivative \[df^{\circ n}(z)/dz\] at \[z_0\] takes the value \[\lambda^n\],
which diverges to infinity as \[n\to\infty\]. (Compare \S1.4.)\QED

The case of an indifferent periodic point is much more difficult. (Compare
\S\S7, 8.) One particularly important case is the following.\smallskip

{\bf Definition.} A periodic point \[f^{\circ n}(z_0)=z_0\] is called
{\bit parabolic~} if the multiplier \[\lambda\] at \[z_0\] is equal to
\[+1\], yet \[f^{\circ n}\] is not the identity map, or more generally
if \[\lambda\] is a root of unity, yet no iterate of \[f\] is the identity
map.

As an example, the two fixed points of the rational map \[f(z)=z/(z-1)\]
both have multiplier equal to \[-1\]. These do not count as parabolic
points since \[f\circ f(z)\] is identically equal to \[z\]. This exclusion
is necessary so that the following assertion will be true.

{\QP{\bf 3.4. Lemma.} \it Every parabolic periodic point belongs to the Julia
set.\medskip}

{\bf Proof.} Let \[w\] be a local uniformizing parameter, with \[w=0\]
corresponding to the periodic point. Then some iterate \[f^{\circ m}\]
corresponds to a local mapping of the \[w$-plane with
power series expansion of the form \[\;w\,\mapsto\, w+a_kw^k+a_{k+1}w^{k+1}+
\cdots\,\], where \[k\ge 2\,,\;\; a_k\ne 0\]. It follows that \[f^{\circ mp}\]
corresponds to a power series \[w\mapsto w+pa_kw^k+\cdots\]. Thus the \[k$-th
derivatives of \[f^{\circ mp}\] diverge as \[p\to\infty\].
It follows from 1.4 that
no subsequence can converge locally uniformly.\QED

Now and for the rest of \S3,
let us specialize to the case of a rational map \[f:\hat\C\to\hat\C\] of
degree \[d\ge 2\].

{\QP{\bf 3.5. Lemma.} {\it If \[f\] is rational of degree two or more, then the
Julia set \[J(f)\] is non-vacuous.}\medskip}

{\bf Proof.} If \[J(f)\] were vacuous, then some sequence of iterates
\[f^{\circ n(j)}\] would converge, uniformly over the entire sphere \[\bar
\C\], to a holomorphic limit \[g:\hat\C\to\hat\C\]. (Here we are using the
fact that normality is a local property: Problem 2-5.) A standard topological
argument would then show that the degree of \[f^{\circ n(j)}\] is equal to the
degree of \[g\] for large \[j\]. But the degree of \[f^{\circ n}\] is equal to
\[d^n\], which diverges to infinity as \[n\to\infty\].\QED

A different, more constructive proof of this Lemma will be given in \S9.5.
\smallskip

{\bf Definition.} A point \[z\in \hat\C\] is called
{\bit grand orbit finite} or (to use the classical terminology)
{\bit exceptional}$\,$ under \[f\] if its grand orbit \[\GO(z,f)\subset\hat
\C\] is a finite set. Using Montel's Theorem, we prove the following.

{\QP{\bf 3.6. Lemma.} {\it If \[f\] is rational of
degree two or more, then the set
\[\TF(f)\] of grand orbit finite points can have at most two elements. These
grand orbit finite points,
if they exist, must always be critical points of \[f\], and must belong
to the Fatou set.}\medskip}

{\bf Proof.} (Compare Problem 3-3.)
Note that \[f\] maps any grand orbit \[\GO(z,f)\] onto
itself. Hence, any finite grand orbit must constitute a single
periodic orbit under \[f\]. Each point \[z\] in this finite orbit must
be critical (and in fact
\[(d-1)$-fold critical, where \[d\] is the degree), since otherwise
\[f(z)\] would have two or more pre-images. Therefore, such an
orbit must be attracting, and hence contained in the Fatou set.

If there were three distinct grand orbit finite points, then the union of the
grand orbits of these points would form a finite set whose complement
\[U\] in \[\hat\C\] would be Hyperbolic, with \[f(U)= U\]. Therefore,
the set of iterates of \[f\] restricted to \[U\] would be normal by Montel's
Theorem. Thus both \[U\] and its complement would be contained in the
Fatou set, contradicting 3.5.\QED

{\QP{\bf 3.7. Lemma.} {\it Let \[z_1\] be any point of the Julia set
\[J(f)\subset\hat\C\]. If \[N\]
is a sufficiently small neighborhood of \[z_1\], then the union \[U\] of the
forward images \[f^{\circ n}(N)\] is precisely equal to the complement
\[\hat\C\ssm\TF(f)\] of the set of grand orbit finite points.}\medskip}

In particular, this union \[U\] contains the Julia set \[J(f)\].
(In \S11.2 we will see that just one forward image \[f^{\circ n}(N)\] actually
contains the entire Julia set, provided that \[n\] is sufficiently large.)

{\bf Proof of 3.7.} Let \[E\] be the complementary set \[\hat\C\ssm U\].
We have \[f(U)\subset U\], or equivalently \[f^{-1}(E)\subset E\], by the
construction. Since \[U\] intersects the Julia set, it follows from Montel's
Theorem that its complement \[E\] has at most two points. Now since \[E\] is
finite and \[f\] is onto, a counting argument shows that \[f^{-1}(E)=E\],
hence \[E\] is contained in the set \[\TF(f)\] of grand orbit
finite points. If the initial neighborhood \[N\]
is small enough to be disjoint from \[\TF(f)\], then it follows that
\[E=\TF(f)\].\QED

{\QP{\bf 3.8. Corollary.} {\it If the Julia set contains an interior point
\[z_1\], then it must be equal to the entire Riemann sphere.}\medskip}

For if we choose a neighborhood \[N\subset J\], then 3.7 shows that the
union \[U\subset J\] of forward images of \[N\] is everywhere dense on
\[\hat\C\].\QED

{\QP{\bf 3.9. Theorem.} {\it
If \[z_0\in J(f)\], then the set of all iterated pre-images
$$	\{z\;:\; f^{\circ n}(z)=z_0\quad{\rm for\;some\quad} n\ge0\}	$$
is everywhere dense in \[J(f)\]. In particular, it follows that the
grand orbit\break \[\GO
(z_0\,,f)\] is everywhere dense in \[J(f)\].}\medskip}

For \[z_0\] is not a grand orbit finite point, so 3.7 shows that every point
\[z_1\in J(f)\] can be approximated arbitrarily closely by points \[z\]
whose forward orbits contain \[z_0\].\QED

This Theorem suggests an algorithm for computing pictures of the
Julia set: Starting with any \[z_0\in J(f)\], first compute all pre-images
\[f(z_1)=z_0\], then compute all pre-images \[f(z_2)=z_1\], and so on, thus
eventually coming arbitrarily close to every point of \[J(f)\]. This method
is most often used in the quadratic case, since quadratic equations are very
easy to solve. The method is very insensitive to round-off errors; since
\[f\] tends to be expanding on its Julia set, so that \[f^{-1}\] tends to be
contracting. (Compare Problem 3-4. ??)

{\QP {\bf 3.10. Corollary. } {\it If \[f\] has degree two or more, then
\[J(f)\] has no isolated points.}\medskip}

{\bf Proof.} Since \[J(f)\] is fully invariant, it follows from 3.5 and 3.6
that \[J(f)\] must be an infinite set. Hence it contains at least one
limit point \[z_0\]. Now the iterated pre-images of \[z_0\] form a dense set
of non-isolated points in \[J(f)\].\QED

A property of a point of \[J\] is said to be true for {\bit generic}
\[z\in J\] if it is true for all points in some countable intersection
of dense open subsets \[U_i\cap J\subset J\]. (Compare \S8. Here the notation
is supposed to indicate that the \[U_i\] are open subsets of \[\hat\C\]
and that the closure of \[U_i\cap J\] is equal to \[J\].) By Baire's
Theorem, any such countable intersection is itself a dense subset of \[J\].

{\QP{\bf 3.11. Corollary.} \it For generic \[z\in J(f)\], the forward orbit
\hfil\break \[\{z\,,\,f(z)\,,\,f^{\circ 2}(z)\,,\,\ldots\;\}\;\] is everywhere
dense in \[J(f)\].\medskip}

{\bf Proof.} Let \[\{B_j\}\] be a countable collection of open sets
forming a basis for the topology of \[\hat\C\]. For each \[B_j\] which
intersects \[J=J(f)\], let \[U_j\] be the union of the iterated pre-images
\[f^{-n}(B_j)\] for \[n\ge 0\]. Then it follows from 3.9 that the closure
of \[U_j\cap J\] is equal to the entire Julia set \[J\], and the conclusion
follows.\QED
\pageinsert
\vfil
\insertRaster 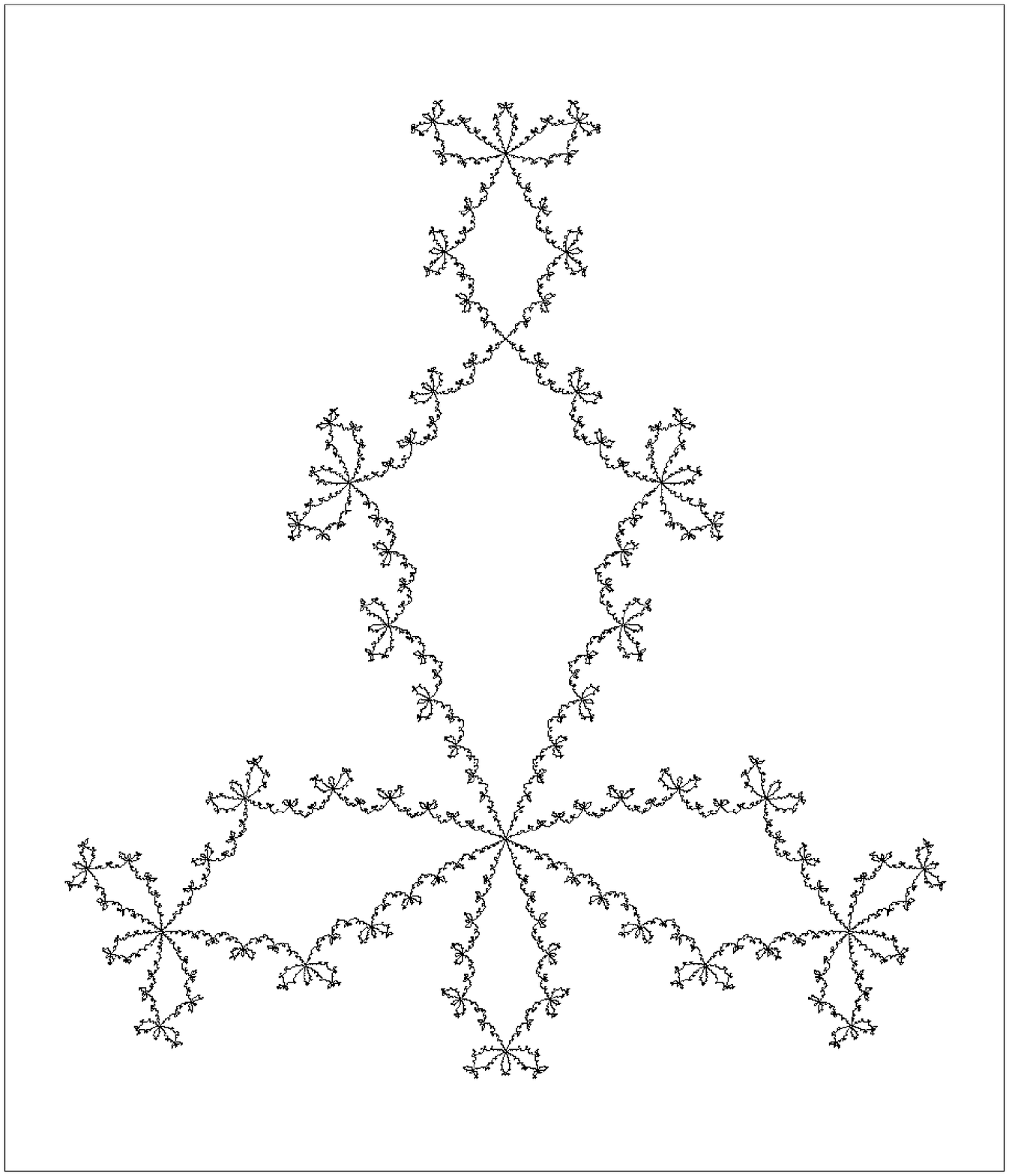 pixels 960 by 1120 scaled 400
\bigskip\bigskip
\centerline{Figure 2. Julia set for \[f(z)=z^3+{12\over 25}z+{116\over 125}i\].}
\vfil
\endinsert

We will continue the study of rational Julia sets in \S11.
%Here are some problems for the reader.
\lln\smallskip

{\bf Concluding problems.}\smallskip

{\bf Problem 3-1.} If \[f:\hat\C\to\hat\C\] is rational of degree \[d=1\],
show that the Julia set \[J(f)\] is either vacuous, or
consists of a a single repelling or parabolic fixed point.

{\bf Problem 3-2.} If \[\Omega\subset\hat\C\]
is the basin of attraction for an attracting
periodic orbit, show that the boundary \[\partial \Omega=\bar\Omega\ssm
\Omega\] is equal to the entire Julia set. (Compare Theorem 3.3.) Here it
is essential that we include all connected components of this basin.\smallskip

{\bf Problem 3-3.} Show that a rational map \[f\] is actually a polynomial
if and only if the point at infinity is a grand orbit finite fixed
point for \[f\]. Show that \[f\] has both zero and
infinity as grand orbit finite points if and only if \[f(z)=\alpha z^n\], where
\[n\] can be any non-zero integer, negative or positive, and where
\[\alpha\ne 0\]. Conclude that
a rational map has grand orbit finite points if and only if it is conjugate,
under some fractional linear change of coordinates, either to a polynomial
or to the map \[z\mapsto 1/z^d\] for some \[d\ge 2\].\smallskip

{\bf Problem 3-4.} Using a hand calculator if necessary, decide what maps to
what in Figure 1d.\smallskip

{\bf Problem 3-5.} If
\[f(z)=z^2-6\], show that \[J(f)\] is a Cantor set contained in the intervals
\[[-3,-\sqrt 3]\cup[\sqrt 3,3]\]. More precisely, show that a point
in \[J(f)\] with\break orbit
\[z_0\mapsto z_1\mapsto\cdots\] is uniquely determined by the sequence of
signs \[z_j/|z_j|=\pm 1\].\smallskip

{\bf Problem 3-6.} For Fatou's function \[f(z)=z^2/(2+z^2)\], show that
the entire completed real axis \[\R\cup\infty\] is contained in the basin
of attraction of the origin. Show that \[J(f)\] is a Cantor set. More
precisely, given any infinite sequence of signs \[\epsilon_0\,,\,\epsilon_1\,
,\,\ldots\] show that there is one and only one point \[z=z(\epsilon_0\,,\,
\epsilon_1\,,\,\ldots)\] which satisfies the condition that \[f^{\circ n}
(z)\] is uniformly bounded away from zero, and
belongs to the half-plane \[\epsilon_n H\] for each \[n\ge 0\]. To
this end, consider the branch \[g(z)=\sqrt{2z/(1-z)}\] of \[f^{-1}\],
mapping \[\hat\C\ssm [0,1]\] onto the upper half-plane \[H\]. Starting with
any \[z_0\not\in[0,1]\], show that the successive images  \[\epsilon_0g(
\epsilon_1g(\cdots \epsilon_ng(z_0)\cdots))\] converge to the required
point \[z\in J(f)\]. Since a Cantor set cannot separate the plane, show that
the basin of attraction of the origin is equal to \[\hat\C\ssm J(f)\].
\smallskip

{\bf Problem 3-7. Self-similarity.} With rare exceptions, any shape which is observed about
one point of the Julia set will be observed infinitely often, throughout the
Julia set. More precisely, for two points \[z\] and \[z'\] of \[J=J(f)\],
let us say that \[(J,z)\] is {\bit locally conformally isomorphic} to
\[(J,z')\]
if there exists a conformal isomorphism from a neighborhood \[N\] of \[z\]
onto a neighborhood \[N'\] of \[z'\] which carries \[z\] to \[z'\] and
\[J\cap N\] onto \[J\cap N'\]. For all but finitely many \[z_0\in J\],
show that the set of \[z\] for which \[(J,z)\] is locally conformally
isomorphic to \[(J,z_0)\] is everywhere dense in \[J\].

As an example, consider the polynomial map \[f(z)=z^3+.48\,z+.928\,i\]
of Figure 2, and explain why the fixed point \[.8\,i\] looks different
from all other points of the Julia set. (See also Example 2 of \S5.)
How many pre-images does this point have?
\vfill

%\rightline{jwm, version of 8-1-91}
\eject

\footline={\hss\tenrm 4-\folio\hss}\pageno=1

\centerline{\bf \S4. Dynamics on Other Riemann Surfaces.}\bigskip

This section will try to say something about the theory of iterated
holomorphic
mappings on an arbitrary Riemann surface. However we will
concentrate on the easy cases, and simply
refer to the literature for the hard cases.
It turns out that there are only three
Riemann surfaces \[S\] for which the study of iterated mappings
is really difficult, namely
the sphere \[\hat\C\], the plane \[\C\], and the punctured plane
\[\C\ssm\{0\}\].\smallskip

The theory of iterated holomorphic mappings from the {\bit Riemann
sphere\/} \[\hat \C\] to itself has been outlined in \S3, and
will form the main goal of all of the subsequent
sections.\smallskip

We can distinguish two different classes of holomorphic maps from the {\bit
complex
plane\/} \[\C\] to itself. A {\bit polynomial map\/} of \[\C\] extends uniquely
over the Riemann sphere \[\hat\C\]. Hence the theory of polynomial
mappings can be subsumed as a special case of the theory of rational maps
of \[\hat\C\]. (See especially \S\S17-18.) On the other hand,
{\bit transcendental mappings\/} from \[\C\] to itself form an essential
distinct and more difficult subject of study. Such mappings have
been studied for more than sixty years by
many authors, starting with Fatou himself. Even iteration of the exponential
map \[\exp:\C\to\C\] provides a number of quite challenging problems.
(See for example Lyubich, and Rees.) A
proof that the Julia set \[J(\exp)\]
is the entire plane \[\C\] is included in Devaney [Dv1].
Further information about iterated transcendental functions may be found
for example in
Baker, Devaney [Dv2], Goldberg \& Keen, and in
Eremenko \& Lyubich.\smallskip

The study of iterated maps from
the {\bit punctured plane\/} \[\C\ssm\{0\}\] to itself is also a difficult and
interesting subject. (See for example Keen.)
\smallskip

For all of the uncountably many other Riemann surfaces, it turns out
that the possible dynamical behavior
is very restricted, and fairly easy to describe.
We must distinguish two very different cases according as the surface
is a {\bit torus\/} or is {\bit Hyperbolic\/}. First consider the case of
a torus \[{\bf T}=\C/\Lambda\], where \[\Lambda\] is a lattice in \[\C\].

{\QP{\bf 4.1. Theorem.} \it Every holomorphic map
\[f:{\bf T}\to{\bf T}\] is an affine map,\break
\[f(z)\equiv \alpha z+c \;\;(\mod\Lambda)\].
The corresponding Julia set \[J(f)\] is either the empty set or the
entire torus according as \[|\alpha|\le 1\] or \[|\alpha|>1\].\medskip}

Proofs will be given at the end of this section.
Here the possible values for the derivative \[\alpha\] are sharply restricted.
(Problem 4-1.) The dynamics of such iterated affine maps on the torus are
of some interest. (See Problem 4-3, as well as Example 3 of \S5.)
However, the possibilities are so limited
that this study cannot be considered very difficult.\smallskip

Finally, suppose that \[S\] is Hyperbolic. Then we will see that
any holomorphic self-map behaves in a rather dull manner under iteration.
In particular, the Julia set is alway vacuous.
Consider first the special case of the
unit disk. The following was proved by Denjoy
in 1926, sharpening an earlier result by Wolff.

{\QP {\bf 4.2. Theorem.} {\it Let \[f:D\to D\] be any
holomorphic map. Then either

$(1)$ \[f\] is a ``rotation'' about some fixed point \[z_0\in D\], or else

$(2)$ the successive iterates \[f^{\circ n}\] converge, uniformly on
compact subsets of \[D\], to a constant function \[z\mapsto c_0\in \bar D\].}
\medskip}

Here the notation \[f^{\circ n}\] stands for the \[n$-fold composition
\[f\circ\cdots\circ f\] mapping \[D\] into itself. Note that the
limiting value \[\lim_{n\to\infty}\,f^{\circ n}(z)\] may belong to the
boundary \[\partial D\]. Similar statements hold for maps of the upper
half-plane \[H\] to itself. As examples, if \[f(w)=2w\] or if \[f(w)=w+i\]
for \[w\in H\], then evidently \[\lim_{n\to\infty}\,f^{\circ n}(w)\] is
equal to
the boundary point \[+\infty\] for every \[w\in H\]. (Here we must measure
convergence with respect to the spherical metric on the compact set
\[\bar H\subset\hat\C\].)\smallskip

The corresponding assertion for an arbitrary Hyperbolic surface
is only slightly more complicated to state.
First some definitions. If \[f:S\to S\], then the sequence of points
\[\;z\mapsto f(z)\mapsto
f^{\circ 2}(z)\mapsto\cdots\;\] is called the {\bit orbit} of the point
\[z\in S\]. A
fixed point \[f(z_0)=z_0\] is said to be {\bit attracting} if the derivative
\[\lambda=f'(z_0)\] satisfies \[|\lambda|<1\]. (To be more precise, we
must first choose some local coordinate in a neighborhood of the fixed point,
and compute this derivative in terms of this local coordinate. We will be
careless about this, since in practice our Riemann surfaces will usually
be open subsets of \[\C\].) For such an attracting fixed point, the
{\bit basin of attraction} \[\Omega=\Omega(z_0)\]
is defined to be the set of all
points \[z\in S\] for which the orbit \[z\mapsto f(z)\mapsto f^{\circ 2}(z)
\cdots\] converges towards \[z_0\]. It is not difficult
to check that \[\Omega\] is an open subset of \[S\] and that \[z_0\in
\Omega\].

{\QP {\bf 4.3. Theorem.} {\it If \[S\] is a Hyperbolic Riemann surface,
then for every holomorphic map \[f:S\to S\] the Julia set \[J(f)\]
is vacuous. Furthermore either:

\noindent \[(a)\] every orbit converges
towards a unique attracting fixed point \[f(z_0)=z_0\],

\noindent\[(b)\] every orbit diverges to infinity with respect to the
\Poincare metric on \[S\],

\noindent \[(c)\] \[f\] is an automorphism of finite order, or

\noindent\[(d)\] \[S\] is isomorphic either to a disk \[D\], a punctured
disk \[D\ssm\{0\}\], or an annulus
$$	{\cal A}_r\;=\;\{z:1<|z|<r\}\,,	$$
and \[f\]
corresponds to an irrational rotation \[z\mapsto e^{2\pi i t}z\] with
\[t\not\in{\bf Q}\].}\medskip}

\noindent
Evidently these four possibilities are mutually exclusive.
Much later, in \S13, we will want to apply this Theorem to the case
where \[S\] in an open subset of the sphere \[\hat\C\] and \[f\] is a
rational map carrying \[S\] into itself. In this case, as in 4.2, it
is convenient to study limiting values which belong to the closure of \[S\].
We can then sharpen the statement in Case (b) as follows.

{\QP{\bf 4.4. Addendum.} {\it Suppose that \[U\] is a Hyperbolic open
subset of \[\hat\C\] and that \[f:U\to U\] extends holomorphically throughout
a neighborhood of the closure \[\bar U\].
Then in Case $(b)$ above, all orbits in \[U\]
must converge within \[\hat\C\] to a single boundary fixed point
$$	f(\hat z)=\hat z\;\in\;\partial U\,.	$$
This convergence is uniform on compact subsets of \[U\].}\medskip}

\noindent (In the special case where \[U\] is the standard unit disk, of course
we have this same result without the extra hypothesis
that \[f\] extends over the boundary.)\medskip

Now let us begin the proofs.\smallskip

{\bf Proof of Lemma 4.1.} To fix our ideas, suppose that \[{\bf T}=\C/\Lambda\]
where the lattice \[\Lambda
\subset\C\] is spanned by the two numbers \[1\] and \[\tau\], and where
\[\tau\not\in\R\]. Any holomorphic map \[f:{\bf T}\to {\bf T}\] lifts to a holomorphic
map \[F:\C\to\C\] on the universal covering surface. Note first that there
exist two lattice elements \[\lambda_1\,,\,\lambda_2\in\Lambda\] so that
$$	F(z+1)\;=\;F(z)+\lambda_1\,,\quad F(z+\tau)\;=\;F(z)+\lambda_2	$$
for every \[z\in\C\]. For example we certainly have \[F(z+1)\equiv F(z)\;\;
(\mod\Lambda)\], and the difference function \[F(z+1)-F(z)\in\Lambda\] must
be constant since \[\C\] is connected and the target space \[\Lambda\] is
discrete.
Now let \[g(z)=F(z)-\lambda_1 z\], so that \[g(z+1)=g(z)\]. Then
$$	g(z+\tau)\;=\;g(z)+(\lambda_2-\lambda_1\tau)\,.	$$
Thus \[g\] gives rise to a map from
the torus \[{\bf T}\] to the infinite cylinder
$$ 	\C/(\lambda_2-\lambda_1\tau)\Z\;\;\cong\;\;\C\ssm\{0\}\;,	$$
or from the torus
to \[\C\] itself if \[\;\lambda_2-\lambda_1\tau=0\]. Using Liouville's
Theorem or the Maximum Modulus Principle, we see easily that \[g\] must be
constant Thus \[g(z)\equiv c\], and\break \[F(z)=\lambda_1 z+c\]
as required. The computation of \[J(f)\] will be left as an
exercise. (See Problems 4-2 and 4-3 below.) \QED\smallskip

We will skip over 4.2 for the moment.
The proof of Theorem 4.3 begins as follows. If we are not
in Case (b), then some orbit \[z_0\mapsto z_1
\mapsto\cdots\] must possess an accumulation point \[\hat z\in S\]. That is, we
can find integers \[n(1)<n(2)<\cdots\] so that the sequence \[\{z_{n(i)}\}\]
converges to \[\hat z\]. Consider the sequence of maps \[f^{\circ(
n(i+1)-n(i))}\], carrying \[z_{n(i)}\] to \[z_{n(i+1)}\].
By Montel's Theorem, in the version 2.6, there exists some subsequence
which converges, uniformly on compact subsets, to a holomorphic map
\[h:S\to S\]. Evidently \[h(\hat z)=\hat z\].\smallskip

{\it First suppose that \[f\] strictly contracts the \Poincare metric.} Then
\[h\] must
also contract the \Poincare metric, hence \[h\] cannot have two distinct
fixed points. But \[f\] and \[h\] commute, since \[h\] is a limit of
iterates of \[f\]; hence \[f\] must map fixed points of \[h\] to fixed
points of \[h\]. Therefore the unique fixed point \[\hat z\] of \[h\]
must also be a fixed point of \[f\].
This fixed point is attracting, since \[f\] contracts the \Poincare metric.
{\it It follows that every orbit under \[f\] converges to \[\hat z\],
so that we are in Case $(a)$ of 4.3}. For otherwise,
if \[z'\] were the closest point which did not belong to the attractive
basin of \[\hat z\], then the \Poincare distance \[\rho(\hat z\,,\,z')\] could
not strictly decrease under the map \[f\].\smallskip

{\it Suppose, on the other hand, that \[f\] preserves the \Poincare metric.}
Still assuming that \[f\] has an orbit with a finite limit point, we will
prove the following.

{\QP\bf 4.5. Lemma. {\it
Under these hypotheses, some sequence of iterates \[f^{\circ m(i)}\]
converges, uniformly on compact subsets, to the identity map of \[S\].
It follows that \[f\] is necessarily an automorphism
of \[S\].}\medskip}

{\bf Proof.} Since \[f\] preserves the \Poincare metric, it
lifts to an isometry from the universal covering \[\tilde S\]
to itself. It follows easily that \[f\] is either an automorphism
or a covering map from \[S\] onto \[S\]. On the other hand,
we know that some sequence of iterates of \[f\] converges to a
map \[h\] which has a fixed point. We can lift \[h\] to a map \[H:
\tilde S\to\tilde S\] which has a corresponding fixed point.
Since \[H\] must also preserve
the \Poincare metric, it follows that \[H\] can only be a rotation about
this fixed point. Therefore, some sequence
of iterates of \[h\] converges to the identity map \[\iota\] of \[S\].
It then follows
easily that some sequence of iterates of \[f\] converges to \[\iota\]. Therefore
\[f\] is one-to-one. For if \[f(z)=f(z')\], then \[f^{\circ m}(z)
=f^{\circ m}(z')\], and passing to the limit we have \[\iota(z)=\iota(z')\],
or in other words \[z=z'\].
Since a one-to-one covering map is clearly a conformal isomorphism,
this proves the Lemma.\QED

To complete the proof of 4.3, we must prove the following.

{\QP{\bf 4.6. Lemma.} \it If an automorphism \[f\] of a Hyperbolic surface
\[S\] has iterates \[f^{\circ m}\] which are
arbitrarily close to the identity map
(uniformly on compact subsets), then either \[f\]
has finite order, or else \[S\] is isomorphic to \[D\] or \[D\ssm\{0\}\]
or to an annulus, and \[f\] corresponds to
an irrational rotation.\medskip}

(For this Lemma, we don't really need the hypothesis that \[S\] is Hyperbolic.
In the non-Hyperbolic case, the corresponding statement would be that
\[S\] is isomorphic to either \[\C\] or \[\C\ssm\{0\}\] or \[\hat\C\],
and that \[f\] corresponds to an irrational rotation.)\smallskip

Before proving 4.6, it
will be convenient to briefly consider the more general situation of
a map \[f:S\to S'\] between different Riemann surfaces, where \[S\cong\tilde S
/\Gamma\] and \[S'\cong\tilde S'/\Gamma'\]. Such a map \[f\] lifts to a map
\[F:\tilde S\to\tilde S'\] which is unique up to composition with elements
of the group \[\Gamma'\] of deck transformations of the target space.
As in the proof of 4.1,
to each deck transformation \[\gamma\in\Gamma\] there corresponds one and only
one deck transformation \[\gamma'\in\Gamma'\] so that the identity
$$	F(\gamma(z))\;=\;\gamma'(F(z))	$$
holds for all \[z\in S\]. In fact the correspondence \[\gamma\mapsto\gamma'\]
is a group homomorphism, which can be identified with the ``induced
homomorphism'' \[f_*:\pi_1(S)\to\pi_1(S')\] between fundamental groups.

We are interested in the special case \[S=S'\], with universal covering
\[\tilde S\] isomorphic to the unit disk \[D\]. It will be convenient to
choose some compact disk \[K\subset S\] and a corresponding compact disk \[K^*\]
in the universal covering surface \[\tilde S\]. If \[f^{\circ m(j)}\] is
uniformly close to the identity map on \[K\], note
that the lifted map \[F^{\circ m(j)}\] may be far from the identity map on
\[K^*\]. However, there must exist some deck transformation
\[\gamma_j\] so that the composition \[F_j=\gamma_j\circ F^{\circ m(j)}\]
is uniformly close to the identity on \[K^*\].\smallskip

We will also need the following.

{\QP{\bf 4.7. Lemma.} {\it Let \[\Gamma\subset G\] be any discrete
subgroup of a topological group \[G\]. Then for each \[\gamma\in\Gamma\]
there exist a neighborhood \[N\] of the identity in
\[G\] so that elements \[g\in N\] satisfy \[\;g\circ\gamma\circ
g^{-1}\in\Gamma\;\] only if they commute with \[\gamma\].}\medskip}

The proof is straightforward.\QED

{\bf Proof of 4.6.} By the discussion above, we have a sequence of
automorphisms\break
\[\;F_j\,=\,\gamma_j\circ F^{\circ m(j)}\,\in\, G(\tilde S)\;\]
which converge locally uniformly to
the identity automorphism. Furthermore, for each fixed \[\gamma_0\in\Gamma\]
and each \[F_j\]
there is a corresponding \[\gamma'_j\] so that
$$	F_j\circ\gamma_0\;=\gamma'_j\circ F_j\,.	$$
By 4.7, it follows that \[F_j\] commutes with \[\gamma_0\] for large \[j\].

Now let us identify \[G(\tilde S)\] with the automorphism group \[G(D)\],
which operates not only on the open disk \[D\] but also
on the closed disk \[\bar D\]. By 1.15 we know that
two non-identity elements of \[G(D)\] commute if and only if they have exactly
the same fixed points in \[\bar D\]. Thus, if we fix some non-identity
\[\gamma\in\Gamma\], and if we exclude the case where some \[F_j\] is actually
equal to the identity automorphism, then we can say that
\[\gamma\] has exactly the same fixed points as
\[F_j\] for \[j\] sufficiently large. {\it This implies that all non-identity
elements of \[\Gamma\] have the same fixed points, so that \[\Gamma\] is a
commutative group.} Evidently, either \[\Gamma\] is trivial and \[S\cong D\],
or \[\Gamma\] is free cyclic and \[S\] is an annulus or punctured disk.
(Compare Problem 2-3.) Further details are straightforward, and will be left
to the reader.\QED

{\bf Proof of Addendum 4.4.}
Starting at some arbitrary point \[z_0\] in the connected open set
\[U\subset\hat\C\], choose a path \[p:[0,1]\to U\] from the point
\[z_0=p(0)\] to \[f(z_0)=p(1)\], and continue
this path inductively for all \[t\ge 0\] by setting \[p(t+1)=f(p(t))\].
By hypothesis, the orbit \[p(0)\mapsto p(1)\mapsto p(2)\mapsto\cdots\]
converges to infinity with respect to the \Poincare metric on \[U\].
Hence it must tend to the boundary of the open set \[U\], with respect to
the spherical metric. Let \[\delta\] be the diameter
of the image \[p[0,1]\]
in the \Poincare metric for \[U\]. Then each successive image \[p[n,\,n+1]\]
must also have diameter less than \[\delta\]. Since these sets converge
to the boundary of \[U\], it follows that the diameter of \[p[n,\,n+1]\] in
the spherical metric must tend to zero as \[n\to\infty\]. (\S2.4.)
Since \[f(p(t))=
p(t+1)\], this implies that every accumulation point of the path \[p(t)\]
as \[t\to \infty\] must be a fixed point of \[f\]. It is not difficult
to check that the set of all such accumulation points must be a connected
subset of the boundary \[\partial U\]. But our hypothesis that \[f\]
continues holomorphically throughout a neighborhood of the closure \[\bar U\]
guarantees that \[f\] can have only finitely many fixed points in \[\bar U\].
This proves that the path \[p(t)\] converges to just one fixed point \[f(\hat z)
=\hat z\in\partial U\]. Thus the orbit of \[z_0\] converges to \[\hat z\]
in the spherical metric. Using 2.2 and 2.4, it follows easily that every
orbit in \[U\] converges to \[\hat z\], and that this convergence is
uniform on compact subsets of \[U\].\QED

{\bf Proof of the Denjoy-Wolff Theorem 4.2.} We
now assume that \[S\] is the unit disk \[D\subset\C\],
but do not assume that \[f\] can be continued outside of the open disk.
The following argument is taken from a
lecture of Beardon, as communicated to me by Shishikura. For any \[\epsilon
>0\], let us approximate \[f\] by the map \[z\mapsto (1-\epsilon)f(z)\]
from \[D\] into a proper subset of itself.
Then it is not difficult to check that there is one and only one fixed
point \[z_\epsilon=(1-\epsilon)f(z_\epsilon)\]. If \[f\] itself has no fixed
point, then these \[z_\epsilon\] must tend to the boundary of the disk
as \[\epsilon\to 0\].
Let \[r(\epsilon)\] be the \Poincare distance of \[z_\epsilon\] from the
origin, and consider the closed neighborhood \[N_{r(\epsilon)}(z_\epsilon
\,,\,\rho_D)\], which contains the origin as a boundary point. By Pick's
Theorem 1.10, this neighborhood is necessarily carried into itself by the map
\[z \mapsto(1-\epsilon)f(z)\]. These neighborhoods are actually round disks
(although with a different center point) with respect to the Euclidean metric.
(Problem 1-8.)
After passing to a subsequence as \[\epsilon
\to 0\], we may assume that these disks \[N_{r(\epsilon)}(z_\epsilon
\,,\,\rho_D)\] tend to a limit disk \[N_0\], bounded  by a circle (known as a
``horocircle'') which is tangent to the boundary of \[D\] at a single
point \[\hat z\]. Now \[f\] must map \[N_0\] into itself, hence the entire
orbit of \[0\] under \[f\] must be contained in \[N_0\]. Therefore, if
this orbit has no limit point in the open disk \[D\], then it must
converge towards the point of tangency
\[\hat z\]. The argument now proceeds as in 4.4.\QED
\smallskip

\lln

We conclude with problems for the reader.
\medskip

{\bf Problem 4-1.} Given some torus \[\T=\C/\Lambda\] and some \[\alpha\in\C\],
show that there exists a holomorphic map \[f(z)\equiv
\alpha z+c\] from \[\T\] to itself with derivative
\[\alpha\] if and only if
\[\alpha\Lambda\subset\Lambda\]. Evidently an arbitrary rational
integer \[\alpha
\in\Z\] will satisfy this condition. Show that there exists such a map with
derivative \[\alpha\not\in\Z\]
if and only the ratio of two generators of \[\Lambda\] satisfies a quadratic
equation with integer coefficients. Such a torus is said to admit
``{\bit complex multiplications\/}''.\smallskip

{\bf Problem 4-2.} If \[|\alpha|=1\] but \[\alpha\ne 1\], show that \[f\] is
an automorphism of finite order, and in fact of order either 2, 3, 4, or 6.
Conclude that \[J(f)\] is vacuous.
Show that all four cases can occur, for suitably chosen \[\Lambda\].
\smallskip

{\bf Problem 4-3.} If \[\alpha\ne 0\], show that any equation of the form
\[f(z)=z_0\] has exactly \[|\alpha|^2\] solutions \[z\in\T\].
If \[\alpha\ne 1\], show that \[f\]
has exactly \[|\alpha-1|^2\] fixed points. (In particular, both
\[|\alpha|^2\] and \[|\alpha-1|^2\] are
necessarily integers.) More generally, if \[|\alpha|>1\]
show that the equation \[f^{\circ n}(z)=z\] has exactly \[|\alpha^n-1|^2\]
solutions in \[\T\]. Show that the periodic points of \[f\]
are everywhere dense in \[\T\] whenever \[\alpha\ne 0\].
Conclude that the Julia set \[J(f)\] is the entire torus whenever
\[|\alpha|>1\].
\smallskip

{\bf Problem 4-4.} On the other hand, show that a map from a Hyperbolic
surface to itself can have at most one periodic point
(necessarily a fixed point), unless every point
is periodic.
\vfill\bigskip
%\rightline{jwm, version of 8-28-91}

\eject

\footline={\hss\tenrm 5-\folio\hss}\pageno=1

\centerline{\bf \S5. Smooth Julia Sets.}\medskip

Most Julia sets tend to be complicated fractal subsets of \[\hat\C\].
(Compare Brohlin [Br]. We will give a very partial explanation for this fact
in \S6.3.) This section, however,
will be devoted to three exceptional Julia sets which are actually smooth
manifolds.\smallskip

{\bf Example 1: the Circle.} As discussed already in \S3,
the unit circle appears as Julia set for the mapping
\[z\mapsto z^{\pm n}\] for any \[n\ge 2\]. Other rational
maps with this same Julia set are described in Problem 5-1. Similarly,
the real axis \[\R\cup\infty\], as the image of the unit circle under
a conformal automorphism, can appear as a Julia set. (Problem 5-2.)\smallskip

{\bf Example 2: the Interval.} Following Ulam and Von Neumann,
consider the map \[f(z)= z^2-2\], which carries the closed interval
\[I=[-2,\,2]\] onto itself. We will show that the Julia set \[J(f)\] is
equal to this interval \[I\], and that every point outside of \[I\]
belongs to the attractive basin \[\Omega(\infty)\] of the point at infinity.
(For higher degree maps with this same property, see Problem 5-3.)
For \[z_0\in I\], it is easy to check
that both solutions of the equation \[f(z)=z_0\] belong to this interval \[I\].
Since \[I\] contains a repelling fixed point \[z=2\], it follows from
Theorem 3.9 that \[I\] contains the entire Julia set \[J(f)\]. On the other
hand, the basin \[\Omega(\infty)\] is a neighborhood of infinity
whose boundary \[\partial\Omega(\infty)\] is contained in \[J(f)\subset I\]
by Theorem 3.3. Hence everything outside of \[J(f)\] must belong to this
basin. Since \[f(I)\subset I\], it follows that \[J(f)=I\].

Here is a more constructive proof that \[J(f)=I\]. We make use of the
substitution\break
\[g(w)=w+w^{-1}\], which carries the unit circle in a two-to-one
manner onto \[I=[-2,\,2]\]. For \[z_0\ne I\], the equation \[g(w)
=z_0\] has two solutions, one of which lies inside the unit circle and one
of which lies outside. Hence \[g\] maps the exterior of the closed unit disk
isomorphically onto the complement \[\C\ssm I\].
Since the squaring map in the \[w$-plane is related to \[f\] by the identity
$$	g(w^2)=g(w)^2-2=f(g(z))\,,	$$
it follows easily that the orbit of \[z\] under \[f\] either remains
bounded or diverges to infinity according as \[z\] does or does not belong to
this interval. Again using 3.3, it follows that that \[J(f)=I\].\QED

{\bf Example 3: all of \[\hat\C\].} The rest of this section will describe
an example constructed by S. Latt\`es, shortly before his death in 1918. Given
any lattice \[\Lambda\subset \C\]
we can form the quotient torus \[\T=\C/\Lambda\], as in \S4.
Thus \[\T\] is a compact
Riemann surface, and is also an additive Lie group. Note that the automorphism
\[z\mapsto -z\] of this surface has just four fixed points. For example, if
\[\Lambda=\Z+\tau\Z\] is the lattice with basis \[1\] and \[\tau\],
where \[\tau\not\in\R\],
then the four fixed points are \[0,\,1/2,\,\tau/2,\] and \[(1+\tau)/2\] modulo
\[\Lambda\].

Now form a new Riemann
surface \[S\] as a quotient of \[\T\] by identifying each \[z\in T\] with
\[-z\]. Evidently \[S\] inherits the structure of a Riemann surface (although
it loses the group structure). In fact we can use \[(z-z_j)^2\] as a local
uniformizing parameter for \[S\] near each of the four fixed points \[z_j\].
Thus the natural map \[\T\to S\] is two-to-one, except at the four ramification
points. To compute the genus of \[S\], we use the following.\eject

{\QP {\bf 5.1. Riemann-Hurwitz Formula.} {\it Let \[\T\to S\] be a branched
covering map from one compact Riemann surface onto another. Then the number of
branch points, counted with multiplicity, is equal to \[\chi(S)d-
\chi(T)\], where \[\chi\] is the Euler characteristic and \[d\] is the
degree.}\medskip}

{\bf Sketch of Proof.} Choose some triangulation of \[S\] which includes
all critical
values (that is all images of ramification points) as vertices; and let
\[a_n(S)\] be the number of\break
\[n$-simplexes, so that \[\chi(S)=a_2(S)-a_1(S)
+a_0(S)\]. In general, each simplex of \[S\] lifts to \[d\] distinct
simplices in \[\T\]. However, if \[v\] is a critical value, then there
are too few pre-images of \[v\]. The number of missing pre-images is precisely
the number of ramification points over \[v\], each counted with an
appropriate multiplicity; and the conclusion follows. %\QED
\rlap{$\sqcup$}$\sqcap$

{\bf Remark.} This proof works also for Riemann surfaces with smooth
boundary. The Formula remains true for proper maps between non-compact
Riemann surfaces, as can be verified, for example, by means of
a direct limit argument.\smallskip

In our example, since \[\T\] is a torus with Euler characteristic
\[\chi(T)=0\],
and since there are exactly four simple branch points, we conclude
that \[2\chi(S)-\chi(T)=4\], or \[\chi(S)=2\]. {\it Using the standard
formula \[\chi=2-2g\], we conclude that \[S\] is a surface of
genus zero, isomorphic to the Riemann sphere.} (Remark: The projection
map from \[\T\] to the sphere \[S\] is, up to a choice of normalization,
known as the {\bit Weierstrass \[\wp$-function}.)\smallskip

Now consider the doubling map \[z\mapsto 2z\] on \[\T\]. This commutes with
multiplication by \[-1\], and hence induces a map \[f:S\to S\].
Since the doubling map has degree four, it follows that \[f\] is a rational
map of degree four. (More generally, in place of the doubling map, we could
use any linear map which carries the lattice \[\Lambda\] into itself,
as in Problem 4-1.)

{\QP {\bf 5.2. Latt\`es Theorem.} {\it The Julia set for this rational map \[f\] is
the entire sphere \[S\].}\medskip}

{\bf Proof.} Evidently the doubling map on \[\T\] has the property that
periodic points are everywhere dense. For example, if \[r\] and \[s\]
are any rational numbers with odd denominator, then \[r+s\tau\] is periodic.
These periodic orbits are all repelling, since the multiplier is a
power of two. Evidently \[f\] inherits the same property, and the conclusion
follows by \S3.3. (Alternatively, given a small open set \[U\subset S\],
it is not difficult to show that \[f^{\circ n}(U)\] is equal to the entire
sphere \[S\] whenever \[n\] is sufficiently large. Hence no sequence of
iterates of \[f\] can converge to a limit on any open set.)\QED

In order to pin down just which rational map \[f\] has these properties,
we must first label the points of \[S\]. The four branch
points on \[\T\] map to four ``ramification points''
on \[S\], which will play a special role. Let us choose a conformal
isomorphism from \[S\] onto \[\hat\C\] which maps the first three of these
points to \[\infty\,,\,0\,,\,1\] respectively. The fourth ramification point
must then map to some \[a\in\C\ssm\{0\,,\,1\}\]. In this way we construct a
projection map \[\pi:T\to\hat\C\] of degree two, which satisfies
\[\pi(-z)=\pi(z)\], and which has critical values
$$	\pi(0)=\infty\,,\;\; \pi(1/2)=0\,,\;\; \pi(\tau/2)=1\,,\;\;
	\pi((1+\tau)/2)=a\,.	 $$
Here \[a\] can be any number distinct from \[0\,,\,1\,,\,\infty\]. In
fact, given \[a\in\C\ssm\{0\,,\,1\}\], it is not difficult to show that
there is one and only one branched covering \[ \T' \to \hat\C\]
of degree two with precisely
\[\{\infty\,,\,0\,,\,1\,,\,a\}\] as ramification points. (Compare Appendix E.)
The Riemann-Hurwitz formula shows that
this branched covering space \[\T'\] is a torus, necessarily isomorphic to
\[\C/(\Z+\tau\Z)\] for some \[\tau\not\in\R\]. The unique deck transformation
which interchanges the two pre-images of any point must preserve the linear
structure, and hence be multiplication by \[-1\].

Now the doubling map on \[\T\] corresponds under \[\pi\] to a specific
rational map\break \[f_a:\hat\C\to\hat\C\], where
$$	f_a(\pi(z))=\pi(2z)\,,	$$
and where \[J(f_a)=\hat\C\] by 5.2.
A precise computation of this map \[f_a\] is described in
Problem 5-5 below.\smallskip

{\bf Remark.} Mary Rees has proved the existence
of many more rational maps with \[J(f)=\hat\C\]. (See also Herman.)
For any degree \[d\ge 2\],
let \[{\rm Rat}(d)\] be the complex manifold consisting of all rational
maps of degree \[d\]. Rees shows that there is a subset of \[{\rm Rat}(d)\]
of positive measure consisting of maps \[f\] which are ``ergodic". By
definition, this means that any measurable subset of \[\hat\C\] which
is fully invariant under \[f\] must either have full measure or measure zero.
It can be shown that any ergodic map must necessarily have \[J(f)=\hat\C\].
\smallskip

We will study these smooth Julia sets further in \S14.7.
\lln

{\bf Concluding problems.}\smallskip

{\bf Problem 5-1.} For any \[a\in D\] the map \[\phi_a(z)=(z-a)/(1-\bar
a z)\] carries the unit disk isomorphically onto itself. (Problem 1-2.)
A finite product of the form
$$	f(z)\;=\; e^{i\theta}
\phi_{a_1}(z) \phi_{a_2}(z)\cdots\phi_{a_n}(z) $$
with \[a_j\in D\] is
called a {\bit Blaschke product} of degree \[n\]. Evidently any such \[f\]
is a rational map which carries \[D\] onto \[D\]. If \[n\ge 2\], and if
one of the factors is \[\phi_0(z)=z\] so that \[f\] has a fixed
point at the origin, show that the Julia set \[J(f)
=J(1/f)\]
is the unit circle.\smallskip

{\bf Problem 5-2.} Newton's method applied to the polynomial equation\break
\[f(z)=z^2+1=0\] yields the rational map
$$	\nu(z)\;=\;z-f(z)/f'(z)\;=\; {\textstyle{1\over 2}}(z-1/z)	$$
from \[\hat\C=\C\cup\infty\] to itself.
Show that \[J(\nu)=\R\cup\infty\].\smallskip

{\bf Problem 5-3.} The {\bit monic Tchebycheff polynomials}
$$	P_1(z)=z\;,\quad P_2(z)=z^2-2\;,\quad P_3(z)=z^3-3z\;,\;\ldots	$$
can be defined inductively
by the formula \[P_{n+1}(z)+P_{n-1}(z)=zP_n(z)\].
Show that \[P_n(w+w^{-1})
=w^n+w^{-n}\], or equivalently that \[P_n(2\cos\theta)\;=\;2\cos(n\theta)\],
and show that \[P_m\circ P_n=P_{mn}\]. For \[n\ge 2\] show that
the Julia set of \[\pm P_n\] is the interval \[[-2,\,2]\].  For \[n\ge
3\] show that \[P_n\] has \[n-1\] distinct critical points but only two
critical values, namely \[\pm 2\].\smallskip

{\bf Problem 5-4.} Show that the Julia sets studied in this section
have the following extraordinary property.
(Compare Problems 14-2, 14-1.) For all but one or two of the
periodic orbits
\[z_0\mapsto z_1\mapsto\cdots\mapsto z_n=z_0\], the multiplier
\[\lambda=f'(z_1)\cdot\cdots\cdot f'(z_n)\] satisfies \[|\lambda|=d^n\] when \[J\]
is 1-dimensional, or \[|\lambda|=d^{n/2}\] when \[J=\hat\C\],
where \[d\] is the degree.\smallskip

{\bf Problem 5-5.} For the torus \[\T=\C/(\Z+\tau\Z)\] of Example 3, show that
the involution \[z\mapsto z+\,\scriptstyle{1/ 2}\] of \[\T\] corresponds
under \[\pi\] to the involution
\[w\mapsto a/w\] of \[\hat\C\], with fixed points \[w=\pm\sqrt{a}\].
Show that the rational map \[f=f_a\] has poles at \[\infty\,,\,0\,,\,1\,,
\,a\] and double zeros at \[\pm\sqrt{a}\].  Show that \[f\] has a fixed
point of multiplier \[\lambda=4\] at infinity, and conclude that
$$	f(w)\;=\;{(w^2-a)^2 \over 4w(w-1)(w-a)}\,\;.	$$
As an example, if \[a=-1\] then
$$	f(w)\;=\;{(w^2+1)^2\over 4w(w^2-1)}\,\;.	$$
Show that the correspondence \[\;\tau\mapsto a=a(\tau)\in\C\ssm\{0,1\}\;\]
satisfies the equations
$$	a(\tau+1)=1/a(\tau)\;,\quad a(-1/\tau)=1-a(\tau)\;,	$$
and also \[\;a(-\bar\tau)=\bar a(\tau)\,\]. Conclude, for example, that
\[\;a(i)=1/2\,\], and that\break \[a((1+i)/2)=-1\]. (This correspondence
\[\tau\mapsto
a(\tau)\] is an example of an ``elliptic modular function'', and provides
an explicit representation of the upper half-plane \[H\] as a universal
covering of the thrice punctured sphere \[\C\ssm\{0,1\}\]. Compare Ahlfors
[A1, pp.269-274].)\smallskip

{\bf Problem 5-6.} For each of the six critical points \[\omega\] of \[f\],
show that \[f(f(\omega))\] is the repelling fixed point at infinity.
(According to \S13.5, the fact that each critical orbit terminates
on a repelling cycle is enough to imply that \[J(f)=\hat\C\].)
% \bigskip\bigskip

\vfill
%\rightline{jwm, version of 8-28-91}
\end

%% file: ras.tex
% This file is supposed to imitate the -ras.tex- file which contained
% the macros -insertRaster- and -RasterBox- for including SUN raster files
% in TeX.  But here we assume that the original raster files have been
% converted to PostScript, so we use -psfig- rather than -sunbitmap-.
% There was a bug in the way -ras.tex- scaled these bitmaps, and this file
% attempts to preserve this bug so that the resulting figures are identical
% to the raster figures in the original paper.

\input psfig
%\magnification=1200

\newif\ifboxfigure      % set to true if you want the figure boxed
\boxfigurefalse

\def\BoxIt#1#2{%        % put a box around #1, leaving a gap of #2
	\vbox{\hrule
	\hbox{\vrule\kern#2\vbox{\kern#2#1\kern#2}\kern#2\vrule}
		   \hrule}}

\def\insertRaster #1 pixels #2  by #3 scaled #4 {
%  #1 is the raster file
%  #2 is the width of the picture in pixels
%  #3 is the height of the picture in pixels
%  #4 is scaling 500 means half size 2000 means double size 1000 is actual  
			\medskip
			 \hbox to \hsize{%

			 \hss
			 \RasterBox {#1} {#2} {#3} {#4}
			 \hss
			 }%
}
\def\RasterBox #1 #2 #3 #4{

% dimen0 and dimen1 are the dimensions of the hbox which are setup in
% the same way as -ras.tex- unfortunately psfig includes images in a 
% different manner from sunbitmap, so it must be scaled using different
% values, hence dimen2 and dimen3 and MAGIC number 65/72

\dimen5=65pt
\divide\dimen5 by 72

\dimen0=#2\dimen5
\divide\dimen0 by 1000
\dimen1=#3\dimen5
\divide\dimen1 by 1000
\dimen2=#3\dimen5
\divide\dimen2 by 1000
\dimen3=#2\dimen5
\divide\dimen3 by 1000

\setbox4=\hbox to #4\dimen0{
 \vbox to #4\dimen1{
 \vss
 \psfig{figure=#1,height=#4\dimen2,width=#4\dimen3}
 }
 \hss
 }
 \ifboxfigure\BoxIt{\box4}{0pt}
 \else\box4
 \fi
 }

% \insertRaster towrrab.pic pixels 500 by 500 scaled 400

%% file: cd6-8.tex
\input header
\input ras
\def\Koenigs{K\oe nigs }
\font\tencyr=wncyr10
\input cyracc.def
\def\cyr{\tencyr\cyracc}

\footline={\hss\tenrm 6-\folio\hss}\pageno=1

\centerline{\bf LOCAL FIXED POINT THEORY}\bigskip

\centerline{\bf \S6. Attracting and Repelling Fixed Points.}\medskip

Consider a function
$$	f(z)\;=\;\lambda z+a_2 z^2+a_3 z^3+\;\cdots	\eqno (1) $$
which is defined and holomorphic in some neighborhood of the origin, with
a fixed point of multiplier \[\lambda\] at \[z=0\]. If
\[|\lambda|\ne 1\], we will show
that \[f\] can be reduced to a simple normal form by a suitable change of
coordinates. First consider the case \[\lambda\ne 0\],
so that the origin is not a critical point. The following was proved by G.
\Koenigs in 1884.

{\QP{\bf 6.1. \Koenigs Linearization Theorem.} {\it If the
 multiplier \[\lambda\] satisfies\break \[|\lambda|
 \ne 0, 1\], then there exists a local holomorphic change of coordinate
 \[w=\phi(z)\], with \[\phi(0)=0\], so that \[\;\phi\circ\! f\!\circ
 \phi^{-1}\,\] is the linear
 map \[w\mapsto\lambda w\] for all \[w\] in some neighborhood of the origin.
 Furthermore, \[\phi\] is unique up to multiplication by a non-zero constant.
 }\medskip}

\noindent (This functional equation \[\,\phi\circ\! f\!\circ
 \phi^{-1}(w)\,=\,\lambda w\,\] had been studied some years earlier by
 Schr\"oder, who believed however that it did not have many solutions.)
\smallskip

 {\bf Proof of uniqueness.} If there were two such maps \[\phi\] and \[\psi\],
 then the composition
 $$	\psi\circ\phi^{-1}(w) \;=\;b_1 w+ b_2 w^2+b_3 w^3+\,\cdots	$$
 would commute with the map \[w\mapsto\lambda w\]. Comparing coefficients
 of the two resulting power series, we see that \[\lambda b_n=b_n\lambda^n\]
 for all \[n\]. Since \[\lambda\] is neither zero nor a root of unity,
 this implies that
 \[b_2=b_3=\cdots=0\]. Thus \[\psi\circ\phi^{-1}(w) \;=\;b_1 w\], or in
 other words \[\psi(z)=b_1\phi(z)\].\smallskip

 {\bf Proof of existence when \[0<|\lambda|<1\].} Choose a constant \[c<1\]
 so that\break
 \[c^2<|\lambda|<c\], and choose a neighborhood \[D_r\] of the origin
 so that \[|f(z)|\le c|z|\] for \[z\in D_r\]. Thus for any starting point
 \[z_0\in D_r\], the orbit \[z_0\mapsto z_1\mapsto\cdots\] under \[f\]
 converges geometrically towards the origin, with \[|z_n|\le rc^n\].
 By Taylor's Theorem,
 $$	|f(z) - \lambda z| \;\le\; k|z^2|	$$
 for some constant \[k\] and for all \[z\in D_r\], hence
 $$	|z_{n+1}-\lambda z_n|\;\le\; kr^2 c^{2n}\,.	$$
 It follows that the numbers \[w_n=z_n/\lambda^n\] satisfy
 $$	|w_{n+1}- w_n|\;\le\; k'(c^2/|\lambda|)^n\,,	$$
 where \[k'=kr^2/|\lambda|\]. These differences converge
 uniformly and geometrically to zero.\break\eject
 \noindent Thus the holomorphic
 functions \[z_0\mapsto w_n(z_0)\] converge, uniformly throughout \[D_r\],
 to a holomorphic limit \[\;\phi(z_0)\,=\,\lim_{n\to\infty}\; z_n/\lambda^n\].
 The required identity \[\;\phi(f(z))\;=\;\lambda\phi(z)\]\break follows
 immediately. A similar argument shows that \[|\phi(z)-z|\] is less than or
 equal to some constant times \[|z^2|\]. Therefore \[\phi\] has derivative
 \[\phi'(0)=1\], and hence is a local conformal diffeomorphism.

 {\bf Proof when \[|\lambda|>1\].} Since \[\lambda\ne 0\], the inverse
 map \[f^{-1}\] is locally well defined and holomorphic, having the origin as
 an attractive fixed point with multiplier \[\lambda^{-1}\]. Applying the
 above argument to \[f^{-1}\], the conclusion follows.\QED

{\bf 6.2. Remark.} More generally, suppose that we consider a
family of maps \[f_\alpha\] of the form (1) which depend holomorphically on
one (or more) complex parameters
\[\alpha\] and have multiplier \[\lambda=\lambda(\alpha)\]
satisfying \[|\lambda(\alpha)|\ne 0,\,1\]. Then
a similar argument shows that {\it the \Koenigs
 function \[\phi(z)=\phi_\alpha(z)\] depends holomorphically on
\[\alpha\].}  This fact will be important in \S8.5.
To prove this statement, let us fix \[0<c<1\]
 and suppose that \[|\lambda(\alpha)|\] varies
through some compact subset of the interval
 \[(c^2,\,c)\]. Then it is easy to check that the convergence in the proof
 of 6.1 is uniform in \[\alpha\]. The conclusion now follows easily.\QED
\smallskip

 {\bf 6.3. Remark.} The \Koenigs Theorem helps us to understand
 why the Julia set \[J(f)\] is so often a complicated ``fractal" set.
 {\it Suppose that there exists a single repelling periodic point \[\hat z\]
 for which
 the multiplier \[\lambda\] is not a real number. Then \[J(f)\] cannot
 be a smooth manifold, unless it is all of \[\hat\C\].} To see this,
 choose any point \[z_0\in J(f)\] which is close to \[\hat z\], and let
 \[w_0=\phi(z_0)\]. Then \[J(f)\] must also contain an infinite sequence
 of points \[z_1\,,\,z_2\,,\,\ldots\] with \Koenigs coordinates \[\;\phi(z_n)
 =w_n/\lambda^n\;\] which lie along a logarithmic spiral and converge to zero.
 Evidently such a set can not lie in any smooth submanifold of \[\C\].
 Furthermore, if we recall that the iterated pre-images of our fixed
 point are everywhere dense in \[J(f)\], then we see that such sequences
 lying on logarithmic spirals are extremely pervasive. Compare Figures 3
 and 4 which show typical examples of such spiral structures, associated with
 repelling points of periods 2 and 1 respectively.\smallskip

 We can restate 6.1 in a more global form as follows.

 {\QP {\bf 6.4. Corollary.} {\it Suppose that \[f:S\to S\] is a holomorphic
 map from a Riemann surface to itself with an attractive fixed point of
 multiplier \[\lambda\ne 0\] at \[\hat z\].
 Let \[\Omega\subset S\] be the basin of attraction,
 consisting of all \[z\in S\] for which \[\lim_{n\to\infty}\;
 f^{\circ n}(z)\] exists and is equal to \[\hat z\].
 Then there is a holomorphic map \[\phi\] from \[\Omega\] onto \[\C\] so that
 the diagram
 $$	\matrix{\Omega& \buildrel f\over\longrightarrow & \Omega\cr
	 \quad\downarrow\phi & & \quad\downarrow\phi\cr
	 \C&\buildrel \lambda{\textstyle\cdot}\over\longrightarrow &\;\C\cr}
	 \eqno (2) $$
 is commutative, and so that \[\phi\] takes a neighborhood of \[\hat z\]
 diffeomorphically onto a neighborhood of zero.
 Furthermore, \[\phi\] is unique up to multiplication
 by a\break constant.}\medskip}

 \noindent In fact, to compute \[\phi(z_0)\] at an arbitrary point of \[\Omega\]
 we must simply follow the orbit
 of \[z_0\] until we reach some point \[z_n\] which is very close to
 \[\hat z\], then evaluate the \Koenigs coordinate
 \[\phi(z_n)\] and multiply by \[\lambda^{-n}\].\QED\eject

 \pageinsert
 \vfil %\endinsert

% \pageinsert\vfil
 \insertRaster 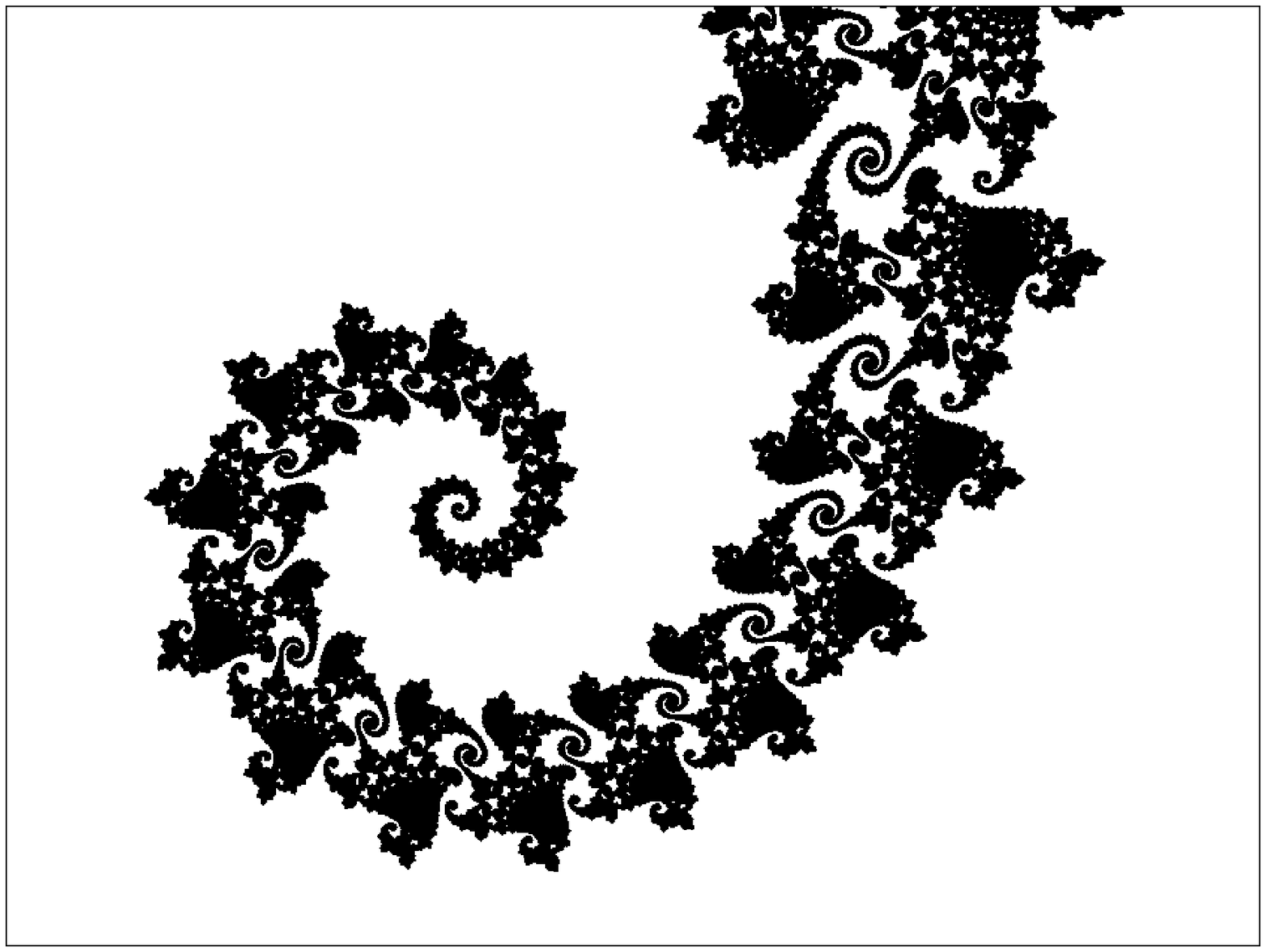 pixels 960 by 720 scaled 310
 \smallskip
 \centerline{Figure 3. Detail of Julia set for \[z\mapsto z^2-.744336+.121198i\].}
 \bigskip
 \insertRaster 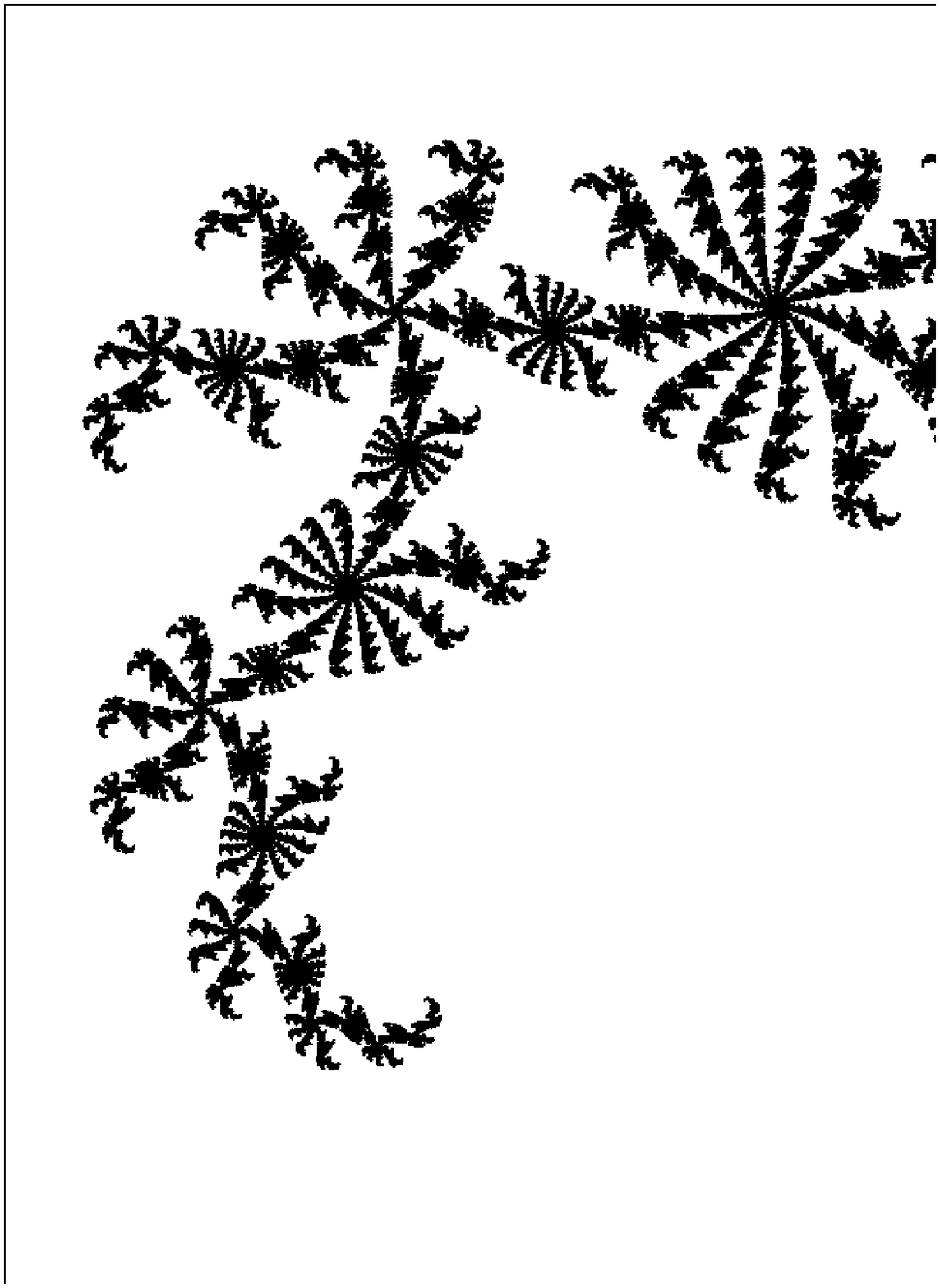 pixels 720 by 990 scaled 320
 \smallskip
 \centerline{Figure 4. Detail of Julia set for
 \[z\mapsto z^2+.424513+.207530i\].}
 \vfil\endinsert

 As an example, Figure 5 illustrates the map
 \[f(z)= z^2+0.7z\]. Here the Julia set \[J\] is the outer Jordan curve,
 bounding the basin \[\Omega\] of the attracting fixed point.
 The critical point \[\omega=-0.35\]
 is at the center of the picture, and the attracting fixed point \[\hat z=0\]
 is directly above it.
 The curves \[\;|\phi(z)|={\rm constant}=|\phi(\omega)/\lambda^n|\;\] have
 been drawn in.
 Note that \[\phi\] has zeros at all iterated preimages of \[\hat z\],
 and critical points at all iterated preimages of the critical point \[\omega\].
 The function \[z\mapsto\phi(z)\] is unbounded, and oscillates wildly as
 \[z\] tends to \[J=\partial\Omega\].\medskip

 \pageinsert\vfil
 \insertRaster 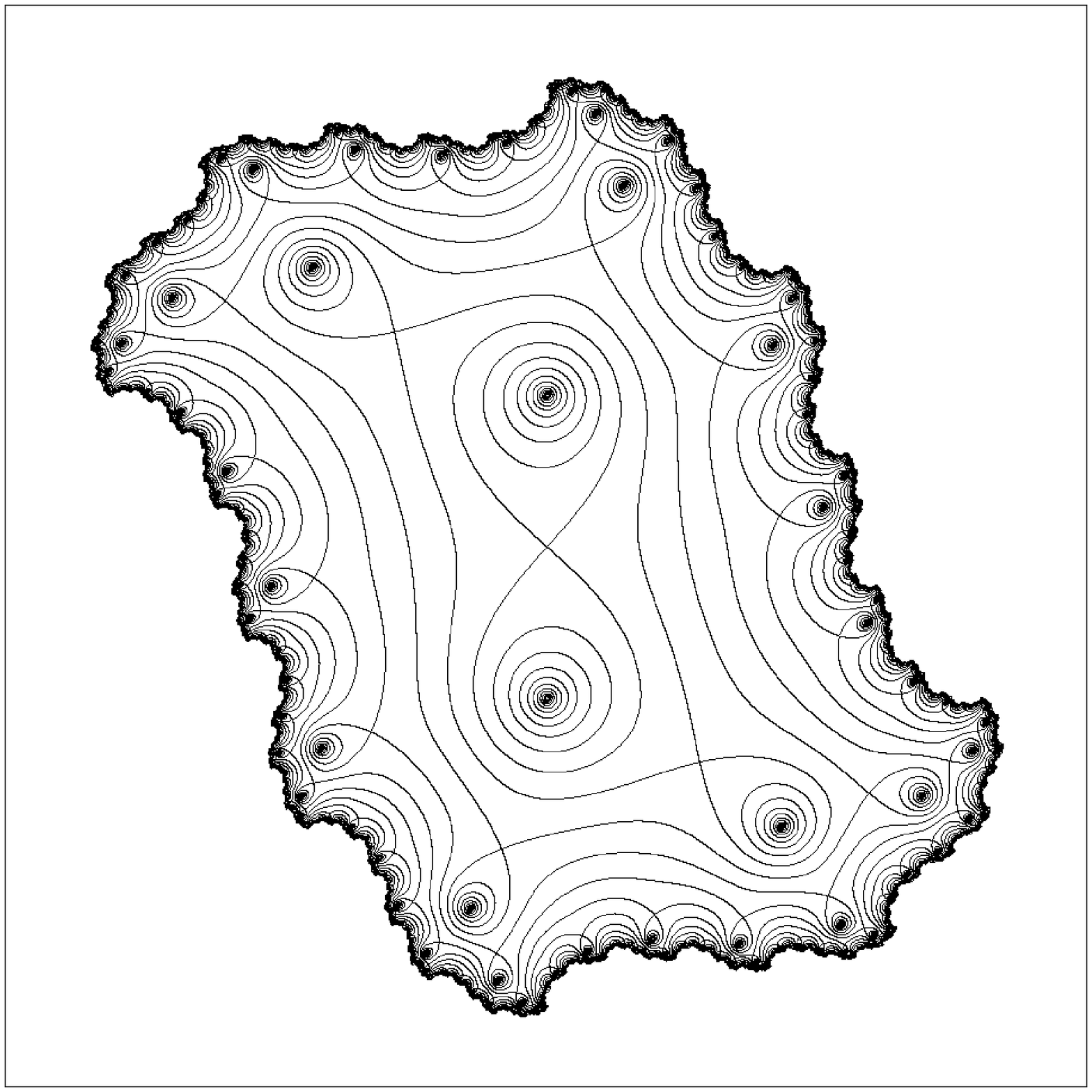 pixels 960 by 960 scaled 400
 \bigskip\bigskip
 \centerline{Figure 5. Julia set for \[z\mapsto z^2+.7z\], with curves
 \[|\phi|={\rm constant}\].}
 \vfil\endinsert

 The statement is the repelling case is somewhat different:

 {\QP{\bf 6.5. Corollary.} {\it If
 \[\hat z\] is a repelling fixed point, then there
 is a holomorphic map \[\psi:\C\to S\] in the opposite direction,
 so that the diagram
 $$      \matrix{S& \buildrel f\over\longrightarrow & S\cr
	 \quad\uparrow\psi & & \quad\uparrow\psi\cr
	 \C&\buildrel \lambda{\textstyle\cdot}\over\longrightarrow &\C\cr}   $$
 is commutative, and so that \[\psi\] maps a neighborhood of zero
 diffeomorphically onto a neighborhood of \[\hat z\]. Here
 \[\psi\] is unique except that it may be replaced by \[w\mapsto
 \psi(cw)\] for any constant \[c\ne 0\].}\medskip}

 To compute \[\psi(w)\] we simply choose some \[\lambda^{-n}w\] which
 is so small that \[\phi^{-1}(\lambda^{-n}w)\] is defined, and then
 apply \[f^{\circ n}\] to the result.\QED\smallskip

 Now suppose that \[f:\hat\C\to\hat\C\] is a rational function with an
 attracting fixed point \[\hat z\].
 By the {\bit immediate basin} \[\Omega_0(\hat z)\]
 we mean the connected component of \[\hat z\] in the
 basin of attraction \[\Omega=\Omega(\hat z)\], or equivalently the
 connected component of \[\hat z\] in the Fatou set \[\hat\C\ssm J\].
 (Compare 4.3.) The following is due to Fatou and Julia.

 {\QP{\bf 6.6. Theorem.} \it If \[f\] has degree two or more,
 then the immediate basin of any attracting fixed point \[\hat z\] of \[f\]
 contains at least one critical point.
 Furthermore, if the multiplier \[\lambda\] is not zero, then there exists a
 unique compact neighborhood \[\bar U\] of \[\hat z\] in \[\Omega_0\] which:

 \noindent $(a)$ maps bijectively onto some round disk \[\bar D_r\] under the
 \Koenigs map \[\phi\], and

 \noindent $(b)$ has at least one critical point on its boundary \[\partial U\].
 \medskip}

 \noindent Evidently \[\bar U\] can be described as
 the largest neighborhood which maps bijectively to
 a round disk centered at the origin. As an example, in Figure 5 the region
 \[U\] is bounded by the top half of the central figure 8 shaped curve.
 \smallskip

 {\bf Proof of 6.6.} If \[\lambda=0\], then \[\hat z\] itself is critical.
 Thus we may
 assume that \[\lambda\ne 0\] and apply 6.1. Evidently some branch
 \[\phi_0^{-1}\] of the inverse map
 can be defined as a single valued holomorphic function over some small disk
 \[D_\epsilon\], with \[\phi_0^{-1}(0)=\hat z\]. Let us try to extend
 \[\phi_0^{-1}\] by analytic continuation along radial lines through the
 origin in \[D_\epsilon\]. It cannot be possible to extend indefinitely far in
 every direction; for then \[\phi_0^{-1}\] would be a non-constant holomorphic
 map from the entire plane \[\C\] into the basin \[\Omega_0(\hat z)\]. This
 is impossible, since this basin is Hyperbolic.
 Thus there must exist some largest
 radius \[r\] so that \[\phi_0^{-1}\] extends analytically throughout the open
 disk \[D_r\]. Let \[U=\phi_0^{-1}(D_r)\]. We must prove that the closure
 \[\bar U\] is a compact subset of the basin
 \[\Omega(\hat z)\], and also that there is at least one critical point of \[f\]
 on the boundary \[\partial U\].
 If \[z_1\in\partial U\] is any boundary point, then
 using the identity \[\phi(f(z))=\lambda\phi(z)\in D_{|\lambda|r}\] for \[z\in U\]
 arbitrarily close to \[z_1\], we see that \[f(z_1)\] belongs to the open
 set \[U\subset \Omega\]. Therefore \[z_1\] also belongs to the basin
 \[\Omega\], with \[|\phi(z_1)|=r\] by continuity.
 Now at least one such \[z_1\] must be a critical
 point of \[f\]. For whenever \[z_1\] is non-critical we can continue
 \[\phi_0^{-1}\] analytically throughout a neighborhood of the image point
 \[\phi(z_1)\in\partial D_r\] simply by chasing around Diagram (2),
 composing the map \[z\mapsto
 \phi_0^{-1}(\lambda z)\] with the branch of \[f^{-1}\] which carries \[f(z_1)\]
 to \[z_1\].\QED

 For further information about the attractive basin \[\Omega_0\] see \S13.4,
 and also \S17.1.\medskip

 Next let us consider the {\bit superattracting} case \[\lambda=0\]. The
 following was proved by B\"ottcher in 1904.\smallskip

{\bf Historical Note:} L. E. B\"ottcher was born in Warsaw in 1872. He
took his doctorate in Leipzig in 1898, working in Iteration Theory, and
then moved to Lvov. He published in Polish and Russian. (The Russian form
of his name is ~{\cyr B\"etkher\cdprime}$\, $.)

 {\QP {\bf 6.7. Theorem of B\"ottcher.} {\it Suppose that
 $$	f(z)\;=\;a_nz^n\,+\,a_{n+1}z^{n+1}\,+\,\cdots\;,	$$
 where \[n\ge 2\], \[a_n\ne 0\]. Then there exists a local
 holomorphic change of coordinate \[w=\phi(z)\] which conjugates \[f\] to the
 \[n$-th power map \[\;w\mapsto w^n\;\]
 throughout some neighborhood
 of \[\phi(0)=0\]. Furthermore, \[\phi\] is unique up to multiplication by
 an \[(n-1)$-st root of unity.}\medskip}

 Thus near any critical fixed point, \[f\] is conjugate to a map of the form
 $$	\phi\circ f\circ\phi^{-1}\,:\,w\;\mapsto\;w^n\,,	$$ 
 with \[n\ne 1\]. This Theorem is most often applied in
 the case of a fixed point at infinity. For example, any polynomial
 map \[p(z) =a_nz^n+a_{n-1}z^{n-1}+\cdots+ a_0\] of degree \[n\ge 2\] has
 a superattractive fixed point at infinity. It follows easily from 6.7
 that there is a
 local holomorphic map \[\psi\] taking infinity to infinity which conjugates
 \[p\] to the map \[w\mapsto w^n\] around \[w=\infty\].\smallskip

 The proof is quite similar to the \Koenigs proof. It is only
 necessary to first make a logarithmic change of coordinates, and to be
 careful since the logarithm is not a single valued function. Suppose, to fix
 our ideas, that we consider a map having a fixed point at infinity,
 with a Laurent series expansion of the form
 $$	f(z)\; =\;a_nz^n+a_{n-1}z^{n-1}+\cdots+a_0+a_{-1}z^{-1}+\cdots $$
 with \[n\ge 2\],
 convergent for \[|z|>r\]. Note first that the linearly conjugate map\break
 \[z\mapsto \alpha f(z/\alpha)\], where \[\alpha^{n-1}=a_n\], has leading
 coefficient equal to \[+1\]. Thus, without loss of generality, we may
 assume that \[f\] itself has leading coefficient \[a_n=1\].
 Then\break \[f(z)=z^n(1+{\cal O}|1/z|)\] for \[|z|\] large,
 where \[{\cal O}|1/z|\] stands for some expression which is bounded by a
 constant times \[|1/z|\].
 Let us make the substitution \[z=e^Z\], where \[Z\]
 ranges over the half-plane \[{\cal R}(Z)>\log(r)\]. Then \[f\] lifts to
 a continuous map
 $$	F(Z)=\log f(e^Z)\;,	$$
 which is uniquely
 defined up to addition of some multiple of \[2\pi i\]. With correct choice of
 this lifting \[F\], it is not hard to check that
 \[\;F(Z)\;= \;nZ+{\cal O}(e^{-{\cal R}(Z)})\;\]
 for \[{\cal R}(Z)\] large. In fact, we will only need the weaker
 statement that
 $$	|F(Z)-nZ|<1	\eqno (3)	$$
 for \[{\cal R}(Z)\] large. Let us choose \[\sigma>1\] to be large enough
 so that the inequality (3) is satisfied for all \[Z\] in the half-plane
 \[{\cal R}(Z)>\sigma\]; note that \[F\] necessarily maps this half-plane
 into itself. Note the identity
 \[F(Z+2\pi i)=F(Z)+2\pi in\], which follows since \[F(Z+2\pi i)-F(Z)\]
 is a multiple of \[2\pi i\] which differs from \[n(Z+2\pi i)-nZ\] by at most
 2. If \[Z_0\mapsto
 Z_1\mapsto\cdots\] is any orbit under \[F\] in this half-plane, then we have
 \[|Z_{k+1}-nZ_k|<1\]. Setting \[W_k=Z_k/n^k\], it follows that
 $$	|W_{k+1}-W_k|<1/n^{k+1}	\,. $$
 Thus the sequence of holomorphic functions \[W_k=W_k(Z_0)\] converges
 uniformly and geometrically as \[k\to\infty\] to a holomorphic limit
 \[\Phi(Z_0)=\lim_{k\to\infty}\; W_k(Z_0)\].
 Evidently this mapping \[\Phi\] satisfies the identity
 $$	\Phi(F(Z))\;=\;n\Phi(Z)\,.	$$
 Note also that \[\Phi(Z+2\pi i)=\Phi(Z)+2\pi i\]. Therefore the mapping
 \[\phi(z)=\exp(\Phi(\log z))\] is well defined near infinity, and
 satisfies the required identity \[\phi(f(z))=\phi(z)^n\].

 To prove uniqueness, it suffices to study mappings \[w\mapsto \eta(w)\]
 near infinity which satisfy \[\eta(w^n)=\eta(w)^n\]. Setting \[
 \;\eta(w)= c_1w+c_0+c_{-1}w^{-1}+\cdots\;,	\]
 this becomes
 $$	c_1w^n+c_0+c_{-1}w^{-n}+\cdots\;=\; (c_1w+c_0+\cdots)^n\;=\;
	 c_1^nw^n+nc_1^{n-1}c_0w^{n-1}+\cdots\,.	$$
 Since \[c_1\ne 0\], \[c_1\] must be an \[(n-1)$-st root of unity,
 and an easy induction shows that the remaining coefficients are zero.\QED
 \smallskip

 {\bf Caution.} In analogy with 6.4, one might hope that the change of
 coordinates\break \[z\mapsto\phi(z)\] extends throughout the entire basin
 of attraction of the superattractive point as a holomorphic mapping.
 (Compare \S17.3.) However, this is not always possible. Such an extension
 involves computing expressions of the form
 $$	z\mapsto \root n \of {\phi(f(z))}\;,	$$
 and this does not work in general since the \[n$-th root cannot be
 defined as a single valued function. For example, there is trouble
 whenever some other point in the basin maps exactly onto the superattractive
 point, or whenever the basin is not simply-connected.\smallskip

 We conclude with a problem.
 \lln

 {\bf Problem 6-1.} What maps to what in Figure 5?
 \vfil
\eject

\footline={\hss\tenrm 7-\folio\hss}\pageno=1

 \centerline{\bf \S7. Parabolic Fixed Points: the Leau-Fatou Flower.}
 \medskip

 Again we consider functions \[\;   f(z)\;=\;\lambda z+a_2 z^2+a_3 z^3+\;
 \cdots    \;\]
 which are defined and holomorphic in some neighborhood of the origin,
 but in this section
 we suppose that the multiplier \[\lambda\] at the fixed point
 is a root of unity, \[\lambda^q=1\]. Such a fixed point
 is said to be {\bit parabolic}, provided that \[f^{\circ q}\] is not the
identity map. (More generally, any periodic orbit
 with \[\lambda\] a root of unity is called parabolic, provided that no
iterate of \[f\] is the identity map.) First consider
 the special case
 \[\lambda=1\]. It will be convenient to write our map as
 $$	f(z)\;=\;z+az^{n+1}+\,({\rm higher\;terms})\,,	\eqno (7.1) $$
 with \[a\ne 0\]. The integer
 \[n+1\ge 2\] is called the {\bit multiplicity}$\,$ of the fixed point.
 (By definition, the `` simple" fixed points with \[\lambda\ne 1\] have
 multiplicity equal to 1.) Choose a neighborhood \[N\]
 of the origin which is small enough so that \[f\] maps \[N\]
 diffeomorphically onto some neighborhood \[N'\] of the origin.\smallskip

 {\bf Definition.} A connected open set \[U\],
 with compact closure \[\bar U\subset N\cap N'\], will be called an
 {\bit attracting petal}$\,$ for \[f\] at the origin if
 $$	f(\bar U)\subset U\cup\{0\}\qquad{\rm and}\qquad
	 \bigcap_{k\ge 0}\;f^{\circ k}(\bar U)=\{0\}\,.	$$
 Similarly, \[U'\subset N\cap N'\] is a {\bit repelling petal}$\,$ for \[f\]
 if \[U'\] is an attracting petal for \[f^{-1}\].\smallskip

 {\QP{\bf 7.2. Leau-Fatou Flower Theorem.} {\it If the origin is a fixed
 point of multiplicity \[n+1\ge 2\], then there exist \[n\] disjoint
 attracting petals \[U_i\] and \[n\] disjoint repelling petals \[U'_i\] so
 that the union of these \[2n\] petals, together with the origin itself,
 forms a neighborhood \[N_0\] of the origin. These petals
 alternate with each other, as illustrated in Figure 6, so
 that each \[U_i\] intersects only \[U'_i\] and \[U'_{i-1}\] $($where
 \[U'_0\] is to be identified with \[U'_n\,)\].}\medskip}

 If \[U_i\] is an attracting petal, then evidently
 the sequence of maps \[f^{\circ k}\] restricted to \[\bar U_i\]
 converges uniformly to zero.
 On the other hand, if \[U_i'\] is a repelling
 petal, then every orbit \[z_0\mapsto z_1\mapsto\cdots\] which starts out in
 \[U_i'\] must eventually leave \[U_i'\], and in fact must leave the union
 \[U_1'\cup\cdots\cup U_n'\]. (However it may later return, perhaps
 even infinitely often.) Here are three immediate consequences of 7.2.

 {\QP{\bf 7.3. Corollary.} {\it There is no periodic orbit,
 other than the fixed point at the origin, which is
 completely contained within the neighborhood \[N_0\].}\medskip}

 \pageinsert
\insertRaster 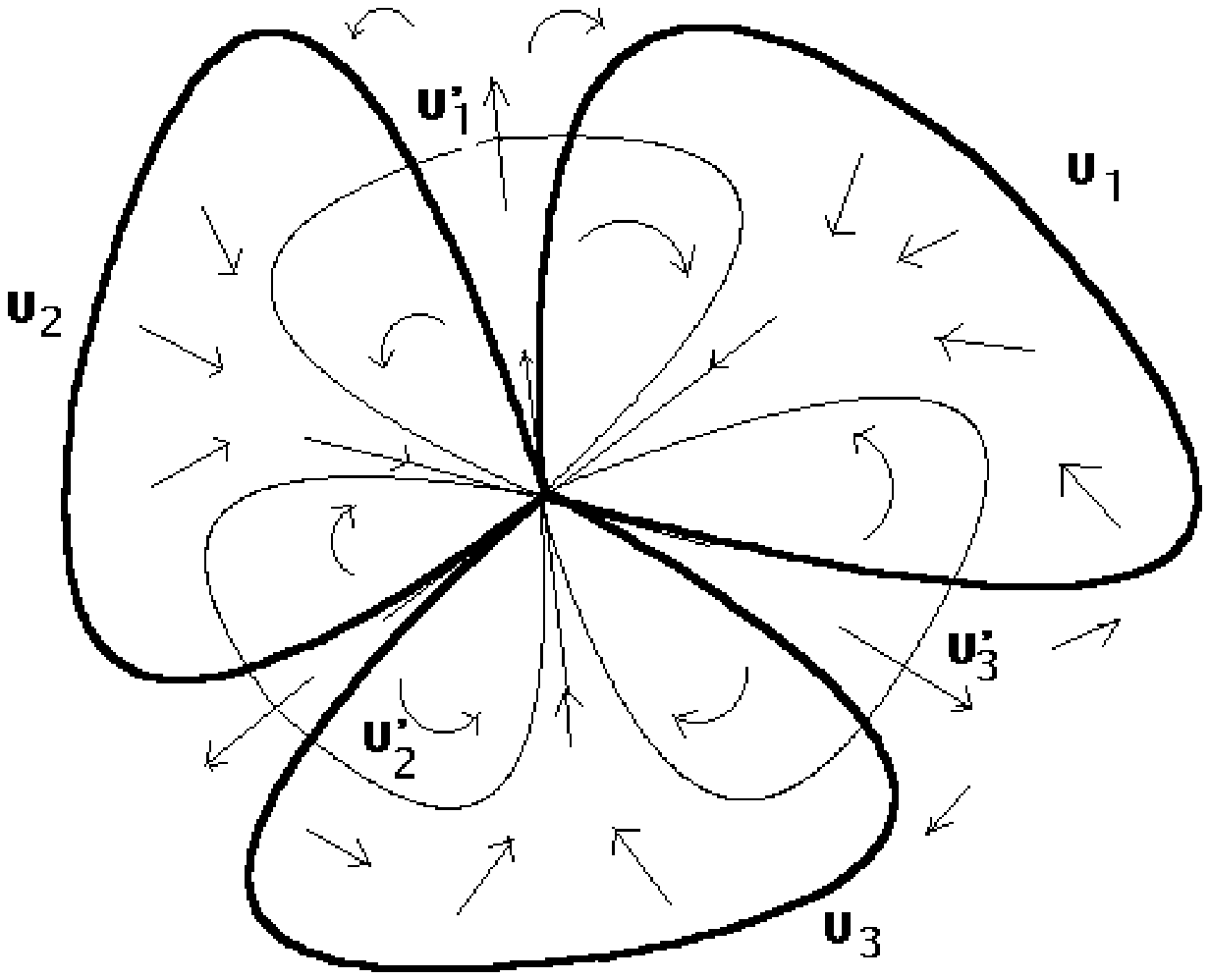 pixels 640 by 480 scaled 500
\smallskip
 \centerline{Figure 6. Leau-Fatou Flower with three attracting}
\centerline{ petals \[U_i\] and three repelling petals \[U'_i\].}
\bigskip
\insertRaster 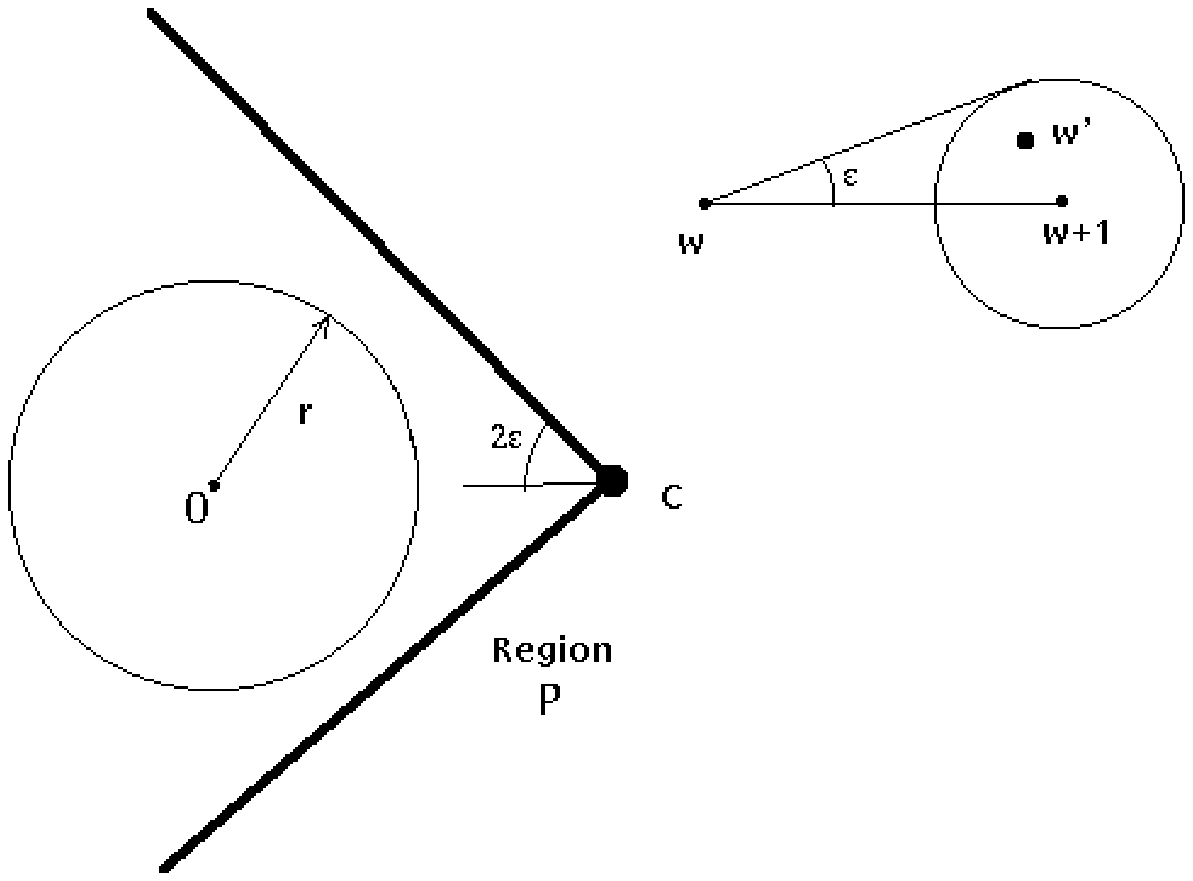 pixels 640 by 480 scaled 550
 \centerline{Figure 7.}\bigskip
 \endinsert

 Now suppose that \[f\] is a globally defined rational function. We continue
 to assume that the origin is a fixed point with \[\lambda=1\].
 Each attracting petal \[U_i\] determines a parabolic {\bit basin of attraction}
 \[\Omega_i\], consists of all \[z_0\] for which the orbit
 \[z_0\mapsto z_1\mapsto\cdots\] eventually lands in the attracting
 petal \[U_i\], and hence converges to the fixed point through \[U_i\].
 Evidently these basins \[\Omega_1\,,\,\ldots\,,\, \Omega_n\] are disjoint open
 sets.

 {\QP {\bf 7.4 Corollary.} {\it If we exclude the case of an orbit
 which exactly hits the fixed point, 
 then an orbit \[z_0\mapsto z_1\mapsto\cdots\] under \[f\] converges to
 the fixed point if and only if it eventually
 lands in one of the attracting petals \[U_i\], and hence belongs to
 the associated basin \[\Omega_i\].}\medskip}
 \eject % ???
 % \advancepageno

 {\QP{\bf 7.5. Corollary.} {\it Each parabolic
 basin \[\Omega_i\] is contained in
 the Fatou set \[\hat\C\ssm J(f)\], but each basin boundary \[\partial \Omega_i\]
 is contained in the Julia set \[J(f)\].
 It follows that each repelling petal \[U'_i\] must intersect
 \[J(f)\].}\medskip}

 \noindent In particular, it is claimed
 that the parabolic fixed point \[z=0\] must
 belong to \[J(f)\].\smallskip

 {\bf Proof of 7.5.} We first show that \[0\in
 J(f)\]. It follows from 7.1 that
 $$	f^{\circ k}(z)\;=\; z+kaz^{n+1}+{\rm (higher\;terms)}\,. $$
 Evidently no sequence of iterates \[f^{\circ k}\] can converge uniformly
 in a neighborhood of the origin, since the corresponding
 \[(n+1)$-st derivatives do not converge. (Compare 1.3.)
 Thus \[0\in J(f)\], and it follows
 that every point in the grand orbit of zero belongs to \[J(f)\].
 If \[z_1\in\partial \Omega_i\] is not in the grand orbit of zero,
 then by 7.4 we can extract
 a subsequence from the orbit of \[z_1\] which remains bounded away from
 zero. Since the sequence of iterates \[f^{\circ k}\] converges to zero
 throughout the open set \[\Omega_i\], it follows that \[\{f^{\circ k}\}\] can not
 be normal in any neighborhood of the boundary point \[z_1\]. The proof
 is now straightforward.\QED\medskip

 {\bf Proof of Theorem 7.2.} We will say that a vector \[v\in\C\] points in
 an {\bit attracting direction} at the fixed point of 7.1 if
 the product \[av^n\] is
 real and negative. If we ignore higher order terms, then these are just the
 directions for which the vector from \[v\] to \[f(v)\approx v(1+av^n)\]
 points straight in towards the origin. Similarly, \[v\] points in
 a {\bit repelling
 direction} if \[av^n\] is real and positive.
 Evidently there are \[n\] equally spaced attracting directions which are
 separated by the \[n\] equally spaced repelling directions.

 We will make use of the substitution \[w=b/z^n\] with inverse
 \[z=\root n \of {b/w}\], where \[b=-1/(na)\].
 Evidently the sector between two repelling
 directions in the \[z$-plane will correspond under this transformation
 to the entire \[w$-plane, slit along the negative real axis. In particular,
 a neighborhood of zero in such a sector will correspond to a neighborhood
 of infinity in such a slit \[w$-plane. Let us write the transformation
 7.1 as
 $$	z\mapsto f(z)=z(1+az^n+o|z^n|)	$$
 as \[|z|\to 0\]. Here \[o|z^n|\] stands for a remainder term whose ratio
 to \[|z^n|\] tends to zero.
 Substituting \[z=(b/w)^{1/n}\], the corresponding
 self-transformation in the \[w$-plane is
 $$	w\mapsto w'\;=\;b/f(z)^n\;=\;(b/z^n)(1+az^n+o|z^n|)^{-n}\;=\;
	 w(1-naz^n+o|z^n|)	\,.	$$
 But \[z^n=b/w\] and \[nab=-1\], so this can be written simply as
 $$	w'\;=\; w(1+w^{-1}+o|w^{-1}|)\;=\; w+1+o(1)	$$
 as \[|w|\to\infty\]. In other words, given any small number, which it will
 be convenient to write as \[\sin\epsilon>0\], we can choose a radius \[r\]
 so that
 $$	|w'-w-1|<\sin\epsilon\quad{\rm for}\quad |w|>r\,.	$$
 It follows that the slope of the vector from \[w\] to \[w'\] satisfies
 \[\;|{\rm slope}|<\tan
 \epsilon\,\], as long as \[|w|>r\]. Now we can construct an ``attracting
 petal for the point at infinity'' in the \[w$-plane as follows.
 Let \[P\]
 consist of all \[w=u+iv\] with \[|w|>r\], and with \[u>c-|v|/\tan 2\epsilon\],
 where the constant \[c\] is large enough so that all points \[w\in P\]
 satisfy \[|w|>r\]. (Figure 7.) Then an easy geometric argument shows that
 the closure \[\bar P\] is mapped into \[P\], and that every backward orbit
 starting in \[\bar P\] must eventually leave \[\bar P\]. Translating
 these statements back to the \[z$-plane, the proof can easily be
 completed.\QED

 Now suppose that the multiplier \[\lambda\] is a \[q$-th root of unity,
 say \[\lambda=\exp(2\pi ip/q)\] where \[p/q\] is a fraction in lowest
 terms. Then we can apply the discussion above to the \[q$-fold iterate
 \[f^{\circ q}\].

 {\QP{\bf 7.6. Lemma.} {\it If the multiplier \[\lambda\] at a fixed point
 \[f(z_0)=z_0\] is a primitive \[q$-th root of unity, then the 
 number \[n\] of attractive petals around \[z_0\] must be a multiple of \[q\].
 In other words, the
 multiplicity \[n+1\] of \[z_0\] as a fixed point of
 \[f^{\circ q}\] must be congruent to 1 modulo \[q\].}\medskip}

 Intuitively, if we perturb \[f\] so as to change \[\lambda\] slightly, then
 the multiple fixed point of \[f^{\circ q}\] will split up into one
 point which is still fixed by \[f\] together with some finite collection
 of orbits which have period \[q\] under \[f\]. This Lemma can be proved
 geometrically by showing that multiplication by \[\lambda=f'(z_0)\] must
 permute the \[n\] attractive directions at \[z_0\]. It can be proved
 by a formal power series computation based on the observation that
 \[f\circ f^{\circ q}\;=\;f^{\circ q}\circ f\]. Details will be left to
 the reader.\QED
 % \eject
 % \pageno=33

 As an example, Figure 8 shows part of the Julia set for the polynomial
 \[z\mapsto z^2+\lambda z\] where \[\lambda\] is a seventh root of unity,
 \[\lambda=e^{2\pi it}\] with \[t=3/7\]. There are seven attractive petals
 about the origin.\medskip

 We can further describe the geometry around a parabolic
 fixed point as follows.
 As in 7.1, we consider a local analytic map with a fixed point
 of multiplier \[\lambda=1\]. Let \[U\]
 be either one of the \[n\] attracting petals or one of the \[n\]
 repelling petals, as described in the Flower Theorem, \S7.2.
 Form an identification space \[U/f\] from \[U\] by identifying \[z\]
 with \[f(z)\] whenever both \[z\] and \[f(z)\] belong to \[U\]. (This
 means that \[z\] is identified with \[f(z)\]
 for every \[z\in U\] in the case of an attracting petal,
 and for every \[z\in U\cap f^{-1}(U)\] in the case of a repelling petal.)
 By definition, a holomorphic map \[\alpha:U\to\C\] is {\bit univalent}
 if distinct points of \[U\] correspond to distinct points of \[\C\].
 The following was proved by Leau and Fatou.
 \midinsert
 \vfil
 \insertRaster 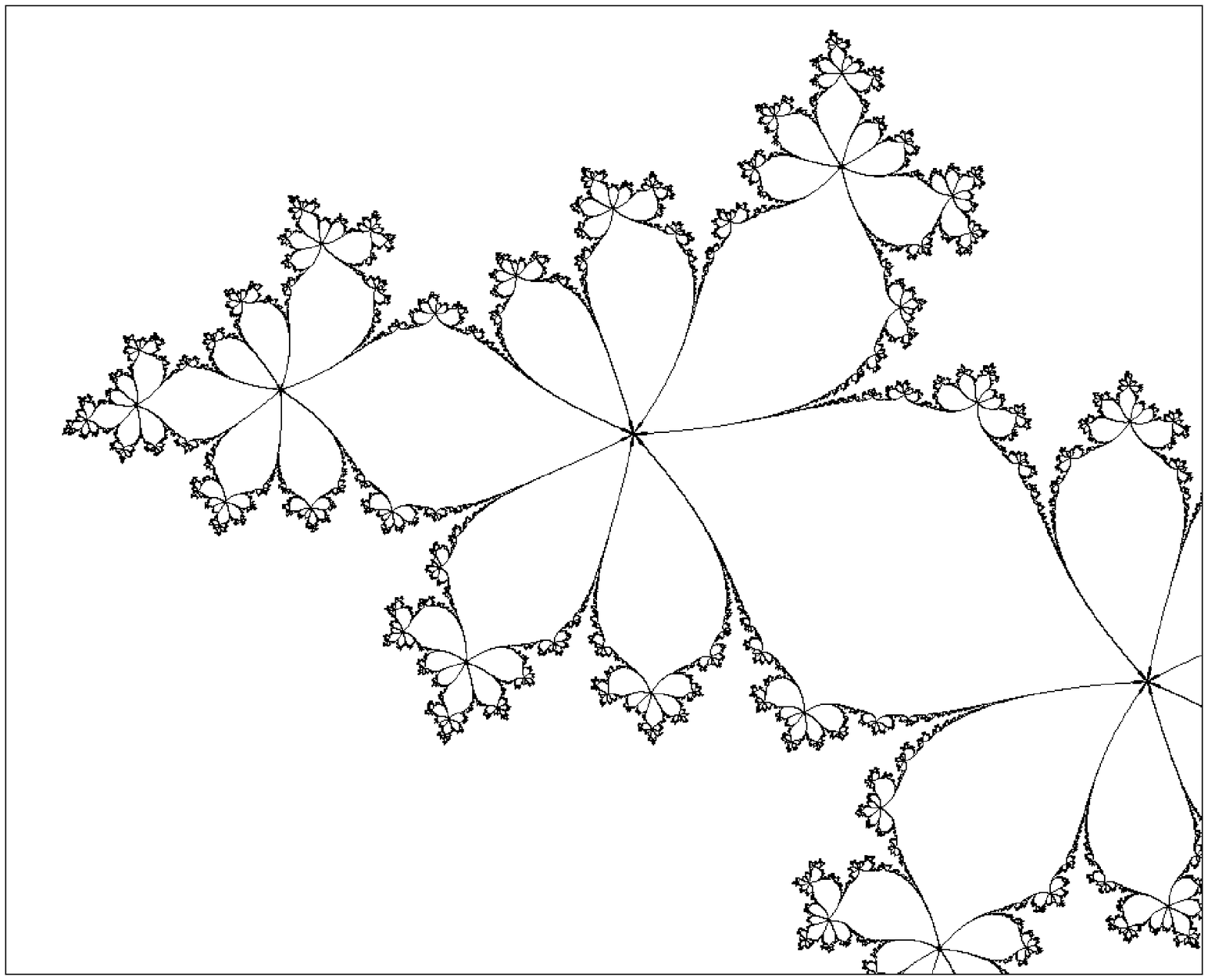 pixels 1008 by 816 scaled 350
 \bigskip\bigskip
 \centerline
 {Figure 8. Julia set for \[z\mapsto z^2+e^{2\pi it}z\] with \[t=3/7\].}
 \vfil\endinsert

 {\QP{\bf 7.7. Theorem.} {\it The quotient manifold \[U/f\] is conformally
 isomorphic to the infinite cylinder \[\C/\Z\]. Hence there is one, and
 up to composition with a translation only one, univalent embedding \[\alpha\]
 from \[U\] into the universal covering space \[\C\] which satisfies the
 {\bit Abel functional equation}
 $$	\alpha(f(z))\;=\;1+\alpha(z)	$$
 for all \[z\in U\cap f^{-1}(U)\]. With suitable choice of \[U\], the image
 \[\alpha(U)\subset \C\] will contain
 some right half-plane \[\{w:{\cal R}(w)>c\}\] in the case of an attracting
 petal, or some left half-plane in the case of a repelling petal.}\medskip}

 By definition, the quotient \[U/f\] is called an {\bit \'Ecalle cylinder}
 for \[U\]. (This term is due to Douady, suggested by the work of \'Ecalle
 on holomorphic maps tangent to the identity.)\smallskip
 % \eject
 % \advancepageno

 The proof of 7.7 begins as follows.
 To fix our ideas, we consider only the case of an
 attracting petal. As in the proof of 7.2, a substitution of the form
 \[w=b/z^n\] will conjugate the map \[f\] of \S7.1 to a map
 which has the form
 $$	 g(w)\;=\; w + 1 + a_1 w^{-1/n} + a_2 w^{-2/n}+\,\cdots\;. $$
 Here \[w\] ranges over a neighborhood of infinity, with the negative
 real axis removed.\smallskip

 {\bf Definition.} Let \[{\cal G}\] be the group
 consisting of all holomorphic maps which are defined and univalent in
 some region of the form
 $$ \{\;u+iv\in\C\;:\; u> c_1-c_2|v|\;\}\;,	\eqno (1)	$$
 and which are
 asymptotic to the identity map as \[|u+iv|\to\infty\]. Evidently our
 map \[w\mapsto g(w)\] belongs to this group \[{\cal G}\]. Our
 object is to show that \[g\] is conjugate to the translation \[w\mapsto w+1\]
 within \[{\cal G}\]. More generally, we will prove the following.

 {\QP{\bf 7.8. Lemma.} \it If a transformation \[g_0\in{\cal G}\] has the form
 \[g_0(w)=w+1+o(1)\] as \[|w|\to\infty\], then \[g_0\] is conjugate
 within \[{\cal G}\] to the translation \[w\mapsto w+1\].\medskip}

 {\bf Proof.} We assume that \[g_0\] can be written as
 \[g_0(w)=w+1+\eta_0(w)\], where \[\eta_0(w)\to 0\]
 as \[|w|\to \infty\] within some region of the form (1).
 We will first make two preliminary transformations to improve the
 error bound. Let
 $$	F_0(w)=\int(1+\eta_0(w))^{-1}dw	$$
 be any indefinite integral of \[\,1/(1+\eta_0)\,\]
 within this region. Then it is not difficult to check that \[F_0\in{\cal G}\].
 Using the Schwarz Lemma (\S1.3), note that \[|\eta_0'(w)|\,=\,o(1/|w|)\]
 and hence
 $$	F_0''(w)\,=\,o(1/|w|)\,,	$$
 within a smaller region of the same form.
 By Taylor's Theorem we have
 $$\eqalign{	F_0\circ g_0(w)\;=\;F_0(w+(1+\eta_0(w)))\;&=\; F_0(w)+F_0'(w)
	 (1+\eta_0(w))+o(1/|w|)\cr
	 &=\;F_0(w)+1+o(1/|w|)\;.\cr}	$$
 In other words, setting \[\;g_1\,=\,F_0\circ g_0\circ F_0^{-1}\,\] we have
 \[\;g_1(w)\,=\,w+1+o(1/|w|)\,\].
 Now repeating exactly this same construction, we see that \[g_1\] is conjugate
 to a map\break \[g_2=F_1\circ g_1\circ F_1^{-1}\] within \[{\cal G}\], where
 \[\;g_2(w)\;=\;w+1+o(1/|w|^2)\;\]
 within some smaller region of the same form. In particular,
 $$	|g_2(w)-w-1|\;\le \;1/|w|^2	\eqno (2) $$
 provided that \[|w|\] is sufficiently large.

 Starting with any \[w_0\] in this region, consider the orbit
 \[w_n=g_2^{\circ n}(w_0)\]. We will prove that the differences \[\{w_n-n\}\]
 form a Cauchy sequence which converges locally uniformly, so that the limit
 $$	w_0\mapsto\phi(w_0)=\lim_{n\to\infty}\,(w_n-n)	$$
 defines a transformation \[\phi\] which belongs
 to the group \[{\cal G}\]. Since \[\phi\circ g_2(w)\] is
 evidently equal to \[\phi(w)+1\], this will prove Lemma 7.8, and hence
 prove 7.7.

 As a preliminary remark, using the weaker inequality \[g_2(w)=w+1+o(1)\],
 we see that for any \[\epsilon>0\] we have \[|w_{n+1}-w_n-1|<
 \epsilon\] for \[|w_0|\] sufficiently large, and hence \[|w_n-w_0-n|
 <n\epsilon\]. In particular, it is not difficult to check that
 \[|w_n|\;\ge\; |w_0+n|/2\], and hence
 $$	|w_{n+1}-w_n-1|\;\le\; 1/|w_n|^2\;\le\;4/|w_0+n|^2\;,	$$
 provided that \[|w|\] is large. For \[m>n\ge 0\], this implies that
 $$	|(w_m-m)-(w_n-n)|\,\;<\;\sum_{n\le j<\infty}\; 4/|w_0+j|^2\;\approx\;
	 4\int_n^\infty dj/|w_0+j|^2\;.	$$
 Setting \[w_0+n=re^{i\theta}\] with \[|\theta|<\pi\], this integral can be
 evaluated as \[\theta/(r\sin\theta)\le c/r\], for some constant \[c\]
 depending on the region. This tends to zero as \[r=|w_0+n|\to\infty\],
 hence the \[w_n-n\] form a Cauchy sequence.
 Further details will be left to the reader.\QED

 {\bf Remark.} Note that this preferred Fatou coordinate system is defined
 only within one of the \[2n\] attracting or repelling petals. In order to
 described a full neighborhood of the parabolic fixed point, we would
 have to describe how these \[2n\] Fatou coordinate systems are to be pasted
 together in pairs by means of univalent mappings. In fact each of the
 \[2n\] required pasting maps has the form \[w\mapsto w+\Upsilon(e^{\pm
 2\pi i w})\], where \[\Upsilon\] is defined and holomorphic in
 some neighborhood of the origin, and where the signs \[\pm\] alternate.
 By studying this construction,
 one sees that there can be no normal form depending on only finitely many
 parameters for a general holomorphic map \[f\] in the neighborhood of
 a parabolic fixed point.\medskip

 Now suppose that \[f:\hat\C\to\hat\C\] is a globally defined rational map.
 Although attracting petals behave much like repelling petals in the local
 theory, they behave quite differently in the large.

 {\QP{\bf 7.9. Corollary.} {\it If \[U\] is an attracting petal, then the
 Fatou map
 $$ 	\alpha:U\to\C	$$
 extends uniquely to a map which is defined
 and holomorphic throughout the attractive basin \[\Omega\] of \[U\], still
 satisfying the Abel equation \[\alpha(f(z))\,=\,1+\alpha(z)\].}\medskip}

 \noindent This extended map \[\Omega\to\C\] is surjective. However, it is no
 longer univalent, but
 rather has critical points whenever some iterate \[f\circ \cdots\circ f\]
 has a critical point. In fact, we have the following.

 {\QP{\bf 7.10. Corollary.} \it For each attracting petal \[U_i\], the
 corresponding immediate basin \[\Omega_i\supset U_i\] contains at least one
 critical point of \[f\]. Furthermore, there exists a unique preferred
 petal \[U^*_i\] for this basin which maps precisely onto a right half-plane
 under \[\alpha\], and which has at least one critical point on its boundary.
 \medskip}

 \noindent The proofs are completely analogous to the corresponding proofs in
 \S6.4 and \S6.6. Thus \[\alpha^{-1}\] can be defined on some right half-plane,
 and if we try to extend leftwards by analytic continuation then we must run
 into some obstruction, which can only be a critical point of \[f\].
 (For an alternative proof that every attracting basin contains a critical
 point, see Milnor \& Thurston, pp. 512-515.)\QED

 As an example, Figure 9 illustrates the map \[f(z)=z^2+z\], with a parabolic
 fixed point of multiplier \[\lambda=1\] at \[z=0\], which is the cusp
 point at the right center of the picture. Here the Julia set
 \[J\] is the outer Jordan curve (the ``cauliflower'') bounding the basin
 of attraction \[\Omega\]. The critical point \[\omega=-1/2\] lies exactly at
 the center of the basin, and all orbits in this basin converge towards
 \[z=0\] to the right. The curves \[{\cal R}(\alpha(z))=
 {\rm constant}\in\Z\] have been drawn in. Thus the preferred petal \[U^*\],
 with the critical point on \[\partial U^*\], is bounded by the right half
 of the central figure \[\infty\] shaped curve.\medskip

 \midinsert
 \vfil
 \insertRaster 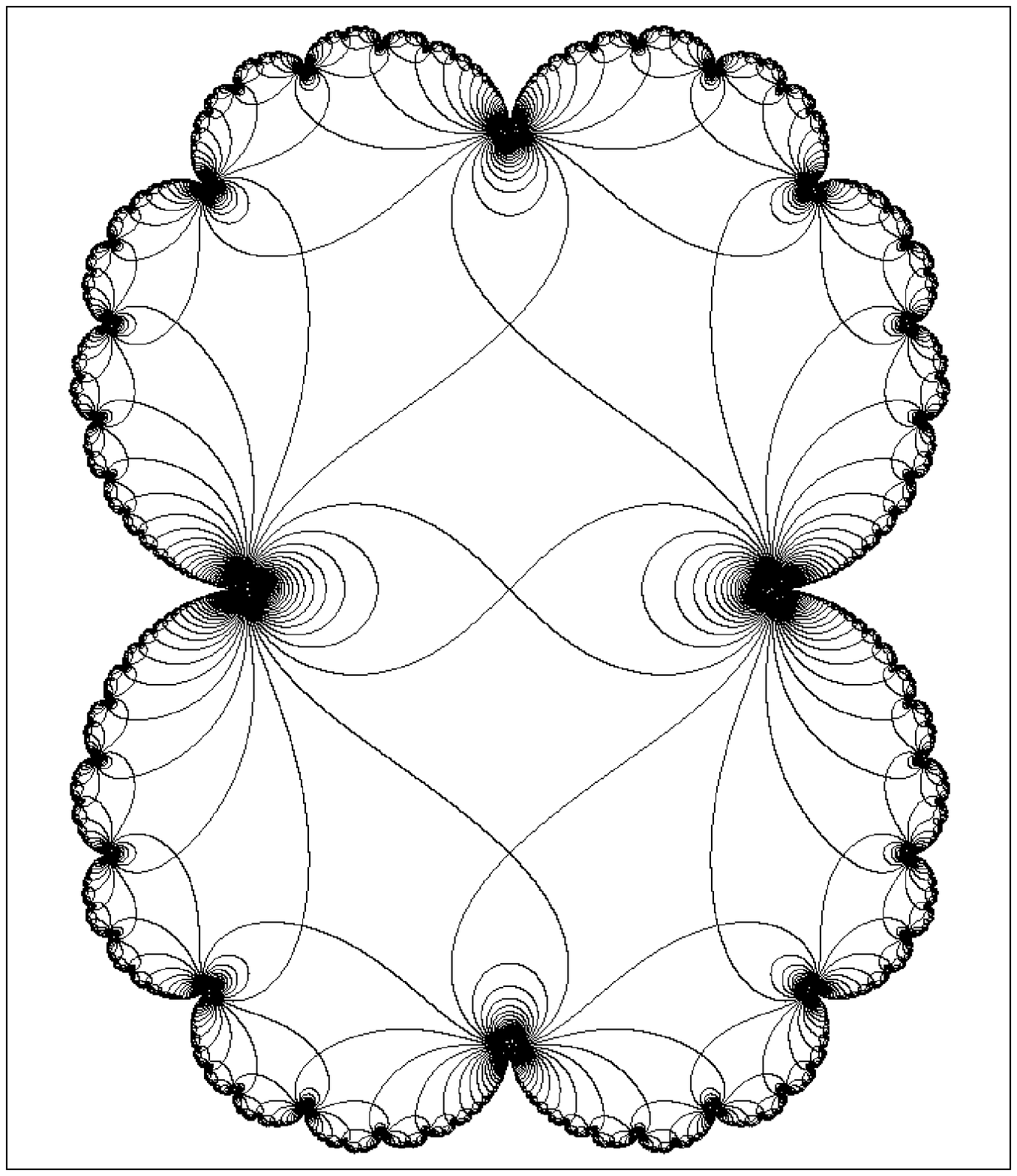 pixels 720 by 834 scaled 400
 \bigskip\bigskip
 \centerline{Figure 9. Julia set for \[\,z\mapsto z^2+z\,\], with the curves
 \[\,\alpha(z)\in\Z+i\R\,\] drawn in.}
 \vfil
 \endinsert

% \eject
 \noindent For a repelling petal, the corresponding statement is the following.

 {\QP{\bf 7.11. Corollary.} {\it If \[U'\] is a repelling
 petal, then the inverse map
 $$	\alpha^{-1}:\alpha(U')\to U'	$$
 extends uniquely to a globally defined holomorphic map \[\beta:\C\to\hat\C\]
 which satisfies the corresponding equation \[f(\beta(w))=\beta(1+w)\]. The
 image \[\beta(\C)\] is equal to the finite plane \[\C\] if \[f\] is a
 polynomial map,
 and is the entire sphere \[\hat\C\] if \[f\] is not conjugate to any
 polynomial.}\medskip}

 \noindent Again the proof is easily supplied. (Compare 6.5, together with 3.7
and Problem 3-3.)
 \QED\lln\bigskip

 {\bf Problem 7-1.} If \[z_0\] belongs to one of the basins
 of attraction \[\Omega_i\] of
 Corollary 7.4, with orbit \[z_0\mapsto z_1\mapsto z_2\mapsto\cdots\],
 show that \[\;\lim_{k\to\infty}\;z_k/|z_k|\;\] exists
 and is a unit vector which points in one of the \[n\] attracting directions.
 \smallskip

 {\bf Problem 7-2.} Define two attracting petals \[U\] and \[V\] for \[f\]
 to be {\bit equivalent} if every orbit for one intersects the other. Show that
 the petals which occur in 7.2 are unique up to equivalence. Show however
 that a petal as defined at the very beginning of \S7 may be too small,
 so that it cannot occur in 7.2, and so that the quotient \[U/f\] is not a
 full cylinder.

 \vfil\eject % \smallskip\bigskip

\footline={\hss\tenrm 8-\folio\hss}\pageno=1

 \centerline{\bf\S8. Cremer Points and Siegel Disks.}\medskip

Once more we consider holomorphic maps of the form
$$	f(z)\,=\,\lambda z+a_2z^2+a_3 z^3+\cdots\;, $$ % \eqno (8.1a)	$$
defined throughout some neighborhood of the origin.
% with a fixed point of multiplier \[\lambda\] at the origin.
% near the origin. 
In \S6 we supposed that \[|\lambda
|\ne 1\], while in \S7 we took \[\lambda\] to be a root of unity. This
section considers the remaining cases where \[|\lambda|=1\] but \[\lambda\]
is not a root of unity. Thus we assume that the multiplier \[\lambda\]
can be written as
$$	\lambda\;=\;e^{2\pi i\xi}\qquad{\rm with}\quad \xi\quad{\rm
real\; and \;irrational}\,.	$$ % \eqno (8.1b)	$$
Briefly, we will say that the origin is an {\bit irrationally indifferent\/}
fixed point. The number \[\xi\in\R/\Z\] may be descibed as the {\bit angle
of rotation\/} in the tangent space at the fixed point.\smallskip

The central question here is whether or not there exists a local change
of coordinate \[z=h(w)\] which conjugates \[f\] to the irrational rotation
\[\;w\mapsto\lambda w\,\], so that
$$	f(h(w))\;=\; h(\lambda w) $$
% $$	h(f(z)\;=\;\lambda h(z) $$
% $$	w\;\mapsto h\circ f\circ h^{-1}(w)\;=\;\lambda w $$
near the origin. (Compare \S6.)
This is the so called ``center problem''. If such a linearization is
possible, then a small disk \[|w|<\epsilon\] in the $w$-plane
corresponds to an open set \[U\] in the $z$-plane which is mapped
bijectively onto itself by \[f\]. Evidently such a neighborhood
\[U\] contains no periodic points of \[f\] other than the fixed
point at zero. If \[f\] is a rational function, note that \[U\] is
contained in its Fatou set \[\hat\C\ssm J\].
\smallskip

This section will first survey what is known about this problem, and then
prove some of the easier results.\smallskip
%was clearly posed in 
At the International Congress in
1912, E. Kasner conjectured that such a linearization
is {\bit always\/} possible. Five years later, G. A. Pfeiffer disproved
this conjecture
by giving a rather complicated description of certain holomorphic functions
for which no linearization is
possible. In 1919 Julia claimed to settle the question completely for
rational functions of degree two or more by showing that such a linearization
is {\bit never\/} possible. His proof was incorrect. H. Cremer
put the situation in much clearer perspective in 1927 with a beautiful
note which proved the following.\smallskip % (Compare 8.? below.)

{\bf Definition.} It will be convenient
to say that a property of a unit complex number
is true for {\bit generic\/} \[\lambda\in S^1\] if the set of \[\lambda\]
for which it is true contains a countable intersection of dense open
subsets of the circle. According to Baire, such a countable intersection
of dense open sets is necessarily dense and uncountably infinite.

{\QP{\bf Cremer Non-linearization Theorem.}
\it For a generic choice of \[\lambda\]
on the unit circle, the following is true. If \[z_0\] is a fixed
point of multiplier \[\lambda\] for a completely
arbitrary rational function of degree two or more, then \[z_0\]
is the limit of an infinite sequence of periodic
points. Hence % \[z_0\] belongs to the Julia set, and 
there is no
linearizing coordinate in a neighborhood of \[z_0\].\medskip}

\noindent (See 8.5 below.)
The question as to whether this statement is actually true for {\bit all\/}
numbers \[\lambda\] on the unit circle
remained open until 1942, when Siegel proved the following. (Compare 8.4
and 8.6.)
\eject

{\QP{\bf Siegel Linearization Theorem.}
\it For almost every \[\lambda\] on the unit
circle (that is for every \[\lambda\] outside of a set with one-dimensional
Lebesgue measure equal to zero) any germ of a holomorphic function
with a fixed point of multiplier \[\lambda\] can be linearized by a local
holomorphic change of coordinate.\medskip}

{\bf Remark.}  Thus there is a total contrast between behavior for
generic \[\lambda\]
and behavior for almost every \[\lambda\].
This contrast is quite startling, but is
not uncommon in dynamics. Here is quite different and equally remarkable
example. Consider the exponential function \[\exp:\C\to\C\] as a dynamical
system. {\it For a generic choice of \[z\in\C\],
the orbit of \[z\] is everywhere dense in
\[\C\]. On the other hand, for almost every \[z\in\C\] the
% $\omega$-limit set, that is the 
set of all accumulation points for the
orbit of \[z\], consists only of those points
$$	0\,,\;1\,,\;e\,,\;e^e\,,\;\ldots$$
which belong to the orbit
of zero.} (See Rees [R1], Lyubich [L2]. It is amusing to test this statement
on a computer:
From a random start, the iterated exponential usually either gets to zero
to within computer accuracy, or else gets too big to compute,
within five to ten iterations.)
In applied dynamics, it usually understood that
behavior which occurs for a set of parameter values of measure zero
has no importance, and can be ignored. However, even in applied dynamics
the study of generic behavior remains an extremely valuable tool.\medskip

{\bf Definition.} We will say that an irrationally indifferent fixed point is
a {\bit Siegel point\/} or a {\bit Cremer point\/} according
as a local linearization is possible or not. (In the classical literature,
Siegel points are called ``centers''.)\smallskip

{\QP{\bf 8.1. Lemma.} \it An irrationally indifferent fixed point of a
rational function is either a Cremer point or a Siegel point according
as it belongs to the Julia set or not. In the case of a Siegel point \[z_0\],
the entire connected component \[U\] of the Fatou set \[\hat\C\ssm J\]
which contains \[z_0\] is conformally isomorphic to the open unit disk
in such a way that the map \[f\] from \[U\] onto itself corresponds to
the irrational rotation \[w\mapsto\lambda w\] of the unit disk.\medskip}

\noindent
By definition, such a component \[U\] is called a {\bit Siegel disk\/},
or a {\bit rotation disk\/}.\smallskip

{\bf Proof of 8.1.} If \[z_0\] is a Siegel point, then the iterates of \[f\]
in a neighborhood correspond to iterated rotations of a small disk, and
hence form a normal family. Thus \[z_0\] belongs to the Fatou set.
Conversely, whenever \[z_0\] belongs to the Fatou set, we see easily from
Theorem 4.3 that \[z_0\] must be a Siegel point.\QED

Both Cremer and Siegel proved theorems which are much sharper than the rough
versions stated above. In order to state these precise results,
and their more recent generalizations, it is convenient to introduce
a number of different classes of irrational numbers, which are related
to each other as indicated in the following schematic diagram.

\centerline{\psfig{figure=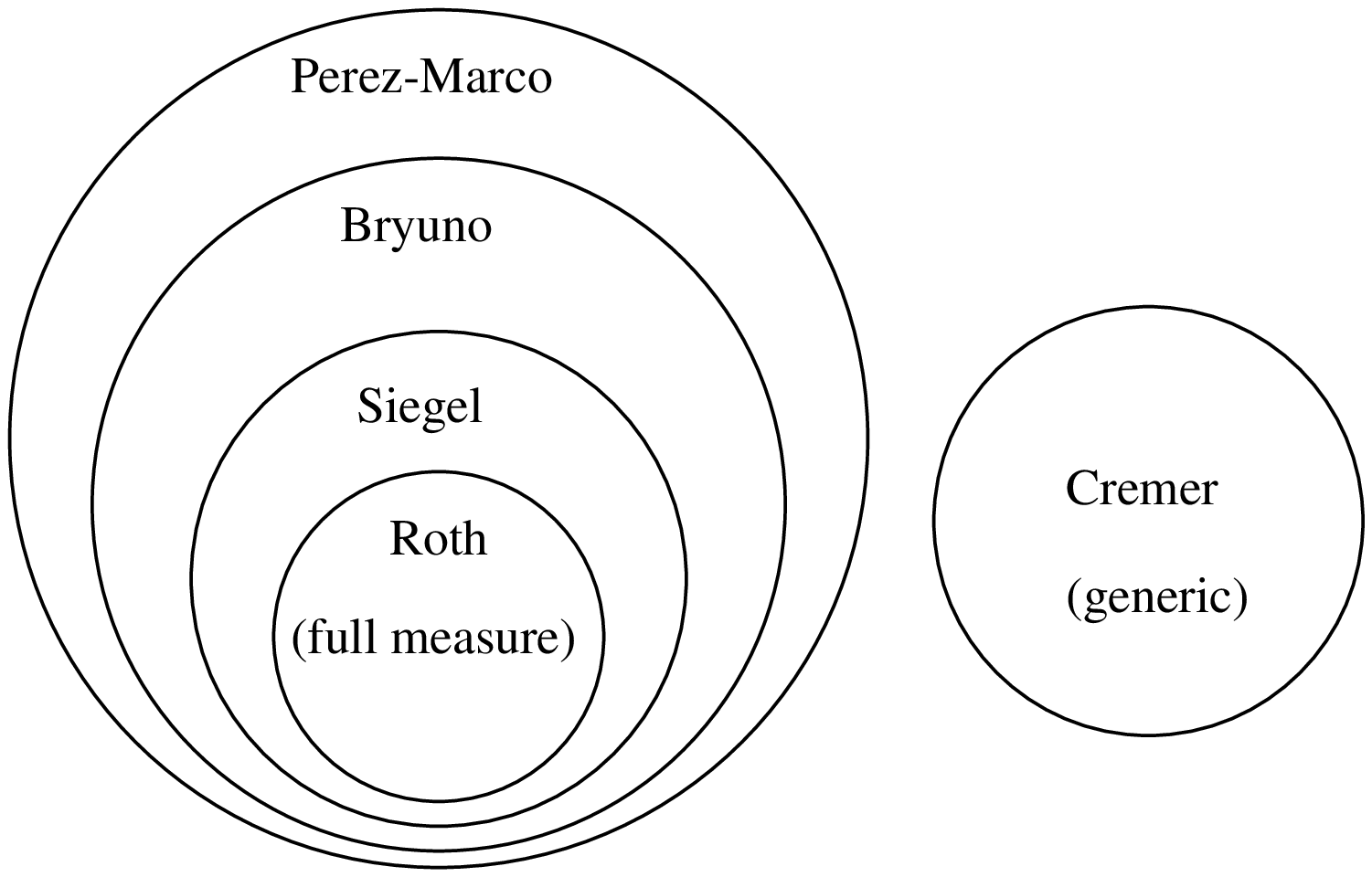,height=2.5in}}\medskip

\noindent
Roughly speaking, the Cremer numbers are those which can be approximated\break
extremely closely by rational numbers, while the Roth numbers are those
which can only be approximated badly by rationals.\smallskip

To give more
precise definitions, given some fixed real number \[\kappa\ge 2\] let us
% first introduce the set \[D_\kappa\] consisting
say that an irrational angle \[\xi\] satisfies a {\bit Diophantine
condition\/} of order \[\kappa\] if there exists some
\[\epsilon=\epsilon(\xi)>0\] so that
$$	\big|\xi-{p\over q}\big|\;>\;{\epsilon\over q^\kappa} $$
for every rational number \[p/q\]. Setting \[\lambda=e^{2\pi i\xi}\] as above,
since
$$	|\lambda^q-1|\;\;=\;\;|e^{2\pi i(q\xi-p)}-1|\;\;\sim\;\;
		% 2\pi |q\xi-p|\;=\;
	2\pi q|\xi-p/q| $$
as \[\;(q\xi-p)\to 0\,\], this
is equivalent to the requirement that
$$	|\lambda^q-1|\;>\;\epsilon'/ q^{\kappa-1} $$
for some \[\epsilon'>0\] which depends on \[\lambda\], and for all positive
integers \[q\]. Let \[D_\kappa\subset\R\ssm\Q\] be the set of all numbers
\[\xi\] which satisfy such a condition. Note that \[D_\kappa\subset
D_\eta\] whenever \[\kappa<\eta\]. % It is not hard to
% check that \[D_\kappa\] is vacuuous for \[\kappa<2\]. (Problem 8-??.)
We define the set \[\rm Si\] of {\bit Siegel numbers} (also called ``{\bit
Diophantine numbers\/}'')
to be the union of the \[D_\kappa\]. We can now make the following
more precise statement.

{\QP{\bf Theorem of Siegel.}
\it If the angle \[\xi\] belongs to this union \[\,{\rm Si}
=\bigcup D_\kappa\,\], then any holomorphic germ with multiplier
\[\lambda=e^{2\pi i\xi}\] is locally linearizable.\medskip}

\noindent
Proofs may be found in Siegel, or Siegel and Moser, or Carleson.\smallskip

A classical theorem of
Liouville asserts that every algebraic number of degree \[d\]\break belongs
to the class \[D_d\]. (Compare Problem 8-1.)
Hence every irrational number outside of the class
\[\rm Si\] must be transcendental. Such numbers in the complement of
\[{\rm Si}\] are often called {\bit Liouville numbers\/}.\smallskip

Define the set of {\bit Roth numbers\/} to be the intersection
$$ {\rm Ro}\;=\;\bigcap_{\kappa>2}\; D_\kappa\;. $$
%of \[D_\kappa\] over all real numbers \[\kappa>2\]. 

\medskip
\midinsert
\insertRaster 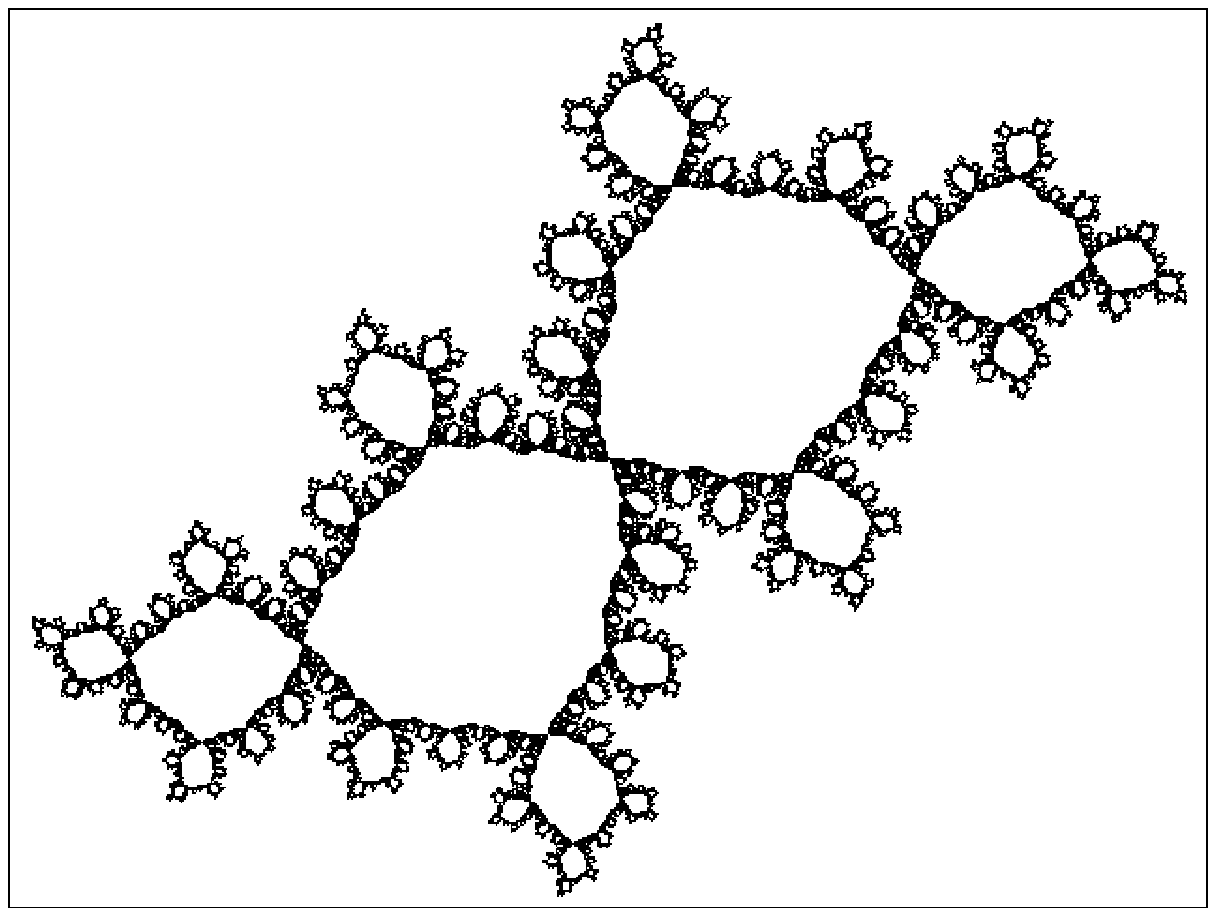 pixels 640 by 480 scaled 350
\smallskip
\centerline{Figure 10a. Julia set for \[z^2+e^{2\pi i\xi}z\]}
\centerline{ with \[\xi=\root 3 \of {1/4}=.62996\cdots\].}
\bigskip
\insertRaster 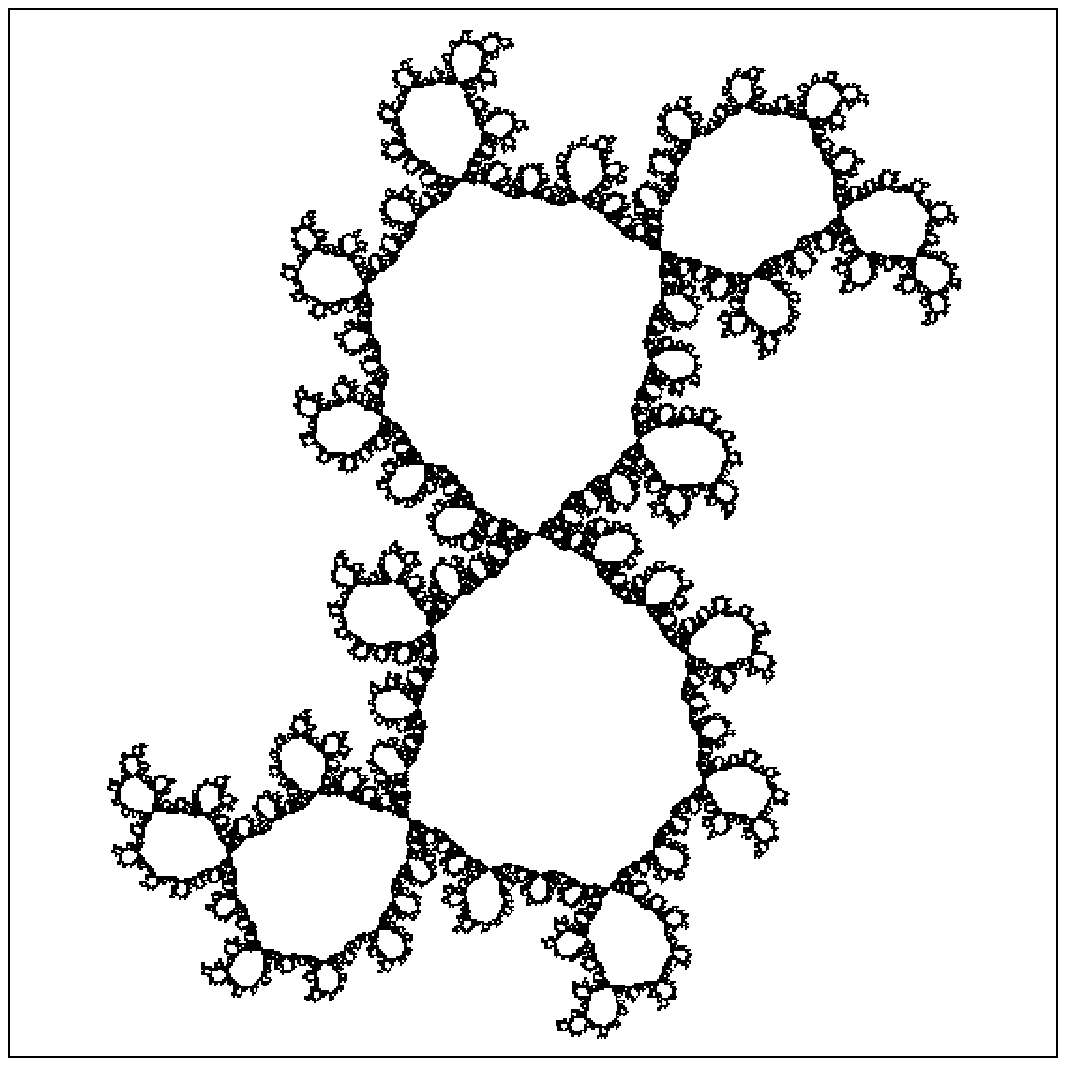 pixels 560 by 560 scaled 350
\smallskip
\centerline{Figure 10b. Corresponding Julia set with a}
\centerline{ randomly chosen angle \[\xi=.78705954039469\].}
\vfil\endinsert

Roth, in 1955, proved the much sharper result that
every algebraic number belongs to this intersection \[\rm Ro\].
It is quite easy to check
that: {\it Almost every real number belongs to \[\rm Ro\]}. (See Problem 8-2.)
\smallskip

{\it Thus if \[\xi\] is a completely arbitrary irrational algebraic number,
then any rational map, such as \[f(z)=z^2+e^{2\pi i\xi}z\],
which has a fixed point of multiplier \[e^{2\pi i\xi}\] must have a
Siegel disk.} Similarly, if \[\xi\] is a randomly chosen real number,
then the same will be true with probability one. Examples illustrating both
cases are shown in Figure 10.\medskip
%\eject

For a more precise analysis of the approximation of an irrational number
\[\xi\in (0\,,\,1)\] by rationals, it is useful
to consider the continued fraction expansion
$$	\xi\;\;=\;% 1/\bigg(a_1+\Big(1/(a_2+1/(a_3+\cdots\;))\Big)\bigg) \;, $$
\;\;{1\hfil\over
   \displaystyle a_1+{\strut 1\hfil\over \displaystyle a_2 +% \;\;{{}\atop
{\strut 1\hfil\over \displaystyle a_3 +\cdots}}} $$
%   \ddots}{{}\atop {{}\atop
%  \displaystyle \;\; +{1\hfil\over\displaystyle a_{n-1}}}}}}    $$
where the \[a_i\] are uniquely defined strictly positive integers. The
rational numbers
$$	{p_n\over q_n}\;\;=\;\;\;{1\hfil\over
   \displaystyle a_1+{\strut 1\hfil\over \displaystyle a_2 + \;\;{{}\atop
   \ddots}{{}\atop {{}\atop
  \displaystyle \;\; +{1\hfil\over\displaystyle a_{n-1}}}}}}	$$
are called the {\bit convergents\/} to \[\xi\]. The denominators \[q_n\]
will play a particularly important role. These denominators always grow at
least exponentially with \[n\]. In fact
$$	q_{n+1}\;>\; q_n\;>\; ((\sqrt 5 +1)/2)^{n-2}\;>\;1 $$
for \[n>2\]. We will need two basic
facts, which are proved in Appendix C. Each \[p_n/q_n\] is the best
approximation to \[\xi\] by rational numbers with denominator at most
\[q_n\]. In fact, setting \[\lambda=e^{2\pi i
\xi}\] as usual, we have the following.

{\QP{\bf 8.2. Assertion.}
\[	|\lambda^k-1|\;>\;|\lambda^{q_n}-1|\qquad{\it for}\qquad
	k= 1\,,\;2\,,\;\ldots\,,\;\;q_n-1\;\].\medskip}

\noindent
Furthermore, the error \[|\lambda^{q_n}-1|\] has the order of magnitude
of \[1/q_{n+1}\]. That is:

{\QP{\bf 8.3. Assertion.} \it There are constants \[\;0<c_1<c_2<\infty\;\]
so that
$$	{c_1\over q_{n+1}}\;\;\le\;\;
	 |\lambda^{q_n}-1|\;\;\le\;\; {c_2\over q_{n+1}}\qquad\qquad{\it
in\; all\; cases.} $$\smallskip}

\noindent For example we can take \[c_1=2\] and \[c_2=2\pi\].
Using these two facts, we can write the Roth condition as
$$\lim_{n\to\infty}\;{\log q_{n+1}\over\log q_n}\;=
	\;1\; ,\leqno \qquad	{\bf Ro:} $$
and the Siegel condition as
$$ \sup \;{\log q_{n+1}\over\log q_n}\;<\;\infty\;. \leqno \qquad{\bf Si:} $$
With these same notations, we now introduce the weaker {\bit Bryuno
condition}
$$	\sum_n\;{\log(q_{n+1})\over q_n}\;<\;\infty\;,\leqno \qquad{\bf Br:} $$
and also the much weaker {\bit Perez-Marco condition}
$$	\sum_n\;{\log\log(q_{n+1})\over q_n}\;<\;\infty\;.\leqno \qquad{\bf PM:} $$

\noindent
It is not difficult to check that \[\;{\bf Ro}\Rightarrow{\bf Si}\Rightarrow
{\bf Br}\Rightarrow{\bf PM}\].~ %\smallskip
Bryuno, in 1972, proved an extremely sharp version of Siegel's Theorem.

{\QP{\bf 8.4. Theorem of Bryuno.} \it If the angle \[\xi\] satisfies
the condition that\break \[\sum\log(q_{n+1})/q_n<
\infty\;,\;\] then any holomorphic germ of the form
$$	f(z)\;=\;e^{2\pi i\xi}z+a_2z^2+\cdots $$
can be linearized by a local holomorphic change of variable.\medskip}

\noindent A proof will be outlined at the end of this section.
Yoccoz, in 1987, showed that this result is best possible.

{\QP{\bf Theorem of Yoccoz.} \it Conversely, if \[\;\sum\log(q_{n+1})
/q_n=\infty\,\], then the quadratic map
$$	f(z)\;=\;z^2+e^{2\pi i\xi} z $$
%there exists a holomorphic germ \[\;z\mapsto e^{2\pi i\xi}z+a_2z^2+\cdots\;\]
has the property that every neighborhood of the origin contains infinitely
many periodic orbits. Hence the origin is a Cremer point.\medskip}
%has the property that the origin can be
%approximated arbitrarily closely by ``small cycles'', and hence is a Cremer
%point.\medskip}
% . Hence no such local linearization is possible.\medskip}
%\eject

\noindent
Briefly, we will sat that the fixed point can be approximated by ``small
cycles''.\break Perez-Marco, in 1990, completely
characterized the multipliers for which such small cycles must appear.
%\eject

{\QP{\bf Theorem of Perez-Marco.} \it If \[\xi\] satisfies the condition
that $$ \sum\log\log(q_{n+1})/q_n<\infty\;, $$
then any non-linearizable germ
with multiplier \[e^{2\pi i\xi}\] contains infinitely many periodic orbits
in every neighborhood of the fixed point. However, whenever \[\sum\log\log
(q_{n+1})/q_n=\infty\] there exists a non-linearizable
germ which has no\break periodic orbit other than the fixed point itself
within some neighborhood of the fixed point.\medskip}

\noindent (In the special case of a rational function, it is not known
whether every Cremer point can necessarily be approximated by small cycles.
Compare 8.5. below.)
\bigskip

The rest of this section will provide a few proofs. We
first prove a slightly sharper form of Cremer's
Theorem, and then a rather weak form of Siegel's Theorem. Finally, we give
a very rough outline proof for Bryuno's Theorem.\smallskip

We begin with Cremer's Theorem. Let us say that an irrational
angle \[\xi\] satisfies a ``{\bit Cremer condition\/}''
of degree \[d\] if the associated \[\lambda=e^{2\pi i\xi}\] satisfies
$$	{\rm lim \;sup}_{q\to\infty}\;{\log\log(1/|\lambda^q-1|)\over q}\;\;
	>\;\;\log d\;.  \leqno \quad{\bf Cr}_d: $$
Thus the error \[|\lambda^q-1|\] must tend to zero extremely rapidly
for suitable large \[q\].
This is equivalent to the hypothesis that
\[\quad {\rm lim\; sup}_{q\to\infty}\; q^{-1}\log\log(1/|\xi-p/q|)\;>
 \;\log d\;\], or to the hypothesis that
\[\quad{\rm lim\;sup}\;(\log\log q_{n+1})/q_n\,>\,\log d\;.\;\]
%We will say that \[\xi\] satisfies \[{\bf Cr}_\infty\] if this ~lim sup~
%is infinite, so that \[\xi\] satisfies \[{\bf Cr}_d\] for every degree \[d\].
{\it It is not difficult to show that a generic real number satisfies this
condition \[{\bf Cr}_d\] for {\bit every} degree \[d\].\/}
% belongs to this intersection \[{\rm Cr}_\infty\].}
(See Problem 8-3.)
% Essentially following Cremer's 1927 paper, we have the following.

{\QP{\bf 8.5. Theorem.}
\it If \[\xi\] satisfies \[{\bf Cr}_d\] with \[d\ge 2\],
then for a completely arbitrary rational function of degree \[d\], any
neighborhood of a fixed point of multiplier\break \[\lambda=e^{2\pi i\xi}\]
must contain infinitely many periodic orbits. Hence no local
linear\-ization is possible.\medskip}

\noindent % Briefly, we will say that such a fixed point {\bit
% can be approximated by small cycles\/}.
In particular, for a generic choice of \[\xi\] this statement will
be true for non-linear
rational functions of arbitrary degree.\smallskip
% , and for any fixed point of multiplier \[e^{2\pi i\xi}\].\smallskip

The proof which follows is nearly all due to Cremer. However, Cremer
used the slightly
weaker hypothesis that \[{\rm lim\;inf}\; |\lambda^q|^{1/d^q}\;=\;0\;\]
and concluded only that the fixed
point is a limit of periodic points, rather than full periodic orbits.
\smallskip

{\bf Proof of 8.5.} First consider a monic polynomial
\[\; f (z)\;=\;z^d+\cdots+\lambda z \;\]
of degree \[d\ge 2\] with a fixed point of multiplier \[\lambda\] at the
origin. Then \[\;f^{\circ q}(z)=z^{d^q}+\cdots+\lambda^q z\,\],\break
so the fixed points of \[f^{\circ q}\] are the roots of the equation
$$	z^{d^q}+\cdots+(\lambda^q-1)z\;=\; 0\;. $$
Therefore, the product of the \[d^q-1\]
non-zero fixed points of \[f^{\circ q}\]
is equal to \[\pm(\lambda^q-1)\]. If \[|\lambda^q-1|<1\], then it
follows that there exists at least one such fixed point \[z_q\] with
$$  0\;<\; |z_q|
	\;<\; |\lambda^q-1|^{1/(d^q-1)}\;<\; |\lambda^q-1|^{1/d^q}\;. $$
By hypothesis, for some \[\epsilon>0\], we can choose \[q\] arbitrarily large
with
$$	q^{-1}\log\log(1/|\lambda^q-1|)\;>\;\log(d)+\epsilon\;, $$
or in other words
$$	|\lambda^q-1|^{1/d^q}\;<\;\exp(-e^{\epsilon q})\;. $$
This tends to zero as \[q\to\infty\], so we certainly have % non-zero
periodic points \[z_q\ne 0\] in\break every neighborhood of zero. By Taylor's
Theorem, if \[\delta>0\] is
sufficiently small, then\break \[\;|f(z)|\,<\,e^\epsilon |z|\]
whenever \[|z|<\delta\]. It follows that
$$	|f^{\circ k}(z)|\;<\;\delta\quad{\rm for}\quad 1\le k\le q
\quad{\rm whenever}\quad |z|\;<\;
		e^{-\epsilon q}\delta\;. $$
Now note that there exist periodic points \[z_q\] which satisfy
the inequality $$ |z_q|\;<\;
\exp(-e^{\epsilon q})\;<\;e^{-\epsilon q}\delta $$
for arbitrarily large values of \[q\]. It follows that the
entire periodic orbit of such a point, with period at most \[q\], is
contained in the \[\delta\] neighborhood of zero.
Since \[\delta\] can be arbitrarily small, this completes
the proof of 8.5 in the polynomial case.\smallskip

In order to extend this argument to the case of a
rational function \[f\], Cremer first notes that \[f\]
must map at least one point \[z_1\ne 0\] to the fixed point \[z=0\]. After
conjugating by a M\"obius transformation which carries \[z_1\] to infinity,
we may assume that  \[f(\infty)=f(0)=0\]. If we set
\[f(z)=P(z)/Q(z)\], this means that \[P\] is a polynomial of degree
strictly less than \[d\], Furthermore, after
a scale change we may assume that \[P(z)=({\rm higher\;terms})+\lambda z\],
and that \[Q(z)=z^d+\cdots+1\] is monic. A brief computation then shows
that \[\;f^{\circ q}(z)=P_q(z)/Q_q(z)\;\] where
\[\;P_q(z)=({\rm higher\;terms})+\lambda^q z\;\]
and where \[Q_q(z)\] has the form \[\;z^{d^q}+\cdots+1\,\].
Thus the equation for fixed points of \[f^{\circ q}\] has the form
$$	0\;=\;zQ(z)-P(z)\;=\; z(z^{d^q}+\cdots+(1-\lambda^q))\;. $$
The proof now proceeds just as in the polynomial case.\QED\smallskip

For further information about Cremer points, see \S11.5 and \S18.6.\medskip

Let us next prove that Siegel disks really exist. We will
describe a proof, due to Yoccoz, of the following special
case of Siegel's Theorem. (Compare Herman [He2] or Douady [D2].) 
% \eject

{\QP{\bf 8.6. Theorem.} \it For Lebesgue almost every angle \[\xi\in\R/\Z\],
taking \[\lambda=e^{2\pi i\xi}\] as usual,
the quadratic map \[f_\lambda(z)=z^2+\lambda z\;\] possesses a
Siegel disk about the origin.\medskip}%\eject

{\bf Remark.}
Somewhat more precisely, we can define the {\bit size} of a
Siegel disk to be the largest number \[\sigma\] such that there exists
a holomorphic embedding \[\psi\] of the disk of radius \[\sigma\]
into the Fatou set \[\hat\C\ssm J(f_\lambda)\] so that \[\psi'(0)=1\], and
so that \[f_\lambda(\psi(w))=\psi(\lambda w)\].
If there is no Siegel disk, then we set \[\sigma=0\]. Using a normal
family argument, it is not difficult to show that this size \[\sigma\] is upper
semicontinuous as a function of \[\lambda\]. (Compare the proof of 8.8
below.) {\it In other words,
for any fixed \[\epsilon>0\] the set of \[\;\lambda=
e^{2\pi i\xi}\;\] with
\[\sigma(\lambda)\ge\epsilon\] is compact.} This set is totally disconnected,
since it contains no roots of unity. As \[\epsilon\to 0\], it
grows larger, and the proof will show that its
measure tends to the measure of the full unit circle.\smallskip

The proof of 8.6 has
three steps. The first two steps will be carried out here, while the
third will be put off to Appendix A. Here is the first step.
Consider the dynamics of \[f_\lambda\] for \[\lambda\] inside the open
disk \[D\]. According to K\oe nigs, for \[\lambda\in D\ssm\{0\}\],
there exists a neighborhood \[U\] of
zero and a holomorphic map \[\phi_\lambda(z)=\lim_{k\to\infty}
f_{\lambda}^{\circ k}(z)/\lambda^k\] which carries \[U\]
diffeomorphically onto some disk \[D_\rho\], so that \[f_\lambda\] on \[U\]
corresponds to multiplication by \[\lambda\] on \[D_\rho\], and so that
\[\phi_\lambda\] has derivative \[+1\] at the origin.

{\QP{\bf 8.7. Lemma.} {\it We can choose the open set
\[U\] so as to map diffeomorphically onto
the disk \[D_\rho\] of radius \[\rho=|\phi_\lambda(-\lambda/2)|\]
under \[\phi_\lambda\]. However
no larger radius is possible. Furthermore, the correspondence
$$	\lambda\;\mapsto\; \eta(\lambda)=\phi_\lambda(-\lambda/2)	$$
is bounded, holomorphic, and non-zero throughout the punctured disk \[D
\ssm\{0\}\].}\medskip}

{\bf Proof.} Note that \[-\lambda/2\] is the unique critical point, that is
the unique point at which the derivative
\[f'_\lambda\] vanishes. Thus the first assertion of 8.7 follows immediately
from \S6.6. Furthermore, it follows from \S6.2 that this correspondence
\[\lambda\mapsto\eta(\lambda)\] is holomorphic. To show that \[\rho=|\eta|\]
is bounded,
note first that \[U\] must be contained in the disk \[D_2\] of radius 2.
For if \[|z|>2\] then an easy estimate shows that \[|f_\lambda(z)|>|z|\],
so the orbit of \[z\] cannot converge to zero. Thus \[\phi^{-1}\] maps
the disk \[D_\rho\] holomorphically onto \[U\subset D_2\] with derivative
1 at the origin; so it follows from the Schwarz Lemma, \S1.3, that \[\rho
\le 2\].\QED

Since this function \[\eta\] is bounded, it follows that \[\eta\] has
a removable
singularity at the origin; that is, it can be extended as a holomorphic
function throughout the disk \[D\]. (See for example Ahlfors, 1966 p. 114,
or 1973 p. 20.)
\smallskip

Now consider
the {\bit radial limit} of \[\eta(r\exp(2\pi it))\] for fixed \[t\], as
\[r\to 1\].

{\QP{\bf 8.8. Lemma.} {\it Suppose, for some fixed \[\lambda=e^{2\pi i\xi}\],
that the quadratic map \[f_\lambda\] does not possess any Siegel disk.
Then the radial limit
$$	\lim_{r\to 1} \eta(r e^{2\pi i\xi})	$$
must exist and be equal to zero.}\medskip}

\noindent
Conversely, if the quantity \[\;\rho\;=\;\lim\,{\rm sup}_{r\to 1}\; |\eta(re^
{2\pi i\xi})|\;\] is strictly positive, the proof will show that
\[f_{\exp(2\pi i\xi)}\] admits a Siegel disk of ``size'' \[\ge\rho\].
\smallskip

{\bf Remark.}
Yoccoz has shown that this estimate is best
possible. That is, there exists a Siegel disk if and only if \[\rho>0\];
and \[\rho\] is precisely the size of the maximal Siegel disk,
as defined in the Remark following 8.6.
\smallskip

{\bf Proof of 8.8.} If the \[{\rm lim\, sup}\]
of \[|\eta(r\exp(2\pi i\xi))|\] as
\[r\to 1\] is equal to \[\rho_0>0\], then for any \[\rho<\rho_0\]
and for some sequence \[\lambda_j\in D\] tending to \[\lambda=\exp(2\pi i\xi)
\in \bar D\], the
inverse diffeomorphism \[\phi_{\lambda_j}^{-1}\] mapping \[D_\rho\] into
\[D_2\] is well defined. By a normal family argument, we can choose a
subsequence which converges, uniformly on compact sets, to a holomorphic
limit \[\psi:D_\rho\to \C\]. It is easy to check that this limit \[\psi\]
satisfies the required equation \[\psi(\lambda w)=f_\lambda(
\psi(w))\], and hence describes a Siegel disk.\QED

Finally, the third step in the proof of 8.6
is a classical theorem by F. and M. Riesz
which asserts that such a radial limit cannot exist and be equal to
zero for a set of \[\xi\] of positive Lebesgue measure. In other words
the quantity
$$	\lim\,{\rm sup}_{r\to 1}\; |\eta(re^{2\pi i\xi})|	$$
must be strictly positive for almost every \[\xi\]. This theorem
will be proved in Appendix A.3. Combining these three steps, we
obtain a proof of the special case 8.6 of Siegel's Theorem.\QED
\bigskip

To conclude this section, here is a very rough outline of a proof of the
Bryuno Theorem, due to Yoccoz. The proof is based on a ``renormalization
construction'' due to Douady and Ghys. Consider first a map \[f_1:D\to\C\]
which is univalent (that is, holomorphic and one-to-one) on the open
unit disk \[D\], with a fixed point of multiplier \[\lambda_1=e^{2\pi i\xi_1}\]
at the origin. We introduce a new coordinate by setting \[z=e^{2\pi iZ}\]
where \[Z=X+iY\] ranges over the half-plane \[Y>0\]. Then \[f_1\] corresponds
to a map of the form
$$	F_1(Z) \;=\; Z+\xi_1+\sum_1^\infty a_n e^{2\pi inZ} $$
which is defined and univalent on this upper half-plane. This map \[F_1\]
commutes with the translation \[T_1(Z)=Z+1\], and is approximately equal to
the translation	\[T_{\xi_1}(Z)=Z+\xi_1\]. More precisely, we have
$$	F_1(Z)\;=\; T_{\xi_1}(Z)+o(1)\;,
	\qquad{\rm uniformly\;in}\quad X\quad{\rm as}\quad Y\to\infty\;,
\; .$$
In fact, the  \[e^{2\pi inZ}\] terms
decrease exponentially fast as \[Y\to\infty\], so that if \[Y\]
is bounded well away from zero then \[F_1\] is extremely close to the
translation \[Z\mapsto Z+\xi_1\]. In particular, we can choose some
height \[h_1\] so that \[F_1\] moves points \[Z=X+iY\] definitely to the
right, and has derivative close to 1, throughout the half-plane \[Y>h_1\].
\smallskip

Construct a new Riemann surface \[S'_1\] as follows. Take a vertical strip
\[S_1\]
in the\break \[Z$-plane which is bounded on the left by the vertical line \[\;L
=\{iY : h_1\le Y<\infty\}\,\], on the right by its image \[F_1(L)\,\],
and from below by the straight line from \[ih_1\] to
\[F_1(ih_1)\]. Now glue the left edge to the right edge by \[F_1\].
The resulting Riemann surface \[S'_1\]
is conformally isomorphic to the punctured
unit disk. Hence it can be parametrized by a variable\break
\[w\in D\ssm\{0\}\].
It is convenient to fill in the puncture point, \[w=0\], which corresponds to
the improper points \[Z=X+i\infty\] in the \[Z$-plane.
We now introduce a holomorphic map \[f_2\] from a neighborhood of zero
in \[S'_1\] into
\[S'_1\] as follows. Starting from any point \[Z\] in the strip \[S_1\]
which is not too close to the bottom,
let us iterate the map \[F_1\] until
we reach some point \[1+Z'\] of the translated strip \[1+S_1\]. The
corrrespondence \[Z\mapsto Z'\] on \[S_1\] now yields the required
holomorphic map \[f_2\] from a neighborhood of zero in the quotient surface
\[S'_1\] to \[S'_1\]. Note that \[w\] is asymptotic to some constant times
\[e^{2\pi iZ/\xi_1}\] as \[Y\to\infty\]. If \[1+Z'=F_1^{\circ a}(Z)\approx
 Z+a\xi_1\] as \[Y\to\infty\], then \[2\pi iZ'/\xi_1\approx 2\pi iZ/\xi_1
+a-1/\xi_1\]. Setting \[1/\xi_1\equiv \xi_2\;(\mod 1)\] with \[0<\xi_2<1\],
it follows that
the corresponding map \[f_2(w)=w'\] in the disk \[S'_1\] is asymptotic to
$$	w\;\mapsto\; w\,e^{-2\pi i\xi_2}\; . $$
{\it Thus this Douady-Ghys construction relates a map \[f_1\]
with rotation angle \[\xi_1\] to a map \[f_2\] with rotation angle
\[\;-\xi_2\equiv -1/\xi_1\;\;(\mod 1)\].}

\bigskip\bigskip
\centerline{***  to be continued ***}
\bigskip\bigskip

See \S11.4 and \S12 for closely related results.
We conclude this section with some problems.\smallskip
\lln

{\bf Problem 8-1 (Liouville).} Let \[f\] be a polynomial of degree \[d\]
with integer
coeficients, and suppose that \[f(\xi)=0\] where \[\xi\] is irrational.
If every other root of
this equation has distance at least \[\epsilon\] from \[\xi\], and
if \[|f'(x)|<K\] throughout the \[\epsilon\] neighborhood of \[\xi\],
show that
$$	K|\xi-p/q|\;\ge\;|f(p/q)| \;\ge\; 1/q^d $$
for every rational number \[p/q\] in the \[\epsilon\] neighborhood
of \[\xi\]. Conclude that \[\xi\in D_d\], and hence that all irrational numbers
in the complement of \[{\rm Si}=\bigcup D_d\] must be transcendental.
\medskip
\eject

{\bf Problem 8-2.} If \[\kappa>2\] and \[\epsilon>0\], show that the set
\[S(\kappa\,,\,\epsilon)\] of numbers
\[\xi\in [0,1]\] which satisfy
\[|\xi-p/q|\,\le\,\epsilon/q^\kappa\] for some rational number
\[p/q\] has measure less than or equal
to \[\epsilon\sum q/q^\kappa<\infty\]. Since this tends to zero as
\[\epsilon\to 0\],
conclude that almost every real\break number belongs to \[D_\kappa\],
and hence that almost every real number belongs to the Roth set
\[{\rm Ro}=\bigcap_{\kappa>2}\; D_\kappa\]. (On the other hand, the
subset \[D_2\] has measure
zero, and \[D_\kappa\] is vacuous for \[\kappa<2\]. Compare Hardy
and Wright. Numbers in \[D_2\] are said to be ``of constant type''.)\medskip

{\bf Problem 8-3 (Cremer).} Given a completely arbitrary function
\[q\mapsto\eta(q)>0\], show that the set \[S_\eta\], consisting
of all irrational numbers \[\xi\] such that
$$	|\xi-{p\over q}|\;<\; \eta(q)\qquad{\rm
for\; infinitely\; many\; rational\; numbers}\; \;\;p/q\;, $$
is a countable intersection
of dense open subsets of \[\R\]. As an example, taking \[\phi(q)=2^{-q!}\]
conclude that a generic real number belongs to the set \[S_\phi\], which is
contained in the Cremer class \[{\rm Cr}_\infty\].\medskip

{\bf Problem 8-4 (Cremer 1938).} If \[f(z)=\lambda z+a_2z^2+a_3z^3+\cdots\],
where \[\lambda\] is not zero and not a root of unity,
show that there is one and only one formal power series of the form
\[h(z)=z+h_2z^2+h_3z^3+\cdots\] which formally satisfies the condition that
\[h(\lambda z)=f(h(z))\]. In fact
$$	h_n\;=\;{a_n+X_n\over \lambda^n-\lambda} $$
for \[n\ge 2\], where \[X_n=X(a_2\,,\,\ldots\,,\, a_{n-1}\,,\, h_2\,,\,
\ldots\,,\, h_{n-1})\]
is a certain polynomial\break expression whose value
can be computed inductively. Now suppose that
we choose the \[a_n\] inductively, always equal to zero
or one, so that \[|a_n+X_n|\ge 1/2\]. If
$$	{\rm lim\;inf}_{q\to\infty}\;\; |\lambda^q-1|^{1/q}\;=\; 0\;, $$
show that the uniquely defined power series \[h(z)\] has radius of
convergence zero. Conclude that \[f(z)\] is a holomorphic germ which is
not locally linearizable. Choosing the \[a_n\] more carefully, show that
we can even choose \[f(z)\] to be an entire function.

\vfill
\end

%% file: cd9-11.tex
\input header
\input ras
\footline={\hss\tenrm 9-\folio\hss}\pageno=1

\centerline{\bf GLOBAL FIXED POINT THEORY}\bigskip

\centerline{\bf \S9. The Holomorphic Fixed Point Formula}\medskip

First note the following.

{\QP{\bf 9.1. Lemma.} \it Every rational map \[f(z)\not\equiv z\]
of degree \[d\] has exactly \[d+1\]
fixed points, counted with multiplicity.\medskip}

Here, by definition, the {\bit multiplicity} of a fixed point
is equal to one whenever the multplier satisfies \[\lambda\ne 1\], and
is strictly greater than one otherwise. (Compare \S7.) As an example, in
the special case of a polynomial map of degree \[d\ge 2\],
there is exactly one fixed point at infinity, and hence
\[d\] finite fixed points counted with multiplicity.\smallskip

{\bf Proof.}
Conjugating \[f\] by a fractional linear automorphism if necessary, we
may assume that the point at infinity is not fixed by \[f\]. Hence we
can write \[f\] as a quotient \[f(z)=p(z)/q(z)\] where
the polynomial \[q(z)\] has degree equal to \[d\] and \[p(z)\] has
degree at most \[d\]. Now the equation
\[f(z)=z\] is equivalent to the polynomial
equation \[p(z)=zq(z)\] of degree \[d+1\], and the assertion follows.\QED
\smallskip

Both Fatou and Julia made use of a ``well known'' relation between the
multipliers at the fixed
points of a rational map. Consider first an isolated
fixed point \[f(z_0)=z_0\] of a holomorphic map in one complex variable.
If \[z_0\ne\infty\], we define the {\bit holomorphic index} of \[f\] at
\[z_0\] to be the residue
$$	\iota(f\,,\,z_0)\;=\;{1\over 2\pi i}\oint {dz\over z-f(z)}	$$
where we integrate in a small loop in the positive direction around \[z_0\].
As usual, if \[z_0\in\hat\C\] happens to be the point at infinity, then
we must first introduce the local uniformizing parameter \[\zeta=\phi(z)=1/z\].
In this case, we define the residue of \[f\] at \[\infty\] to be equal to
the residue of \[\phi\circ f\circ \phi^{-1}\] at the origin.

{\QP{\bf 9.2. Theorem.} {\it For any rational map \[\;f:\hat\C\to\hat\C\;\]
with \[f(z)\] not identically equal to \[z\], we have the relation
$$	\sum_{f(z)=z}\,\iota(f,z)\;=\; 1\,,	$$
to be summed over all fixed points.
In the case of a simple fixed point, where the multiplier
\[\lambda=f'(z_0)\] satisfies \[\lambda\ne 1\], the index is given by
$$	\iota(f,z_0)\;=\;{1\over 1-\lambda}\;.	$$
In any case, the index at \[z_0\] is a local analytic invariant. That
is, if \[\,g=\phi\circ f\circ\phi^{-1}\]
where \[\phi\] is a local holomorphic change of coordinate, then
\[\iota(f,z_0)=\iota(g,\phi(z_0))\].}\medskip}

\noindent
As an application, if \[f\] has one fixed point with multiplier very close
to 1, and hence with \[|\iota|\] large, then it must have at least one
other fixed point with \[|\iota|\] large and hence with \[\lambda\]
close to (or equal to) 1.\smallskip

{\bf Proof of 9.2.} Conjugating \[f\] by a linear fractional automorphism
if necessary, we may assume that the point at infinity is not fixed by \[f\].
Then \[f(z)\] remains bounded as \[|z|\to\infty\], and there is
a Laurent series expansion of the form
$$	\big( z-f(z)\big)^{-1}\;=\; z^{-1}+c_2z^{-2}+c_3z^{-3}+\cdots	$$
for large \[|z|\]. It follows that
the integral of \[\;{1\over 2\pi i}\cdot {dz\over z-f(z)}\;\] in a large loop
around the origin is equal to \[+1\]. Evidently this integral is equal to
the sum of the residues \[\iota(f,z_j)\] at the fixed points of \[f\];
hence the summation formula. The
computation of \[\iota(f,z_j)\] at a fixed point with multiplier \[\lambda_j
\ne 1\] is similar, using the Laurent series expansion of \[\big(z-f(z)
\big)^{-1}\] in a neighborhood of \[z_j\].

The proof that \[f\] is a local analytic
invariant, even at a multiple fixed point where \[\lambda=1\],
can be sketched as follows. Choose a 1-parameter family of perturbations
\[f_t\] of the given map \[f_0\]
so that \[f_t\] has distinct fixed points, all with multiplier
different from 1, for all small \[t\ne 0\]. For example we can set \[f_t(z) 
=f(z)+t\]. Thus a fixed point \[z_0\]
of \[f=f_0\] will split up into a cluster of nearby simple fixed points for
\[t\ne 0\].  Since the integral in a fixed loop around \[z_0\] varies
continuously with \[t\], and since this integral is a sum of residues which
are evidently local analytic invariants when \[t\ne 0\], it follows that
\[\iota(f_0\,,\,z_0)\] is also a local analytic invariant.\QED\smallskip

For generalizations of this formula, the reader is referred to
Atiyah and Bott.\medskip

{\bf Examples.} A rational map \[f(z)=c\] of degree zero has just one fixed
point, with multiplier zero and hence with index \[\iota(f,c)=1\]. A
rational map of degree one
usually has two distinct fixed points, and the relation
$$	{1\over 1-\lambda_1}+{1\over 1-\lambda_2}\;=\;1	$$
simplifies to \[\lambda_1\lambda_2=1\]. Any polynomial map \[p(z)\]
of degree two or more has a ``superattractive'' fixed point at infinity, with
multiplier zero and hence with index \[\iota(p,\infty)=1\]. Thus
the sum of the indices
of the {\it finite} fixed points of a polynomial map of degree \[\ge 2\]
is always zero.
For a polynomial of degree exactly two, the relation
$$      {1\over 1-\lambda_1}+{1\over 1-\lambda_2}\;=\;0 $$
for the finite fixed points simplifies to \[\;{1\over 2}
(\lambda_1+\lambda_2)=1\,\].\medskip

{\QP{\bf 9.3. Lemma.} {\it A fixed point with multiplier \[\lambda\ne 1\]
is attracting if and only if its index \[\iota\]
has real part \[{\cal R}(\iota)>{1\over 2}\].}\bigskip}

Geometrically, this is proved by noting that inversion carries the disk
\[1+D\] having the origin as boundary point to the appropriate half-plane.
Computationally, it can be proved by noting that \[{1\over 2}<{\cal R}
({1\over 1-\lambda})\] if and only if
$$	1\;<\;{1\over 1-\lambda}+{1\over 1-\bar\lambda}\;.	$$
Clearing denominators, this reduces easily to the required inequality
\[\lambda\bar\lambda<1\].\QED

One important consequence is the following.\medskip

{\QP{\bf 9.4. Corollary.} {\it Every rational map of degree two or more
must have either a repelling fixed point, or a fixed point
with \[\lambda= 1\], or both.}\bigskip}

{\bf Proof.} If the \[d+1\] fixed points were distinct and were all
attracting or neutral, then each index would have real part \[{\cal R}(\iota)
\ge{1\over 2}\], hence the sum would have real part greater than or equal to
\[{d+1\over 2}>1\]; but this would contradict the Fixed Point Formula.\QED

Since repelling points and parabolic points both belong to the Julia set,
this yields another proof of the following. (Compare \S4.4, as well as
\S4.3 and \S7.5.)

{\QP{\bf 9.5. Corollary.} \it The Julia set, for any rational map of
degree two or more, is always non-vacuous.\smallskip}
\lln\medskip

{\bf Problem 9-1.} If \[\;f(z)=z+\alpha z^2+\beta z^3+({\rm higher\;terms})
\,\],
with \[\alpha\ne 0\], show that the holomorphic index is given by
\[\iota(f,0)=\beta/\alpha^2\]. As an example,
consider the one-parameter family of cubic maps
$$	f_\alpha(z)\;=\;z+\alpha z^2+z^3	$$
with a double fixed point at the origin. Show that the remaining finite
fixed point of \[f_\alpha\]
is attracting if and only if \[\alpha^2\] lies within a unit disk centered
at \[-1\], or if and only if \[\alpha\] lies within a figure 8 shaped region
bounded by a lemniscate. Show that \[f_\alpha\] can be perturbed so that the
double fixed point at the origin splits up into two fixed points which
are {\it both} attractive
if and only if \[\alpha^2\] lies inside the disk of radius \[1/2\] centered
at \[1/2\], or if and only if \[\alpha\] lies within a region bounded
by a lemniscate shaped like the symbol \[\infty\].\smallskip

{\bf Problem 9-2.} Any fixed point \[z_0\] for \[f\] is evidently also a
fixed point for \[f^{\circ n}\]. If \[z_0\] is attracting [or repelling],
show that \[\iota(f^{\circ n},z_0)\] tends to the limit 1 [or 0] as
\[n\to\infty\]. If \[\;f(z)=z+\alpha z^k+({\rm higher\;terms})\,\], show that
\[\iota(f^{\circ n},0)\] tends to the limit \[k/2\].

\vfil\eject % \bigskip\bigskip
\footline={\hss\tenrm 10-\folio\hss}\pageno=1

\centerline{\bf \S10. Most Periodic Orbits Repel.}\bigskip

This section will prove the following theorem of Fatou. By a {\bit cycle}
we will mean simply a periodic orbit of \[f\]. Recall that a cycle is
called {\bit attracting}, {\bit neutral}, or {\bit repelling} according as its
multiplier \[\lambda\] satisfies \[|\lambda|<1\,,\;|\lambda|=1\], or
\[|\lambda|>1\].

{\QP{\bf 10.1. Theorem.} {\it Let \[f:\C\to\C\] be a rational map
of degree two or more. Then \[f\] has at most a finite number of
cycles which are attracting or neutral.}\medskip}

We will see in \S11 that there always exist infinitely many repelling
cycles. Shishikura has given
the sharp upper bound of \[2d-2\] for the number of attracting or neutral
cycles, using methods of quasi-conformal surgery. However the classical
proof, which is given here, shows only that this number is less than
or equal to \[6d-6\].\smallskip

First consider the case of an attracting cycle \[\,\hat z_0\mapsto \hat z_1\mapsto\cdots
\mapsto \hat z_m=\hat z_0\,\].
Each \[\hat z_j\] is an attracting fixed point for the \[m$-fold composition
\[f^{\circ m}\]. If \[\Omega_j\] is the immediate basin for \[\hat z_j\] under
\[f^{\circ m}\], then the union \[\Omega_0\cup\cdots\cup\Omega_{m-1}\] will be
called the {\bit immediate basin}$\,$ for our \[m$-cycle. It can be described
as the union of those components of the Fatou set which intersect the given
attracting cycle. We continue to assume that \[f\] has degree two or more.

{\QP{\bf 10.2. Lemma.} \it The immediate basin \[\Omega_0\cup\cdots\cup
\Omega_{m-1}\] for any attracting cycle
of \[f\] must contain at least one critical point of \[f\].\medskip}

{\bf Proof.} In the special case of an attracting fixed point, this was proved
in \S6.6. Applying this result to the \[m$-fold iterate of \[f\], we see
that the immediate basin \[\Omega_0\] for the fixed point \[\hat z_0\]
of \[f^{\circ m}\] must contain at least one critical point \[z_0\] of
\[f^{\circ m}\]. Let
\[\,z_0\mapsto z_1\mapsto\cdots\,\] be the orbit of this critical point
under \[f\]. By the chain rule,
at least one the the \[m\] points \[z_0\,,\,z_1\,,\,
\ldots\,,\,z_{m-1}\] must be a critical point for \[f\]. Since \[z_j
\in\Omega_j\], this proves the Lemma.\QED

Next consider a parabolic cycle \[\,\hat z_0\mapsto \hat z_1\mapsto\cdots
\mapsto \hat z_m=\hat z_0\,\], with multiplier
\[\lambda\] equal to a \[q$-th root of unity. Then the \[qm$-fold
iterate \[f^{\circ qm}\]
maps each \[\hat z_j\] to itself with multiplier \[\lambda^q\] equal to \[1\].
Recall that the number \[n\ge 1\] of attracting petals around \[\hat z_j\]
must be some multiple of \[q\]. (Compare 7.2 and 7.6.)\smallskip

{\bf Definition.} The union of
those components of the Fatou set which contain one of these \[n\] attracting
petals around one of the \[m\] points of the cycle will be called
the {\bit immediate basin} of this parabolic
cycle. Thus the number \[nm\] of connected components of
this immediate basin is some multiple of \[qm\].

{\QP{\bf 10.3. Lemma.} \it If \[f\] has degree \[d\ge 2\], then
the immediate basin for any parabolic cycle must also
contain at least one critical point.\medskip}

{\bf Proof.} In the case of a fixed point with multiplier \[\lambda=1\],
this was proved in \S7.10. The general case follows easily, just as in the
argument above.\QED

In fact, it is not difficult to check that there must be at least
\[n/q\] distinct critical points in such an immediate basin.
Combining 10.2 and 10.3 we obtain the following.

{\QP{\bf 10.4. Lemma.} {\it A rational map of degree \[d\ge 2\] can
have at most \[2d-2\] cycles which are attracting or parabolic.
Similarly, a polynomial map of degree \[d\ge 2\] can
have at most \[d-1\] cycles in the finite plane
which are attracting or parabolic.}\medskip}

{\bf Proof.} Since the immediate basins for distinct cycles are distinct
by definition, it follows from 10.2 and 10.3 that
the number of attracting or parabolic cycles is less
than or equal to the number of critical points. In the case of a polynomial
of degree \[d\], the number of finite critical points is clearly at
most \[d-1\]. In the case of a rational function of degree \[d\], by
the Riemann-Hurwitz formula \S5.1 the number of critical
ponts, counted with multiplicity, is equal to \[2d-2\]. Hence
the number of distinct critical points is at most \[2d-2\]. In either
case, the conclusion follows.\QED

{\bf 10.5. Remark.} This Lemma gives a practical computational procedure
for finding all attracting or parabolic cycles. We must simply follow the
orbits of all critical points, and test for convergence.

{\QP{\bf 10.6. Lemma.} {\it For a rational map of degree \[d\ge 2\],
the number of neutral cycles which have multiplier \[\lambda\ne 1\]
is at most \[4d-4\].}\medskip}

Evidently 10.4 and 10.6 together imply Theorem 10.1.
The proof of 10.6 begins as follows. Following Fatou, we perturb the given
map \[f\] and then applying 10.4. It will be convenient
to compare \[f\] with the map \[z\mapsto z^d\] of the same degree.
Note that this model map has no neutral cycles. It has
superattractive fixed points at zero and infinity,
but otherwise all of its periodic points are strictly repelling, with
\[|\lambda|=d^m>1\] where \[m\] is the period. Let
\[f(z)=p(z)/q(z)\] where \[p(z)\] and \[q(z)\] are polynomials, at least
one of which has degree \[d\]. Consider the one-parameter family of
rational maps
$$	f_t(z)\;= \;{(1-t)p(z)\,+\,tz^d\over (1-t)q(z)\,+\,t\;\;\;}	$$
with \[f_0(z)=f(z)\] and \[f_1(z)=z^d\]. There will be some finite set
of exceptional values of \[t\] for which the numerator and denominator
of this fraction have a common divisor. Algebraically this condition
is expressed by setting the ``resultant'' of the numerator and denominator
equal to zero.  Geometrically, it means that
a zero and pole of \[f_t\] crash together. If we exclude this finite set of
bad values of \[t\], then clearly \[f_t(z)\] is holomorphic as a function
of two variables, and each \[f_t\] has degree \[d\].

Suppose that \[f=f_0\] had \[4d-3\] distinct neutral cycles
with multipliers \[\lambda_j\ne 1\].
By the Implicit Function Theorem, we could follow each
of these cycles under a small deformation of \[f_0\]. Thus, for small
values of \[|t|\], the map \[f_t\] would have corresponding cycles
with multipliers \[\lambda_j(t)\] which depend holomorphically on \[t\],
with \[|\lambda_j(0)|=1\].

{\QP{\bf 10.7. Sub-Lemma.} {\it None of these functions \[t\mapsto
\lambda_j(t)\] can be constant throughout a neighborhood of \[t=0\].}\medskip}

{\bf Proof.} Suppose that for some \[j\] the function \[t\mapsto
\lambda_j(t)\] were constant throughout a neighborhood of \[t=0\]. Then
we will show that it is possible to continue analytically along any path
from 0 to 1 in the \[t$-plane which misses the finitely many
exceptional values of \[t\].
To prove this, we must check that the set of \[t\] for which we can continue
is both open and closed. But it is closed since any limit point of periodic
points with fixed multiplier \[\lambda_j\ne 1\] is itself a periodic point
with this same multiplier; and it is open since any such cycle varies
smoothly with \[t\], throughout some open neighborhood
in the \[t$-plane, by the Implicit Function Theorem.
Now continuing analytically to \[t=1\], we see that the map \[z\mapsto z^d\]
must also have  a cycle with mutiplier equal to \[\lambda_j\], with \[|
\lambda_j|=1\]. But this is known to be false, which proves 10.7.\QED

The proof of 10.6 continues as follows. We can express each of our \[4d-3\]
multipliers as a locally convergent power series
$$	\lambda_j(t)/\lambda_j(0)\;=\;1+a_j t^{n(j)}+({\rm higher\;terms})\,,
	 $$
where \[a_j\ne 0\] and \[n(j)\ge 1\].
Hence \[\log|\lambda_j(t)|\] is equal to the real part of\break
\[\;a_j t^{n(j)}+
{\rm (higher\;terms)}\,\]. If we ignore the higher order terms, this means
that we can divide the \[t$-plane into \[n(j)\] sectors for which
\[|\lambda_j(t)|>1\] and \[n(j)\] congruent sectors for which
\[|\lambda_j(t)|<1\]. Taking account of the higher order terms, we have
the following statement. Note that \[{\rm sgn}(
\log|\lambda|)\] is equal to \[+1\] or \[-1\] according as \[|\lambda|>1\] or
\[|\lambda|<1\].

{\QP{\bf Assertion.} \it The step function
$$	\theta\;\mapsto\;\sigma_j(\theta)\;\;=\;\;\lim_{r\to 0}
	\;{\rm sgn}\,\log|\lambda_j(r e^{i\theta})|	$$
is well defined with value \[\pm 1\]
except at finitely many jump discontinuities; and has average equal to
zero.\smallskip}

\noindent Therefore the sum \[\;\sigma_1(\theta)+\cdots+\sigma_{4d-3}
(\theta)\;\] is also a well defined step function with average zero.
Since this sum takes odd
values almost everywhere, we can choose some \[\theta\] for which
\[\sigma_1(\theta)+\cdots+\sigma_{4d-3}(\theta)\le -1\].
If we choose \[r\] sufficiently small and set \[t=re^{i\theta}\],
this means that \[f_t\] has at least \[2d-1\] distinct cycles
with multiplier satisfying \[|\lambda_j|<1\].
But this is impossible by 10.4, which proves 10.6 and 10.1.\QED

\vfil\eject % \bigskip\bigskip
\footline={\hss\tenrm 11-\folio\hss}\pageno=1

\centerline{\bf \S11. Repelling Cycles are Dense in \[J\].}\bigskip

We saw in \S4.3 that every repelling cycle is contained in the
Julia set. The following much sharper statement was proved
by Fatou and by Julia. (Compare 11.9.)

{\QP{\bf 11.1. Theorem.} {\it The Julia set for any rational map
of degree \[\ge 2\] is equal
to the closure of its set of repelling periodic points.}\medskip}

\noindent Since the proofs by Julia and by Fatou are interesting and
different, we will give both.\medskip

{\bf Proof following Julia.} We will use
the Holomorphic Fixed Point Formula of \S9.
Recall from 9.4 that every rational map \[f\] of degree two or more has either
a repelling fixed point, or a fixed point with \[\lambda=1\]. In either
case, this fixed point belongs to the Julia set \[J(f)\].
(Compare \S4.3 and \S7.5.)\smallskip

Thus we can start with a fixed point \[z_0\] in the Julia set. Let \[U
\subset\hat\C\] be any open set, disjoint from \[z_0\],
which intersects \[J(f)\]. The next step
is to construct a special orbit \[\;\cdots\mapsto z_2\mapsto z_1
\mapsto z_0\;\] which passes through \[U\] and terminates at this fixed point
\[z_0\]. By definition, such an orbit is called {\bit homoclinic} if
the backwards limit \[\lim_{j\to\infty} \;z_j\] exists and is equal to
the terminal point \[z_0\]. To construct a homoclinic orbit, we will
appeal to Theorem 4.8 which says that there exists an integer \[r>0\] and
a point \[z_r\in J(f)\cap U\] so that the \[r$-th forward image \[f^{\circ
r}(z_r)\] is equal to \[z_0\]. Given any neighborhood \[N_0\] of \[z_0\], we
can repeat this argument and conclude that there exists an integer \[q>r\]
and a point \[z_q\in N_0\] so that \[f^{\circ(q-r)}(z_q)=z_r\]. (Figure 11.)
\vskip .5in
\midinsert
\insertRaster 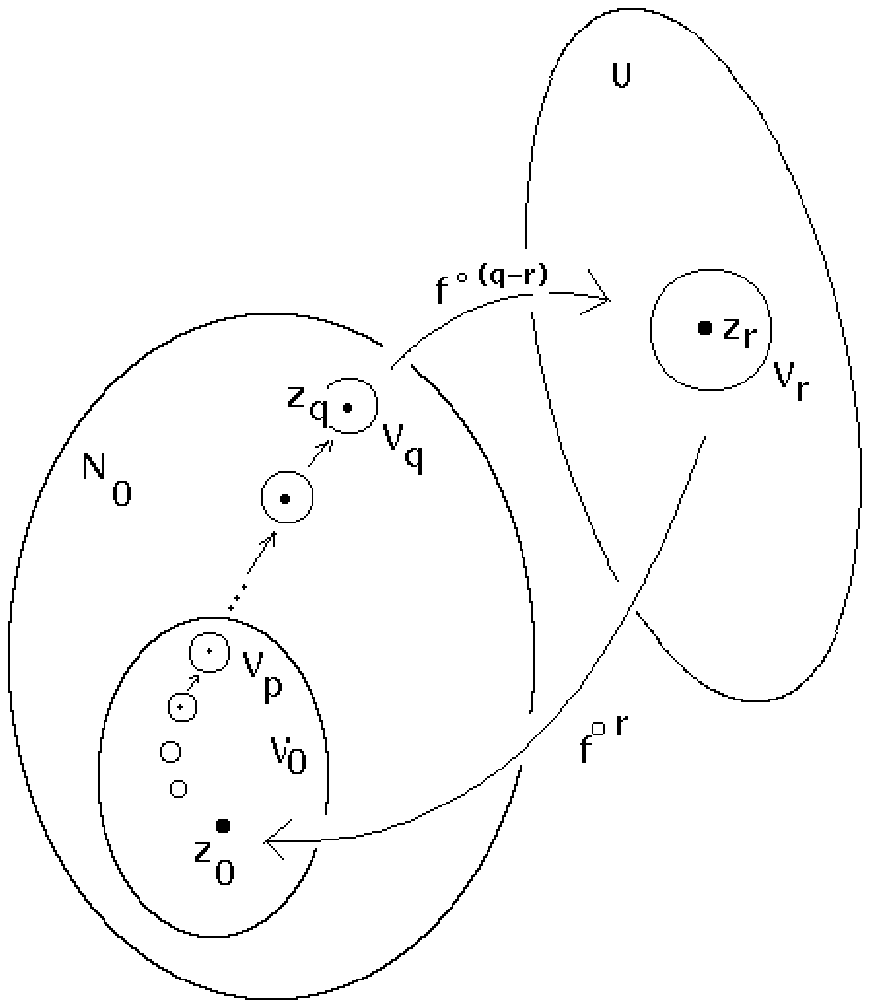 pixels 480 by 640 scaled 500
\vskip -.3in
\centerline{Figure 11. A homoclinic orbit.}
\vfil
\endinsert

To be more explicit, in the case where \[z_0\] is a repelling fixed point
we choose \[N_0\] to be a linearizing neighborhood, as in the
Koenigs Theorem 6.1. In the parabolic case, we choose \[N_0\]
to be a flower neighborhood, as in 7.2. In either case, we choose \[N_0\]
small enough to be disjoint from \[z_r\]. It then follows that we can
inductively choose preimages \[\cdots\mapsto
z_j\mapsto z_{j-1}\mapsto\cdots\mapsto
z_q\], all inside of the neighborhood \[N_0\]. {\it These preimages \[z_j\]
will automatically converge to \[z_0\] as \[j\to\infty\].} If \[z_0\] is
repelling, this is clear. In the parabolic case, \[z_q\]
cannot belong to an attracting petal; hence it must belong to a repelling
petal, and again this statement is clear.\smallskip

First suppose that none of the points \[\;\cdots\mapsto z_j\mapsto\cdots
\mapsto z_0\;\] in this homoclinic orbit are critical points of \[f\].
Then a sufficiently small disk neighborhood \[V_q\] of\break
\[z_q\in N_0\] will map
diffeomorphically under \[f^{\circ q}\] onto a neighborhood \[V_0\] of
\[z_0\]. Pulling this neighborhood \[V_q\] back under iterates of \[f^{-1}\],
we obtain neighborhoods \[z_j\in V_j\] for all \[j\], shrinking down towards
the limit point \[z_0\] as \[j\to\infty\]. In particular, if we choose
\[p\] sufficiently large, then \[\bar V_p\subset V_0\]. Now \[f^{-p}\] maps
the simply connected open set \[V_0\] holomorphically into a compact subset
of itself. Hence it contracts the Poincar\'e metric on \[V_0\] by a
factor \[c<1\], and therefore must have a attractive fixed point \[z'\]
within \[V_p\]. Evidently this point \[z'\in V_p\] is a repelling periodic
point of period \[p\] under the map \[f\]. Since the orbit of \[z'\]
under \[f\] intersects the
required open set \[U\], the conclusion follows.
\bigskip

If our homoclinic orbit contains critical points, then this argument must
be modified very slightly as follows. We can still choose simply
connected neighborhoods \[V_j\] of the \[z_j\] so that \[\bar V_p\subset
V_0\] for some large \[p\], and so that \[f\] maps each \[V_j\] onto
\[V_{j-1}\]. However, some finite number of these mappings will be branched.
Choose a slit \[S\] in \[V_0\] from the boundary to the midpoint \[z_0\]
so as to be disjoint from
\[\bar V_p\], and choose some sector in \[V_p\] which maps isomorphically
onto \[V_0-S\] under \[f^{\circ p}\]. The proof now proceeds just as
before.\QED\medskip
% \advancepageno

{\bf Proof of 11.1 following Fatou.} In this case, the main idea is an easy
application of Montel's Theorem (\S2.5). However, we must use Theorem 10.1 to
finish the argument.\smallskip

To begin the proof, recall from 4.9 that the Julia set \[J(f)\] has
no isolated points. Hence we can exclude finitely many points of \[J(f)\]
without affecting the argument. Let \[z_0\] be any point of \[J(f)\]
which is not a fixed point, and not a critical value. In other words, we
assume that there are \[d\] preimages \[z_1\,,\,\ldots\,,\,z_d\], which
are distinct from each other and from \[z_0\], where \[d\ge 2\] is the degree.
By the Inverse Function Theorem, we can find \[d\] holomorphic functions
\[z\mapsto \varphi_j(z)\] which are defined throughout some neighborhood \[N\] of
\[z_0\], and which satisfy \[f(\varphi_j(z))=z\], with \[\varphi_j(z_0)=z_j\].
We claim that for some \[n>0\] and for some \[z\in N\] the function
\[f^{\circ n}(z)\] must take one of the three values \[z\,,\,
\varphi_1(z)\] or \[\varphi_2(z)\]. For otherwise the family of holomorphic
functions
$$	g_n(z)\;=\;	{\bigl(f^{\circ n}(z)-\varphi_1(z)\bigr)\,\bigl(
	z-\varphi_2(z)\bigr)
	\over \bigl(f^{\circ n}(z) -\varphi_2(z)\bigr)\,\bigl
	(z-\varphi_1(z)\bigr)}	$$
on \[N\] would avoid the three values \[0\,,\,1\] and \[\infty\], and hence
be a normal family. It would then follow easily that \[\{f^{\circ n}|N\}\]
was also a normal family, contradicting the hypothesis that \[N\]
intersects the Julia set. Thus we can find \[z\in N\] so as to satisfy either
\[f^{\circ n}(z)=z\] or \[f^{\circ n}(z)=\varphi_j(z)\]. Clearly it
follows that \[z\] is a periodic point of period \[n\] or \[n+1\]
respectively.\smallskip

This shows that every point in \[J(f)\] can be approximated arbitrarily
closely by periodic points. Since all but finitely many of these periodic
points must repel, this completes the proof.\QED
\smallskip

There are a number of interesting corollaries.

{\QP{\bf 11.2. Corollary.} {\it If \[U\] is an open set which intersects the
Julia set \[J\] of \[f\], then for \[n\] sufficiently large
the image \[\;f^{\circ n}(U\cap J)\;\] is equal
to the entire Julia\break set \[J\].}\medskip}

{\bf Proof.} We know that \[U\] contains a repelling periodic point \[z_0\] of
period say \[p\]. Thus \[z_0\] is fixed by the iterate \[g=f^{\circ p}\].
Choose a
small neighborhood \[V\subset U\] of \[z_0\] with the property that
\[V\subset g(V)\]. Then clearly \[V\subset g(V)\subset g^{\circ 2}
(V)\subset\cdots\,\]. But it follows from 4.6 or 4.8 that the union of
the open sets \[g^{\circ n}(V)\] contains the entire Julia set \[J=
J(f)=J(g)\].
Since \[J\] is compact, this implies that \[J\subset g^{\circ n}(V)\subset
g^{\circ n}(U)\] for \[n\] sufficiently large, and the corresponding
statement for \[f\] follows.\QED

{\QP{\bf 11.3. Corollary.} {\it If a Julia set \[J\] is not connected,
then it has uncountably many distinct connected components.}\medskip}

{\bf Proof.} Suppose that \[J\] is the union \[J_0\cup J_1\] of two disjoint
non-vacuous compact subsets. After replacing \[f\] by some iterate
\[g=f^{\circ n}\], we may assume by 11.2 that\break
\[g(J_0)=J\] and \[g(J_1)=J\].
Now to each point \[z\in J\] we can assign an infinite sequence of symbols
$$	\epsilon_0(z)\,,\,\epsilon_1(z)\,,\,\epsilon_2(z)\,,\,\ldots\;\;\in\;
	\{0,1\}		$$
by setting
\[\epsilon_k(z)\] equal to zero or one according as \[g^{\circ k}(z)\]
belongs to \[J_0\] or \[J_1\]. It is not difficult to check that
points with different symbol sequences must belong to different connected
components of \[J\], and that all possible symbol sequences actually occur.
\QED

{\bf Remark.} Evidently only countably many of the connected components of
\[J\] can have positive Lebesgue measure. I don't know whether
all but countably many of the components of \[J\] must be single points.
\smallskip

Siegel disks and Cremer points do not seem directly related to critical
points, yet there is a connection as follows. By a {\bit critical orbit}
we will mean the forward orbit of some critical point.

{\QP\bf 11.4. Corollary. {\it Every boundary point of a maximal
Siegel disk \[U\]
belongs to the closure of some critical orbit.}\medskip}

{\bf Proof.} Otherwise, we could construct a small disk \[V\] around
the given point \[z_0\in\partial U\] so that the forward orbits of all
critical points avoid \[V\]. This would mean that
every branch of the \[n$-fold iterated inverse function \[f^{-n}\]
could be defined as a single valued holomorphic function \[f^{-n}:V\to
\hat\C\]. Let us choose that particular branch which carries the intersection
\[U\cap V\] into \[U\] by a
``rotation'' of \[U\]. Since the rotation number is irrational, we can choose
some subsequence of these inverse maps which converges to the identity
map on \[U\cap V\]. This is evidently a normal family, since it avoids
the central
part of \[U\]. Hence there is a sub-sub-sequence \[\{f^{-n(i)}\}\] which
converges on all of \[V\], necessarily to the identity map of \[V\].
It follows easily that the corresponding
sequence of forward iterates \[f^{\circ n(i)}\] also converges to
the identity on \[V\]. But this contradicts 11.2.\QED

{\QP\bf 11.5. Corollary. {\it Every Cremer point
is a non-isolated point in the closure of some critical orbit.}\medskip}

{\bf Proof.} (This proof applies also to parabolic cycles. Compare \S10.3.) 
Otherwise we could choose a small disk \[V\]
around the given point \[z_0\] so that no critical orbit intersects the
punctured disk \[V\ssm\{z_0\}\]. Replacing \[f\] by some iterate if necessary,
we may assume that \[z_0\] is fixed by \[f\].
Arguing as above, there exists a unique
holomorphic branch  \[f^{-n}:V\to\hat\C\] of the \[n$-fold iterated inverse
function which fixes the point \[z_0\]. These inverse maps form a normal
family since, for example, they avoid any periodic orbit which is disjoint
from \[V\]. Thus we can choose a subsequence \[f^{-n(i)}\]
converging locally uniformly to some holomorphic map \[h:V\to\hat\C\].
By the Inverse function Theorem, since \[|h'(z_0)|=1\], this \[h\] must map some
small neighborhood of \[z_0\] isomorphically onto a neighborhood \[V'\] of \[z_0\].
It follows that the corresponding forward maps \[f^{\circ n(i)}\] converge
on \[V'\] to the inverse map \[h^{-1}:V'\to V\]. Again,
this contradicts 11.2.\QED
\smallskip

By definition, the rational map \[f\] is called {\bit post critically finite~}
(or is called a {\bit Thurston map\/}) if
it has the property that every critical orbit is finite, or in other words
is either periodic or eventually periodic. According to Thurston, such a
map can be uniquely specified by a finite topological description.
(Compare Douady \& Hubbard [DH1].)

{\QP\bf 11.6. Corollary. {\it If \[f\] is post critically finite,
then every periodic
orbit of \[f\] is either repelling or superattracting.}\smallskip}

{\QP\bf 11.7. Corollary. {\it More generally, suppose that \[f\] has the
property that every critical orbit either is finite, or converges to an
attracting periodic orbit. Then every periodic orbit of \[f\] is either
attracting or repelling; there are no parabolic cycles, Cremer cycles, 
or Siegel cycles.}\medskip}

\noindent (Compare 13.5, 14.4)
The proofs are immediate, and will be left to the reader.\QED

By definition, \[z_0\] belongs to the Julia set if {\it some} sequence of
iterates of \[f\] has no subsequence which converges throughout a neighborhood
of \[z_0\]. However, a priori, there could be other sequences which do
converge.

{\QP{\bf 11.8. Corollary.} \it If \[z_0\] belongs to the Julia set, then
no sequence of iterates of \[f\] can converge uniformly
throughout a neighborhood of \[z_0\].\medskip}

For if \[\{f^{\circ n(i)}\}\] converges uniformly,
with \[n(i)\to\infty\], then according to Weierstrass (\S1.4) the
sequence of derivatives \[df^{\circ n(i)}(z)/dz\] must also converge. But
if \[z_1\] is a repelling periodic point, then evidently the sequence of
derivatives at \[z_1\] diverges to infinity.\QED
\medskip

{\bf 11.9. Concluding Remark.} The statement that the Julia set is equal to the
closure of the set of repelling periodic points is actually true for an
arbitrary holomorphic map of an arbitrary Riemann surface, providing that we
exclude just one exceptional case. For transcendental functions this
was proved by Baker [Ba1], and for maps of a torus or a Hyperbolic surface
it follows easily from \S4. The unique exceptional case occurs for degree one
maps of \[\hat\C\]
which have just one parabolic fixed point --- for example the map \[f(z)=z+1\]
with \[J(f)=\{\infty\}\].

On the other hand, for a holomorphic map of a complex 1-dimensional manifold
with two or more connected components (for example where both components
map into one), the Julia set can clearly be much bigger than the repelling
orbit closure.

\vfill
\bigskip

%\rightline{jwm, version of 9-3-91}

\end

%% file: cd12-14.tex
\input header
\input ras
\def\rot{{\rm rot}}
\def\Rot{{\rm Rot}}

\footline={\hss\tenrm 12-\folio\hss}\pageno=1

\centerline{\bf STRUCTURE OF THE FATOU SET}\bigskip

\centerline{\bf \S12. Herman Rings.}\medskip

The next two sections will be surveys only, with no proofs for several
major statements. This section
will describe a close relative of the Siegel disk.\smallskip

{\bf Definition.} A component \[U\] of the Fatou set \[\hat\C\ssm J(f)\] is
called a
{\bit Herman ring} if \[U\] is conformally isomorphic to some annulus
\[{\cal A}_r=\{z:1<|z|<r\}\], and if \[f\] (or some itterate of \[f\])
corresponds to an irrational rotation of this annulus.

{\bf Remark.} These objects are sometimes called
``Arnold-Herman rings'', since the existence of some examples would follow 
easily from Arnold's work in 1965.
Siegel disks and Herman rings are often collectively called
``{\bit rotation domains}''.

There are two known methods for constructing Herman rings. The original
method, due to Herman, is based on a careful analysis of real analytic
diffeomorphisms of the circle, as first studied by Arnold. An alternative
method, due to Shishikura, uses quasiconformal
surgery to cut and paste together two Siegel disks in order to fabricate
such a ring.

The original method can be outlined as follows. (Compare [S1], [He1].)
First a
number of definitions. If \[f:\R/\Z\to \R/\Z\] is an orientation
preserving homeomorphism, then we can lift to a homeomorphism \[F:\R\to\R\]
which satisfies the identity \[F(t+1)=F(t)+1\], and
is uniquely defined up to addition of an integer constant.\smallskip

{\bf Definition.} Following Poincar\'e, the {\bit rotation number}
of the lifted map \[F\] is defined to be the real number
$$	\Rot(F)\;=\;\lim_{n\to\infty}\; {F^{\circ n}(t_0)\over n}	$$
for any constant \[t_0\],
while the rotation number \[\;\rot(f)\in\R/\Z\;\] of the circle map \[f\]
is the residue class of \[\Rot(F)\] modulo 1.\smallskip

It is well known that this construction is well defined, and invariant
under orientation preserving topological conjugacy, and that it has
the following properties. (Compare Coddington and Levinson, or de Melo.)

{\QP{\bf 12.1. Lemma.} {\it The homeomorphism \[f\]
has a periodic point with period \[q\] if and only if its rotation number
is rational with denominator \[q\].}\smallskip}

{\QP{\bf 12.2. Denjoy's Theorem.} {\it If \[f\] is smooth of class \[C^2\],
and if the rotation number \[\rho=\rot(f)\]
is irrational, then \[f\] is topologically conjugate to the rotation
\[t\;\mapsto\; t+\rho\;\;(\mod 1)\].}\smallskip}

{\QP{\bf 12.3. Lemma.} {\it Consider a one-parameter family of lifted
maps of the form
$$	F_\alpha(t)\;=\;F_0(t)+\alpha\,.	$$
Then the rotation number \[\Rot(F_\alpha)\] increases
continuously and monotonically with \[\alpha\], increasing by \[+1\] as
\[\alpha\] increases by \[+1\,.\;\;(\]However, this dependence is not strictly
monotone. Rather, there is an interval of constancy corresponding to each
rational value of \[\Rot(F_\alpha)\,.\,)\]}\medskip}
\eject

In the real analytic case, Denjoy's Theorem has an analog which
can be stated as follows. Recall from \S8 that a real number \[\xi\]
is said to be {\bit Diophantine\/} if there exist a (large) number \[n\]
and a (small) number \[\epsilon\] so that the distance of \[\xi\] from
every rational number \[p/q\] satisfies \[|\xi-p/q|>
\epsilon/q^n\]. The following
was proved in a local version by Arnold, and
sharpened first by Herman and then by Yoccoz.

{\QP{\bf 12.4. Herman-Yoccoz Theorem.} \it If \[f\] is a real analytic
diffeomorphism of \[\R/\Z\] and if the rotation number \[\rho\] is
Diophantine, then \[f\] is real analytically conjugate
to the rotation \[\;t\;\mapsto\; t+\rho\;\;(\mod 1)\].\medskip}

I will not attempt to give a proof. (In the \[C^\infty\] case, Herman
and Yoccoz prove a\break corresponding
if and only if statement: {\it Every \[C^\infty\] diffeomorphism
with rotation number \[\rho\] is
\[C^\infty$-conjugate to a rotation if and only if \[\rho\] is
Diophantine.}) % However, as far as I know,
%no such sharp result is known in the real analytic case.
\smallskip

Next we will need the concept of a {\bit Blaschke product}. (Compare Problems
1-2 and 5-1.) Given any constant \[a\in\hat\C\]
with \[|a|\ne 1\], it is not difficult to show that there is one and only
one fractional linear transformation \[z\mapsto\beta_a(z)\] which maps the
unit circle \[\partial D\] onto itself fixing the base point \[z=1\],
and which maps \[a\] to \[\beta_a(a)=0\]. For example \[\beta_0(z)=z\],
\[\beta_\infty(z)=1/z\], and in general
$$	\beta_a(z)\;=\; {1-\bar a\over 1-a}\cdot{z-a\over 1-\bar a z}	$$
whenever \[a\ne\infty\].
If \[|a|<1\], then \[\beta_a\] preserves
orientation on the circle, and maps the unit disk into itself. On the other
hand, if \[|a|>1\], then \[\beta_a\] reverses orientation on \[\partial D\]
and maps \[D\] to its complement.

{\QP{\bf 12.5. Lemma.} \it A rational map of degree \[d\]
carries the unit circle into itself if and only if it
can be written as a {\bit ``Blaschke product''}
$$	f(z)\;=\;e^{2\pi it}\beta_{a_1}(z)\cdots\beta_{a_d}(z)	\eqno (*) $$
for some constants \[e^{2\pi it}\in\partial D\] and \[a_1\,,\,\ldots\,,\,a_d
\in\hat\C\ssm \partial D\].\medskip}

Here the \[a_i\] must satisfy the conditions that \[a_j\bar a_k\ne 1\]
for all \[j\] and \[k\]. For if \[a\bar b=1\], then a brief computation shows
that \[\beta_a(z)\beta_b(z)\equiv 1\]. Evidently the expression in 12.5 is
unique, since the constants \[\;e^{2\pi it}=f(1)\;\] and \[\;\{a_1\,,\,\ldots
\,,\,a_d\}\,=\,f^{-1}(0)\;\] are uniquely determined by \[f\]. The proof of
12.5 is not difficult: Given \[f\], one simply chooses any solution to the
equation \[f(a)=0\], then divides \[f(z)\] by \[\beta_a(z)\] to obtain
a rational map of lower degree, and continues inductively.\QED

Such a Blaschke product carries the unit disk into itself if
and only if all of the \[a_j\] satisfy \[|a_j|<1\]. (Compare Problems 5-1,
12-3.) However, we will rather be
interested in the mixed case, where some of the \[a_j\] are inside the unit
disk and some are outside.

{\QP\bf 12.6. Theorem. {\it For any odd degree \[d\ge 3\] we can choose a Blaschke
product \[f\] of degree \[d\] which carries the unit circle \[\partial D\] into
itself by a diffeomorphism with any desired rotation number \[\rho\]. If
this rotation number \[\rho\] is Diophantine, then \[f\] possesses a Herman
ring.}\medskip}\break

{\bf Proof Outline.} Let 
\[d=2n+1\], and choose the \[a_j\] so that \[n+1\] of them are close to
zero while
the remaining \[n\] are close to \[\infty\]. Then it is easy to check that
the Blaschke product \[z\mapsto\beta_{a_1}(z)\cdots\beta_{a_d}(z)\] is
\[C^1$-close to the identity map on the unit circle \[\partial D\].
In particular, it induces an orientation
preserving diffeomorphism of \[\partial D\]. Now multiplying by
\[e^{2\pi it}\] and using 12.3, we can adjust the rotation number to be any
desired constant. If this rotation number \[\rho\]
is Diophantine, then there
is a real analytic diffeomorphism \[h\] of \[\partial D\] which conjugates
\[f\] to the rotation \[z\mapsto e^{2\pi i\rho}z\]. Since \[h\] is real
analytic, it extends to a complex analytic diffeomorphism on some small
neighborhood of \[\partial D\], and the conclusion follows.\QED

As an example, Figure 12 shows the Julia set for the cubic rational map\break
\[f(z)\,=\,e^{2\pi it}z^2(z-4)/(1-4z)\;\] with zeros at \[0,\,0,\,4\], where
the constant
\[t=.6151732\cdots\] is adjusted so that the rotation number will be equal to
\[(\sqrt 5 -1)/2\].
There is a critical point near the center of this picture, with a
Herman ring to its left, surrounding the
superattractive basin about the origin in the left center.
This is the simplest kind
of example one can find, since Shishikura has shown that such a ring can
exist only if the degree \[d\] is at
least three, and since it is easy to check that a polynomial map cannot
have any Herman ring. (Problem 12-1 or \S17.1.)
\midinsert
\vfill
\insertRaster 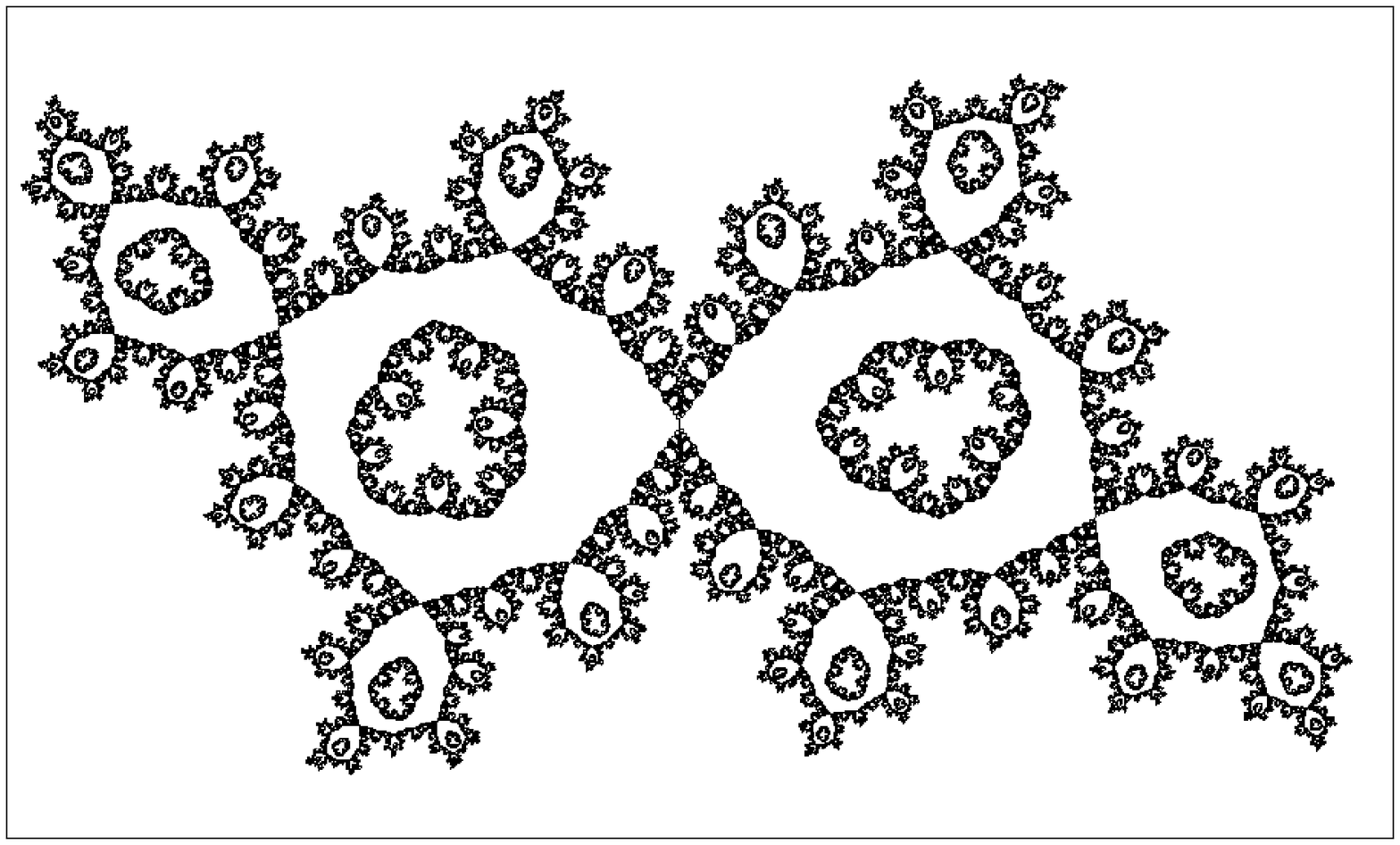 pixels 960 by 576 scaled 400
\bigskip\bigskip
\centerline{Figure 12. Julia set for a cubic rational map possessing
a Herman ring.}
\vfill
\endinsert

\eject

The rings constructed in this way are very special in that they are symmetric
about the unit circle, with \[f(1/\bar z)=1/\bar f(z)\]. However
Herman's original construction, based on work of Helson and Sarason,
was more flexible.
% Furthermore, the
%rotation number is necessarily Diophantine. 
Shishikura's more general
construction also avoids the need for symmetry. Furthermore
Shishikura's construction makes it clear that the possible rotation
numbers for Herman rings are exactly the same as the possible
rotation numbers for Siegel disks. In particular, any number satisfying
the Bryuno condition of \S8.4 can occur. The idea is roughly
that one starts with two rational maps having
Siegel disks with rotation numbers \[+\rho\] and \[-\rho\] respectively.
One cuts out a small concentric disk from each, and then glues the resulting
boundaries together. After making corresponding modifications at each of the
infinitely many iterated pre-images of each of the Siegel disks, Shishikura
applies
the Morrey-Ahlfors-Bers Measurable Riemann Mapping Theorem in order to
conjugate the resulting topological picture to an actual rational map.\medskip

Although Herman rings do not contain any critical
points, none-the-less they are closely associated with critical points.

{\QP\bf 12.7. Lemma. {\it If \[U\] is a
Herman ring, then every boundary point of \[U\]
belongs to the closure of the orbit of some critical point. The boundary
\[\partial U\]
has two connected components, each of which is an infinite set.}\medskip}

\noindent The proof is almost identical to the proof of 11.4.\QED % \medskip
\lln

{\bf Problem 12-1.} Using the maximum modulus principle, show that no
polynomial map can have a Herman ring.\smallskip

{\bf Problem 12-2.} For any Blaschke product \[f:\hat\C\to\hat\C\] show that
\[z\] is a critical point of \[f\] if and only if \[1/\bar z\] is a critical
point,
and show that \[z\] is a zero of \[f\] if and only if \[1/\bar z\] is a pole.
\smallskip

{\bf Problem 12-3.} A holomorphic map \[f:D\to D\] is said to be {\bit
proper} if the inverse image of any compact subset of \[D\] is compact.
Show that any proper holomorphic map from \[D\] onto itself can be expressed
uniquely as a Blaschke product $(*)$, with \[a_j\in D\].\smallskip

{\bf Problem 12-4.} Show that the rotation number \[\rot(f)\], as well as its
continued fraction expansion, can be deduced directly from the cyclic order
relations on a single orbit, without passing to the universal covering, as
follows. Let \[0=t_0\mapsto t_1\mapsto t_2\mapsto\cdots\] be the orbit of
zero, where we can choose representatives modulo 1 so that \[t_1\le t_j<
t_1+1\] for all \[j\]. Let us define \[t_j\] to be a ``closest return on the
left'' if \[t_j<0\], and if none of the numbers \[t_k\] with \[k<j\]
belong to the open interval \[(t_j\,,\,0)\]. Similarly, \[t_j\] is a ``closest
return on the right'' if the interval \[(0,t_j)\] does not contain any
\[t_k\] with \[k<j\]. If the sequence \[t_1\,,\,t_2\,,\,\ldots\] contains
first \[n_1\] closest returns on the left then \[n_2\] closest returns on
the right, and so on, show that \[\;	\rot(f)\;=\;1/(n_1+(1/n_2+1/(n_3
+\cdots)))\]. (Compare Appendix C.)\smallskip
\vfil\eject % \bigskip\bigskip

\footline={\hss\tenrm 13-\folio\hss}\pageno=1

\centerline{\bf \S13. The Sullivan Classification of Fatou Components.}
\medskip

The results in this section are due in part to Fatou and Julia, but
with very major contributions by Sullivan.\smallskip

By a {\bit Fatou component}
we will mean any connected component of the Fatou set\break \[\hat\C\ssm J(f)\].
Evidently \[f\] carries each Fatou component \[U\] onto some Fatou
component \[U'\] by a proper holomorphic map. First consider the special case
\[U=U'\].

{\QP\bf 13.1. Theorem. {\it If \[f\] maps the Fatou component \[U\] onto
itself, then there are just four possibilities, as follows. Either \[U\]
is the immediate attractive basin of an attracting fixed point, or of one
petal of a parabolic fixed point, or else \[U\] is a Siegel disk or
Herman ring.}\medskip}

Here we are lumping together the case of a superattracting fixed point, with
multiplier \[\lambda=0\], and the case of an ordinary attracting fixed
point, with \[\lambda\ne 0\]. Note that immediate attractive basins
always contain critical points by 10.2 and 10.3, while rotation domains
(that is Siegel disks and Herman rings)
evidently cannot contain critical points.\smallskip

Much of the proof of 13.1 has already been carried out in \S4. In fact,
according to 4.3 and 4.4, a priori there are just four possibilities.
Either:

(a) \[U\] contains an attractive fixed point;

(b) all orbits in \[U\] converge to a boundary fixed point;

(c) \[f\] is an automorphism of finite order; or

(d) \[f\] is conjugate to an irrational rotation of a disk, punctured
disk, or annulus.

\noindent In Case (a) we are done. Case (c) cannot occur, since our
standing hypothesis that the degree is two or more guarantees that there
are only countably many periodic points. In Case (d) we cannot have a
punctured disk, since the puncture point would have to be a fixed point
belonging to the Fatou set, so that we would just have a Siegel disk
with its center point incorrectly removed. Thus, in order to prove 13.1,
we need only show that the boundary fixed point  in Case (b)
must be parabolic. But
this boundary fixed point certainly cannot be an attracting point or a
Siegel point, since it belongs
to the Julia set. Furthermore, it cannot be repelling, since it attracts
all orbits in \[U\]. The only other possibility, which we must exclude,
is that it might be a Cremer point.\smallskip

The proof will be based on the following statement, which is
due to Douady and Sullivan. (Compare Sullivan [S1], or Douady-Hubbard [DH2,
p. 70]. For a more classical alternative, see Lyubich [L1, p. 72].)
Let
$$	f(z)\;=\;\lambda z+a_2z^2+a_3z^3+\cdots	$$
be a map
which is defined and holomorphic in some neighborhood \[U\] of the origin,
and which has
a fixed point with multiplier \[\lambda\] at \[z=0\]. By a {\bit path
converging to the origin} in \[U\] we will mean a continuous map \[p:(0,\infty)
\to U\ssm\{0\}\], where \[(0,\infty)\] is the open interval consisting
of all positive real numbers, satisfying the condition that \[p(t)\]
tends to zero as \[t\to\infty\].
\eject

{\QP{\bf 13.2. Snail Lemma.} {\it Let \[\;p:(0,\infty)\to U\ssm\{0\}\;\]
be a path which converges to the origin, and suppose that \[f\] maps \[p\]
into itself so that \[f(p(t))=p(t+1)\]. Then either the origin is
either an attracting or a parabolic fixed point. More precisely, the
multiplier \[\lambda\] satisfies either \[|\lambda|<1\] or \[\lambda=1\].}
\medskip}

\noindent Replacing \[f\] by \[f^{-1}\], we obtain the following
completely equivalent formulation, which will be useful in \S18.

{\QP{\bf 13.3. Corollary.} {\it If \[\;p:(0,\infty)\to U\ssm\{0\}\;\]
is a path which
converges to the origin, with \[f(p(t))=p(t-1)\] for \[t>1\], then either
\[|\lambda|>1\] or \[\lambda=1\].}\medskip}

{\bf Proof of 13.2.}
By hypothesis, the orbit \[p(0)\mapsto p(1)\mapsto
p(2)\mapsto\cdots\;\] in \[\;U\ssm\{0\}\;\] converges towards the origin. Thus
the origin certainly cannot be a repelling fixed point: we must have
\[|\lambda|\le 1\]. Let us assume that \[|\lambda|=1\] with
\[\lambda\ne 1\], and show that this hypothesis leads to a contradiction.

The intuitive idea can be described as follows.
As the path \[t\mapsto p(t)\] winds closer and closer to the origin,
the behavior of the map \[f\] on \[p(t)\] is more and more dominated by the
linear term \[z\mapsto\lambda z\]. Thus \[p(t+1)\approx \lambda p(t)\]
for large \[t\],
and the image must resemble a very tight spiral as shown in Figure 13.
Draw a radial segment \[E\] joining two turns of this spiral, as shown.
Then the region \[V\] bounded by \[E\] together with a segment of the
spiral will be mapped strictly into itself by \[f\]. Therefore, by the
Schwarz Lemma, the fixed point of \[f\] at the point \[0\in V\] must be
strictly attracting; which contradicts the hypothesis that \[|\lambda|=1\].
\midinsert
\hbox to \hsize{
\hfil
\hbox{\RasterBox {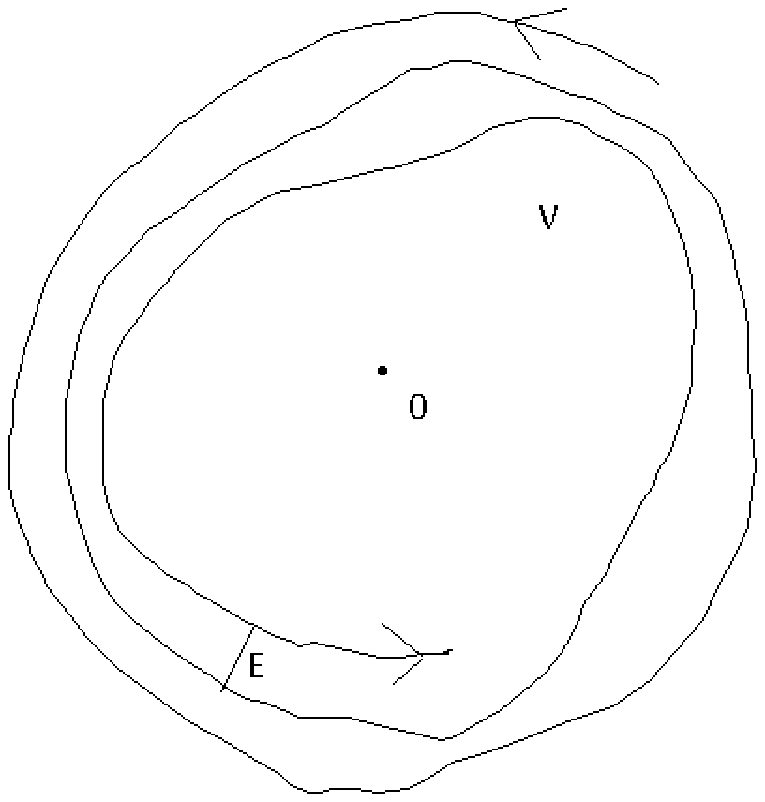} {640} {480} {350}}
\hskip -1in
\hbox{\RasterBox {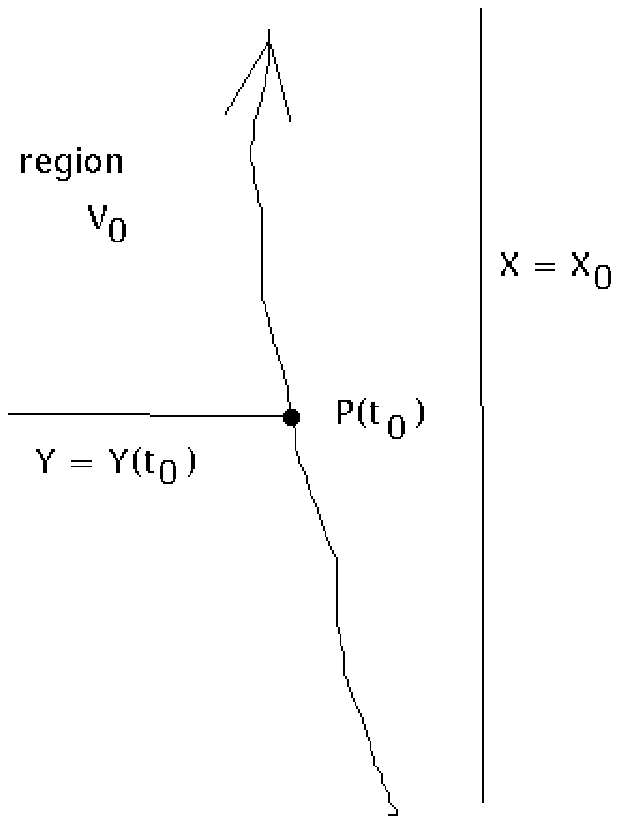} {640} {480} {350}}
\hfil
}
\smallskip
\centerline{Figure 13. Diagram in the \[z$-plane (left) and in the
\[Z=\log(z)\] plane (right).}
\endinsert % \vfil\eject

In order to fill in the details of this argument, it is convenient
to set \[Z=\log z\], and to lift \[p\] to a continuous path
\[Z=P(t)=X(t)+iY(t)\], with \[e^{P(t)}=p(t)\]. Then evidently as \[t\to\infty\] we have
\[X(t)\to -\infty\] and
$$	P(t+1)\;=\;P(t)+ic+o(1)	\;.$$
Here \[ic\] is a pure imaginary constant with
\[e^{ic}=\lambda\], and \[o(1)\] stands for a remainder term which tends to
zero. Similarly, we can lift \[f\] to a map \[Z\mapsto F(Z)\]
which is defined on some left half-plane  \[X<{\rm constant}\], and which
satisfies
$$	F(Z)\;=\;Z+ic+o(1)	$$
as the real part \[X={\cal R}(Z)\] tends to \[-\infty\], with the same
non-zero constant \[\,ic\]. Suppose, to fix our ideas, that \[c>0\]. Then,
for \[Z=X+iY\] within some half-plane \[X<X_0\], we may assume that \[F\]
is univalent, and that the
imaginary part \[Y\] increases by at least \[c/2\] under each iteration
of \[F\]. Choose \[t_0\] large enough so that \[X(t)<X_0\] for \[t\ge t_0\].
Then the path \[t\mapsto P(t)\] cuts the half-plane \[Y>Y(t_0)\] into two
or more connected components. Let \[V_0\] be that component whose intersection
with each horizontal line is unbounded to the left. Thus \[V_0\] is
contained in the quarter-plane \[X<X_0\,,\;Y>Y(t_0)\]. Evidently \[F\] maps
\[V_0\] diffeomorphically into itself. Furthermore, the imaginary part of
\[X+iY\] increases by at least \[c/2\] under each iteration, hence the real
part must decrease towards \[-\infty\] under iteration.

Let \[V=\exp(V_0)\cup\{0\}\] be the corresponding open set in the \[z$-plane.
Then it follows that \[f\] maps \[V\] into itself,
and that every orbit of \[f\] in \[V\] converges towards the origin.
Since the open set \[V\] is Hyperbolic, it follows by the
Schwarz Lemma that the origin is an attractive fixed point; which
contradicts our hypothesis that \[|\lambda|=1\].\QED

To prove 13.3, we simply note that the orbit
$$	\cdots\mapsto p(2)\mapsto p(1)\mapsto p(0)	$$
is repelled by the origin, so the multiplier \[\lambda\] cannot be zero.
Hence \[f^{-1}\] is defined and holomorphic near the origin. Applying
13.2 to \[f^{-1}\], the conclusion follows.\QED

{\bf Proof of 13.1.} Let \[U\] be a Fatou component which is mapped into
itself by \[f\] in such a way that all orbits converge to a boundary fixed
point \[w_0\]. Choose any base point \[z_0\] in \[U\], and choose any
path \[p:[0,1]\to U\] from \[p(0)=z_0\] to \[p(1)=f(z_0)\]. Extending for
all \[t\ge 0\] by setting \[p(t+1)=f(p(t)\], we obtain a path in \[U\]
which converges to the boundary point \[w_0\] as \[t\to\infty\]. Therefore,
according to 13.2, the fixed point \[w_0\] must be either parabolic or
attracting. But \[w_0\] belongs to the Julia set, and hence cannot be
attracting.\QED\smallskip
%  \advancepageno

Thus we have classified the Fatou components which are mapped onto themselves
by \[f\]. There is a completely analogous description of Fatou components
which cycle periodically under \[f\]. These are just the Fatou components
which are fixed by some iterate of \[f\]. Hence each such component is
either

(1) the immediate attractive basin for some attracting periodic  point,

(2) the immediate attractive basin for some petal of a parabolic periodic
point,

(3) one member of a cycle of Siegel disks, or

(4) one member of a cycle of Herman rings.

\noindent In Cases (3) and (4), the topological type of the domain \[U\]
is specified by this description.
In Cases (1) and (2) it can be described as follows.

{\QP{\bf 13.4. Lemma.} \it Every immediate attractive basin is either simply
connected or infinitely connected.\medskip}

{\bf Proof.} If \[U\] were a non-simply-connected region of finite
connectivity, then it would have a finite Euler characteristic \[\chi(U)
\le 0\]. By 10.2 or 10.3, f maps \[U\] onto itself by a branched covering
map with at least
one branch point, and hence with degree \[n\ge 2\]. But the precise
number of branch points, counted with multiplicity, is equal to
\[(n-1)\chi(U)\le 0\] by the Riemann-Hurwitz formula of \S5.1. This
contradiction completes the proof.\QED	

Sullivan showed that there can be at most a finite number of such periodic
Fatou components. {\it In fact, according to Shishikura, there can be at most 2d-2
distinct cycles of Fatou components.}\smallskip

In order to complete the picture, we need the fundamental theorem
that there are no ``wandering'' Fatou components. (Compare [S2], [C].)

{\QP{\bf 13.5. Theorem of Sullivan.} \it Every Fatou component \[U\] is
eventually periodic. That is, there necessarily exist integers \[n\ge 0\]
and \[p\ge 1\] so that the \[n$-th forward image \[f^{\circ n}(U)\] is
mapped onto itself by \[f^{\circ p}\].\medskip}

Thus every Fatou component is a preimage, under some iterate of \[f\],
of one of the four types described above.
The proof, by quasiconformal deformation, is beyond the scope of these notes.
Roughly speaking, if a wandering Fatou component were to exist, then one
could construct an infinite dimensional space of deformations, all of which
would have to be rational maps of the same degree. But the space of
rational maps of fixed degree is finite dimensional.\smallskip

Recall from \S11.6
that \[f\] is {\bit post critically finite\/} (or a {\bit Thurston map\/})
if every critical orbit is finite. Combining 13.1 and
13.4 with 11.6 and 12.7, we easily obtain the following. (Compare Problem 5-7.)

{\QP{\bf 13.6. Corollary.} \it If a post critically finite
rational map has no superattractive
periodic orbit, then its Julia set is the entire sphere \[\hat\C\].\medskip}

\noindent We will give a different proof for this statement in 14.6.\lln
% \eject
We can sharpen the defining property of the Fatou set as follows.\smallskip

{\bf Problem 13-1.} If \[V\] is a connected open subset of the Fatou
set \[\hat\C\ssm J\], show that the set of all limits of
successive iterates \[f^{\circ n}|V\] as \[n\to\infty\] is either
(1) a finite set of constant maps from \[V\] into an attracting or
parabolic periodic orbit, or (2) a one-parameter family of maps,
consisting of all compositions
\[R_\theta\circ f^{\circ k}|V\] where \[f^{\circ k}\] is some fixed
iterate with values in a rotation domain and \[R_\theta\] is the rotation
of this domain through angle \[\theta\].
\bigskip

\vfil\eject
\footline={\hss\tenrm 14-\folio\hss}\pageno=1

\centerline{\bf \S14. Sub-hyperbolic and hyperbolic Maps.}\bigskip

This section will discuss two important classes of rational maps. The
exposition will be based on Douady-Hubbard [DH2].\smallskip

First some standard definitions from the theory of smooth dynamical systems.
Let \[f:M\to M\] be a \[C^1$-smooth map from a smooth
Riemannian manifold to itself, and let\break \[X\subset M\] be a compact
\[f$-invariant
compact subset, \[f(X)\subset X\]. Let\break \[Df_x:T_xM\to T_{f(x)}M\],
or briefly \[Df\], be the {\bit first
derivative map} at \[x\], that is the induced linear map on the tangent
space at
\[x\], and let \[\|v\|\] be the Riemannian norm of a vector \[v\in T_xM\].
By definition, the map \[f\] is {\bit expanding}
on \[X\] if the length of any tangent vector at a point of \[X\]
expands exponentially under iteration of \[Df\], that is, if there
exist constants \[c>0\] and \[k>1\] so that
$$	\|Df^{\circ n}(v)\|\;\ge\; ck^n\|v\|	$$
for every \[x\in X\] and \[v\in T_xM\], and every \[n\ge 0\]. Since \[X\]
is compact, a completely
equivalent requirement is that there exists some fixed \[n\ge 1\] so that
$$	\|Df^{\circ n}(v)\|\;>\;\|v\|	$$
for all non-zero tangent vectors \[v\] over \[X\]. Similarly, \[f\] is {\bit
contracting} on \[X\] if there are constants
\[c>0\] and \[k<1\] so that \[\|Df^{\circ n}(v)\|\;\le\;ck^{n}\|v\|\,\].
It is not difficult to check that
these conditions do not depend on the particular choice of Riemannian metric.
\smallskip

In the higher dimensional case, a map is called {\bit hyperbolic} on \[X\]
if the tangent space of \[M\]
restricted to \[X\] splits as the Whitney sum of two sub-bundles, each
invariant under \[Df\], so that \[Df\] is expanding on one sub-bundle and
contracting on the other. In the
one-dimensional case, this concept simplifies as follows.

{\bf Definition:} Let \[f\]
be a holomorphic map from a Riemann surface to itself, and let \[X\] be a
compact \[f$-invariant subset. The map \[f\] is
{\bit hyperbolic} on \[X\] if \[f\] is expanding on \[X\] or contracting
on \[X\], or if \[X\] is the disjoint union
of a compact \[f$-invariant subset \[X^+\] on which \[f\] is expanding
and a compact \[f$-invariant subset \[X^-\] on which \[f\] is contracting.

{\QP{\bf 14.1. Theorem.} {\it For a rational map \[f:\hat\C\to\hat\C\] of
degree \[d\ge 2\], the following two conditions are equivalent:

$(1)$ \[f\] is expanding and hence hyperbolic on its Julia set \[J\].

$(2)$ The forward orbit of each critical point of \[f\] converges
towards some attracting periodic orbit.}\medskip}
%$(4)$ The {\bit postcritical set} \[P\], that is the union of the forward
%orbits of the critical points, has closure \[\bar P\] which is disjoint
%from the Julia set \[J\].

Map satisfying these conditions are briefly called {\bit hyperbolic}
rational maps.  As examples, Figures 1a, 1b, 1d and 5 show the Julia sets of
hyperbolic maps, but the other figures do not.
{\bf Caution:} This use of the word ``hyperbolic'' has
nothing to do with the concept of ``Hyperbolic Riemann surface'',
which we write with a capital {\bit H}. (Compare \S2.)\eject % \smallskip

The proof that \[(1)\Rightarrow (2)\] proceeds as follows. We suppose
given some smooth Riemannian metric on \[\hat\C\] so that the lengths of
tangent vectors to \[\hat\C\] at points of \[J\] increase exponentially under
iteration of \[Df\]. After replacing \[f\] by some fixed iterate \[f^
{\circ p}\] (or alternatively after carefully modifying the metric), we may
assume that \[\|Df(v)\|\ge k\|v\|\] for all tangent vectors at points
within some neighborhood \[V\] of
\[J\], where \[k>1\]. In particular, there can be no critical points within
this neighborhood \[V\].

{\QP{\bf 14.2. Lemma.} \it If there exists such an expanding metric in a
neighborhood of the\break Julia set, then every orbit outside of the Julia set
converges towards an attracting periodic orbit.\medskip}

{\bf Proof.} We first show that there exists a constant
\[\epsilon>0\] so that the neighborhood \[N_\epsilon\], consisting of all
points at Riemannian distance less than \[\epsilon\] from \[J\], has the
following property:

{\QP \it If we iterate the mapping \[f\] starting at any point which
belongs to the neighborhood \[N_\epsilon\] but not to \[J\], then the
distance from \[J\] will be multiplied by
\[k\] or more at each iteration until we leave the set \[N_\epsilon\].
Furthermore, no point outside of \[N_\epsilon\] maps into \[N_\epsilon\].
\smallskip}

\noindent
In fact we need only choose \[\epsilon\] small enough so that the neighborhood
\[N_\epsilon\] is contained in \[V\], but disjoint from the compact set
\[f(\hat\C\ssm V)\]. It then follows that \[f^{-1}(N_\epsilon)\subset V\]. In
fact we will prove the stronger statement that \[f^{-1}(N_\epsilon)\subset
N_{\epsilon/k}\subset N_\epsilon\]. For if \[f(z)\in N_\epsilon\], then
\[f(z)\] can be joined
to \[J\] within \[N_\epsilon\] by a geodesic of length equal to the
Riemannian distance \[\rho(f(z),J)<\epsilon\]. We can then lift this back to
a smooth curve which joins \[z\] to \[J\] within \[V\]. The length of
this pulled back curve is evidently \[\le\rho(f(z),J)/k<\epsilon/k\].

Since \[f^{-1}(N_\epsilon)\subset N_\epsilon\], it follows that the
complementary compact set \[L=\hat\C\ssm N_\epsilon\] satisfies \[f(L)
\subset L\].
We have shown that every orbit which starts outside of \[J\] must eventually
be absorbed by this invariant set \[L\].
Since \[L\] is a compact subset of the Fatou set
\[\hat\C\ssm J\], it must be covered by finitely many of the connected components
of \[\hat\C\ssm J\]. Let \[U_1\,,\,\cdots\,,\, U_p\] be those components
of \[\hat\C\ssm J\] which intersect \[L\]. Since \[f(L)\subset
L\], it follows that \[f\] carries each one of these \[U_i\] onto some
\[U_j\]. As we follows such an orbit, we must eventually reach a component
\[U_j\] which is mapped onto itself by some iterate \[f^{\circ q}\].

Now let us examine the four possibilities for such a map \[f^{\circ q}:U_j\to
U_j\], as listed in Theorem 3.1.
This map cannot be a rotation of \[U_j\], or an automorphism
of finite order, since all orbits in \[U_j\]
are eventually absorbed by the compact subset \[L\cap U_j\]. Similarly,
orbits cannot diverge to infinity with
respect to the \Poincare metric on \[U_j\]. Thus the only
possibility is that all orbits of \[f^{\circ q}\] within \[U_j\]
converge to an attracting fixed point of \[f^{\circ q}\]. This proves
the Lemma, and hence completes the proof that
\[(1)\Rightarrow (2)\].\smallskip

{\bf Proof that \[(2)\Rightarrow (1)\].} Let \[P\] be the {\bit post-critical
set}, that is the union of the forward orbits of the critical points.
Then Condition (2) clearly implies that the closure \[\bar P\] is disjoint
from the Julia set. Conversely, if \[\bar P\cap J=\emptyset\], then we
will construct an expanding metric near the Julia set. Let
\[U\] be the open set \[\hat\C\ssm \bar P\]. Since
\[f(\bar P)\subset\bar P\], we have \[f^{-1}(U)\subset U\].
There are no critical points in \[U\], so it follows that \[f^{-1}\]
lifts to a single valued analytic function from the universal covering
surface \[\tilde U\] into itself. But \[U\] contains the Julia set, which
contains at least one repelling fixed point. Therefore this inverse map
on \[\tilde U\] can be chosen so as to contain an attracting fixed point.
If \[\tilde U\] is isomorphic to the unit disk, then this implies that
\[f^{-1}\] strictly contracts the \Poincare metric on \[\tilde U\], so that
\[f\] must strictly increase the \Poincare metric on \[U\]. This proves that
\[f\] is expanding on the compact subset \[J\subset U\].
On the other hand, if \[\tilde U\] is not isomorphic to \[D\], then the
complement \[\bar P=\hat\C\ssm U\] can contain at most two points. It is then
easy to check that \[f\]
must be conjugate to a map of the form \[z\mapsto z^{\pm d}\].
In this case, clearly \[f\] is expanding with respect to the Euclidian metric
on \[U\]. Thus \[(2)\Rightarrow (1)\], which completes the proof.\QED

One invariant set of particular interest is the ``non-wandering set''
of \[f\]. By definition, a point is {\bit wandering} if it has a
neighborhood \[U\] so that the forward images \[f(U)\,,\, f^{\circ 2}(U)\,,\,
\ldots\] are all disjoint from \[U\]. Otherwise, it is {\bit non-wandering}.
The closed set consisting of all non-wandering points is called
the {\bit non-wandering set} \[\Omega=\Omega(f)\]. The map \[f\] is said to
satisfy Smale's {\bit Axiom A} if \[f\] is
hyperbolic on \[\Omega\], and if \[\Omega\] is precisely equal to the
closure of the set of periodic points.

{\QP{\bf 14.3. Corollary.} \it A rational map satisfies Axiom A
if and only if it is hyperbolic.\medskip}

The proof is straightforward,
since we know by \S11 that the repelling periodic points are dense in the
Julia set. Details will be left to the reader. (Compare Problem 14-1.)
\QED

{\bf Remark.} These hyperbolic maps have other extremely important properties.
If \[f\] is hyperbolic, then every nearby map is also hyperbolic. Furthermore,
according to Ma\~n\'e, Sad, Sullivan, and also Lyubich, the Julia set \[J(f)\]
varies continuously under a deformation of \[f\] through hyperbolic maps.
In the non-hyperbolic case, a small change in \[f\] may well lead to a drastic
alteration of \[J(f)\]. It is generally conjectured that every rational map
can be approximated arbitrarily closely by a hyperbolic map.
\medskip

Douady and Hubbard, using ideas of Thurston, also consider a wider class
of mapping which they call {\bit sub-hyperbolic}. The only change in the
definition is that the Riemannian metric is now allowed to have a finite
number of relatively mild singularities,
like that of the metric \[|d\root n\of z |\] at the origin. Such singularities
have the property that the Riemannian distance between points still tends to
zero as the points approach each other.\smallskip

{\bf Definition.} A conformal metric on a Riemann surface, with the
expression  \[\gamma(z)|dz|\] in terms of a local uniformizing parameter
\[z\], will be called an {\bit orbifold metric} if the function \[\gamma(z)\]
is smooth and non-zero except at a locally finite collection of points
\[a_1\,,\,a_2\,,\,\ldots\;\]
where it blows up in the following special way. There should be
integers \[\nu_j\ge 2\] called the {\bit ramification indices} at the points
\[a_j\] with the following property. If we take a local branched covering
by setting \[z(w)=a_j+w^{\nu_j}\], then the induced metric \[\gamma(z(w))
|(dz/dw)dw|\] on the \[w$-plane should be smooth and non-singular throughout
some neighborhood of the origin.

The rational map \[f\] is {\bit
sub-hyperbolic} if it is expanding with respect to some orbifold metric
on a neighborhood of its Julia set.
\eject

{\QP{\bf 14.4. Theorem.} \it A rational map is sub-hyperbolic if and only if:

$(\alpha)$ every critical orbit in the Julia set is eventually periodic,
and

$(\beta)$ every critical orbit outside of the Julia set converges to an
attracting periodic orbit.\medskip}

(Compare 11.7.) As examples, Figures 1-5 show sub-hyperbolic maps,
but 8-10, 12, 17 do not.

The proof in one direction is completely analogous to
the proof above. If \[f\] is expanding with respect to some orbifold metric,
defined near the Julia set, then, just as in 14.2 above, every orbit
outside the Julia
set \[J\] will be pushed away from \[J\], and can only converge to an
attractive periodic orbit. On the other hand, if \[\omega\] is a critical
point inside the Julia set, then every post-critical
point \[f^{\circ n}(\omega)\,,\; n>0\], must be one of the singular
points  \[a_j\] for our
orbifold metric, since the map \[f^{\circ n}\] is critical, and yet is required
to be expanding, at the point \[\omega\]. Thus the forward
orbit of \[\omega\] cannot have any limit points, and hence must be finite.
\smallskip

For the proof in the other direction, we must introduce the concept of the
``universal branched covering'' of an orbifold. Compare Appendix E.
For our purposes, an
{\bit orbifold} \[(S,\,\nu)\] will just mean a Riemann surface
\[S\], together with a locally finite collection of marked points \[a_j\],
each of which is assigned
a {\bit ramification index} \[\nu_j=\nu(a_j)\ge 2\] as above. For any
point \[z\] which is not one of the \[a_j\] we set \[\nu(z)=1\].\smallskip

{\bf Definition.} To every rational
map \[f\] which satisfies Conditions ($\alpha$) and ($\beta$) of Theorem
14.4, we assign
the {\bit canonical orbifold} \[(S,\,\nu)\] as follows. As underlying
Riemann surface \[S\] we take the Riemann sphere \[\hat\C\] with all
attracting periodic orbits removed. As ramified points \[a_j\] we take
all strictly {\bit post-critical} points, that is, all points which have
the form \[a_j=f^{\circ n}(\omega)\] for some \[n>0\], where \[\omega\] is
a critical point of \[f\]. Using ($\alpha$) and ($\beta$), we see easily that this
collection of points \[a_j\] is locally finite in \[S\] (although perhaps
not in \[\hat\C\]). In order to specify the corresponding ramification indices
\[\nu_j=\nu(a_j)\], we will need another definition.
If \[f(z_1)=z_2\], with local power series development
$$	f(z)=z_2+c(z-z_1)^n+({\rm higher\;terms})\,,	 $$
where \[c\ne 0\] and \[n\ge 1\],
then the integer \[n=n(f\,,\,z_1)\] is called the {\bit local degree} or the
{\bit branch index} of \[f\] at \[z_1\].
{\it Now choose the \[\nu(a_j)\ge 2\] to be the smallest integers which
satisfy the following:}

{\QP{\bf Condition ($\ast$).} For any \[z\in S\],
the ramification index \[\nu(f(z))\] at the image point must be a multiple
of the product \[n(f\,,\,z)\nu(z)\]. \medskip}

\noindent
In fact, to construct these integers \[\nu(a_j)\]
we simply consider all critical pre-images
$$	f^{\circ m}(\omega)=a_j\,,\quad f'(\omega)=0\,,	$$
and let \[\nu(a_j)\] be the least common multiple of the corresponding
branch indices\break \[n(f^{\circ m}\,,\,\omega)\]. Since we have removed all
attracting periodic orbits, it is not difficult to
check that this least common multiple is finite,
and that it provides a minimal solution to the required condition $(\ast)\].
\smallskip

As in Appendix E, we consider the {\bit universal covering surface}
$$	\tilde S_\nu\;\to\; (S\,,\,\nu)	$$
for this canonical orbifold, that is the unique regular branched
covering of \[S\] which has the given \[\nu\] as ramification function,
and which is simply connected. In order
to determine the geometry of \[\tilde S_\nu\], we introduce
the {\bit orbifold Euler characteristic}
$$	\chi(S\,,\,\nu)\;=\;\chi(S)+\sum({1\over\nu(a_j)}-1)\;=\;2-k+
	\sum({1\over\nu(a_j)}-1)\,,	$$
to be summed over all ramified points, where \[k\] is the number of
attracting periodic points. According to Lemma E.1 in the Appendix,
such a universal covering nearly always exists; and the few pathological
cases where it does not exist necessarily have Euler characteristic \[\chi
(S\,,\,\nu)>0\].
Furthermore, by E.4, the surface \[\tilde S_\nu\] is either Spherical,
Hyperbolic, or Euclidean according as \[\chi(S\,,\,\nu)\] is positive,
negative, or zero.\smallskip

Still assuming Conditions ($\alpha$) and ($\beta$) of 14.4, we will
prove the following.

{\QP{\bf 14.5. Lemma.} \it This canonical orbifold has Euler
characteristic
\[\chi(S\,,\,\nu)\le 0\].\break Hence its universal covering \[\tilde
S_\nu\]
has either a unique \Poincare ~metric or a\break
Euclidean metric which is unique
up to a scale change. In either case, there is an induced orbifold metric
on \[S\], and \[f\] is expanding with respect to this orbifold metric
throughout any compact subset of \[S\].\medskip}

Evidently this will complete the proof of Theorem 14.4.
In the Euclidean case, the proof yields much more. Let \[z\] be a uniformizing
parameter on the covering surface \[\tilde S_\nu\cong\C\].

{\QP{\bf 14.6. Theorem.} \it If \[\chi(S\,,\,\nu)=0\], then \[f\] induces
a linear isomorphism\break \[z\mapsto az+b\] from the Euclidean covering space
\[\tilde S_\nu\cong \C\] onto itself. In this case, the Julia set is
either a circle or line segment, or the entire Riemann sphere. The
coefficient of expansion \[|a|\] satisfies \[|a|=d\] in the first two
cases, and \[|a|=\sqrt d\] in the last case, where \[d\] is the degree
of \[f\].\medskip}

(Caution: The coefficient \[a\] itself is not uniquely determined; for
the lifting of \[f\] to the covering surface is determined only up to
composition with a deck transformation. The deck transformations
may well have fixed points, since we are dealing with a branched covering,
but necessarily
have the form \[z\mapsto \omega z+c\] where \[\omega\] is some root of unity.)
\smallskip

{\bf Proof of 14.5.} First note that Condition ($\ast$) above is exactly
what is needed in order to assert that \[f^{-1}\] lifts, locally at least,
to a holomorphic map from the universal covering \[\tilde S_\nu\] to
itself. But since \[\tilde S_\nu\] is simply connected, there is then no
obstruction to constructing a globally lifted inverse
$$	\tilde f^{-1}\,:\, \tilde S_\nu\;\to\;\tilde S_\nu\,.	$$
The proof now divides into three cases.\smallskip

{\bf Spherical Case.} Note that the composition
$$	\tilde S_\nu\;\buildrel \tilde f^{-1}\over\longrightarrow\;
	\tilde S_\nu\;\buildrel p\over\longrightarrow\;S\;\buildrel
	 f\over\longrightarrow\;S $$
always coincides with the projection map \[p:\tilde S_\nu\to S\].
If \[\tilde S_\nu\] were Spherical, then this composition
would have a well defined
degree, which would have to be at least \[d\] times the degree of \[p\]. This
is impossible, hence the Spherical case cannot occur.\smallskip

{\bf Hyperbolic Case.} If \[\tilde S_\nu\] is Hyperbolic, then it has a
uniquely defined Poincar\'e metric, and the proof proceeds just as in the
proof that \[(2)\Rightarrow (1)\] in 14.1 above.\smallskip

{\bf Euclidean Case.} Suppose that \[\tilde S_\nu\] is Euclidean, or
equivalently that \[\chi(S\,,\,\nu)=0\].
Construct a new orbifold \[(S'\,,\,\mu)\] as follows. Let \[S'\subset S\]
be the open set \[S\] with all immediate pre-images of attracting periodic
points removed, and define \[\mu=f^*(\nu)\] by the formula
$$	\mu(z)\;=\; \nu(f(z))/n(f\,,\,z)	$$
where \[n(f\,,\,z)\] is the branch index.
Note that \[\mu(z)\ge\nu(z)\] for all \[z\in S'\]. Evidently it follows
that
$$	\chi(S'\,,\,\mu)\;\le\; \chi(S\,,\,\nu)	$$
with equality only if \[S'=S\] and \[\mu=\nu\]. But by Lemma E.2, since the
map\break \[f:(S'\,,\,\mu)\to(S\,,\,\nu)\] is a ``$d$-fold covering'' of
orbifolds, we conclude that \[f\] induces an isomorphism \[\tilde S'_\mu
\to\tilde S_\nu\] of universal
covering surfaces, and also that the Riemann-Hurwitz formula takes the form
$$	\chi(S'\,,\,\mu)\;=\;\chi(S\,,\,\nu)d\,.	$$
Thus \[\chi(S'\,,\,\mu)\] is also zero. We conclude that
\[S'=S\] and \[\mu=\nu\], so that \[f\] must lift to a (necessarily linear)
isomorphism from the Euclidean covering space \[\tilde S_\nu\] to itself.
Further details will be left to the reader. This completes the proof
of 14.4 through 14.6.\QED

As a corollary, we obtain another proof of 13.6.
We continue to assume that every
critical orbit either converges to an attracting orbit or is eventually
periodic, according as it belongs to the Fatou set or the Julia set.

{\QP{\bf 14.7. Corollary.} \it If \[f\] is sub-hyperbolic with
no attracting periodic orbits,
so that \[S\] is the entire Riemann sphere, then \[f\] is expanding on the
entire sphere, and it follows that \[J(f)\] is the entire sphere.
\medskip}

\lln

{\bf Problem 14-1.} Using the results of \S13, show that the non-wandering
set for a rational map \[f\] is the (disjoint) union of its Julia set,
its rotation domains (if any), and
its set of attracting periodic points.\smallskip

{\bf Problem 14-2.} Show that the Julia set for the rational map \[z\mapsto
(1-2/z)^{2n}\] is the entire Riemann sphere. Show that the orbifold
metric for this example is Euclidean when \[n=1\], but is Hyperbolic
for \[n>1\].\smallskip

{\bf Problem 14-3.} For any sub-hyperbolic map whose
canonical orbifold metric is\break Euclidean,
show that every periodic orbit outside of the finite post-critical set has
multiplier \[\lambda\] satisfying \[|\lambda|=d^{p/\delta}\] where
\[d\] is the degree,
\[p\] is the period, and \[\delta\] is the dimension (1 or 2) of the Julia
set.\smallskip

{\bf Problem 14-4.} A rational map \[f\] is said to be {\bit expansive}
in a neighborhood of its Julia set if there exists \[\epsilon>0\] so that,
for any two points \[x\ne y\] whose orbits remain in the neighborhood forever,
there exists some \[n\ge 0\] so that \[f^{\circ n}(x)\] and \[f^{\circ n}(y)\]
have distance greater than \[\epsilon\]. Using Sullivan's results from \S13,
show that this condition is satisfied if and only if \[f\] is hyperbolic.
(However, a map with a parabolic fixed point may be expansive on the Julia set
itself.)

\vfill
%\rightline{jwm, version of 9-2-91}
\end

%% file: cd15-16.tex
\input header
\input ras.tex
\def\caratheodory{Carath\'eodory }
\footline={\hss\tenrm 15-\folio\hss}\pageno=1

\centerline{\bf CARATH\'EODORY THEORY}\bigskip

\centerline{\bf \S15. Prime Ends.}\medskip

Let \[U\] be a simply connected subset of \[\hat\C\] such that the complement
\[\hat\C-U\] is infinite. Then the Riemann Mapping Theorem asserts that
there is a conformal isomorphism
$$	\phi:U \buildrel \approx\over\longrightarrow D\,.	$$
In some cases, \[\phi\] will extend to a homeomorphism from the closure
\[\bar U\] onto the closed disk \[\bar D\]. (Compare Figures 1a and 14a,
together with \S16.7.)
However, this is certainly not true in general, since
the boundary \[\partial U\] may be an extremely complicated object.
As an example, Figure 14b shows a region \[U\] such that one point of
\[\partial U\] corresponds to a Cantor set of distinct points of the circle
\[\partial D\], and Figures 14c, 14d show examples
for which an entire interval
of points of \[\partial U\] corresponds to a single point of the circle.
An effective analysis of the relationship between the compact
set \[\partial U\] and the boundary circle \[\partial D\] was carried out
by \caratheodory in 1913, and will be described here.\smallskip

\midinsert
\hbox to \hsize{
\hfil
\hbox{\RasterBox {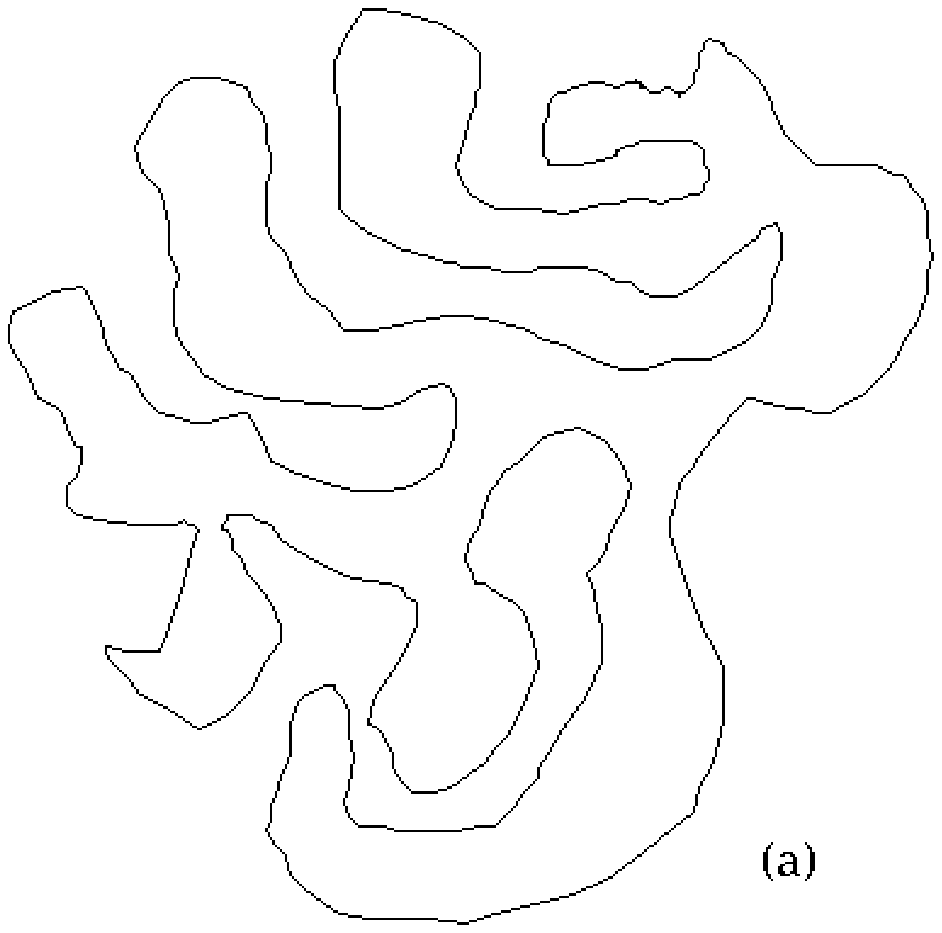} {480} {480} {400}}
\hbox{\RasterBox {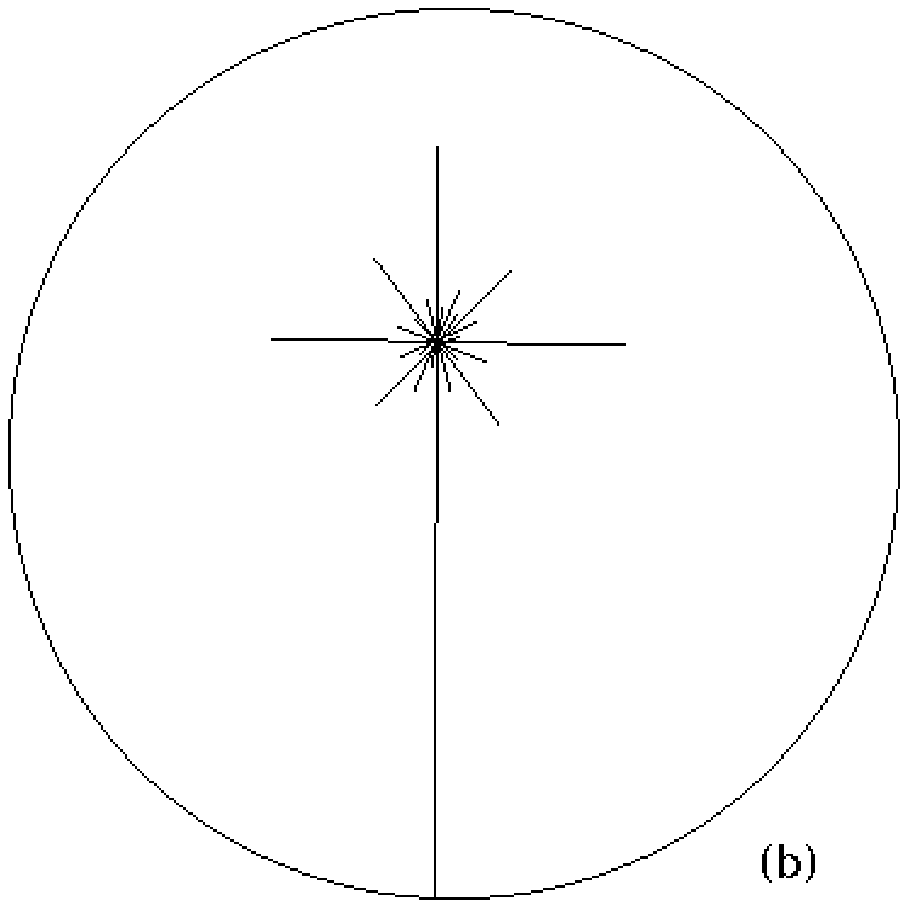} {480} {480} {400}}
\hfil
}
\vskip -.5in
\hbox to \hsize{
\hfil
\hbox{\RasterBox {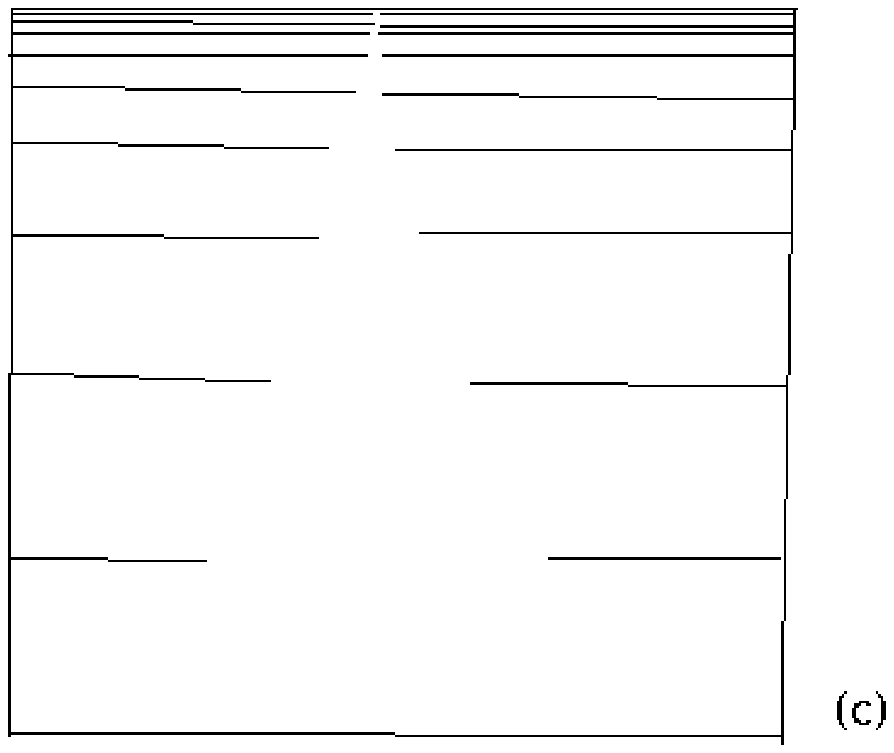} {480} {480} {400}}
\hbox{\RasterBox {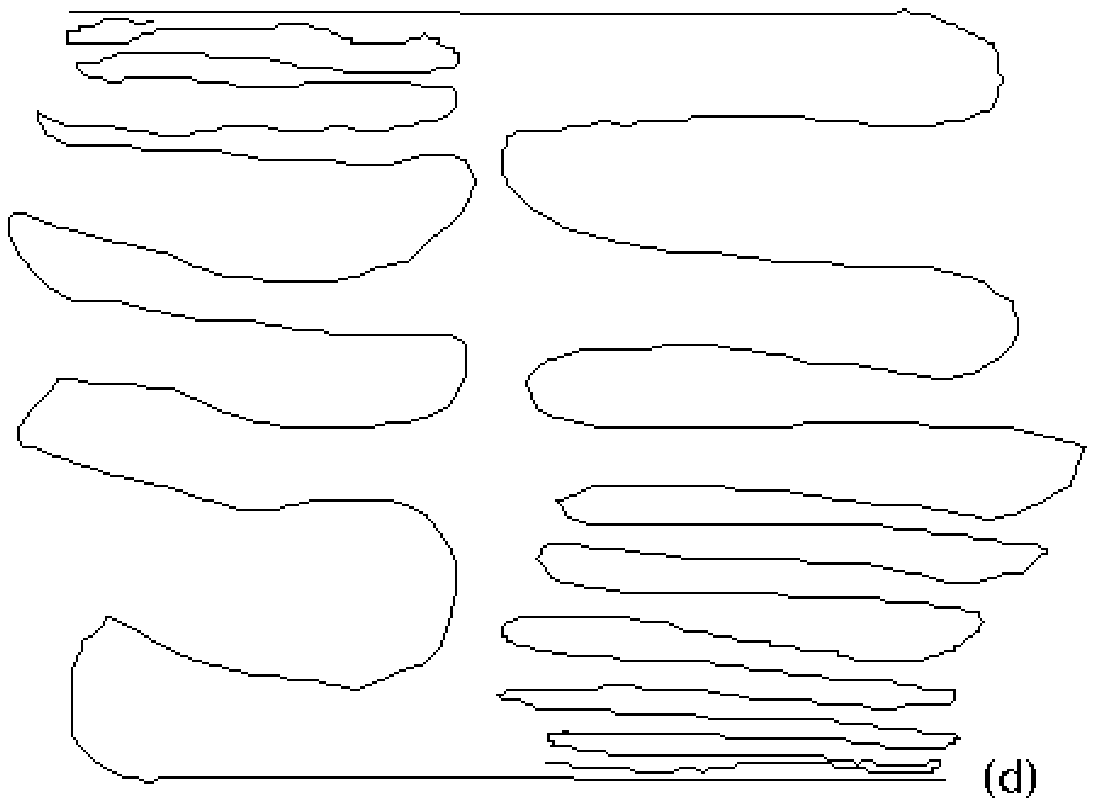} {480} {480} {400}}
\hfil
}

\centerline{Figure 14. The boundaries of
four simply connected regions in \[\C\].}
\vfil
\endinsert

The final construction will be purely topological, but we must first
use analytic methods, as in Appendix B,
to prove several lemmas about the existence of short arcs.
Let \[I=(0,\delta)\] be an open interval of real numbers, and let \[I^2\subset
\C\] be the open square, consisting of all \[z=x+iy\] with \[x,\,y\in I\].
We will always use the spherical metric on \[\hat\C\] when measuring
either arclengths or areas. Recall from \S1(6) that this metric has the form
\[ds=\eta(z)|dz|\] with \[\eta>0\]. (In fact
\[\eta(z)=2/(1+|z|^2)\].) Setting \[z=x+iy\], the corresponding spherical
area element on \[\hat\C\] can be written as \[dA=\eta^2(z)
dxdy\], with total area \[\int\!\!\int dA=4\pi<\infty\].

{\QP{\bf 15.1. Lemma.} {\it If \[f:I^2\to V\] is a conformal isomorphism
from the open square onto an open subset of \[\hat\C\], then for
Lebesgue almost every \[x\in I\] the arc \[f(x\times I)\]
has finite spherical arclength.}\medskip}

{\bf Proof.} (Compare B.1 in the Appendix.)
The length \[L(x)\] of this image \[f(x\times I)\]
can be expressed as \[\int_I\, |\eta f'|dy\]
where the function \[\eta\] is to be evaluated at \[f(x+iy)\],
while the area of \[V\] can be expressed as \[A=\int_I\, a(x)dx\] with
\[a(x)=\int_I\, |\eta f'|^2dy\]. By the Schwarz inequality we have
$$      \bigl(\int_I 1\cdot |\eta f'|\;dy\,\bigr)^2\;\le \;\bigl(
	\int_I 1^2dy\bigr)\cdot\bigl(\int_I |\eta f'|^2 dy\bigr)\,, $$
or in other words \[L(x)^2\le \delta a(x)\]. Since the integral of \[a(x)\] is
finite, it follows that \[a(x)\] is finite almost everywhere, hence
\[L(x)\] is finite almost everywhere.\QED

Here is a more quantitative version of this statement. Again let
\[A\] be the spherical area of \[f(I^2)=V\].

{\QP{\bf 15.2. Lemma.} {\it With \[f:I^2\buildrel \approx\over\to
V\] as above, the length \[L(x)\] of the curve \[f(x\times I)\] satisfies
\[L(x)<2\sqrt A\] for more than half of the points \[x\in I\]. Similarly, the
length of \[f(I\times
y)\] is less than \[2\sqrt A\] for more than half of the points
\[y\in I\].}\medskip}

{\bf Proof.} We must show that the Lebesgue measure of
\[S=\{x\in I:L(x)\ge 2\sqrt{A}\}\] satisfies \[\ell(S)<{1\over 2}
\ell(I)={1\over 2}\delta\]. It is clear that \[4A\ell(S)=(2\sqrt A)^2
\ell(S)\le\int_I L(x)^2dx\].
On the other hand, the
discussion above shows that \[\int_I L(x)^2dx\le\delta A\]. Combining these
two inequalities and dividing by \[4A\], we obtain \[\ell(S)\le\delta/4\],
as required.\QED

Now consider a simply connected open set \[U\subset\hat \C\] and some
choice of Riemann map \[\phi:U\to D\], with inverse \[\psi:D\to U\].

{\QP{\bf 15.3. Theorem.} {\it For almost every point \[e^{i\theta}\] of
the circle \[\partial D\] the radial line \[r\mapsto re^{i\theta}\]
maps under \[\psi\] to a curve of finite spherical
length in \[U\]. In particular, the radial limit
$$	\lim_{r\to 1}\; \psi(re^{i\theta})\in\partial U	$$
exists for Lebesgue almost every \[\theta\]. Furthermore, if we fix
\[\theta\], then for Lebesgue almost every \[\theta'\] the radial limit
of \[\psi(re^{i\theta'})\] is different from the radial limit of \[\psi(r
e^{i\theta})\].}\medskip}

\eject
We will say briefly that almost every image curve \[r\mapsto\psi
(re^{i\theta})\]
in \[U\] {\bit lands} at some single point of \[\partial U\], and that
different values of \[\theta\] almost always correspond to distinct
landing points.

{\bf Remark.} The first part of this Lemma is the univalent version of a
basic {\bit Theorem of Fatou}, which
Carath\'eodory used as the starting point of his theory. Fatou's Theorem
says that any bounded holomorphic function on \[D\] has radial limits
in almost all directions, whether or not it is univalent. (See for example
Hoffman, p. 38.) However the univalent case is all that we will need,
and is easier to prove than the general theorem.\smallskip

{\bf Proof of 15.3.}
Map the upper half-plane onto \[D-\{0\}\] by the exponential
map \[w\mapsto e^{iw}\]. In terms of polar coordinates \[r,\,\theta\] in
\[D\], this means that we set \[w=\theta-i\log r\], so that \[e^{iw}=r
e^{i\theta}\]. Applying the argument of 15.1
to the composition \[w\mapsto \psi(e^{iw})\], we see that \[\psi\] has
radial limits in almost all directions \[e^{i\theta}\].

If this map \[\psi\] is bounded, then a
Theorem of F. and M. Riesz asserts that any
given radial limit can occur only for a set of directions \[e^{i\theta}\]
of measure zero. A proof may be found in Appendix A, Theorem A.3. Now for any
univalent \[\psi\], we can reduce to the bounded case in two steps,
as follows. First suppose that
the image \[\psi(D)=U\] omits an entire neighborhood of some point \[z_0\]
of \[\hat\C\]. Then by composing \[\psi\] with a fractional linear
transformation which carries \[z_0\] to \[\infty\], we reduce to the
bounded case. In general, \[\psi(D)\] must omit at least two values,
which we may take to be \[0\] and \[\infty\]. Then
\[\sqrt \psi\] can be defined as a single valued function which omits an
entire open set of points; and we are reduced to the previous case.\QED

Similarly, using the change of coordinates \[z=e^{iw}\] together
with 15.2, we obtain the following. Let \[\psi:D\buildrel\approx\over
\to U\] as above.

{\QP{\bf 15.4. Lemma.} {\it Any point of \[\partial D\] has a nested
sequence of neighborhoods\break
\[{\cal N}_1\supset{\cal N}_2\supset\cdots\] in \[\bar D\] with the following
properties:

\ref $(1)\] each boundary \[\overline{{\cal N}_k}\cap\overline{D-{\cal N}_k}\]
consists of an open arc \[{\cal A}_k\]
in the interior of \[D\] together with two end points in \[\partial D\],

\ref $(2)\] these boundaries \[\bar{\cal A}_k\] are pairwise disjoint,

\ref
$(3)\] both \[{\cal A}_k\] and its image \[\psi({\cal A}_k)\subset U\] have
finite length which tends to zero as \[k\to\infty\], and

\ref $(4)\] each image arc \[\psi({\cal A}_k)\] has two distinct boundary
points in \[\partial U\], and the\break
\hbox{\kern .8em} closures \[\overline{\psi({\cal A}_k)}\]
are pairwise disjoint.}\medskip}

{\bf Proof.} We choose a neighborhood bounded by two short radial line
segments going out to the boundary of \[D\] together with an arc of a
circle \[|z|={\rm constant}\]. Using 15.2 we see that the image of such
a curve \[{\cal A}\] can have length less than any given epsilon, and using
the last part of 15.3 we see that \[{\cal A}\] can be chosen so that the
two endpoints of \[\psi({\cal A})\] in \[\partial U\] are distinct,
and so that successive arcs have distinct endpoints.\QED

Conversely, let us start with an embedded arc \[\alpha:[0,1)\to U\] which
lands at a point \[z_0\] of \[\partial U\]. That is, we suppose
that \[\lim_{t\to 1}\, \alpha(t)=z_0\].
\eject

{\QP{\bf 15.5. Theorem.} {\it Any arc in \[U\] which lands at one point
\[z_0\] of \[\partial U\] corresponds, under the Riemann map,
to an arc in \[D\] which lands at one point of \[\partial D\].
Furthermore, arcs which land at distinct points of \[\partial U\] necessarily
correspond to arcs which land at distinct points of \[\partial D\].}\medskip}

{\bf Proof.}
If the image arc \[t\mapsto \phi(\alpha(t))\] did not land at a single
point of \[D\], then it would have to oscillate back and forth (or round
and round), accumulating
towards some connected set of limit points on \[\partial D\]. It would then
follow, for some set of \[e^{i\theta}\in \partial D\]
of positive measure, that the radial
limit of \[\psi(re^{i\theta})\] could only be \[z_0\], which is impossible
by the Riesz Theorem.

If arcs landing at two distinct points
of \[\partial U\] corresponded to arcs landing at a single point of \[\partial
D\]. Then, using 15.4, we could cut both of these arcs by the image of an arc
which is arbitrarily short and arbitrarily close to the
boundary both in \[D\] and in \[U\]. Evidently this is impossible.\QED

Using these simple facts, \caratheodory showed that we can reconstruct the
topology of the closed disk \[\bar D\] from the topology of the embedding
of \[U\] into \[\bar U\]. There are a number of possible variations on his
basic definitions. (Compare Ahlfors [1973], Epstein, Ohtsuka.)
We will make use of a variation which is fairly close to the original
construction.\smallskip

A {\bit transverse arc} (or ``crosscut'') in \[U\] is a set
\[A\subset\bar U\] which is homeomorphic to the closed interval \[[0,1]\],
and which intersects the boundary \[\partial U\] only in the two end
points of the interval.\smallskip

Note that it is very easy to construct examples of transverse arcs.
For example we can start with any short line segment inside \[U\] and extend
in both directions until it first hits the boundary.

{\QP {\bf 15.6. Lemma.} {\it Any transverse arc \[A\] cuts \[U\] into two
connected components.}\medskip}

{\bf Proof.} The quotient space \[\bar U/\partial U\], in which the boundary
is identified to a point, is evidently homeomorphic to the 2-sphere. Since
\[A\] corresponds to a Jordan curve in this quotient 2-sphere, the
conclusion follows from the Jordan Curve Theorem. (See for example
Munkres.)\QED

It will be convenient to distinguish these two components of \[U-A\]
by choosing some base point \[b_0\in U\].
Then for any transverse arc \[A\] which is disjoint from \[b_0\] we
define the  neighborhood \[N(A)\] which is ``cut off'' by \[A\] to be the
component of \[U-A\] which does not contain \[b_0\].\smallskip

{\bf Definition.}
A {\bit fundamental chain} is a infinite sequence \[A_1\,,\, A_2\,,\,\ldots\]
of disjoint transverse arcs with the property that the corresponding
neighborhoods are nested:
$$	 N(A_1)\supset N(A_2)\supset N(A_3)\supset\cdots\,,	$$
and with the property that the diameter of \[A_i\] tends to zero as
\[i\to\infty\].\smallskip

Note however that the neighborhoods \[N(A_i)\] are not required to become
small as \[i\to\infty\]. (Compare Figure 14c, d.)
By the {\bit impression} (or {\bit support})
of a fundamental chain \[\{A_i\}\] we mean
the intersection of the closures \[\bar N(A_i)\]. It is not difficult to
check that this impression is always a compact connected subset
of \[\partial U\].
Evidently the impression consists of a single point if and only if the
diameter of \[N(A_i)\] tends to zero as \[i\to\infty\].\smallskip

Two such fundamental chains \[\{A_i\}\] and \[\{A'_j\}\] are said to be
{\bit equivalent} if each \[N(A_i)\] contains some
\[N(A'_j)\] and each \[N(A'_j)\] contains some \[N(A_i)\]. An equivalence
class of fundamental chains is also called a {\bit prime end} \[{\cal E}\]
in \[U\], or a point in the {\bit \caratheodory boundary} of \[U\].

Two fundamental chains \[\{A_i\}\] and \[\{A'_j\}\] are said to be
{\bit disjoint} if
\[	N(A_i)\cap N(A'_j)=\emptyset	\]
for some \[i\] and \[j\].

{\QP{\bf 15.7. Lemma.} {\it Any two fundamental chains are either equivalent
or disjoint.}\medskip}

{\bf Proof.} Given \[A_i\] we can choose \[j\] large enough so that the
diameter of \[A'_j\] is less than the distance between \[A_i\] and
\[A_{i+1}\]. It then follows easily that \[N(A'_j)\] must be either
contained in \[N(A_i)\] or disjoint from \[N(A_{i+1})\].\QED

We can now define the {\bit \caratheodory completion} \[\hat U\] of \[U\]
to be the disjoint union of the set \[U\] and the set of all prime ends of
\[U\], topologized as follows. For any transverse arc \[A\subset \bar U-
\{b_0\}\] we define the neighborhood \[{\cal N}_A\] to be
the neighborhood \[N(A)\subset U\] together with the set consisting of all
prime ends \[{\cal E}\] for which some representative fundamental chain
\[\{A_i\}\] satisfies \[A_1=A\]. These neighborhoods \[{\cal N}_A\],
together with the open subsets of \[U\], form a basis for the required
topology.

{\QP{\bf 15.8. Lemma.} {\it In the case of the unit disk \[D\], the
identity map of \[D\] extends uniquely to a homeomorphism between the closure
\[\bar D\subset\C\] and the \caratheodory completion \[\hat D\].}
\medskip}

(Compare 16.7.) In particular, the prime ends of \[D\]
are in one-to-one correspondence with the points of \[\partial D\].
More precisely, the impression of any prime end of \[D\] is a single point of
\[\partial D\], and
each point of \[\partial D\] is the impression of one and only one
prime end. The proof is not difficult.\QED

For an arbitrary simply connected and Hyperbolic
\[U\subset\hat\C\] we now have the following.

{\QP{\bf 15.9. Theorem.} {\it The Riemann map \[\phi:U\to D\] extends
uniquely to a homeomorphism between the completion \[\hat U\] and the
completion \[\hat D\cong\bar D\].}\medskip}

{\bf Proof.} It follows immediately from 15.5 that every transverse arc
in \[U\] corresponds to a transverse arc in \[D\], and that disjoint
transverse arcs correspond to disjoint transverse arcs. Thus every end
in \[U\] gives rise to an end in \[D\]. On the other hand,
given any point of \[\partial D\], we can find a nested sequence of
neighborhoods in \[\bar D\], each bounded by a figure consisting of two
radial segments and an arc which is completely inside \[D\]. As in 15.4,
we can choose these figures so as to correspond to transverse arcs in
\[\bar U\] which are pairwise disjoint, with lengths converging to zero.
The resulting fundamental chain in \[U\] determines the required prime end.
Further details will be left to the reader.\QED
\eject
\footline={\hss\tenrm 16-\folio\hss}\pageno=1

\centerline{\bf \S16. Local Connectivity.}\bigskip

This will be a continuation of the preceeding section, describing
Carath\'eodory's theory in the locally connected case. First a number
of definitions and well known lemmas.
We will always assume that our topological spaces are Hausdorff.
By definition, a topological space \[X\] is {\bit locally connected} at
the point \[x\in X\] if there exist arbitrarily small connected neighborhoods
of \[x\] in \[X\].

Here our ``neighborhoods'' are not required to be open subsets
of \[X\], but only to be sets which contain a (relatively) open set
containing \[x\]. Unfortunately, some authors make a slightly different
definition by considering only open neighborhoods. Let us say
that \[X\] is {\bit openly locally connected} at a point \[x\] if there exist
arbitrarily small connected open neighborhoods of \[x\] in \[X\]. This
is a stronger condition (compare Problem 16-1); however we have the following.

{\QP{\bf 16.1. Lemma.} \it The space \[X\] is locally connected at every point
\[x\in X\] if and only if every open subset of \[X\] is a union of
connected open subsets of \[X\].\medskip}

\noindent If this conditions is satisfied, then evidently
\[X\] is openly locally connected at every point.
Such spaces are simply said to be {\bit locally connected}.
The proof is straightforward, and will be left to the reader.\QED

{\bf 16.2. Remark.} If \[X\] is compact metric, with distance
function \[\rho(x,x')\], then an equivalent condition is that:

{\narrower\smallskip\noindent\it
For any \[\epsilon>0\] there exists \[\delta>0\] so that any two points
with distance\break \[\rho(x,x')<\delta\] are contained in a connected set of
diameter less than \[\epsilon\].\smallskip}

\noindent Again the proof is straightforward, and will be omitted.\smallskip

The space \[X\] is {\bit path connected}$\, $ if there exists a continuous
map from the unit interval \[[0,1]\] into \[X\] which joins any two given
points. It is {\bit arcwise connected}$\,$ if there is a topological embedding
of \[[0,1]\] into \[X\] which joins any two given distinct points.

{\QP{\bf 16.3. Lemma.} {\it A space is path connected if and only if it
is arcwise connected.}\medskip}

\noindent (Recall that all spaces must be Hausdorff.)
Proofs will be deferred until the end of this section.
There is a corresponding concept of {\bit ``locally path connected''} or
equivalently {\bit ``locally arcwise connected''}.
In general, a connected space need not be path connected (Figure 14d). However:

{\QP{\bf 16.4. Lemma.} {\it If a compact metric space \[X\] is locally
connected, then it is locally path connected. Hence every connected
component of \[X\] is actually path connected.}\medskip}

\noindent We will also need the following statement.

{\QP{\bf 16.5. Lemma.} {\it Any continuous image of a compact locally
connected space is compact and locally connected.}\medskip}

\noindent The following is the principal result of this section. Let \[U\subset
\hat\C\] be an open set which is conformally isomorphic to the unit disk.
\eject

{\QP{\bf 16.6. Theorem of Carath\'eodory.} {\it The inverse Riemann map
\[\psi:D\buildrel
\approx\over\to U\] extends continuously to a map from the
closed disk \[\bar D\] onto \[\bar U\] if and only if the boundary \[\partial
U\] is locally connected, or if and only if the complement \[\hat\C-U\]
is locally connected.}\medskip}

\noindent \caratheodory also proved the following.

{\QP{\bf 16.7. Theorem.} {\it If the boundary of \[U\] is a Jordan curve, then
the Riemann map extends to a homeomorphism from the closure \[\bar U\]
onto the closed disk \[\bar D\].}\bigskip}

\noindent The proofs begin as follows.\smallskip

{\bf Proof of 16.3.} Let \[f=f_0:[0,1]\to X\] be any continuous path with\break
\[f(0)\ne f(1)\]. We must construct an embedded arc \[A\subset X\] from
\[f(0)\] to \[f(1)\].
Choose a closed subinterval \[I_1=[a_1\,,\,b_1]\subset [0,1]\] whose
length \[0\le\ell(I_1)=b_1-a_1<1\] is as large as possible, subject to the
condition that \[f(a_1)=f(b_1)\]. Now,
among all subintervals of \[[0,1]\] which are disjoint from \[I_1\],
choose an interval \[I_2=[a_2\,,\,b_2]\]  of maximal length subject to
the condition \[f(a_2)=f(b_2)\].
Continue this process inductively, constructing disjoint subintervals
of maximal lengths \[\;\ell(I_1)\ge\ell(I_2)\ge\cdots\ge 0\;\] subject to
the condition that \[f\]
is constant on the boundary of each \[ I_j\].

Let \[\alpha:[0,1]\to X\] be
the unique map which is constant on each closed interval \[I_j\], and which
coincides with \[f\] outside of these subintervals. Thus \[\alpha(t)\] must
coincide with \[f(t)\] for \[t\in\partial I_j\].
Then it is easy to check
that \[\alpha\] is continuous, and that for each point\break
\[x\in \alpha([0,1])\] the
preimage \[\;\alpha^{-1}(x)\subset[0,1]\;\] is a (possibly degenerate) closed
interval of real numbers. {\it We will show that these conditions,
with \[\alpha\] non-constant,
imply that the image \[A=\alpha([0,1])\subset X\;\] is an embedded arc
joining \[\alpha(0)\] to \[\alpha(1)\].}

Choose a countable dense subset \[\{t_1\,,\, t_2\,,\,\ldots\}\subset[0,1]\].
Define a total ordering of the image \[A\] by specifying
that \[\alpha(s)<\alpha(t)\] if and only if \[\alpha(s)\] and \[\alpha(t)\]
are distinct points with \[s<t\].
An order preserving homeomorphism \[h:[0,1]\buildrel\approx\over
\to A\] can now be constructed inductively as follows.
We must set \[h(0)=\alpha(0)\] and \[h(1)=\alpha(1)\]. For each
dyadic fraction \[0<m/2^k<1\] with \[m\] odd, let us assume inductively that
\[h((m-1)/2^k)\] and \[h((m+1)/2^k)\] have already been defined. Then choose
the smallest index \[j\] so that
$$	h\big({m-1\over 2^k}\big)\;<\;\alpha(t_j)\;<\;h\big({m+1\over2^k}\big)	$$
and set \[h(m/2^k)=\alpha(t_j)\]. It is not difficult to check that the \[h\]
constructed in this way on dyadic fractions extends uniquely to an order
preserving map from \[[0,1]\] to \[A\], which is necessarily a
homeomorphism.\QED

{\bf Proof of 16.4.} Let \[X\] be compact metric and locally connected.
Given \[\epsilon>0\], it follows from 16.2 that we can choose a
sequence of numbers \[\delta_n>0\] so
that any two points with distance \[\rho(x,x')<\delta_n\] are contained in
a connected set of diameter less than \[\epsilon/2^n\]. {\it We will prove
that any two points
\[x(0)\] and \[x(1)\] with distance \[\rho(x(0),x(1))
<\delta_0\] can be joined by a path of diameter at most
\[4\epsilon\].} The plan of attack is as follows. We will choose a sequence
of denominators \[1=k_0<k_1<k_2<\cdots\], each of which divides the next.
Also, for each fraction of the form \[i/k_n\] between 0 and 1 we will
choose an intermediate point \[x(i/k_n)\]. These are to satisfy two
conditions:

{\narrower\smallskip\noindent
(1) Any two consecutive fractions \[i/k_n\] and \[(i+1)/k_n\] must correspond
to points \[x(i/k_n)\] and \[x((i+1)/k_n)\] which have distance less
than \[\delta_n\] from each other.\smallskip}

\noindent
Hence any two such points are contained in a connected set \[C(i,k_n)\]
which has diameter less than \[\epsilon/2^n\]. The second condition is that:

{\QP
(2) Each point of the form \[x(j/k_{n+1})\], where \[j/k_{n+1}\] lies between
\[i/k_n\] and \[(i+1)/k_n\], must belong to this set \[C(i,k_n)\], and
hence have distance less than \[\epsilon/2^n\] from \[x(i/k_n)\]
and from \[x((i+1)/k_n)\].\smallskip}

\noindent The construction of such denominators \[k_n\] and intermediate
points \[x(i/k_n)\], by induction on \[n\], is completely straightforward,
and will be left to the reader. Thus we may assume that \[x(r)\] has been
defined for a dense set of rational numbers \[r\] in the unit interval.

Next we will prove that this densely defined correspondence \[\;r\mapsto x(r)
\;\] is uniformly continuous. Let \[r\] and \[r'\] be any two rational numbers
for which \[x(r)\] and \[x(r')\] are defined. Suppose
that \[|r-r'|<1/k_n\]. Then choosing \[i/k_n\] as close as possible to both
\[r\] and \[r'\],
it is easy to show that \[x(r)\] and \[x(r')\] have distance less
than \[2\epsilon/2^n\] from \[x(i/k_n)\]. Hence they have distance less than
\[4\epsilon/2^n\] from each other. This proves uniform continuity; and it
follows that there is a unique continuous extension \[t\mapsto x(t)\]
which is defined for all \[t\in[0,1]\]. In this way, we have constructed the
required path of diameter \[\le 4\epsilon\] from \[x(0)\] to \[x(1)\].
Thus \[X\] is locally path connected; and the rest of the
argument is straightforward.\QED

{\bf Proof of 16.5.} Let \[f(X)=Y\] where \[X\] is compact and locally
connected. Then \[Y\] is clearly compact, and we must show that it is
locally.connected. Given any point \[y\in Y\] and open neighborhood \[V\subset Y\],
choose a smaller neighborhood \[U\] with \[\bar U\subset V\].
Then \[f^{-1}(V)\] is a union of disjoint connected open sets, and finitely
many of these connected open sets, say \[V_1\,,\,\ldots\,,\,V_q\],
will suffice to cover the compact subset \[K=f^{-1}(\bar U)\]. Let \[N\]
be the union of those \[f(V_i)\] which contain the given point \[y\].
Evidently \[N\] is connected and contained in \[V\]. In order to prove
that \[N\] is a neighborhood of \[y\], let \[L\] be the union of those compact
images \[f(K\cap V_i)\] which do not contain \[y\]. Then \[U-L\] is
an open neighborhood of \[y\] which is contained in \[N\]. Thus \[N\subset V
\] is a connected neighborhood of \[y\], which proves that
\[Y\] is locally connected.\QED

{\bf Proof of Theorem 16.6.} We will always use the spherical metric on
\[\hat\C\]. If either \[\partial U\] or \[\hat\C-U\] is
locally connected, then given \[\epsilon\] we can choose \[\delta\] so
that any two points of distance less than \[\delta\] in \[\partial U\] are
joined by an arc \[A\] of diameter less than \[\epsilon\] in \[\hat\C-U\].
Now suppose that we start with some point \[w_0\in\partial D\]. According to
Lemma 15.4, we can find a neighborhood of diameter less than \[\delta\] in
\[\bar D\] which is bounded by an arc which maps under \[\psi:D\buildrel
\approx\over\to U\] to an arc \[A'\subset U\] which has length less than \[\delta\].
Furthermore, we can choose this arc so that its two endpoints in \[\partial
U\] are distinct. Let \[A\subset\hat\C-U\] be an arc in the complement which
has diameter less than \[\epsilon\], and has the same endpoints. Then the union \[A\cup A'\] is a Jordan
curve of diameter less than \[\epsilon+\delta\] in \[\hat\C\]. Hence it
bounds a region of diameter less than \[\epsilon+\delta\]. Now as \[\epsilon
\to 0\] this region must converge to a point in \[\partial U\], which
we define to be \[\psi(w_0)\]. Thus we have constructed a function from the
closed disk \[\bar D\] to \[\bar U\] which extends \[\psi:D\to U\]. It
is not difficult to check that this extension is continuous.

Conversely, if \[\psi\] extends continuously over \[\bar D\], the the boundary
\[\partial U\] is a continuous image of the circle \[\partial D\]. Hence
it follows from Lemma 16.5 that \[\partial U\] is locally connected. To
see that \[\hat\C-U\] must also be locally connected, let \[N\] be an
arbitrarily small connected neighborhood of a point \[z\] within \[\partial U\].
Choose \[\delta>0\] so that the \[\delta$-neighborhood of \[z\] within
\[\partial U\] is contained in \[N\]. Then the \[\delta$-neighborhood of \[z\]
within \[\hat\C-U\], together with \[N\], is clearly also connected. Therefore
\[\hat\C-U\] is locally connected.\QED

{\bf Proof of Theorem 16.7.} If \[\partial U\] is a Jordan curve, that is
a homeomorphic image of the circle, then we certainly have a continuous
extension \[\psi:\bar D\to \bar U\] by the preceeding theorem. Suppose that
two distinct points of the circle \[\partial D\] mapped to the same point \[z_0\]
of \[\partial U\]. Then the straight
line segment \[A\] joining these two points in \[\bar
D\] would map to a Jordan curve \[\Gamma\subset\hat\C\] under \[\psi\].
Such a Jordan curve must separate \[\hat\C\] into two components, of
which \[\Gamma\] is the common boundary. Note that
the two components of \[D-A\] must map into the two distinct components
of \[\hat\C-\Gamma\] under \[\psi\]. Since
\[\Gamma\] intersects the Jordan curve \[\partial U\] in a single point
\[z_0\], it follows that the entire connected set \[\partial U-\{z_0\}\]
must be contained in just one of the two connected components of
\[\hat\C-\Gamma\].
But this means that one of the two components of \[D-A\] must
map into the other component of \[\hat\C-\Gamma\], which is disjoint from
\[\partial U\]. Hence all of the corresponding  boundary points in
\[\partial D\] can only map to the single point
\[z_0\in\Gamma\cap\partial U\] under \[\psi\]. This is impossible by
the Riesz Theorem A.3 of Appendix A.\QED\lln

{\bf Problem 16-1.} Let \[X\subset\C\] be the compact connected set which
is obtained from the unit interval \[[0,1]\] by drawing line segments from
\[1\] to the points \[{1\over 2}(1+i/n)\] for\break
\[n=1\,,\,2\,,\,3\,,\,\ldots\]
and then adjoining the successive images of this configuration under the map
\[z\mapsto z/2\]. (Figure 15.)
Show that \[X\] is locally connected at the origin,
but not openly locally connected.\vfil
\midinsert
\insertRaster 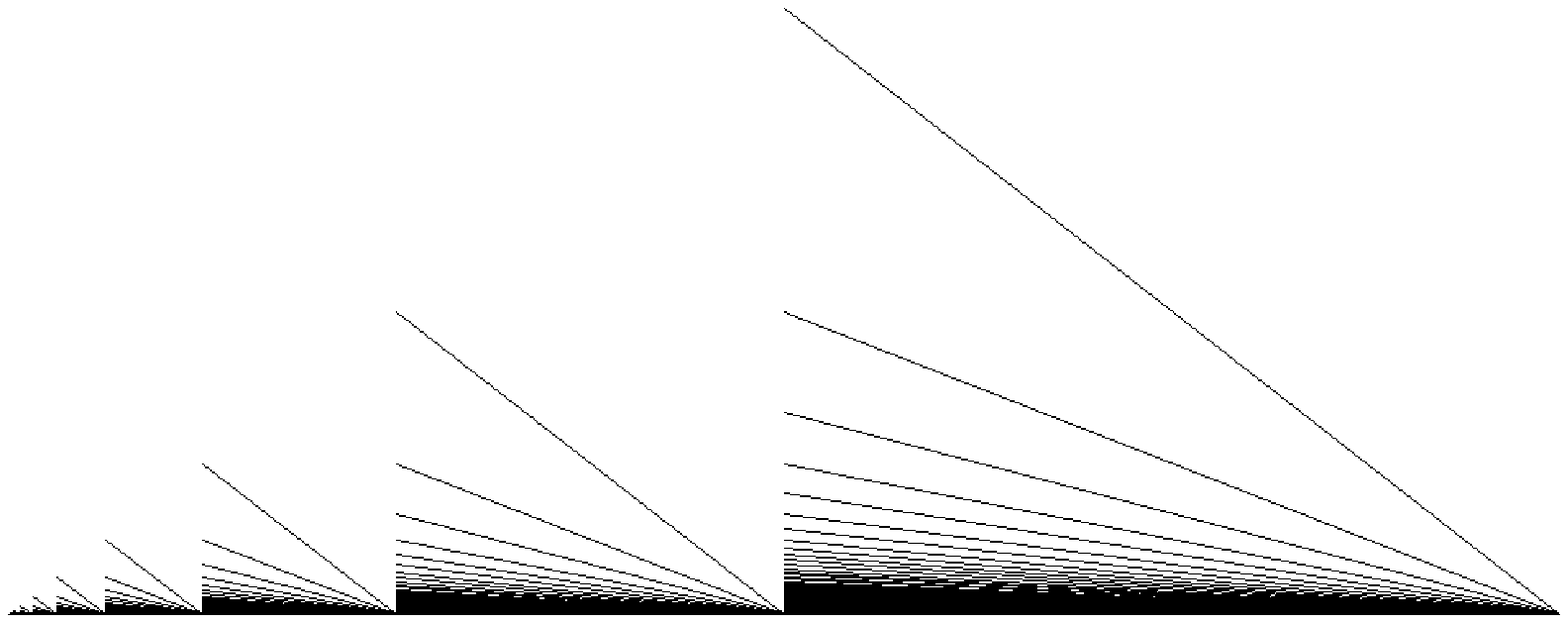 pixels 1152 by 450 scaled 300
\bigskip
\centerline{Figure 15. The witch's broom.}
\endinsert
\vfil

% \rightline{Milnor, Conformal Dynamics Lectures \S15,16 draft, 1-21-90}
\end

%% file: cd17-18.tex
\input header
\input ras.tex
\footline={\hss\tenrm 17-\folio\hss}\pageno=1

\centerline{\bf POLYNOMIAL MAPS}\bigskip

\centerline{\bf \S17. The Filled Julia Set.}\bigskip

Let \[f:\hat\C\to\hat\C\] be a polynomial map of degree \[d\ge 2\], say
$$	f(z)=a_dz^d+\cdots+a_1z+a_0	$$
with \[a_d\ne 0\]. Then \[f\] has a superattracting fixed point at infinity.
In particular, there exists a constant \[k=k_f\] so that any point \[z\]
with \[|z|>k_f\] belongs to the basin of attraction of infinity.

{\bf Definition.} The complement of the attractive basin of infinity, that
is the set of all points \[z\in\C\] with bounded forward orbit under \[f\],
is called the {\bit filled Julia set} \[K=K(f)\].

{\QP{\bf 17.1. Lemma.} {\it This filled Julia set \[K\] is a compact set
consisting of Julia set \[J\]
together with all of the bounded components of the complement \[\C\ssm J\].
These bounded components $($if any$)$ are all simply connected.
The Julia set \[J\] is equal to the topological boundary \[\partial K\].
}\medskip}

{\bf Proof.} (Compare Problem 4-1.)
If \[z\] is a boundary point of \[K\], then \[z\] itself has
bounded orbit, but points arbitrarily close have orbits which converge to
the point at infinity. Thus the family of iterates of \[f\] cannot be normal
in any neighborhood, hence \[z\in J\]. Conversely, if \[z\in J\], then
\[z\] has bounded orbit, but it follows from 4.6 that points arbitrarily
close to \[f\] have unbounded orbit, which therefore converges to infinity.
Thus \[J=\partial K\subset K\], and it follows that the unbounded component
of the complement \[\C\ssm J\] consists of points attracted to infinity. But if
\[z\] belongs to a bounded component of \[\C\ssm J\], then it
follows from the Maximum Modulus Principle that every point in the orbit
of \[z\] belongs to the closed disk of radius \[k_f\], hence \[z\in K\].
Finally, any bounded component \[U\] of \[\C\ssm J\] must be simply connected,
since for any Jordan
curve \[\Gamma\subset U\] the entire bounded component of \[\C\ssm \Gamma\] must
also be contained in \[U\], again by the Maximum Modulus Principle.
\QED

It follows from the B\"ottcher
Theorem \S6.7 that exists a new coordinate \[w=\phi(z)\]
which is defined throughout some neighborhood \[|z|>c_f\] of infinity, with
\[w\to\infty\] as \[z\to\infty\], and which satisfies
$$	\phi\circ f\circ\phi^{-1}(w)\;=\; w^d	\,. $$
Further, \[\phi\] is unique up to multiplication with a \[(d-1)$-st root
of unity.\smallskip

The function \[G(z)=\log|\phi(z)|\] is called the canonical {\bit potential
function} or {\bit Green's function}$\,$ for \[K(f)\]. Note the identity
$$	G(z)\;=\;G(f(z))/d\,.	\eqno (*)	$$
Although we have so far defined \[G\] only in a neighborhood of infinity,
there is one and only one extension to all of \[\C\] which is continuous
and satisfies (*). In fact we define \[G(z)=0\] for \[z\in K(f)\], and
\[G(z)=G(f^{\circ n}(z))/d^n\] otherwise, where \[n\] is large enough so
that \[|f^{\circ n}(z)|>c_f\]. Alternatively, we can simply set
$$	G(z)\;=\;\lim_{n\to\infty}\, \log^+|f^{\circ n}(z)|/d^n\,, $$
where \[\log^+(x)=\log(x)\] for \[x\ge 1\] and \[\log^+(x)=0\] for
\[0\le x\le 1\].
Note that \[G\] is a smooth real analytic function outside of the Julia set
\[J\], but is only continuous at points of \[J\].
The verification of these facts will be left to the reader.

The lines
\[|\phi(z)|={\rm constant}>0\] are called {\bit equipotential curves}.
Their orthogonal trajectories, that is the images under \[\phi^{-1}\]
of the radial lines extending out towards infinity from the unit disk
are called {\bit external rays} for \[f\]. Both of these families of
curves are clearly
smooth except at critical points of \[G\]. {\it Note that a point \[z\in
\C\ssm K\] is a critical point of \[G\] if and only if it is a {\bit pre-critical
point} of \[f\], that is a critical point of some iterate \[f\circ\cdots
\circ f\].} Again proofs are easily supplied.

For each \[r>1\] let \[V(r)\] be the compact region with boundary consisting
of all \[z\] with \[G(z)\le\log(r)\]. Evidently the boundary \[\partial
V(r)\] is just the equipotential curve \[\{z:G(z)=\log r\}\].

{\QP{\bf 17.2. Lemma.} \it Suppose that there are no critical points of \[G\]
$($or pre-critical points of \[f\,)\] on the
boundary \[\partial V(r)\]. If \[m\ge 0\] is the algebraic number of critical
points of \[G\] in the complement \[\C\ssm V(r)\], then \[V(r)\] is a disjoint
union of \[m+1\] closed topological disks, each of which intersects the
Julia set \[J(f)\].\medskip}

{\bf Proof.} If \[\partial V(r)\] contains no critical points (or equivalently
if \[\log r\] is a regular value of the map \[G\]), then evidently \[V(r)\]
is a smooth compact manifold with boundary. Each component must be a
closed topological disk. For otherwise, the complement \[\C\ssm V(r)\] would
have a bounded component, which is impossible since the function \[G\]
would have to take its maximum on the boundary of such a component. Similarly,
the minimum value of \[G\] on each component of \[V(r)\] must occur at
a point of \[K(f)\]. Hence each such component must intersect the Julia set
\[\partial K(f)\].

Given \[r>1\] as above, let us choose \[n\] so large that the region
\[V(r^{d^n})\] is a single topological disk. It will be convenient to
set \[s=r^{d^n}\]. Then the composition \[f^{\circ n}
\] maps \[V(r)\] onto \[V(s)\] by a \[d^n$-fold branched covering map.
According to the Riemann-Hurwitz formula, \S5.1, the algebraic number
of branch points of \[f^{\circ n}\] in \[V(r)\] is equal to \[d^n\chi(V(
s))-\chi(V(r))=d^n-\chi(V(r))\],
where in our case \[\chi\] is just the number of
connected components. If \[m\] is the number of branch points of \[f^{\circ n}
\] outside of \[V(r)\], then, since the total number of branch points is
\[d^n-1\], we can write this formula as
$$	d^n-1-m\;=\; d^n-\chi(V(r))\,,	$$
or in other words \[\chi(V(r))=m+1\]. This completes the proof,
since critical points of \[G\] are the same as critical points of
\[f^{\circ n}\] for large \[n\].\QED

Following Brown, a compact set in Euclidean space is said to be$\,$ {\bit
cellular} $\,$ if it is equal to the intersection of a nested sequence
of closed embedded topological disks, each of which contains the next
in its interior. The following result is essentially due to Fatou and Julia.

{\QP{\bf 17.3. Theorem.} {\it For a polynomial map \[f\] of degree \[d\ge 2\],
there are just two mutually exclusive possibilities: If the
filled Julia set \[K\] contains all of the finite critical points of \[f\],
then \[K\] and \[J=\partial K\] are connected, and in fact \[K\] is a cellular
set. Furthermore, the B\"ottcher map near infinity extends to a conformal isomorphism
$$	\phi\;:\;\C\ssm K\;\;\buildrel\approx\over\longrightarrow\;\; \C\ssm\bar D\,.	$$
On the other hand, if \[f\]
has at least one critical point in \[\C\ssm K\], then both \[J\] and \[K\]
have uncountably many connected components.}\medskip}

{\bf Proof of 17.3.} If \[K\] contains all of the finite critical points,
then it follows from 17.2 that \[K\] is the intersection of a strictly nested
family of embedded disks \[V(r)\]. Thus \[K\] is cellular, and hence connected.
Furthermore, the proof shows that each \[\C\ssm V(r)\] is a \[d^n$-fold unramified
covering of some \[\C\ssm V(s)\], which is isomorphic to \[\C\ssm\bar D_s\] under
the B\"ottcher map. We can then lift to a B\"ottcher isomorphism
\[\C\ssm V(r)\buildrel\approx\over\to\C\ssm\bar D_r\] which is compatible near
infinity and satisfies
$$	\phi(f(z))=\phi(z)^d\,.	$$
Passing to the limit as \[r\to 1\], we see that \[\phi\] can be defined
throughout \[\C\ssm K\]. Finally, note that the Julia set
\[J\] can be expressed as the intersection of the closures of the bounded
connected annuli \[\phi^{-1}(\{w:1<|w|<r\})\]; hence \[J\] is also connected.

Conversely, if \[f\] has at least one critical point in \[\C\ssm K\], then the
argument above shows that some \[V(r)\] has two or more connected components,
and that each of these components
intersects the subsets \[J\subset K\subset V(r)\]. Thus \[J\] has uncountably
many components by 11.3. The analogous proof for \[K\] will be left to the
reader.\QED

For information on the structure of \[K\] when it is not connected,
see Branner \& Hubbard, Blanchard.
As an immediate consequence of Carath\'eodory's Theorem 16.6
we have the following.
(Note that a compact set which has infinitely many components clearly cannot
be locally connected.)

{\QP{\bf 17.4. Corollary.} {\it With \[f\] as above, the following three conditions
are equivalent:

$(a)$ The set \[J(f)\] is locally connected.

$(b)$ The set \[K(f)\] is locally connected.

$(c)$ K(f) is connected, and the inverse Riemann mapping \[\psi:\hat\C\ssm\bar
D\buildrel\approx\over\to\hat\C\ssm K\] extends continuously over the boundary,
yielding a continuous map from the unit circle onto \[J\] satisfying
the identity \[f(\psi(w))=\psi(w^d)\].}\medskip}

\noindent Following Douady and Hubbard, we have the following.

{\QP{\bf 17.5. Theorem.} {\it If the polynomial map \[f\] is hyperbolic or
subhyperbolic,
with \[K(f)\] connected, then \[K(f)\] is locally connected.}\medskip}

For non locally connected examples, see 18.6.\smallskip

{\bf Proof of 17.5.} We may assume that there is a smooth Riemannian metric
or orbifold metric
whose restriction to some small neighborhood \[U\] of \[J\] is strictly
increasing under \[f\],
say by a factor of at least \[k>1\]. Choose \[r>1\] small enough so that
the neighborhood \[V(r)\] of \S17.2 is contained in \[K\cup U\], and let
\[s=\root d\of r\]. Consider the (half-closed) annulus
$$	A_0\;=\; V(r)\ssm V(s)\,,	$$
which is isomorphic under \[\phi\] to the ``round'' annulus \[\bar D_r\ssm\bar
D_s\]. Then \[V(r)\ssm K\] can be expressed as the countable
union \[A_0\cup A_1\cup
\cdots\], where \[A_n\] is the annulus \[f^{-n}(A_0)\]. Similarly,
\[\bar D_r\ssm\bar D_1\] is the union of a corresponding family of round
annuli.

Each radial line segment in \[\bar D_r\ssm\bar D_s\] corresponds
to a smooth curve in \[A_0\], which is a segment of some external ray.
Let \[\ell_0\] be the
maximum of the lengths of these curve segments in \[A_0\], and let \[\ell_n\] be the
maximum length of the corresponding curve segments in \[A_n\]. Since \[f^{\circ n}\]
maps \[A_n\] onto \[A_0\], preserving this foliation by external rays,
and stretching all distances by at least \[k^n\], it follows that
\[\ell_n\le\ell_0/k^n\]. Therefore, the total length of any external ray
intersected with \[V(r)\] is uniformly bounded by the quantity\break \[\;\ell_0
(1+k^{-1}+k^{-2}+\cdots)\;=\;\ell_0
/(k-1)\;<\;\infty\,\]. It then follows easily that the maps
$$	e^{i\theta}\;\mapsto\;\psi(re^{i\theta})	$$
from the unit circle onto \[\partial V(r)\]
converge uniformly as \[r\to 1\] to the required limit mapping, which
carries the unit circle continuously onto the Julia set.\QED
\vfil\eject % \bigskip\bigskip

\footline={\hss\tenrm 18-\folio\hss}\pageno=1

\centerline{\bf \S18. External Rays and Periodic Points.}\bigskip

To fix our ideas and simplify the discussion, let us assume
that the filled Julia set \[K\] is connected. (For the general case, see
[DH2].) By 17.3 the complement \[\C\ssm K\]
is isomorphic to the complement \[\C\ssm\bar D\] under a conformal isomorphism
$$	\phi:\C\ssm K\;\buildrel\approx\over\longrightarrow\;\C\ssm\bar D	$$
which is compatible with the dynamics:
$$	\phi(f(z))\;=\;\phi(z)^d\,.	$$
Here \[\phi\] is unique up to multiplication by a \[(d-1)$-st root of unity.
In practice, it is customary to make a linear change of coordinates so that
the polynomial \[f\] is {\bit monic}, ie., has leading coefficient equal to
\[+1\]. There is then a preferred choice of B\"ottcher coordinate \[w=
\phi(z)\], determined by the requirement that \[\;w/z\to 1\] as \[|z|\to
\infty\].\smallskip

{\bf Definition.} Let \[\psi:\C\ssm\bar D\buildrel\approx\over\to\C\ssm K\] be the inverse map.
The image under \[\psi\] of a radial line
$$	\{re^{2\pi it}:r>1\}	$$
in \[\C\ssm\bar D\] is called the {\bit external ray} \[\;R_t\]
at {\bit angle} \[t\] in \[\C\ssm K\]. Note that our {\it angles} are elements of
\[\R/\Z\], measured in fractions of a full turn and not in radians. Evidently
$$	f(R_t)\;=\;R_{dt}	\,.	$$
That is: {\it \[f\] maps the external ray at angle \[t\] to the external ray
at angle \[dt\].} In particular, if \[t\] is a fraction of the form \[m/
(d-1)\] then \[f(R_t)=R_t\].\smallskip

By definition, the external ray \[R_t\] {\bit lands} at some point
\[z_t\], which necessarily belongs to the Julia set \[J\], if the points
\[\psi(re^{2\pi it})\in \C\] tend to a well defined limit \[z_t\] as \[r\to
1\]. If \[J\], or equivalently \[K\], is locally connected, then according to
17.4 {\it every} external ray lands, at a point \[z_t\] which depends
continuously on \[t\]. Douady and Hubbard call the resulting parametrization
\[t\mapsto z_t\] the {\bit Carath\'eodory loop} on the Julia set.
However, in this section we will usually not assume that \[J\] is locally
connected.\smallskip

{\bf Remark.} A priori, it is possible that every external ray may land
even if \[K\]
is not locally connected. An example of a compact set (but not a filled Julia
set) with this property is shown in Figure 16. Evidently, in such a case,
the correspondence \[t\mapsto z_t\] cannot be continuous.\smallskip
\midinsert
\insertRaster 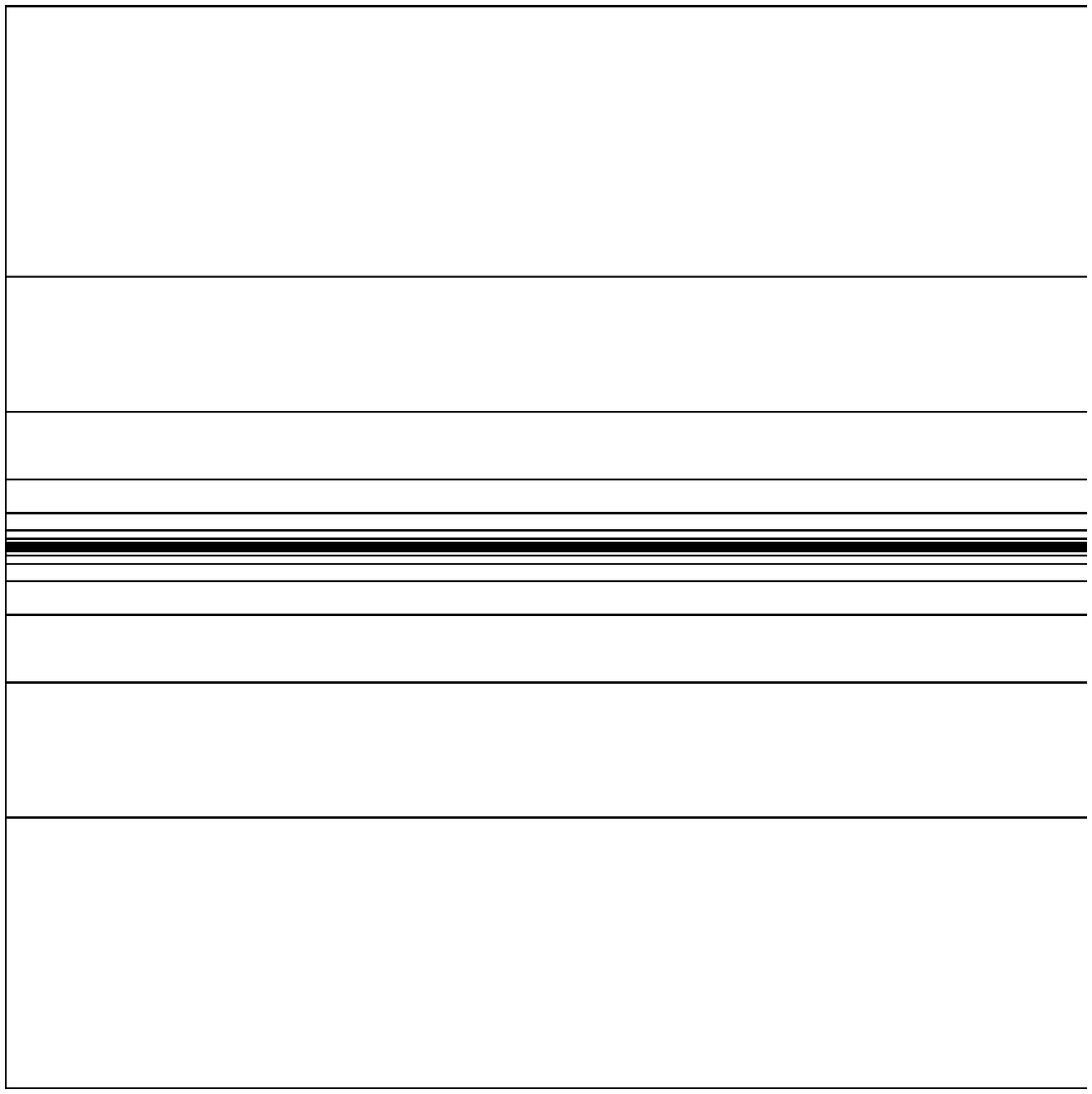 pixels 512 by 513 scaled 225
\bigskip
\centerline{Figure 16. The double comb \[\;[0,2]\times\{0\,,\,\pm 2^{-n}\}\;
\cup\;0\times[-1,1]\]}
\medskip\vfil
\endinsert
%\eject

{\bf Definition.} An external ray \[R_t\] will be called {\bit rational}
if its angle \[t\in\R/\Z\] is rational, and {\bit periodic} if
\[t\] is periodic under multiplication by the degree \[d\],
so that \[\;d^nt\equiv t \;(\mod 1)\] for some \[n\ge 1\].\smallskip

Note that \[R_t\] is eventually periodic
under  multiplication by \[d\] if and only if \[t\] is rational, and
is periodic if and only if the number \[t\] is
rational with denominator relatively prime to \[d\].
(If \[t\] is rational with denominator \[m\], then
the successive images of
\[R_t\] under \[f\] have angles \[dt\,,\,d^2t\,,\,d^3t\,,\,\ldots\;
{\rm (mod\;1
)}\;\;\] with denominators dividing \[m\]. Since there are only finitely
many such fractions modulo 1, this sequence must eventually repeat.\break
% Using ???, it follows that the corresponding sequence
% $$	z_t\buildrel f\over\mapsto z_{dt}\buildrel f\over\mapsto z_{d^2t}
%	\buildrel f\over\mapsto\cdots	$$
%is also eventually periodic. 
In the special case where \[m\] is relatively
prime to \[d\], the fractions with denominator \[m\] are permuted under
multiplication by \[d\] modulo 1, so the landing point \[z_t\] is
actually periodic.)\smallskip

%(Conversely, if \[z_t\] is periodic, then the denominator \[m\]
%must be relatively prime to \[d\]. For otherwise we could find
%two distinct fractions of the form \[d^nt\] modulo 1 which have the same
%image under multiplication by \[d\] modulo 1. Hence we could find
%two distinct external rays landing at
%points on the orbit of \[z_t\], with a common image under \[f\]. But
%this is impossible: Since \[z_t\] is in the Julia set, it cannot belong to
%a superattractive orbit. Therefore, \[f\] is a local diffeomorphism near the
%orbit of \[z_t\].)\smallskip

The following result is due to Sullivan, Douady and Hubbard.
We do not assume local connectivity.

{\QP{\bf 18.1. Theorem.} \it If \[K(f)\] is connected, then every periodic
external ray lands at a periodic point which is either repelling or
parabolic.\medskip}

\noindent The converse result, due to Douady and Yoccoz, is more difficult:

{\QP{\bf 18.2. Theorem.} \it Still assuming that \[K(f)\] is connected,
every repelling or parabolic periodic point is the landing point of at
least one external ray, which is necessarily periodic.\medskip}

In fact we will only prove the parabolic case. For the repelling case, which
is quite similar, the reader is referred to Petersen.
Here is a complementary result.

{\QP{\bf 18.3. Lemma.} \it If a periodic ray lands at \[z_t\], then only
finitely many external rays, all periodic of the same period, can land
at \[z_t\].\medskip}

\noindent The following is an immediate consequence of 18.1.

\midinsert
\vfil
\insertRaster 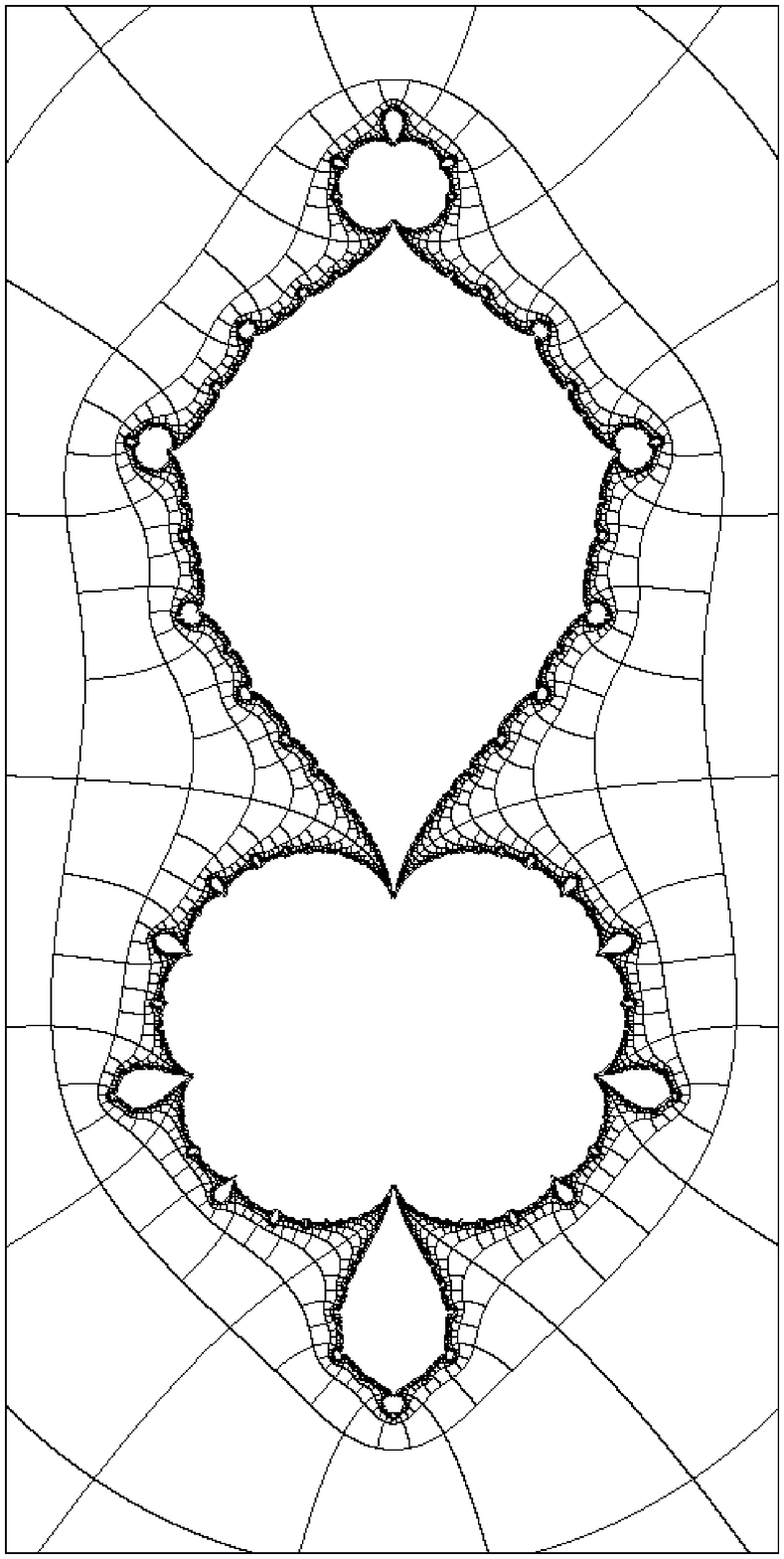 pixels 480 by 960 scaled 350
\bigskip
\centerline{Figure 17. Julia set for \[z\mapsto z^3-iz^2+z\]}
\centerline{with some external rays and equipotentials indicated.}
\vfil
\bigskip
\endinsert

{\QP{\bf 18.4. Corollary.} \it If \[t\] is rational but not periodic,
then \[R_t\] lands at a point \[z_t\] which is eventually periodic
but not periodic.\medskip}

As an example, Figure 17 shows the Julia set for the cubic map \[f(z)=z^3-
iz^2+z\].
Here the rays \[R_0\] and \[R_{1/2}\] are mapped into themselves by \[f\].
Both must land at the parabolic fixed point \[z=0\], since the only other
fixed point, at \[z=i\], is superattracting and hence does not belong
to the Julia set. On the other hand,
the \[1/6\,,\,1/3\,,\, 2/3\] and \[5/6\] rays have denominator divisible
by 3. These four rays land at the two disjoint pre-images of zero.\smallskip

The following result of Sullivan and Douady is closely related.
(Sullivan 1983. For an independent proof, see Lyubich, p. 85.)

{\QP\bf 18.5. Theorem. {\it If the Julia set \[J\] of a polynomial map \[f\]
is locally connected, then every periodic point in \[J\] is either
repelling or parabolic.}\medskip}

Recall from \S8
that a {\bit Cremer point} can be characterized as a periodic point which
belongs to the Julia set but is neither repelling nor parabolic. Thus
the following is a completely equivalent formulation.

{\QP\bf 18.6. Corollary. {\it If \[f\] is a polynomial map with a Cremer
point, then the Julia set \[J(f)\] is not locally connected.}\medskip}

\noindent The proofs begin as follows.

{\QP{\bf 18.7. Lemma.} {\it The ray \[R_t\] lands at a single point \[z_t\]
of the Julia set if and only if \[R_{dt}\] lands at a single point
\[z_{dt}\]. Furthermore, the image \[f(z_t)\] is necessarily equal to
\[z_{dt}\].}\medskip}
\eject

{\bf Proof.} The set \[L(t)\] consisting of all limit points of
\[\psi(re^{2\pi it})\] as \[r\] tends to 1 from above
is certainly compact and connected, with \[f(L(t))=L(dt)
\].  It follows easily that \[L(t)\] is a single point if and only if
\[L(dt)\] is a single point.\QED

Thus if the ray \[R_t\] is periodic, and if \[R_t\] lands at a point \[z_t\],
then it follows that \[z_t\] is a periodic point of \[f\].\smallskip

{\bf Proof of 18.3.}
First consider the special case \[t=0\]. Then \[f\] maps the ray \[R_0\] to
itself, hence \[f(z_0)=z_0\]. If another ray \[R_s\] also lands at \[z_0\],
then we will prove that \[s\] is a rational number of the form \[j/(d-1)\].
This will show that there are
at most \[d-1\] external rays landing at \[z_0\], and that they
are all periodic (and in fact fixed) under multiplication by \[d\], since
the numbers \[s=j/(d-1)\] are exactly the solutions to the congruence
\[s\equiv ds\;({\rm mod}\; 1)\] where \[d\] is the degree.

If a ray \[R_s\] which is not of
this form lands at \[z_0\], then we have \[s\not\equiv ds\;({\rm mod}\; 1)\].
Define a sequence of angles \[0\le s_n<1\] by the condition that \[s_0\equiv
s\] and \[s_{n+1}\equiv ds_n\] modulo 1. By hypothesis, the numbers \[0\,,
\,s_0\] and \[s_1\] are distinct.
 Suppose, to fix our ideas, that
\[0<s_0<s_1<1\]. Since \[f\] is a local diffeomorphism near \[z_0\], it must
preserve the cyclic order of the rays landing at \[z_0\]. Thus the images of
\[0\,,\,s_0\] and \[s_1\] must satisfy \[0<s_1<s_2<1\]. Continuing
inductively, it follows that \[0<s_0<s_1<s_2<\cdots<1\]. Thus the \[s_n\]
must converge to some angle \[s_\infty\], which is necessarily a fixed point
of the map \[s\mapsto ds\;({\rm mod}\; 1)\]. But this is impossible, since
this map has only strictly repelling fixed points.

More generally, if \[R_t\] is any ray which is mapped into itself by \[f\],
so that \[t=j/(d-1)\], then we can translate all angles by \[-t\] and apply
the argument above.
The proof in the general case now follows easily. We can
replace \[f\] by some iterate \[g=f^{\circ n}\], where \[d^ns\equiv s\;({\rm
mod}\; 1)\] so that \[g(R_s)=R_s\]. The argument then proceeds as above.
This proves 18.3.\QED

The proof of 18.1 begins as follows. We will make
use of the \Poincare metric on \[\C\ssm K\]. {\it Note first that the map \[f\] is
a local isometry for this metric.} In fact the
universal covering of \[\C\ssm\bar D\] is isomorphic to the right half-plane
\[\{W=U+iV:U>0\}\] under the exponential map, and the map \[f\] on \[\C\ssm K\]
corresponds to the \[d$-th power map on \[\C\ssm\bar D\], which corresponds
to the automorphism \[W\mapsto dW\] on this right half-plane.\smallskip

We suppose that \[t\] has denominator prime to \[d\], so that some
iterate \[g=f^{\circ n}\] maps \[R_t\] to itself.
This ray \[R_t\] can be expressed as a union
$$	R_t\;=\;\cdots\cup S(-1)\cup S(0)\cup S(1)\cup\cdots	$$
of finite segments, where \[g\] maps
each \[S(j)\] isomorphically onto \[S(j+1)\]. In fact let \[S(j)\] correspond
under \[\phi\] to the set of \[re^{2\pi i t}\] with \[d^{jn}\le\log r\le
d^{(j+1)n}\]. Thus the \[S(j)\] all have the same \Poincare length. But as
\[j\to-\infty\] the segments \[S(j)\] converge towards the boundary \[J\]
of \[\C\ssm K\]. Hence, according to \S2.4, the Euclidean length of \[S(j)\]
must converge to zero. This means that any limit point of the sequence
of sets \[S(j)\], as \[j\to-\infty\], must be a fixed point of \[g\].
But, as noted in the proof of 18.7,
the set of all limit points must be connected.
Since \[g\] has only finitely many fixed points, this proves that the ray
\[R_t\] lands at a single fixed point \[z_t\]
of the map \[g=f^{\circ n}\].\smallskip

Now let us apply the Snail Lemma of \S13.
After translating coordinates so that
our fixed landing point lies at the origin, the map \[g\] and the path
\[x\mapsto \psi(\exp(d^{nx}+2\pi it))\] will satisfies the hypothesis of 13.3.
Hence it follows that this fixed point \[z_t\] of \[g\] is either repelling,
or parabolic with multiplier \[g'(z_t)=1\]. In terms of the original map
\[f\], it follows that \[z_t\] is a periodic point, and that its periodic
orbit is either repelling or has a root of unity as multiplier. This proves
18.1.\QED

{\bf Proof of 18.4.} Consider a rational angle \[t\] with denominator which
is not relatively prime to \[d\]. In this case the multiple \[dt\] will
have smaller denominator,
but the sequence of angles \[d^nt\] modulo 1 must certainly be eventually
periodic. Thus some forward image \[f^{\circ n}R_t\] must have a
well defined landing point. But this implies that \[R_t\] itself has
a well defined landing point by 18.7.\QED

The proof of 18.5 and 18.6 will be based on 18.1, together with the
following. Let \[n\ge 2\] be an integer.

{\QP\bf 18.8. Lemma. {\it Let \[L\subset \R/\Z\] be a compact set which is
mapped homeomorphically onto itself by the map \[t\mapsto nt (\mod 1)\].
Then \[L\] is finite.}\medskip}

{\bf Proof.} In fact we will prove the following more general
statement. Let \[X\] be a compact metric space with distance function
\[\rho(x,y)\], and let \[h:X\to X\] be a homeomorphism which is expanding
in the following sense: There should exist numbers \[\epsilon>0\] and \[k>1\]
so that
$$	\rho(h(x),h(y))\;\ge\;k\rho(x,y)	$$
whenever \[\rho(x,y)<\epsilon\]. Then we will show that \[X\] is finite.
Evidently this hypothesis is satisfied with \[X=L\subset\R/\Z\], so
this argument will prove 18.8.

Since \[h^{-1}\] is uniformly continuous, we can choose \[\delta>0\] so that
\[\rho(x,y)<\epsilon\] whenever \[\rho(h(x),h(y))<\delta\]. But this
implies that \[\rho(x,y)<\delta/k\]. Since \[X\] is compact, we can choose
some finite number, say \[n\], of balls of radius \[\delta\] which cover
\[X\]. Applying \[h^{-p}\], we obtain \[n\] balls of radius \[\delta/k^p\]
which cover \[X\]. Since \[p\] can be arbitrarily large, this proves that
\[X\] can have at most \[n\] distinct points.\QED

{\bf Proof of 18.5.}
Let \[z_0\] be any periodic point. After
replacing the given polynomial map by some iterate, we may assume that
\[z_0\] is fixed by \[f\]. Furthermore, after a linear
change of coordinates, we may assume that \[f\] is monic.
Since \[J(f)\] is locally connected, according to 17.4 the inverse Riemann map
$$	\psi\,:\,\hat\C\ssm\bar D\;\buildrel\approx\over\longrightarrow
	\; \hat\C\ssm K(f)	$$
extends continuously over the boundary, yielding a map from the unit
circle onto \[J(f)\]
which semi-conjugates  the \[d$-th power map on the circle to the map \[f\]
on \[J(f)\]. In other words, identifying the unit circle with \[\R/\Z\],
we have a map
$$	\Psi:\R/\Z\;\to\; J(f)	$$
which satisfies \[\Psi(dt)=f(\Psi(t))\], where \[d\] is the degree. Let
\[X\] be the compact set \[\Psi^{-1}(z_0)\subset\R/\Z\]. Since \[z_0\]
is fixed by \[f\], it follows immediately that the map\break \[t\mapsto td
(\mod 1)\] carries \[X\] homeomorphically onto itself. Hence \[X\] is
finite by 18.8, and the conclusion follows.\QED

We will prove the parabolic case of Theorem 18.2 in the following slightly
sharper form.
For the proof in the repelling case, the reader is referred to Petersen.

{\QP{\bf 18.9. Theorem.} \it Let \[f\] be a polynomial of degree \[d\ge 2\]
with  \[K(f)\]  connected,
and let  \[z_0\]  be a parabolic fixed point
with multiplier \[\lambda=1\]. Then for each repelling petal
\[P\]  at  \[z_0\]  there exists at least one
external ray \[R_t\] which lands at \[z_0\] through \[P\], 
where the angle \[t\] is a rational number of the form \[m/(d-1)\],
so that \[R_t\] is mapped onto itself by \[f\].\medskip}

Very roughly, we will show that
the basin \[\C\ssm K\] of the point at infinity corresponds to a union
of one or more annuli in the \'Ecalle cylinder \[P/f\]. Each of these
annuli contains a unique simple closed geodesic \[\Gamma\] (think of
placing a rubber band around a napkin ring); and \[\Gamma\]
will lift to the required external ray in \[\C\ssm K\]. Before giving
details, let us look at two examples.
\smallskip

{\bf Example 1.} Consider the map \[f(z)=z^2+\exp(2\pi i\cdot 3/7)z\] of
Figure 8. Here the seventh iterate is a map of degree \[2^7=128\] with
seven repelling petals about the origin. Hence there must be at least
seven external rays landing at the origin, and their angles must be
rational numbers of the form \[m/127\], so as to map into themselves
under multiplication by 128 modulo 1. In fact a little experimentation
shows that only the ray with angle \[21/127\] and its successive iterates
under doubling modulo 1 will fit in the right order around the origin.
(Compare Goldberg.) Thus there are just seven rays which land at zero,
one in each repelling petal. The numerators of the corresponding angles
are \[21\,,\,37\,,\, 41\,,\, 42\,,\,74\,,\,
82\,,\,84\].\smallskip

{\bf Example 2.} Now consider the cubic map \[f(z)
=z^3-iz^2+z\] of Figure 17. Evidently there is only one repelling petal
at the parabolic fixed point \[z=0\]. Yet we saw above (following 18.4)
that the two rays \[R_0\]
and \[R_{1/2}\] must both land at this point.\smallskip

The proof of 18.9 begins as follows.
According to \S7.7, there is an essentially unique map \[\beta=\alpha^{-1}\]
which is defined and univalent in some left half-plane \[{\cal R}(w) < c\],
taking
values in the repelling petal \[ P\], and which satisfies the Abel equation
\[\beta(w+1) = f(\beta(w))\] or equivalently:
$$	\beta(w) = f(\beta(w-1)) \,. $$
As in 7.11, we can extend to an analytic map from  \[\C\]  to  \[\C\]
satisfying this same equation
by setting  \[\;\beta(w) = f^{\circ n}(\beta(w-n))\;\] for large  \[n\].
{\bf Definition:} Let  \[K_P\]  be the inverse
image of  \[K = K(f)\]  under  \[\beta\]. Since \[K\] is closed and
\[f$-invariant, it follows that \[K_P\]  is closed and periodic:
\[K_P = K_P+1\]. Furthermore, it is not hard to show that \[K_P\] contains
all points \[w=u+iv\] with \[|v|\] large. For points which are
far from the real axis must correspond, under \[\beta\], to points
which belong to one of the two neighboring attracting petals of \[f\].

{\QP{\bf 18.10. Yoccoz Lemma.} \it Each component \[U\] of \[C\ssm K_P\]
is a universal covering of \[C\ssm K\] with projection map \[\beta\].\medskip}
\eject

{\bf Proof.} Consider the topological space \[E\] consisting of all
left infinite sequences \[{\bf w}=(\,\ldots\,,\,w_{-2}\,,\,w_{-1}\,,\,w_0)\] of points in the
Fatou set \[\C\ssm J\] satisfying
$$	\cdots\buildrel f\over\mapsto w_{-2}\buildrel f\over\mapsto w_{-1}
	\buildrel f\over\mapsto w_0	$$
and \[\;\lim_{\,n\to-\infty}\, w_n\,=\,z_0\]. We can give each component of
\[E\] the structure of a Riemann surface in such a way that the projection \[
\pi:(\,\ldots\,,\, w_{-1}\,,\,w_0)\mapsto w_0\] is a local
conformal isomorphism.
In fact, for any such \[{\bf w}\], we can choose
\[N<0\] so that \[w_n\] belongs to some fixed repelling petal \[P\] for \[n
\le N\]. The coordinate \[w_N\in P\] then serves as a local uniformizing
parameter for \[E\], and \[\pi\] corresponds to the map \[f^{|N|}:w_N\mapsto
w_0\]. (It is important here that there are no critical points
in the finite part of the Fatou set.)

To show that \[\pi:E\to \C\ssm K\] is actually a covering map, we
start with any simply connected neighborhood \[V\] of \[\pi({\bf w})=w_0\] in
\[\C\ssm K\]. Then we can find a unique single valued branch \[f_0^{-n}\]
of \[f^{-n}\] on \[V\] which maps \[w_0\] to \[w_{-n}\]. These maps
\[f_0^{-n}\] must constitute a normal family on \[V\], since they take
values in \[\C\ssm K\]. Since
the sequence \[f_0^{-n}\] converges uniformly to \[z_0\] throughout some
small neighborhood of \[w_0\], it follows easily that this sequence
converges to \[z_0\] throughout the given neighborhood \[V\]. Thus for each
\[w\in V\] the sequence \[\;{\bf w}(w)=\big(\,\cdots\,,\, f_0^{-2}(w)\,,\,
f_0^{-1}(w)\,,\,w\big)\;\] belongs to \[E\], hence \[V\] is evenly covered
under the projection \[\pi\].

It is not difficult to check
that \[\C\ssm K_P\] embeds naturally as an open and closed subset of
our covering space \[E\]. But each component \[U\] of \[\C\ssm K_P\] is simply
connected. For otherwise, there would be a compact component of \[\C\ssm U\]
which, after left translation by a large integer, would embed into a compact
of \[K\] in the petal \[P\]. This is impossible, since \[K\] is connected.
Thus each such \[U\] is a universal covering surface of \[\C\ssm K\].\QED

{\bf Remark.} Although \[U\] is a universal covering of \[\C\ssm K\], the proof
does not tell us just how the free cyclic group of deck transformations acts
on \[U\]. We will show in 18.12
that there is also a natural action of the free cyclic
group \[\Z\] by integer translations of \[U\], but that action
is quite different, and has a quotient surface \[U/\Z\] which is an annulus,
naturally embedded in the \'Ecalle cylinder \[P/f\].\smallskip

Let \[\hat G=G\circ\beta\] be the potential function on \[\C\ssm K\] lifted to \[\C\ssm K_P\], so
that\break
\[\hat G(w+1) = \hat G(w)d\] where \[d\] is the degree. Let \[\hat G_0\] be the maximum
value of \[\hat G\] on the imaginary axis, attained say at \[w_0\]. Then
\[\hat G(w)\le G_0\] for \[{\cal R}(w)\le 0\] by the maximum modulus principle
for harmonic maps, hence
$$	 G(w) \le d^n\hat G_0 \qquad{\rm  whenever}\qquad   {\cal R}(w)
 \le n .\eqno	(*)$$

{\QP{\bf 18.11. Lemma.} \it The real part of \[w\] must take values tending
to \[\infty\] within each component \[U\].\medskip}

\noindent This follows, since \[\hat G\] must be unbounded on each \[U\]
by 18.10.\QED

{\QP{\bf 18.12. Main Lemma.} \it Each such component \[U\] is periodic:
\[\;U = U+1\].\medskip}

{\bf Proof.}
Otherwise the translates \[U-n\] would have to be pairwise disjoint. First
consider the component \[U\] on which the maximum \[\hat G_0\] is attained.
Let \[k_n\] be the vertical width of the intersection of \[U-n\] with the
imaginary axis. Then the sum of \[k_n\] is finite, so \[k_n\] tends to
zero as \[n \to +\infty\]. But \[k_n\] is also the vertical width of the
intersection of \[U\] with the line \[\;{\cal R}(w) = n\,\]. {\it Since these
numbers tend to zero, the \Poincare distance of the imaginary axis from
\[\;{\cal R}(w) = n\;\] within \[U\] must grow more than linearly with \[n\].}
\smallskip

On the other hand, an orthogonal trajectory of the curves \[\;\hat G
= {\rm constant}\,\],
of \Poincare arclength \[n \log(d)\] within \[U\], will get from \[w_0\] on the
imaginary axis to a point \[w_n\] with
\[\hat G(w_n) = (d^n)\hat G_0\]. By $(*)$ above, the real
part of \[w_n\] is \[\ge n\]. Thus this distance grows at most linearly
with \[n\]. This contradiction proves the Main Lemma 18.12 when \[w_0\] belongs
to \[U\]. To get a proof for arbitrary \[U\], we must simply use the union of
translates \[U-n\] in place of the entire open set \[C\ssm K_P\] in the
argument above. This proves 18.12.\QED

{\bf Proof of Theorem 18.9.}
We know that each component \[U\] of \[\C\ssm K_P\]
is a universal covering of \[\C\ssm K\], and that the unit translation \[w\mapsto
w+1\] of \[U\] corresponds to the map \[f\] on \[\C\ssm K\].
Furthermore, it follows easily from Corollary 7.5 that the complement
\[\C\ssm K_P\] is non-vacuous, so that there exists at least one such component
\[U\]. Since the imaginary
coordinate is bounded throughout \[U\], the quotient \[U/\Z\] is clearly
an annulus. Hence it has a unique simple closed geodesic with respect to
the Poincar\'e metric. (Compare Problem 2-3.
We can model the annulus as the upper half-plane
\[H\] modulo the identification \[w\sim t_0w\] for some fixed \[t_0>1\]. The
imaginary axis in \[H\] then covers the unique closed geodesic.)\smallskip

This closed geodesic \[\Gamma\subset U/\Z\] lifts to an infinite
\[\Z$-invariant geodesic \[\tilde\Gamma\subset U\], and hence to an
infinite \[f$-invariant geodesic \[\Gamma'=\beta(\tilde\Gamma)\] in
\[\C\ssm K\]. In the negative direction, this geodesic \[\Gamma'\]
converges to the parabolic
fixed point \[z_0\]. For \[\beta\] was constructed in such a way that for
any compact subset \[L\subset U\] the images \[\beta(L-n)\] converge to \[z
_0\] as \[n\to\infty\]. On the other hand, in the positive direction
\[\Gamma'\] converges to the point at infinity, since every compact subset
\[L\subset U\] has the property that the images \[\beta(L+n)=f^{\circ n}
\circ\beta(L)\] converge to \[\infty\] as \[n\to\infty\]. But a \Poincare
geodesic in \[\C\ssm K\] which leads to the point at infinity is necessarily
an external ray.\QED
\vfill

%\rightline{jwm, version of 9-4-91}
\end

%% file: cda-c.tex
\input header.tex
\footline={\hss\tenrm A-\folio\hss}\pageno=1

\centerline{\bf Appendix A. Theorems from Classical Analysis.}\bigskip

This appendix will describe some miscellaneous theorems from classical
complex variable theory. We first complete the arguments from \S8.6
and \S15.3
by proving Jensen's inequality and the Riesz brothers' theorem.
We then describe results from the theory of univalent functions,
due to Gronwall and Bieberbach, in order to prove
the Koebe Quarter Theorem for use in Appendix G. (By definition, a function
of one complex variable is called {\bit univalent} if it is holomorphic
and injective.)\smallskip

We begin with a discussion of Jensen's inequality. (J. L. W. V. Jensen was
the president of the Danish Telephone Company, and a noted amateur
mathematician.) Let \[f:D\to\C\] be a
holomorphic function on the open disk which is not identically zero.
Given any radius \[0<r<1\], we can form the average
$$	A(f,r)\;=\; {1\over 2\pi}\int_0^{2\pi}\; \log|f(re^{i\theta})| d\theta
	$$
of the quantity \[\log|f(z)|\] over the circle \[|z|=r\].

{\QP{\bf A.1. Jensen's Inequality.} {\it This average \[A(f,r)\] is monotone
increasing $($and also convex upwards$)$ as a function of \[r\]. Hence
\[A(f,r)\]
either converges to a finite limit or diverges to \[+\infty\] as \[r\to 1\].
}\medskip}

\noindent In fact the proof will show something much more precise.

{\QP {\bf A.2. Lemma.} {\it If we consider \[A(f,r)\] as a function of
\[\log r\], then it is piecewise linear, with slope \[dA(f,r)/d\log r\]
equal to the number of roots of \[f\] inside the disk \[D_r\] of radius
\[r\], where each root is to be counted with its appropriate multiplicity.}
\medskip}

In particular, the function \[A(f,r)\] is determined, up to an additive
constant, by the location of the roots of \[f\]. To prove this Lemma,
note first that we can write \[d\theta=dz/iz\] around any loop \[|z|=r\].
Consider an annulus \[{\cal A}\;=\;\{z: r_0< |z|< r_1\}\] which
contains no zeros of \[f\]. According to the ``Argument Principle'',
the integral
$$	n\;=\; {1\over 2\pi i}\oint_{|z|=r}d\log f(z)\;=
	\;{1\over 2\pi i}\oint{f'(z)dz\over f(z)}	$$
measures the number of zeros of \[f\] inside the disk \[D_r\]. It follows
that the difference \[\log f(z)-\log z^n\] can be defined as a single valued
holomorphic function throughout this annulus \[{\cal A}\]. Therefore,
the integral of \[(\log f(z) -\log z^n)dz/iz\]
around a loop \[|z|=r\] must be independent of \[r\], as long
as \[r_0<r<r_1\]. Taking the real part, it follows that the difference
\[A(f,r)-A(z^n,r)\] is a constant, independent of \[r\]. Since \[A(z^n,r)
=n\log r\], this proves that the function \[\;\log r\mapsto A(f,r)\;\]
is linear with slope \[n\] for \[r_0<r<r_1\].

Finally, note that the average \[A(f,r)\] takes a well defined finite
value even when \[f\] has one or more zeros on the circle \[|z|=r\], since
the singularity of \[\log |f(z)|\] at a zero of \[f\] is relatively mild.
Continuity of \[A(f,r)\] as \[r\] varies through such a singularity is
not difficult, and will be left to the reader.\QED

{\QP{\bf A.3. Theorem of F. and M. Riesz.} {\it Suppose that \[f:D\to\C\] is
bounded and holomorphic on the open unit disk. If the radial limit
$$	\lim_{r\to 1}\; f(re^{i\theta})	$$
exists and is equal to zero for \[\theta\] belonging to a set
\[E\subset [0,\,2\pi]\] of positive Lebesgue measure, then
\[f\] must be identically zero.}\medskip}

(Compare the discussion of Fatou's Theorem in \S15.3.)\smallskip

{\bf Proof.} Let \[E(\epsilon,\delta)\] be the measurable set consisting
of all \[\theta\in E\] such that
$$	|f(re^{i\theta})|<\epsilon\qquad {\rm whenever}\qquad 1-\delta<r<1
	\;. $$
Evidently, for each fixed \[\epsilon\], the union of the
nested family of sets \[E(\epsilon, \delta)\] contains \[E\]. Therefore,
the Lebesgue measure \[\ell(E(\epsilon,\delta))\] must tend to a
limit which is \[\ge\ell(E)\]
as \[\delta\to 0\]. In particular, given \[\epsilon\] we can choose \[\delta\]
so that \[\ell(E(\epsilon,\delta))>\ell(E)/2\]. Now consider the average
\[A(f,r)\] of Jensen's Inequality, where \[r>1-\delta\]. Multiplying \[f\]
by a constant if necessary, we may assume that \[|f(z)|<1\] for all \[z\in
D\]. Thus
the expression \[\log |f(re^{i\theta})|\] is less than or equal to zero
everywhere, and less than or equal to \[\log\epsilon\] throughout a set of
measure at least \[\ell(E)/2\]. This proves that
$$	2\pi A(f,r)\;<\; \log(\epsilon)\,\ell(E)/2	$$
whenever \[r\] is sufficiently close to 1. Since \[\epsilon\] can be
arbitrarily small, this implies that \[\lim_{r\to 1} A(f,r)\;=\;-\infty\],
which contradicts A.1 unless \[f\] is identically zero.\QED\medskip

Now consider the following situation. Let \[K\] be a compact connected
subset of \[\C\], and suppose that the complement \[\C\ssm K\] is conformally
diffeomorphic to the complement \[\C\ssm \bar D\].

{\QP {\bf A.4 Gronwall Area Inequality.} {\it Let \[\;\phi:\C\ssm \bar D \;\to
\;\C\ssm K\;\] be a conformal isomorphism and let
$$	\phi(w)\;=\;b_1w+b_0+b_{-1}/w+b_{-2}/w^2+\cdots	$$
be its Laurent expansion. Then \[|b_1|\ge|b_{-1}|\], with equality if and
only if \[K\] is a straight line segment.}\medskip}

{\bf Proof.} For any \[r>1\] consider the image under \[\phi\] of the circle
\[|w|=r\]. This will be some embedded circle in \[\C\] which encloses a
region of area say \[A(r)\]. We can compute this area by Green's Theorem,
as follows. Let \[\phi(re^{i\theta})=z=x+iy\]. Then
$$	A(r)\;=\;\oint xdy\;=\;-\oint ydx\;=\;{\textstyle{1\over 2i}}
	\oint \bar z dz\,, $$
to be integrated around the image of \[|w|=r\].
Substituting the Laurent series\break \[z\,=\,\sum_{n\le 1} b_nw^n\], with
\[w=r^ne^{ni\theta}\], this yields
\eject
$$	A(r)\;=\;{\textstyle{1\over 2}}\sum_{m,n\le 1}\; n\bar b_m b_n
	r^{m+n}\oint e^{(n-m)i\theta}d\theta	\,.	$$
Since the integral equals \[2\pi\] if \[m=n\], and is zero otherwise,
we obtain
$$	A(r)\;=\;\pi\sum_{n\le 1} n|b_n|^2r^{2n}\,.	$$
Therefore, taking the limit as \[r\to 1\], we obtain the simple formula
$$	A(1)\;=\; |b_1|^2-|b_{-1}|^2-2|b_{-2}|^2-3|b_{-3}|^2-\,\cdots
	\eqno (A.5) $$
for the area (that is the 2-dimensional Lebesgue measure) of the compact
set \[K\]. Evidently it follows that \[|b_1|\ge|b_{-1}|\]. Furthermore,
if equality holds then all of the remaining coefficients must be zero:
\[b_{-2}=b_{-3}=\cdots=0\]. After a rotation of the \[w\] coordinate
and a linear transformation of the \[z\] coordinate, the Laurent series
will reduce to
the simple formula \[\;z=w+w^{-1}\]. As noted in \S5, Example 2,
this transformation
carries \[\C\ssm \bar D\] diffeomorphically onto the complement of the interval
\[[-2,2]\].\QED\smallskip

Now consider an open set \[U\subset\C\] which contains the origin and
is conformally isomorphic to the open disk.

{\QP{\bf A.6. Bieberbach Theorem.} {\it If \[\psi: D\to U\] is a
conformal isomorphism with power series expansion \[\psi(\eta)\;=\;\sum_{n\ge 1}
a_n\eta^n\], then \[|a_2|\le 2|a_1|\], with equality if and only if \[\C\ssm U\]
is a closed half-line pointing towards the origin.}\medskip}

{\bf Remark.} The Bieberbach conjecture, recently proved by DeBrange,
asserts that \[|a_n|\le n|a_1|\] for all \[n\ge 2\]. Again, equality
holds if \[\C\ssm U\] is a closed half-line pointing towards the origin.\smallskip

{\bf Proof of A.6.} After composing \[\psi\] with a linear transformation,
we may assume that \[a_1=1\]. Let us set \[\eta=1/w^2\],
so that each point \[\eta\ne 0\] in \[D\] corresponds to two points \[w\in
\C\ssm \bar D\]. Similarly, set \[\psi(\eta)=\zeta=1/z^2\], so that each \[\zeta\ne 0\]
in \[U\] corresponds to two points \[z\] in some centrally symmetric
neighborhood \[N\]
of infinity. A brief computation shows that \[\psi\] corresponds to a
Laurent series
$$	w\mapsto z=1/\sqrt{\psi(1/w^2)}\;=\;w-{\textstyle{1\over 2}}a_2 /w
	+({\rm higher\; terms})	$$
which maps \[\C\ssm \bar D\] diffeomorphically onto \[N\].
Thus \[|a_2|\le 2\] by Gronwall's Inequality, with equality if and
only if \[N\] is the complement of a line segment, necessarily centered
at the origin. Expressing this condition
on \[N\] in terms of the coordinate \[\zeta=1/z^2\], we see that equality
holds if and only if \[U\] is the complement of a half-line pointing
towards the origin.\QED
\eject

{\QP {\bf A.7. Koebe-Bieberbach Quarter Theorem.} {\it
Again suppose that the map
$$	\eta\;\mapsto\;\psi(\eta)\;=\;a_1\eta+a_2\eta^2+\cdots $$
carries the unit disk \[D\] diffeomorphically onto an open set \[U\subset\C\].
Then the distance \[r\] between the origin and the boundary of \[U\]
satisfies
$$	{\textstyle{1\over 4}}|a_1|\;\le \;r\;\le|a_1|\,.	$$
Here the first equality holds if and only if \[\C\ssm U\]
is a half-line pointing towards the origin, and the second equality
holds if and only if \[U\] is a disk centered at the
origin.}\medskip}

\noindent
In particular, in the special case \[a_1=1\] the open set \[U\]
necessarily contains the disk of radius \[1/ 4\] centered at the
origin. The left hand inequality was conjectured and partially proved
by Koebe, and later completely proved by Bieberbach. The right hand
inequality is an easy consequence of the Schwarz Lemma.

Here is an interesting restatement of the Quarter Theorem. Let \[ds=\gamma(z)
|dz|\] be the \Poincare metric on the open set \[U\], and let \[r=r(z)\]
be the distance from \[z\] to the boundary of \[U\].

{\QP{\bf A.8. Corollary.} {\it If \[U\subset\C\] is simply connected, then
the \Poincare metric \[ds=\gamma(z)|dz|\] on \[U\] agrees with the
metric \[|dz|/r(z)\] up to a factor of two in either direction.
That is
$$	{1\over 2r(z)}\;\le\;\gamma(z)\;\le\;{2\over r(z)}	$$
for all \[z\in U\]. Again, the left equality holds if and only if \[\C\ssm U\]
is a half-line pointing towards the point \[z\in U\], and the right equality
holds if and only if \[U\] is a round disk centered at \[z\].}\medskip}

\noindent As an example, if \[U\] is a half-plane, then the \Poincare
metric precisely agrees with the \[1/r\] metric \[|dz|/r\].\smallskip

{\bf Proof of A.7.} Without loss of generality, we may assume that \[a_1=1\].
If \[z_0\in \partial U\] be a boundary point with
minimal distance \[r\] from the origin, then we must prove that \[{1\over 4}
\le r\le 1\].  We will compose \[\psi\] with the
linear fractional transformation \[z\mapsto z/(1-z/z_0)\] which maps \[z_0\]
to infinity. Then the composition  has the form
$$	\eta\;\mapsto\; \psi(\eta)/(1-\psi(\eta)/z_0)\;=\;
	\eta+(a_2+1/z_0)\eta^2+\cdots\,.	$$
By Bieberbach's Theorem we have \[|a_2|\le 2\] and \[|a_2+1/z_0|\le 2\],
hence \[|1/z_0|=1/r\le 4\], or \[r\ge 1/4\]. Here equality holds only
if \[|a_2|=2\] and \[|a_2+1/z_0|=2\]. The exact description of \[U\] then
follows easily.

Now suppose that \[r\ge 1\]. Then the inverse mapping \[\psi^{-1}\] is
defined and holomorphic throughout the unit disk \[D\], and takes values in
\[D\]. Since its derivative at zero is \[1\], it follows from the
Schwarz Inequality that \[\psi\] is the identity map.\QED
\vfil\eject % \bigskip\bigskip

\footline={\hss\tenrm B-\folio\hss}\pageno=1
\def\md{{\rm mod}}
\def\ar{{\rm area}}
\def\L{{\bf L}}
\def\T{{\bf T}}

\centerline{\bf Appendix B. Length-Area-Modulus Inequalities.}\bigskip

The most basic length-area inequality is the following. Let \[I^2\subset\C\]
be the open unit square consisting of all \[z=x+iy\] with \[0<x<1\] and
\[0<y<1\]. By a {\bit conformal metric~} on \[I^2\] we mean a metric
of the form
$$	ds\= \rho(z)|dz| $$
where \[z\mapsto \rho(z)>0\] is any strictly positive continuous real
valued function on the open square. In terms of such a metric, the
{\bit length~} of a smooth curve \[\gamma:(a,b)\to I^2\] is defined to be the
integral
$$	{\bf L}_\rho(\gamma)\= \int_a^b \rho(\gamma(t))|d\gamma(t)|\;, $$
and the {\bit area~} of a region \[U\subset I^2\] is defined to be
$$	\ar_\rho(U)\= \int\!\int_U\, \rho(x+iy)^2dx\,dy\;. $$
In the special case of the Euclidean metric \[ds=|dz|\], with \[\rho(z)\]
identically equal to 1, the subscript \[\rho\] will be omitted.

{\QP{\bf Theorem B.1.} \it If the integral \[\ar_\rho(I^2)\] over the entire
square is finite, then for
Lebesgue almost every \[y\in(0,1)\] the length \[\L_\rho(\gamma_y)\] of the
horizontal line \[\gamma_y:t\mapsto
(t,y)\] at height \[y\] is finite. Furthermore, there exists
\[y\] so that
$$	\L_\rho(\gamma_y)^2\;\le\;\ar_\rho(I^2) \;.\eqno (1) $$
In fact, the set consisting of all \[y\in(0,1)\] for which this inequality
is satisfied has positive Lebesgue measure.\smallskip}

{\bf Remark 1.} Evidently this inequality is best posible. For in the
case of the Euclidean metric \[ds=|dz|\] we have
$$	{\bf L}(\gamma_y)^2\;=\;\ar(I^2)\;=\;1\;. $$

{\bf Remark 2.} It is essential here that we use a square, rather than
a rectangle. If we consider instead a rectangle \[R\] with base
\[\Delta x\] and height \[\Delta y\], then the corresponding inequality
would be
$$	\L_\rho(\gamma_y)^2\;\le\;{\Delta x\over\Delta y}\ar_\rho(R)\eqno (2) $$
for a set of \[y\] with positive measure.\smallskip

{\bf Proof of B.1.} 
We use the Schwarz inequality
$$	\left(\int_a^b f(x)g(x)\,dx\right)^2\;\le\;\left(\int_a^b f(x)^2\,
dx\right)\cdot\left(\int_a^b g(x)^2\,dx\right)\;, $$
which says (after taking a square root) that the inner product of any
two vectors in the Euclidean vector space
of square integrable real functions on an interval is less than or
equal to the product of their norms.
We may as well consider the more general case of a
rectangle\break \[R=(0,\,\Delta x)\times(0,\,\Delta y)\].
Taking \[f(x)\equiv 1\] and \[g(x)=
\rho(x,y)\] for some fixed \[y\], we obtain
$$	\bigg(\int_0^{\Delta x} \rho(x,y)\,dx\bigg)^2\;\le\; \Delta x
 \int_0^{\Delta x} \rho(x,y)^2\,dx\;, $$
or in other words
%for each \[y\in(0\,,\,\Delta y)\]. In other words
$$	L_\rho(\gamma_y)^2\;\le \; \Delta x\int_0^{\Delta x} \rho(x,y)^2dx\;, $$
%where \[\gamma_y\] is the horizontal curve \[x\mapsto(x,y)\].
for each constant height \[y\].
Integrating this inequality over the
interval \[0<y<{\Delta y}\] and then dividing by \[\Delta y\], we get
$$	{1\over\Delta y}\int_0^{\Delta y} L_\rho(\gamma_y)^2dy\;\le\;
	{\Delta x\over\Delta y}\; {\rm area}_\rho(A)\;.\eqno (3) $$
In other words, the {\bit average~} over all \[y\] in the interval
\[(0,\,\Delta y)\] of \[L_\rho(\gamma_y)^2\] is less than or equal
to \[{\Delta x\over\Delta y}\; {\rm area}_\rho(A)\,\]. Further details of the
proof are straightforward.\QED\smallskip

Now let us form a cylinder \[C\] of circumference \[\Delta x\] and height
\[\Delta y\] by gluing the left and right edges of our
rectangle together. More precisely, let \[C\] by the quotient space which
is obtained from the infinitely wide strip \[0<y<\Delta y\] in the \[z$-plane
by identifying each point \[z=x+iy\] with its translate \[z+\Delta x\].
Define the {\bit modulus\/} \[\md(C)\] of such a cylinder to be the ratio
\[\Delta y/\Delta x\] of height to circumference.
By the {\bit winding number\/} of a closed curve \[\gamma\] in \[C\] we
mean the integer
$$	w\={1\over\Delta x}\;\oint_\gamma dx\;. $$

{\QP{\bf Theorem B.2 (Length-Area Inequality for Cylinders).}
\it For any conformal metric
\[\rho(z)|dz|\] on the cylinder \[C\] there exists some simple closed
curve \[\gamma\] with winding number \[+1\] whose length
\[\;\L_\rho(\gamma)=\oint_\gamma \rho(z)|dz|\;\] satisfies the inequality
$$	\L_\rho(\gamma)^2\;\le \; \ar_\rho(A)/\md(A)\;. \eqno (4) $$
Furthermore, this result is best possible: If we use the Euclidean metric
\[|dz|\] then
$$      \L(\gamma)^2\;\ge\;\ar(A)/\md(A) \eqno (5)
$$
for {\bit every} such curve \[\gamma\].\smallskip}

\noindent {\bf Proof.} Just as in the proof of B.1, we find a horizontal curve
\[\gamma_y\] with
$$ \L_\rho(\gamma_y)^2\;\le\;{\Delta x\over\Delta y}\ar_\rho(C)\;=\;
{\ar_\rho(C)\over\md(C)}\;. $$
On the other hand in the Euclidean case, for any
closed curve \[\gamma\] of winding number one we have
$$ \L(\gamma)\;=\;\oint_\gamma |dz|\;\ge\;\oint_\gamma dx\;=\;\Delta x\;, $$
hence \[\L(\gamma)^2\;\ge\; (\Delta x)^2\;=\;\ar(C)/\md(C)\].\QED

{\bf Definitions.}
A Riemann surface \[A\] is said to be an {\bit annulus\/} if it is conformally
isomorphic to some cylinder. An embedded annulus \[A\subset C\] is said
to be {\bit essentially embedded~} if it contains a curve which has winding
number one around \[C\].

\noindent Here is an important consequence
of Theorem B.2.

{\QP{\bf Corollary B.3 (An Area-Modulus Inequality).} \it Let \[A\subset C\]
be an essentially embedded annulus in the cylinder \[C\], and suppose that
\[A\] is conformally isomorphic to a cylinder \[C'\]. Then
$$	\md(C')\;\le\;{\ar(A)\over\ar(C)}\,\md(C)\;. \eqno (6)$$
In particular:
$$	\md(C')\;\le\;\md(C)\;. \eqno (7) $$
\smallskip}

{\bf Proof.} Let \[\zeta\mapsto z\] be the embedding of \[C'\] onto
\[A\subset C\]. The Euclidean metric \[|dz|\] on \[C\], restricted to \[A\],
pulls back to some conformal metric \[\rho(\zeta)|d\zeta|\]
on \[C'\], where \[\rho(\zeta)=|dz/d\zeta|\]. According to B.2, there
exists a curve \[\gamma'\] with winding number \[1\] about \[C'\] whose length
satisfies
$$	L_\rho(\gamma')^2\le\ar_\rho(C')/\md(C') \;. $$
%=\ar(A)/\md(C')\;. $$
This length coincides with the Euclidean length \[\L(\gamma)\] of the
corresponding curve \[\gamma\] in \[A\subset C\], and \[\ar_\rho(C')\]
is equal to the Euclidean area
\[\ar(A)\], so we can write this inequality as
$$	L(\gamma)^2\le\ar(A)/\md(C') \;. $$
But according to (5) we have
$$	\ar(C)/\md(C)\;\le\;L(\gamma)^2\;. $$
Combining these two inequalities, we obtain
$$ \ar(C)/\md(C)\;\le\;\ar(A)/\md(C')\;, $$
which is equivalent to the required inequality (6).\QED

{\QP{\bf Corollary B.4.} \it The modulus of a cylinder is a well defined
conformal invariant.\smallskip}

%\noindent Hence the {\bit modulus~} of an annulus can be defined as the modulus
%of any conformally isomorphic cylinder.\smallskip

{\bf Proof.} If \[C'\] is conformally isomorphic to \[C\]
then (7) asserts
that\break \[\md(C')\le\md(C)\], and similarly \[\md(C)\le\md(C')\].\QED

It follows that the {\bit modulus~} of an annulus \[A\]
can be defined as the modulus of any conformally isomorphic cylinder.
Furthermore, if \[A\] is essentially embedded in some other annulus
\[A'\], then \[\md(A)\le\md(A')\].

{\QP{\bf Corollary B.5 (Gr\"otzsch Inequality).} \it Suppose that \[A'\subset A\]
and \[A''\subset A\] are two disjoint annuli, each essentailly embedded
in \[A\]. Then 
$$ \md(A')+\md(A'')\;\le\;\md(A)\;. $$
\smallskip}

{\bf Proof.} We may assume that \[A\] is a cylinder \[C\].
According to (6) we have
$$	\md(A')\;\le \;{\ar(A')\over\ar(C)}\,\md(C)\;,\qquad
      \md(A'')\;\le \;{\ar(A'')\over\ar(C)}\,\md(C)\;. $$
where all areas are Euclidean. Using the inequality
$$	\ar(A')+\ar(A'')\le \ar(C)\;, $$
the conclusion follows.\QED\smallskip

Now consider a {\bit flat torus~} \[{\bf T}=\C/\Lambda\]. Here \[\Lambda
 \subset \C\] is to be a 2-dimensional {\bit lattice\/}, that is
an additive subgroup of the complex numbers, spanned by
two elements \[\lambda_1\] and \[\lambda_2\] where \[\lambda_1/\lambda_2
\not\in\R\]. Let \[A\subset\T\] be an embedded annulus.

By the ``{\bit winding number\/}'' of \[A\] in \[\T\] we will mean the lattice
element \[w\in\Lambda\] which is constructed as follows. Under the
universal covering map \[\C\to\T\], the central curve of \[A\] lifts
to a curve segment which joins some point \[z_0\in\C\] to a translate
\[z_0+w\] by the required lattice element. We say that \[A\subset\T\] is an
{\bit essentially embedded annulus~} if \[w \ne 0\].

{\QP{\bf Corollary B.6 (Bers Inequality).} \it If the annulus \[A\]
is embedded in the flat torus \[\T=\C/\Lambda\] with winding number
\[w \in\Lambda\], then
$$	\md(A)\;\le\;{\ar(T)\over|w |^2}\;. \eqno (8) $$
\smallskip}

\noindent Roughly speaking, if \[A\] winds many times around the torus,
so that \[|w |\] is large, then \[A\] must be very skinny. A slightly
sharper version of this inequality is given in Problem B-2 below.
\smallskip

{\bf Proof.} Choose a cylinder \[C'\] which is conformally isomorphic
to \[A\]. The Euclidean metric \[|dz|\] on \[A\subset\T\] corresponds to some
metric \[\rho(\zeta)|d\zeta|\] on \[C'\], with
$$ \ar_\rho(C')\= \ar(A)\;.$$
By B.2 we can choose a curve \[\gamma'\] of winding number one on \[C'\],
or a corresponding curve \[\gamma\] on \[A\subset\T\], with
$$	\L(\gamma)^2\=\L_\rho(\gamma')^2\;\le\; {\ar_\rho(C')\over\md(C')}
 \={\ar(A)\over\md(A)}\;\le\;{ \ar(T)\over\md(A)}\;. $$
Now if we lift \[\gamma\] to the universal covering space \[\C\] then
it will join some point \[z_0\] to \[z_0+w \]. Hence its Euclidean
length \[\L(\gamma)\] must satisfy \[\L(\gamma)\ge
|w|\]. Thus
$$ |w|^2\;\le\; {\ar(T)\over\md(A)}\;, $$
which is equivalent to the required inequality (8).\QED\smallskip
\eject

Now consider the following situation. Let \[U\subset\C\] be a bounded
simply connected open set, and let \[K\subset U\] be a compact subset
so that the difference \[A=U\ssm K\] is a topological annulus. As noted in
\S2, such an annulus must be conformally isomorphic to a finite or
infinite cylinder. By definition an infinite cylinder, that is a cylinder of
infinite height, has modulus zero. (Such an infinite cylinder may be either
one-sided infinite, conformally isomophic to a punctured disk, or two-sided
infinite, conformally isomorphic to the punctured plane.)

{\QP{\bf Corollary B.7.} \it Suppose that \[K\subset U\] as described
above. Then \[K\]
reduces to a single point if and only if the annulus \[A=U\ssm K\] has infinite
modulus. Furthermore, the diameter of \[K\] is bounded by the inequality
$$	4\;{\rm diam}(K)^2\;\le\;{\ar(A)\over\md(A)}\;\le\;{\ar(U)
 \over\md(A)}\;. \eqno (9)  $$
\smallskip}

{\bf Proof.} According to B.2, there exists a curve with winding number
one about \[A\] whose
length satisfies \[L^2\le \ar(A)/\md(A)\]. Since \[K\] is enclosed
within this curve, it follows easily that \[{\rm diam}(K)\le\L/2\],
and the inequality (9) follows. Conversely, if \[K\] is a single point then
using (7) we see easily that \[\md(A)=\infty\].\QED

The following ideas are due to McMullen. (Compare [BH, II].)
The {\bit isoperimetric inequality~} asserts that the area enclosed by a plane
curve of length \[\L\] is at most \[\L^2/(4\pi)\], with equality if
and only if the curve is a round circle. (See for example [CR].) Combining
this with the argument above, we see that
$$	\ar(K)\;\le\;{\L^2\over 4\pi}\;\le{\ar(A)\over 4\pi\,\md(A)}\,. $$
Writing this inequality
as \[\;4\pi\,\md(A)\le \ar(A)/\ar(K)\;\] and adding \[+1\]
to both sides we obtain the completely equivalent inequality
\[	1+4\pi\,\md(A)\le\ar(U)/\ar(K)\;, \]
or in other words
$$	\ar(K)\;\le\;{\ar(U)\over 1+4\pi\,\md(A)}\;. \eqno (10) $$
This can be sharpened as follows:

{\QP{\bf Corollary B.8 (McMullen Inequality).} \it If
\[ A=U\ssm K\] as above, then
$$	\ar(K)\;\le\;{\ar(U)/ e^{4\pi\,\md(A)}}\;\;. $$
\smallskip}

{\bf Proof.} Cut the annulus \[A\] up into \[n\]
concentric annuli \[A_i\], each of modulus equal to \[\md(A)/n\]. Let \[K_i\]
be the bounded component of the complement of \[A_i\], and assume that
these annuli are nested so that \[A_i\cup K_i=K_{i+1}\] with \[K_1=K\],
and let \[K_{n+1}=A\cup K=U\]. Then \[\;\ar(K_{i+1})/
\ar(K_i)\,\ge\,1+4\pi\,\md(A)/n\;\] by (10), hence
$$	\ar(U)/\ar(K)\;\ge\;(\,1\,+\,4\pi\,\md(A)/n\,)^n\;, $$
where the right hand side converges to \[\;e^{4\pi\,\md(A)}\;\]
as \[n\to\infty\].\QED
\eject

% \medskip\lln\medskip

{\bf Concluding Problems:}\medskip

{\bf Problem B-1.} In the situation of Theorem
B.1, show that more than half of the
horizontal curves \[\gamma_y\] have length \[L_\rho(\gamma_y)\le
\sqrt{2\,\ar_\rho(I^2)}\]. (Here ``more than half'' is to be interpreted
in the sense of Lebesgue measure.)\smallskip

{\bf Problem B-2 (Sharper Bers Inequality).} If the flat torus
\[\T=\C/\Lambda\] contains several disjoint
annuli \[A_i\], all with the same ``winding number'' \[w\in\Lambda \],
show that
$$	\sum\md(A_i)\;\le\; \ar(\T)/|w |^2\;. $$
If two essentially
embedded annuli are disjoint, show that they necessarily
have the same winding number.

{\bf Problem B-3 (Branner-Hubbard).} Let \[K_1\supset K_2\supset K_3
\supset\cdots\] be compact subsets of \[\C\] with each \[K_{n+1}\]
contained in the interior of \[K_n\]. Suppose further that each interior
\[K_n^o\] is simply connected, and that each difference \[A_n=
K_n^o\ssm K_{n+1}\] is an annulus. If \[\sum_1^\infty\md(A_n)\] is infinite,
show that the intersection
\[\bigcap K_n\] reduces to a single point. Show that the
converse statement is false. (For example, do this by showing that a closed
disk \[\overline{D}'\] of radius \[1/2\] can be embedded in the open unit
disk \[D\] so that the complementary annulus \[A=D\ssm\overline{D}'\] has
modulus arbitrarily close to zero.)

\bigskip\bigskip
%\rightline{JWM, version of 7-29-91}

%\ref {\bf References:}

%\ref [CR] R. Courant and H. Robbins, What is Mathematics?, Oxford U. Press 1941.

%\ref B. Branner and J. H. Hubbard, The iteration of cubic polynomials II,
%patterns and parapatterns, Acta Math., to appear.
\vfil\eject

\footline={\hss\tenrm C-\folio\hss}\pageno=1

\def\={\;=\;}
\centerline{\bf Appendix C. Continued Fractions.}\bigskip

Suppose that we start with a real number \[r_1\] in the open interval
\[(0,1)\]. Then \[1/r_1\] can be written uniquely as the sum \[a_1+r_2\]
of an integer \[a_1\ge 1\] and a remainder term \[0\le r_2<1\].
If \[r_2>0\], then similarly we can set \[1/r_2=a_2+r_3\],
and so on, so that
 $$
	1/r_n\=a_n+r_{n+1}\;,\qquad 0\le r_{n+1}<1\;,  \eqno (1)
 $$
where each \[a_n\] is a strictly positive integer.
If \[r_1\] is rational, then \[r_2\] is rational with smaller
denominator, so this construction must stop with \[r_{n+1}=0\] after
finitely many steps. But if \[r_1\] is irrational then the construction
continues indefinitely:
 $$
  r_1\;\;=\;\;{1\hfil\over a_1+r_2}\;\;=\;\;{1\hfil\over \displaystyle a_1+
{\strut 1\hfil\over a_2+r_3}}\;\;=\;\;\cdots\;\;=\;\;{1\hfil\over
 \displaystyle a_1+{\strut 1\hfil\over \displaystyle a_2 + \;\;{{}\atop \ddots}{{}\atop {{}\atop
\displaystyle \;\; +{1\hfil\over\displaystyle a_{n-1}+r_n}}}}}\;\=\;\cdots\,\cdot
 $$
Closely related is the {\bit Euclidean
algorithm\/}, a procedure for
finding the greatest common divisor of two real numbers, when it exists.
Suppose that we start with two numbers \[\xi_0>\xi_1>0\]. Then we can express
\[\xi_0\] as an integral multiple of \[\xi_1\] plus a strictly smaller
remainder term,
 $$
  \xi_0 =a_1\xi_1 +\xi_2\;,\quad{\rm with}\quad \xi_1>\xi_2\ge 0\;,
 $$
where \[a_1\] is a strictly positive
integer. If \[\xi_2>0\], then similarly we can set
\[\xi_1=a_2\xi_2+\xi_3\] with \[\xi_2>\xi_3\ge 0\], and so on.
If the ratio \[r_1=\xi_1/\xi_0\] is rational, then this procedure
terminates with
 $$
	\xi_n={\rm greatest\;common\;divisor}\,,\qquad\xi_{n+1}=0
 $$
after finitely many steps. However, we will rather assume that \[r_1\]
is irrational so that the procedure can be continued indefinitely,
yielding an infinite sequence of numbers
 $$
  \xi_0 > \xi_1 > \xi_2 > \cdots > 0\qquad{\rm with}\qquad
  \xi_{n-1} \=a_n\xi_n +\xi_{n+1}\;\;, \eqno (2)
 $$
where each \[a_n\] is a strictly positive integer.
If we set \[r_n=\xi_n/\xi_{n-1}\in(0,1)\], then dividing equation (2)
by \[\xi_n\] we obtain \[1/r_n=a_n+r_{n+1}\], as above.\medskip

%Note that
% $$
%	\xi_{n+2}\;<\;{\xi_n\over 2}\quad{\rm or\; equivalently}
%	\quad r_nr_{n+1}\;<\;1/2\;.
% $$
%{\it Thus the sequence \[\{\xi_n\}\] always converges to zero, with \[\;
%\;\xi_n<{\rm constant}/2^{n/2}\,\]}. In fact we must have either
%\[\;\xi_{n+1}\le \xi_n/2\,\], in which case the assertion is clear,
%or \[\;\xi_{n+1}>\xi_n/2\;\] hence\break \[\xi_{n+2}\le \xi_n-\xi_{n+1}< \xi_n/2\].

In practice, it will be convenient to set  \[\;x_n = (-1)^{n+1} \xi_n\;\]
so that the \[x_n\] alternate in sign with \[\;|x_n|=\xi_n\;\]. Then
 $$
	x_{n+1}\; = \;x_{n-1} + a_n x_n\;\;.  \eqno (3)
 $$
In matrix notation we can write this as
 $$
	\left[\matrix{x_n \cr x_{n+1}}\right]\;\;=\;\;\left[\matrix{0 & 1\cr 1
	& a_n}\right]\left[\matrix{x_{n-1}\cr x_n}\right]
 $$
and therefore
 $$
        \left[\matrix{x_n \cr x_{n+1}}\right]\;\;=\;\;
	\left[\matrix{0 & 1\cr 1 & a_n}\right]
	\left[\matrix{0 & 1\cr 1 & a_{n-1}}\right]\cdots
	\left[\matrix{0 & 1\cr 1 & a_1}\right]
	\left[\matrix{x_0\cr x_1}\right]\;\=\;
	\left[\matrix{p_n & q_n\cr p_{n+1} & q_{n+1}} \right]
	\left[\matrix{x_0\cr x_1}\right]  \eqno (4)
 $$
where
 $$
	\left[\matrix{p_n & q_n\cr p_{n+1} & q_{n+1}} \right]\;\=\;
	\left[\matrix{0 & 1\cr 1 & a_n}\right]
	\left[\matrix{0 & 1\cr 1 & a_{n-1}}\right]\cdots
	\left[\matrix{0 & 1\cr 1 & a_1}\right]\;.	\eqno (5)
 $$
It follows inductively that
  $$\eqalign{
	 p_0 = 1\,,\qquad  p_1 = 0\,,\qquad &p_2=1\,,\qquad\,
	\;\ldots\,,\qquad  p_{n+1} = p_{n-1} + a_n p_n\;,  \cr
	 q_0 = 0\,,\qquad  q_1 = 1\,,\qquad &q_2=a_1\,,\qquad
	\ldots\,,\qquad q_{n+1} = q_{n-1} + a_n q_n \;. \cr} \eqno (6)
  $$
Here the integers \[0=p_1<p_2<p_3<\cdots\] and \[0=q_0<q_1<q_2<\cdots\]
grow at least exponentially with \[n\], since
  $$
	q_{n+2}\;\ge\; q_n+q_{n+1}\;\ge\; q_n+ (q_{n-1}+q_n)\;\ge\; 2q_n\,,
  $$
and similarly \[p_{n+2}\ge 2p_n\]. 
\medskip

Suppose, to fix our ideas, that we start with some irrational number
\[\xi=r_1\in(0,1)\] and set \[\xi_0=1\,,\;\;\xi_1=\xi\].
Then we can write equation (4) as
  $$
	x_n \;=\; p_n x_0 + q_n x_1\;=\; -p_n+q_n \xi\,,
  $$
or in other words
  $$
	\xi\;=\;{p_n\over q_n}\,+\,{x_n\over q_n}\;.	\eqno (7)
  $$
{\it Thus the irrational number \[\,\xi\,\] is equal to the rational number
\[\,p_n/q_n\,\] plus a remainder term \[\,x_n/q_n\,\].} In the literature,
these numbers \[p_n/q_n\] are called the {\bit convergents\/} to the
continued fraction expansion of \[\xi\]. (Compare Problem C-5.)
The numbers \[|x_n|\] converge at least
geometrically to zero as \[n\to\infty\] (see
Problem C-1 or inequality (8) below),
so the remainder term \[x_n/q_n\] in
equation (7) tends to zero quite rapidly as \[n\to\infty\].\smallskip

It follows from (7) that the successive convergents \[p_n/q_n\] to \[\xi\]
are ordered as follows:
  $$
	0={p_1\over q_1}\;<\; {p_3\over q_3}\;<\;\cdots\;<\;\; \xi\;\;<\;
	\cdots\;<\; {p_4\over q_4}\;<\;{p_2\over q_2}={1\over a_1}\;\le\; 1\;.
  $$
In order to estimate the error term \[x_n/q_n\] in terms of the
\[q_i\], note that the product matrix ~(5)~ has determinant \[(-1)^n\].
Therefore
  $$
  {p_n\over q_n} -{p_{n+1}\over q_{n+1}}\;\=\;{(-1)^n\over q_nq_{n+1}}\;\,.
  $$
Since \[\xi\] lies between \[p_n/q_n\] and \[p_{n+1}/q_{n+1}\] but is closer
to \[p_{n+1}/q_{n+1}\], it follows that the absolute value of the
error term \[x_n/q_n\] lies between \[1/(2q_nq_{n+1})\] and \[1/(q_nq_{n+1})\].
In other words:
 $$
	{1\over 2q_{n+1}}\;<\;|x_n|\;<\;{1\over q_{n+1}}\,\;.	\eqno (8)
 $$\bigskip
\smallskip

Now let us set \[\;\lambda=e^{2\pi i\xi}\,\] on the unit circle
\[\;S^1\subset\C\,\]. Using the discussion above, we
will study the orbit \[\;\;1\mapsto\lambda\mapsto\lambda^2\mapsto\cdots\;\;\]
under the rotation \[\;z\mapsto\lambda z\;\] of this circle.\medskip

{\bf Definition.} We say that a point \[\lambda^q\] on this orbit
is a {\bit closest return\/} to 1 if
 $$
	|\lambda^q-1|\;<\;|\lambda^m-1|
 $$
for every \[m\] with \[0<m<q\], so that
\[\lambda^q\] is closer to 1 than any preceding point on the orbit.
We will prove the following.

{\QP{\bf Lemma C.1.} \it The point \[\;\lambda^q=e^{2\pi i\xi q}\;\] is a closest
return to 1 along the orbit
 $$
	1\;\mapsto\;\lambda\;\mapsto\;\lambda^2\;\mapsto\;\cdots
 $$
if and only if \[q\] is one of the denominators \[\;1=q_1
\le q_2<q_3<q_4<\cdots\;\]
in the continued fraction approximations to \[\xi\]. Furthermore, if \[q=q_n\]
with \[n\ge 2\]
then the order of magnitude of the distance  \[|\lambda^q-1|\] is given by
 $$
	{2\over q_{n+1}}\;\;<\;\;|\lambda^{q_n}-1|\;\;<\;\;
	{2\pi\over q_{n+1}}\;\;. \eqno (9)
 $$
\medskip}

%q_n\xi\equiv x_n\in\R/\Z\] with \[n\ge 1\]
%are exactly the closest returns to zero along the orbit \[0\mapsto \xi
%\mapsto 2\xi\cdots\] in \[\R/\Z\]. Equivalently, the points \[\lambda^{q_n}\]
%are the closest returns to \[1\] along the orbit \[1\mapsto\lambda\mapsto
%\lambda^2\mapsto\cdots\] in the unit circle.\medskip}

The proof will be based on the following.
Instead of studying the multiplicative group \[S^1\subset\C\],
we can equally well work with the additive group \[\R/\Z\] and the orbit
\[0\mapsto\xi\mapsto 2\xi\mapsto\cdots\] under the translation \[t\mapsto
t+\xi\;(\mod\Z)\]. It will be convenient to introduce the following special
notation. For each real number \[x\] let
 $$
\ll x\gg\;=\;{\rm dist}(x\,,\,\Z)\={\rm Min}\,\{\;|x+n|\;:\;n\in\Z\;\}
% \;\le\;{\textstyle{1\over 2}}
 $$
be the distance to the nearest integer.
% For each \[m\xi\] there is a unique smallest representative
% \[\;\eta\equiv m\xi\;\mod\Z\], with \[\;-{1\over 2}< \eta\le{1\over 2}\,\].
Then an easy geometric argument shows that
 $$
	|\lambda^m-1|\;=\;|e^{2\pi im\xi}-1|\=
	2\,\sin\big(\pi\ll m\xi\gg\big)\;.
 $$
Since \[\;4<{2\sin(\pi t)/ t}<2\pi\;\] for \[t\in(0\,,\;{1/ 2})\],
it follows that
 $$
	4\;<\;{|\lambda^m-1|\over\ll m\xi\gg}\;<\;2\,\pi\;\,.\eqno (10)
 $$
%so that\[\;|\lambda^m-1|\;\] is small if
%and only if \[\;|\,m\xi\;\mod\Z|\;\] is small.
In the special case \[m=q_n\] note that
$$ q_n\,\xi\;\;\equiv\;\;x_n\quad(\mod\Z)	$$
by equation (7). If \[n\ge 2\], then \[|x_n|<{1/ 2}\]. (Compare Problem C-1.)
Hence it follows that \[\;|x_n|\;\] {\it is equal to \[\ll q_n\xi\gg\].}
The inequality (9) now follows from (8) and (10).\smallskip

Suppose, to fix our ideas, that \[n\] is odd, so that
 $$
	x_{n-1}\;<\;0\;<\; x_n\;<\; -x_{n-1}\;\;.
 $$
(The case \[n\] even is completely analogous, but with all inequalities
reversed.) Then it is easy to check that
 $$
	x_{n-1}\;<\; x_{n-1}+x_n\;<\;x_{n-1}+2x_n\;<\;\cdots\;<\;
	x_{n-1}+a_nx_n\;<\;0\;<\;x_n		\eqno (11)
 $$
where \[x_{n-1}+a_nx_n\] is equal to \[x_{n+1}\]. For example:\bigskip

\centerline{\psfig{figure=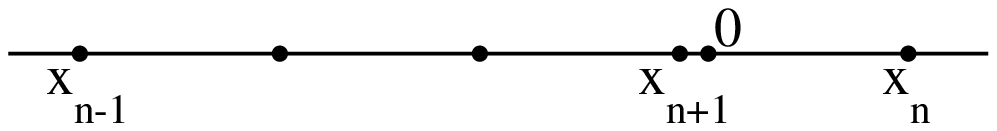,height=.4in}}

\noindent Note that each \[x_{n-1}
+jx_n\] is a representative \[\mod\Z\] for a point \[m\xi\] on our orbit,
where \[m\,=\,q_{n-1}+jq_n\]. For each \[n\ge 1\] we claim the following.

{\QP\bf Assertion \[{\bf A}_n\]. {\it No point \[m\xi\] with \[0<m<q_n\] has a
representative \[\mod\Z\] which lies strictly
between \[x_{n-1}\] and \[x_n\].}\medskip}

\noindent As part of the proof, we will also show the
following sharper statement.

{\QP{\bf Assertion \[{\bf B}_n\].} { \it For \[m\] in
the range \[0< m\le q_{n+1}\], the only points \[m\xi\]
which have a representative \[\eta_m\;\mod\Z\] which lies strictly
between \[x_{n-1}\] and \[x_n\] are the points \[x_{n-1}+jx_n\] which are
listed in formula \[(11)\].}\medskip}

\noindent (In each case, these numbers must be ordered appropriately, according
as \[n\] is even or odd.)
The proof will be by a double induction on \[n\] and \[m\].
The assertion \[{\bf A}_1\] is trivially true, since \[q_1=1\].
We will show that \[{\bf A}_n\Rightarrow{\bf B}_n\]. Since
it is easy to check that \[{\bf B}_n\Rightarrow {\bf A}_{n+1}\],
this will prove
inductively that the assertions \[\;{\bf A}_1\Rightarrow{\bf B}_1
\Rightarrow {\bf A}_2\Rightarrow{\bf B}_2\Rightarrow \cdots\;\]
are all true.

Suppose then that \[{\bf A}_n\] is true, and suppose, for some \[m\] between
\[0\] and \[q_{n+1}\], that \[m\xi\] has a representative \[\eta_m\;
\mod\Z\] which lies strictly between \[x_{n-1}\] and \[x_n\]. We will show
by induction on \[m\] that \[m\] must have the form \[q_{n-1}+jq_n\].
According to \[{\bf A}_n\] we must have \[m>q_n\]. % Let us set \[m=q_n+m'\].
We divide the proof into two cases according as \[\eta_m\] lies between
\[x_{n-1}\] and \[x_{n-1}+x_n\] or between \[x_{n-1}+x_n\] and \[x_n\].
In the former case, we see that \[(m-q_{n-1})\,\xi\] has a representative
\[\eta_m-x_{n-1}\;\mod\Z\] which lies between \[0\] and \[x_n\], thus
contradicting our induction hypothesis on \[m\]. In the later case, we
see that \[(m-q_n)\,\xi\] has a representative \[\;\eta_m-x_n\;\;\mod\Z\;\]
which lies between \[x_{n-1}\] and \[0\]. Therefore, by induction,
\[m-q_n\] has the form \[q_{n-1}+jq_n\], hence \[m\] itself also has this form.
This completes the double induction, proving \[{\bf A}_n\] and \[{\bf B}_n\].
Evidently Lemma C.1 follows easily from the Assertions \[{\bf A}_n\].\QED
%\eject

 \lln\medskip

{\bf Concluding Problems.}\medskip

{\bf Problem C-1.} With \[r_n\] as in equation
(1) or (3), show that the product
$$	r_nr_{n+1}\=\xi_{n+1}/\xi_{n-1}	$$
is always smaller than \[1/2\]. (If both \[r_n\] and \[r_{n+1}\]
are greater than \[1/2\] then\break \[r_nr_{n+1}=1-r_{n+1}<1/2\].)
Conclude that the numbers
$$	\xi_{2n}\;<\; \xi_0/2^n $$
converge rapidly to zero as \[n\to\infty\].\medskip

{\bf Problem C-2.} 
In the simplest possible case \[\;a_1 = a_2 =\cdots = 1\;\], show that
 $$
	\{q_n\} = \{p_{n+1}\} = \{ 0\,,\; 1\,,\; 1\,,\; 2\,,\; 3\,,\; 5\,,\;
	 8\,,\; 13\,,\; 21\,,\;\ldots\; \}\;\;,
 $$
yielding the sequence
of Fibonacci numbers. Prove the asymptotic formula
$q_n\sim \gamma^n/\sqrt 5$ as \[n\to\infty\], and hence \[p_n/q_n\to 1/\gamma\]
as  $n\to\infty$, where $\gamma=(\sqrt 5 + 1)/2$.
Show that this special case corresponds to the slowest possible
growth for the coefficients \[p_n\] and \[q_n\].\medskip

{\bf Problem C-3.} Using Lemma C.1, show for
\[n\ge 2\] that the convergent \[p_n/q_n\]
is the {\bit best} rational
approximation to \[\xi\] among all fractions with
denominator \[m\le q_n\]. In particular,
 $$
	\left|\xi-{p_n\over q_n}\right| \;<\; 
	\left|\xi-{i\over m}\right|
 $$
for any integers \[i\] and \[m\] with \[0<m<q_n\].\medskip

{\bf Problem C-4.}
Define a polynomial \[{\cal P}(x_1\,,\,\ldots\,,\,x_n)\] in \[n\] variables
inductively by the formula
 $$	
   {\cal P}(x_1\,,\,\ldots\,,\,x_{n})\={\cal P}(x_1\,,\,\ldots\,,\,x_{n-2})
  \,+\,{\cal P}(x_1\,,\,\ldots\,,\, x_{n-1})\,x_n\;\;,
 $$
starting with
 $$
	{\cal P}(\,)\=1\;,\quad{\cal P}(x)\=x\;,\quad{\cal P}(x,y)\=1+xy\;,
	\quad{\cal P}(x,y,z)\=	x+z+xyz\;,\;\;\ldots\,.
 $$
Show that the \[p_n\] and \[q_n\]
can be expressed as polynomial functions of the $a_i$
by the formulas
 $$
	p_{n+1}={\cal P}(a_2\,,\,\ldots\,,\,a_n)\;,\qquad
	q_{n+1}={\cal P}(a_1\,,\,\ldots\,,\,a_n) \;\;.
 $$
Using the inverse of the matrix equation (5), show that
 $$
  {\cal P}(a_1\,,\,\ldots\,,\,a_n)\;=\;{\cal P}(a_n\,,\,\ldots\,,\,a_1)\;\;.
 $$
Conclude that we can also write
 $$	
   {\cal P}(x_1\,,\,\ldots\,,\,x_{n})\= x_{1}\,{\cal P}(x_2\,,\,\ldots\,,\,
   x_n)\,+\,{\cal P}(x_3\,,\,\ldots\,,\,x_{n})\;\;. \eqno (12)
 $$
Show that \[{\cal P}(x_1\,,\,\ldots\,,\,x_n)\] is equal to the sum of all
distinct monomials which can be formed out of the product $\;x_1\cdots x_n$
by striking out any number of consecutive pairs. Show that the number of such
monomials is equal to the $n$-th Fibonacci number.\bigskip

{\bf Problem C-5.} Using (12), give an inductive proof of
the finite continued fraction equation
  $$
  {p_{n+1}\over q_{n+1}}\;\;\=\;\;{{\cal P}(a_2\,,\,\ldots\,,\,a_n)\over{\cal P}
  (a_1\,,\,\ldots\,,\,a_n)}\;\;\=\;\; {1\hfil\over
  \displaystyle a_1+{\strut 1\hfil\over \displaystyle a_2 + \;\;{{}\atop
  \ddots}{{}\atop {{}\atop
  \displaystyle \;\; +{1\hfil\over\displaystyle a_{n-1}+{1\over a_n}}}}}}
  \;\;\; .\qquad
 $$
%by induction on \[n\].

\vfill
%\rightline{jwm, version of 7-28-91}

\end

%% file: cdd-g.tex
\input header
\footline={\hss\tenrm D-\folio\hss}\pageno=1

{\bf Appendix D. Remarks Concerning Two Complex Variables.}\bigskip

Many of the arguments in these lectures are strictly
one-dimen\-sional. In fact several of our underlying principles
break down completely in the two variable case.\smallskip

In order to illustrate these differences, it is useful to consider the
family of (generalized)
{\bit H\'enon maps}, which can be described as follows. Choose a complex
constant
\[\delta\ne 0\] and a polynomial map \[f:\C\to\C\] of degree \[d\ge 2\],
and consider doubly infinite sequences of complex numbers \[\;\ldots\,,\,
z_{-1}\,,\,z_0\,,\,z_1\,,\,z_2\,,\,\ldots\;\] satisfying the recurrence
relation
$$	z_{n+1}-f(z_n)+\delta\,z_{n-1}\;=\;0\,.	$$
Evidently we can solve for \[(z_n\,,\,z_{n+1})\] as a polynomial function
\[\;F(z_{n-1}\,,\,z_n)\,\], where the transformation
\[\;\;F(z_{n-1}\,,\,z_n)\,=\,(z_n\,,\, f(z_n)-\delta z_{n-1})\;\;\] has Jacobian matrix
$$	\left[\matrix{0 & 1\cr -\delta & f'(z_n)\cr}\right] $$
with constant determinant \[\delta\], and with trace \[f'(z_n)\]. Similarly
we can solve for\break \[(z_{n-1}\,,\,z_n)\] as a polynomial function
\[F^{-1}(z_n\,,\,z_{n+1})\]. As an example, consider the quad\-ratic polynomial
\[f(z)=z^2+\lambda z\], having a
fixed point of multiplier \[\lambda\] at the origin. Then
$$	F(x\,,\,y)\;=\; (y\,,\; y^2+\lambda y-\delta x)\,,
	\eqno (*) $$
where \[x\] and \[y\] are complex variables. This map
has a fixed point at \[(0,0)\], where the eigenvalues \[\lambda_1\] and
\[\lambda_2\] of the Jacobian matrix satisfy \[\lambda_1+\lambda_2=\lambda\]
and \[\lambda_1\lambda_2=\delta\]. Evidently, by the appropriate choice
of \[\lambda\] and \[\delta\], we can realize any desired non-zero
\[\lambda_1\] and \[\lambda_2\]. If \[\lambda_1\ne\lambda_2\], then we
can diagonalize this Jacobian matrix by a linear change of coordinates.
\smallskip

{\QP{\bf Lemma D.1.} \it Consider any holomorphic transformation
\[F(x,y)=(x',y')\] in two complex variables, with
$$	x'=\lambda_1 x +{\cal O}
	(|x^2|+|y^2|)\,,\quad y'=\lambda_2 y+{\cal O}(|x^2|+|y^2|)	$$
as \[(x,y)\to(0,0)\].
If the eigenvalues \[\lambda_1\] and \[\lambda_2\] of the derivative
at the\break origin satisfy \[1>|\lambda_1|\ge
|\lambda_2|>|\lambda_1^2|\], then \[F\] is conjugate, under a local\break
holomorphic change of coordinates, to the linear map \[L(u,v)=(\lambda_1
u\,,\,\lambda_2v)\].\medskip}

{\bf Proof.} We must show that there exists a change of coordinates \[(u,v)=\phi(x,y)\],
defined and holomorphic throughout a neighborhood of the origin,
so that \[\phi\circ F\circ\phi^{-1}=L\]. As in the proof of the Ko$\!$enigs
Theorem, \S6.1, we first choose a constant \[c\] so that\break
\[1>c>|\lambda_1|\ge|\lambda_2|>c^2\]. To any orbit
$$	(x_0\,,\,y_0)\buildrel F\over\mapsto(x_1\,,\,y_1)
	\buildrel F\over\mapsto\cdots	$$
near the origin, we associate the sequence of points
$$	(u_n\,,\,v_n)\;=\;L^{-n}(x_n\,,\,y_n)\;=\;(x_n/\lambda_1^n\,,\,y_n
	/\lambda_2^n)	$$
and show, using Taylor's Theorem, that it converges geometrically to the
required limit \[\phi(x_0\,,\,y_0)\], with
successive differences bounded by a constant times \[(c^2/\lambda_2)^n\].
Details will be left to the reader.\QED

{\bf Remarks.} Some such restriction on the eigenvalues is essential.
As an example, for the map
$$	f(x,y)\;=\;(\lambda x,\,\lambda^2 y+x^2)\,,	$$
with eigenvalues \[\lambda\] and \[\lambda^2\],
there is no such local holomorphic change of coordinates. (Problem D-1.)
For a much sharper statement as to when linearization is possible,
even when one or both eigenvalues lie on the unit circle, compare Zehnder.

\smallskip

Now consider a H\'enon map \[F:\C^2\buildrel \approx\over\to\C^2\] as
in $(*)$ above,
with eigenvalues satisfying the conditions of D.1. Let \[\Omega\] be the
attractive basin of the origin. We claim that \[\phi\] extends
to a global diffeomorphism \[\Phi:\Omega\buildrel \approx\over\to\C^2\].
In fact, for any \[(x,y)\in\Omega\] we set \[\Phi(x,y)=L^{-n}\circ\phi\circ
F^{\circ n}(x,y)\]. If \[n\] is sufficiently large, then \[F^n(x,y)\] is
close to the origin, so that this expression is defined. Similarly,
\[\Phi^{-1}(u,v)=F^{-n}\circ\phi^{-1}\circ L^{\circ n}\] is well defined
for large \[n\]. This shows that
\[\Phi\] is a holomorphic diffeomorphism with holomorphic inverse.\smallskip

%***

%For suitable choice of \[\lambda_1\] and \[\lambda_2\] in the unit disk,
%it can shown ??? that \[F\] is locally conjugate to the linear map
%\[\;L:(w_1\,,\,w_2)\mapsto(\lambda_1w_1\,,\,\lambda_2w_2)\,\]. {\it Since
%\[F\] is a
%global diffeomorphism, it follows easily that the attractive basin \[B\]
%of the origin is analytically diffeomorphic to the entire coordinate space
%\[\C^2\] under a diffeomorphism \[\phi:B \buildrel \approx\over\longrightarrow
%\C^2\] which satisfies \[\phi\circ F=L\circ\phi\].} (Compare \S6.2.)
%\smallskip

%***

Note that this basin \[\Omega\] is not the entire space \[\C^2\]. For example, if \[
|z_1|\] is sufficiently large compared with \[|z_0|\], and if
$$	(z_0\,,\,z_1)\buildrel F\over\mapsto (z_1\,,\,z_2)\buildrel
	 F\over\mapsto (z_2\,,\,z_3)\buildrel F\over\mapsto\cdots	$$
then it is not difficult to check that \[\;|z_1|<|z_2|<|z_3|<\cdots\,\],
so that \[(z_0\,,\,z_1)\] is not in \[\Omega\].
{\it Thus we have constructed:\smallskip

$(1)$ a proper subset \[\Omega\subset\C^2\] which is analytically diffeomorphic to
all of \[\C^2\], and\smallskip

$(2)$ a non-linear map with an attractive basin which contains no critical
points.}\smallskip

\noindent Evidently neither phenomenon can occur in one complex variable.
Open sets satisfying (1) are called {\bit Fatou-Bieberbach domains} since they
were first constructed by Fatou, by an easy argument similar to that given
here, and then later independently by Bieberbach, who had a much more
difficult construction.\smallskip

The proof that there are only finitely many attracting cycles
also breaks down in two variables. Compare Newhouse.
\lln

{\bf Problem D-1.} For the map \[F(x,y)=(\lambda x\,,\,\lambda^2 y+x^2)\],
where \[\lambda\ne 0\,,\,1\], show that there is only one smooth
\[F$-invariant curve through the origin, namely \[x=0\]. By way of contrast,
for the associated linear map \[L(x,y)=(\lambda x\,,\,\lambda^2 y)\]
note that there are infinitely many \[F$-invariant curves \[y=cx^2\].
Conclude that \[F\] is not locally holomorphically conjugate to a linear map.

\vfil\eject % \bigskip\bigskip
\footline={\hss\tenrm E-\folio\hss}\pageno=1

\centerline{\bf Appendix E. Branched Coverings and Orbifolds.}\bigskip

This will be an outline of definitions and results, without proofs.
(See the list of references at the end.)
We will use ``branch point'' as a synonym for ``critical point''
and ``ramified point'' as a synonym for ``critical value''. Thus if
\[f(z_0)=w_0\] with \[f'(z_0)=0\], then \[z_0\] is called a {\bit branch
point} and the image \[f(z_0)=w_0\] is called a {\bit ramified point}.
More precisely, if
$$	f(z)\;=\;w_0+c(z-z_0)^n+({\rm higher\; terms})\,,	$$
with \[n\ge 1\] and
\[c\ne 0\], then the integer \[n=n(z_0)\] is called the {\bit branch index}
or the {\bit local degree} of \[f\] at
the point \[z_0\]. Thus \[n(z)\ge 2\] if \[z\] is a branch point,
and \[n(z)=1\] otherwise.

A holomorphic map \[p:S'\to S\] between Riemann surfaces is called a
{\bit covering map} if each point of \[S\] has a connected
neighborhood \[U\] which is {\bit evenly covered}, in that each connected
component of \[p^{-1}(U)\subset S'\] maps onto \[U\] by a conformal
isomorphism. A map \[p:S'\to S\] is {\bit proper} if the inverse image
\[p^{-1}(K)\] of any compact subset of \[S\] is a compact subset of \[S'\].
Note that every proper map is finite-to-one, and has a well
defined finite degree \[d\ge 1\]. Such a map may also be called a
\[d$-{\bit fold branched covering}. On the other hand, a covering
map may well be infinite-to-one. Combining these two concepts,
we obtain the following more general concept.\smallskip

{\bf Definition.} A holomorphic map \[p:S'\to S\] between Riemann surfaces
will be called a {\bit branched covering map} if every point of \[S\] has
a connected neighborhood \[U\] so that each connected component of \[p^{-1}
(U)\] maps onto \[U\] by a proper map.\smallskip

Such a branched covering is said to be {\bit regular} or
{\bit normal} if there exists a group \[\Gamma\] of conformal automorphisms
of \[S'\], so that two points \[z_1\] and \[z_2\] of \[S'\] have the same
image in \[S\] if and only if there is a group element \[\gamma\] with
\[\gamma(z_1)=z_2\]. In this case we can identify \[S\] with the quotient
manifold \[S'/\Gamma\]. In fact it is not difficult to check that the
conformal structure of such a quotient manifold is uniquely determined.
This \[\Gamma\] is called the group of {\bit deck transformations} of the
covering.

Regular branched covering maps have several special properties. For
example, each ramified point is isolated,
so that the set of all ramified points is a discrete subset of \[S\].
Furthermore, the branch index \[n(z)\] depends only on the target point
\[f(z)\], that is, \[n(z_1)=n(z_2)\]
whenever \[f(z_1)=f(z_2)\]. Thus we can define the {\bit ramification
function} \[\nu:S\to\{1,\,2,\,3,\,\ldots\}\] by setting \[\nu(w)\] equal
to the common value of \[n(z)\] for all points \[z\] in the pre-image
\[f^{-1}(w)\]. By definition, \[\nu(w)\ge 2\] if \[w\] is a ramified
point, and \[\nu(w)=1\] otherwise.\smallskip

{\bf Definition.} A pair \[(S,\nu)\] consisting of a Riemann surface \[S\]
and a ``ramification function'' \[\nu:S\to\{1,\,2,\,3,\,\ldots\}\] which
takes the value \[\nu(w)=1\] except at isolated points
will be called a Riemann surface {\bit orbifold}.

(Remark: Thurston's general concept of orbifold
involves a structure which is locally modeled on the quotient of a coordinate
space by a finite group. However, in the Riemann surface case only
cyclic groups can occur, so a simpler definition can be used.)

{\bf Definition.}  If \[S'\] is simply connected, then a regular branched
covering \[p:S'\to S\] with ramification function \[\nu\]
will be called the {\bit universal covering}
for the orbifold \[(S,\,\nu)\]. We will use the notation \[\tilde S_\nu\to S\]
for this universal branched covering. The associated
group \[\Gamma\] of deck transformations is called the {\bit fundamental
group} \[\pi_1(S,\,\nu)\] of the orbifold.

{\QP{\bf Lemma E.1.} \it With the following exceptions, every
Riemann surface orbifold
\[(S,\,\nu)\] has a universal covering surface \[\tilde S_\nu\]
which is unique up to conformal isomorphism over \[S\]. The only exceptions
are given by:

(1) a surface \[S\approx\hat\C\]
with just one ramified point, or

(2) a surface \[S\approx\hat\C\] with two ramified points
for which \[\nu(w_1)\ne\nu(w_2)\].

\noindent
In these exceptional cases, no such universal covering exists.\medskip}

By definition,
the {\bit Euler characteristic} of an orbifold \[(S,\,\nu)\] is
the rational number
$$      \chi(S,\,\nu)\;=\;\chi(S)+\sum \big({1\over\nu(w_j)}-1\big)\,, $$
to be summed over all ramified points, where \[\chi(S)\] is the usual
Euler characteristic of \[S\]. Intuitively speaking,
each ramified point \[w_j\] makes a contribution of \[+1\] to the usual
Euler characteristic \[\chi(S)\], but a smaller contribution of \[1/\nu(w_j)\]
to the orbifold Euler characteristic. Thus \[\chi(S,\,\nu)<\chi(S)\le 2\],
or more precisely
$$	\chi(S)-r\;<\; \chi(S,\,\nu)\;\le\;\chi(S)-r/2	$$
where \[r\] is the number of ramified points. {\it As an example, if
\[\chi(S,\,\nu)
\ge 0\], with at least one ramified point, then it follows that \[\chi(S)>0\],
so the base surface \[S\] can only be \[D,\,\C\] or \[\hat \C\], up
to isomorphism.} Compare E.5 below.

If there are infinitely
many ramified points, note that we must set \[\chi(S,\,\nu)=
-\infty\]. Similarly, if \[S\] is
a surface which is not of finite type, then we must set\break
\[\chi(S,\,\nu)=\chi(S)=-\infty\].\smallskip

If \[S'\] and \[S\] are provided with ramification functions \[\mu\] and
\[\nu\] respectively, then a branched covering map \[f:S'\to S\] is said
to yield
a {\bit covering map \[(S'\,,\,\mu)\to(S\,,\,\nu)\] between orbifolds}
if the identity
$$      n(z)\mu(z)\;=\;\nu(f(z))       $$ 
is satisfied for all \[z\in S'\],
% if \[n(z)\nu'(z)\] is a multiple of \[\nu(f(z))\] for each \[z\in S'\],
where \[n(z)\] is the branch index.
% Similarly, we say that \[f^{-1}\]
%is {\bit locally holomorphic} between orbifolds if \[\nu(f(z)\] is a
% multiple of \[n(z)\nu'(z)\], and that \[f\] is a {\bit covering map} of
As an example, the universal covering map \[\tilde S_\nu\to
(S\,,\,\nu)\] is always a covering map of orbifolds, where \[\tilde S_\nu\]
is provided with the trivial ramification function \[\mu\equiv 1\].

{\QP\bf Lemma E.2. {\it \[f:(S'\,,\,\mu)\to(S\,,\,\nu)\] is a covering map
between orbifolds if and only if it lifts to a conformal isomorphism from the
universal covering \[\tilde S'_\mu\] onto \[\tilde S_\nu\].
If \[f\] is a covering in this sense, and has finite degree \[d\],
then the Riemann-Hurwitz formula of \S5.1 takes the form
$$	\chi(S'\,,\,\mu)\;=\; \chi(S\,,\,\nu)d\,.		$$
In particular, if the universal covering of \[(S\,,\,\nu)\] is a covering
of finite degree \[d\], then \[\;\chi(\tilde S_\nu)=\chi(S,\,\nu)d\].}
\medskip}

\noindent The fundamental group and the Euler characteristic
% of any Riemann surface orbifold 
are related to each other as follows.

\smallskip

{\QP{\bf Lemma E.3.} \it Let \[(S,\,\nu)\] be any Riemann surface orbifold
which possesses a universal covering:

If \[\chi(S,\,\nu)>0\], then the fundamental
group \[\pi_1(S,\,\nu)\] is finite.

If \[\chi(S,\,\nu)=0\], then the fundamental group contains either \[\Z\]
or \[\Z\oplus\Z\] as a subgroup of finite index.

If \[\chi(S,\,\nu)<0\], then the fundamental
group contains a non-abelian free product \[\Z *\Z\], and hence cannot contain
any abelian subgroup of finite index.\medskip}

\noindent
The Euler characteristic and the geometry of \[\tilde S_\nu\] are related
as follows.

{\QP{\bf Lemma E.4.} \it Suppose that \[S\] is obtained from a compact
Riemann surface by
removing at most a finite number of points. Then \[\tilde S_\nu\] is
either spherical, hyperbolic or Euclidean according as \[\chi(S\,,\,\nu)\]
is positive, negative or zero.\medskip}

{\bf Remark.} In all other cases, that is whenever \[S\] is
not isomorphic to a finitely punctured compact surface, \[S\] is
hyperbolic, and it follows that the universal covering
\[\tilde S_\nu\] must also be hyperbolic.

{\bf Examples.} If \[S=\C\] with two ramified
points \[\nu(1)=\nu(-1)=2\], then the map\break
\[z\mapsto\cos(2\pi z)\] provides
a universal covering \[\C\to (\C,\,\nu)\]. The Euler characteristic
\[\chi(\C,\,\nu)\] is zero, and the fundamental group \[\pi_1(S,\,\nu)\]
consists of all transformations of the form \[\gamma:z\mapsto  n\pm z\].

For \[S=\hat\C\] with three ramified points \[\nu(0)=\nu(1)=\nu(
\infty)=2\],
the rational map \[\pi(z)=-4z^2/(z^2-1)^2\] provides a
universal covering \[\hat\C\to(\hat\C,\,\nu)\]. In this case we have
\[\chi(\hat\C,\,\nu)=1/2\], and
the degree is equal to \[\chi(\hat\C)/\chi(\hat\C,\,\nu)=4\].
The fundamental group consists of all transformations \[\gamma:z\mapsto
\pm z^{\pm 1}\].

For \[S=\hat\C\] with four ramified points of index \[\nu(w_j)=2\], then the
torus \[T\] described in \S5 provides a regular 2-fold branched covering.
Its universal covering \[\tilde T\] can be identified with the
universal covering of \[(\hat\C,\,\nu)\]. In this case, the Euler
characteristic \[\chi(\hat\C,\,\nu)\] is zero.\smallskip

{\bf Remark E.5.} There are relatively few cases in which \[\chi(S\,,\,
\nu)\ge 0\]. In fact all of these cases can be listed quite explicitly
as follows. The unramified cases are very well known, namely
the sphere, plane or disk with \[\chi>0\]; and the punctured
plane or disk, and the infinite families consisting of annuli and tori,
all with \[\chi=0\].

 By the ``ramification indices'' we will
mean the list of values of the ramification function at the \[r\]
ramified points,
ordered for example so that \[\nu(w_1)\le\cdots\le\nu(w_r)\].

If \[\chi(\hat\C,\,\nu)>0\] with \[r>0\],
then the ramification indices must be either
\[(n,n)\] or \[(2,2,n)\] for some \[n\ge 2\], or \[(2,3,3)\,,\, (2,3,4)\]
or \[(2,3,5)\]. These five possibilities correspond to the five types of
finite rotation groups of the 2-sphere; namely to the cyclic,
dihedral, tetrahedral, octahedral and icosahedral groups respectively.

If \[\chi(\hat\C,\,\nu)=0\], then the ramification indices must be either
\[(2,4,4)\], \[(2,3,6)\], \[(3,3,3)\] or \[(2,2,2,2)\].
These correspond to the automorphism groups of the tilings of \[\C\] by
squares, equilateral triangles, alternately colored equilateral triangles,
and parallelograms respectively. In the parallelogram case, note that there
is actually a one
complex parameter family of distinct possible shapes, corresponding to the
cross-ratio of the four ramified points.

Similarly, if \[\chi(\C,\,\nu)\] or \[\chi(D,\,\nu)\] is strictly positive,
then we must have \[r\le1\]; while if \[\chi(\C,\,\nu)\] or \[\chi(D,\,\nu)\]
is zero, then we must have \[r=2\] with ramification  indices \[(2,2)\].
This is the complete list.

\vfill
\eject % \bigskip\bigskip

\footline={\hss\tenrm F-\folio\hss}\pageno=1

\centerline{\bf Appendix F. Parameter Space.}\bigskip

A very important part of complex dynamics, which has barely been mentioned in
these notes, is the study of parametrized families of mappings. As an example,
consider the family of all quadratic polynomial maps. A priori, a quadratic
polynomial is specified by three complex parameters; however
any such polynomial can be put into the unique normal form
$$	f(z)\;=\;z^2+c \eqno (1)	$$
by an affine change of coordinates. (Other normal forms which have been used
are\break \[	\omega\;\mapsto\;\omega^2+\lambda \omega\],
with a preferred fixed point of multiplier \[\lambda\] at the origin, or
$$	w\; \mapsto\;\lambda w (1-w) \eqno (2) $$
which is more or less equivalent provided that
\[\lambda\ne 0\]. Here \[4c=\lambda(2-\lambda)\] and\break
 \[\omega=z-\lambda/2=-\lambda w\].) Using such a normal form, we can make
a computer picture in the
{\bit parameter space} consisting of all complex constants \[c\] or \[\lambda
\]. Each pixel in such a picture, corresponding to a small square in the
parameter space, is to be assigned some color, perhaps only black or white,
which depends on the dynamics of the corresponding quadratic map.

The first crude pictures of this type were made by Brooks and Matelski,
as part of a study of Kleinian groups. They used the normal form (1), and
introduced the open set consisting
of all points of the \[c$-plane for which the corresponding quadratic map
has an attracting periodic orbit in the finite plane.
I will use the notation \[{\cal H}\] (standing
for ``hyperbolic'') for this Brooks-Matelski set.
At about the same time,
Hubbard (unpublished) made much better pictures of a quite different
parameter space arising from Newton's method for cubics. Mandelbrot,
perhaps inspired by Hubbard, made corresponding pictures for quadratic
polynomials, using the normal form (2) and also a variant of (1).
In order to avoid confusion, let me translate all of Mandelbrot's
definitions to the normal form (1).
He introduced two different sets, which I will call \[Q\subset
M\]. Mandelbrot did not give these sets
different names, since he believed that they were identical. By definition,
a parameter value \[c\] belongs to \[Q\] if the corresponding filled Julia set
contains an interior point, and belongs to \[M\] if its filled Julia
set contains the critical point \[z=0\] (or equivalently is connected).
The Brooks-Matelski set
\[{\cal H}\] is a subset of \[Q\subset M\]. Mandelbrot made quite good computer
pictures, which seemed to show a number of isolated ``islands''. Therefore,
he conjectured that \[Q\] [or \[M\]] has many distinct connected components.
(The editors of the journal thought that his islands were specks of dirt,
and carefully removed them from the pictures.) Mandelbrot also described
an important smaller set \[{\cal L}\subset Q\] which he believed to be the
principal connected component of \[Q\]. This set \[{\cal L}\]
consists of a central cardioid
\[{\cal L}_0\] with some (but not all) boundary points included,
together with countably many smaller nearly-round disks which are pasted on
inductively, in an explicitly described pattern.\smallskip

Although Mandelbrot's statements in this first paper
were not completely right, he deserves a great
deal of credit for being the first to point out the extremely complicated
geometry associated with the parameter space for quadratic maps.
His major achievement has been to demonstrate
to a very wide audience that such complicated ``fractal''
objects play an important role in a number of mathematical sciences.
\smallskip

The first real mathematical breakthrough came with Douady and Hubbard's work in
1982. They introduced the name {\bit Mandelbrot set} for the compact set
\[M\] described above, and
provided a firm foundation for its mathematical study,
proving for example that \[M\] is connected with connected complement.
(Meanwhile,
Mandelbrot had decided empirically that his isolated islands were actually
connected to the mainland by very thin filaments.)
Already in this first paper, they showed that each hyperbolic component
of the interior of \[M\] can be canonically
parametrized, and showed that the boundary \[\partial M\] can be profitably
studied by following external rays.\smallskip

It may be of interest to compare the three sets \[{\cal H}\subset Q\subset M\]
in parameter space. They are certainly different since \[{\cal H}\] is open,
\[M\] is compact, and \[Q\] is neither.
Using Sullivan's work (\S13), we can say that
\[Q\] consists of \[{\cal H}\] together with a very sparse set of boundary
points, namely those for which the corresponding map has
either a parabolic orbit or
a Siegel disk. Quite likely, there is no difference between these three sets
as far as computer graphics are concerned, since it is widely
conjectured{\bf *} that the Brooks-Matelski set \[{\cal H}\]
is equal to the interior of \[M\], and that
\[M\] is equal to the closure
of \[{\cal H}\]. However, as far as practical computing is concerned, it should
be noted that it is quite easy to
test (at least roughly)
whether a parameter value belongs to \[M\], but somewhat harder to decide
whether it belongs to \[{\cal H}\] (compare \S10.5), and very difficult
to decide
whether it belongs to \[Q\]. (Compare Appendix G. Here I am only speaking of
approximate tests. To decide precisely whether a given point belongs to
\[M\] or to whether a given point belongs to \[{\cal H}\] may be very
difficult. As a specific example, I have no idea whether or not the point
\[c=-1.5\] belongs to \[{\cal H}\].)\smallskip

Another important development came in 1983, with
the work of Ma\~n\'e, Sad and Sullivan on stability of the Julia set
\[J(f)\] under deformation of \[f\]. (Compare the discussion following \S14.3.)
These results were obtained independently by Lyubich.
The study of parameter space for
higher degree polynomials began some five years later with the work of
Branner and Hubbard. Using the normal form
$$	f(z)\;=\;z^3-3a^2z+b\,, $$
with the two critical points at \[z=\pm a\], they proved that the {\bit cubic
connectedness locus}, consisting of all parameter pairs \[(a,b)\] for which
\[J(f)\] is connected, is a {\bit cellular set}. (Compare \S17.3.)
In particular, this set is compact and connected. A corresponding result for
polynomials of higher degree has recently been obtained by Lavaurs.
Parameter space studies for rational maps are more awkward,
since there is no obvious normal form. However, an important beginning
has been made by Rees. All of these studies are very new, and much remains
to be done.\smallskip

Here are two problems for the reader.
\vfil
\noindent ------------

\noindent {\bf *} Douady and Hubbard have shown that these conjectures are true
if the set \[M\] is locally connected. Recent work of Yoccoz lends very
strong support to the belief that \[M\] is indeed locally connected.
\eject

{\bf Problem F-1.} Show that every polynomial map of degree \[d\ge 2\]
is conjugate, under an
affine change of coordinates, to one in the ``Fatou normal form''
$$	f(z)\;=\; z^d + a_{d-2}z^{d-2}+\cdots+a_1z+a_0\,.	$$
Let \[P(d)\cong\C^{d-1}\] be the space of all such maps.
Show that the cyclic group \[Z(d-1)\] of \[(d-1)$-st roots of unity acts on
\[P(d)\] by linear conjugation, replacing \[f(z)\] by \[f(\omega z)/\omega\],
and show that the quotient \[P(d)/Z(d-1)\] can be identified with the
``moduli space'' of degree \[d\] polynomials up to affine conjugation.
If \[d\ge 4\], show that this moduli space is not a manifold.\smallskip

{\bf Problem F-2.} Show that every quadratic rational map is conjugate, under
a fractional linear change of coordinates, for example to one in the form
$$	f(z)\;=\;1+(z^2-1)/(az^2-b) $$
where \[a\ne b\],
with critical points at zero and \[\infty\] and a fixed point at \[z=1\].
(In general this form is not unique, since such a
map actually has three fixed points, and since the roles of zero and infinity
can be interchanged.)% Show that the moduli space of quadratic rational maps
%up to fractional linear change of coordinates has just one non-manifold point,
%namely \[z\mapsto 1/z^2\].\smallskip

\vfil\eject
\footline={\hss\tenrm G-\folio\hss}\pageno=1

\centerline{\bf Appendix G. Remarks on Computer Graphics.}\medskip

In order to make a computer picture of some complicated compact set
\[L\subset\C\], for example a Julia set,
we must compute a matrix of small
integers, where the \[(i,j)$-th entry describes the color (perhaps
only black or white) which is assigned to the \[(i,j)$-th ``pixel'' of
the computer screen. Each pixel represents a small square in the complex plane,
and the color which is assigned must tell
us something about intersection of \[L\] with this square.

In the case of a quadratic Julia set \[J(f)\], one very fast method
involves following iterates of the inverse map \[f^{-1}\],
taking all possible branches. (Compare \S3.9.) As Mandelbrot points out,
this method yields an excellent
picture of the outer parts of the \[J(f)\] but shows very little detail in
the inner parts. If we think of the electrostatic field produced by
an electric charge on \[J(f)\], this method will emphasize only those parts of
the Julia set at which ``lines of force'' (or ``external rays'' in the language
of \S18) tend to land.

A slower but much better procedure for plotting \[J(f)\] involves
iterating the map \[f\] for some large number (perhaps 50 to 50000) of times,
starting at the midpoint of each pixel. If the orbit ``escapes'' from
a large disk after \[n\] iterations, then the corresponding pixel is
assigned a color which depends on \[n\]. In more refined versions of
this method, one computes not only the value of the \[n$-th iterate of \[f\]
but also the absolute value of its derivative. Compare the discussion
below. Similar remarks apply to the {\bit Mandelbrot set} \[M\], as defined
by Douady and Hubbard. (Compare
Appendix F). In this case one takes the quadratic map
corresponding to the midpoint of the square and follows the orbit of its
critical point.

{\bf Remark.}
In order to understand some of the limitations of this method, consider
the situation near a fixed point \[f(z_0)=z_0\] in the Julia set.
First consider the repelling
case, with multiplier say of absolute value two. If we start at at point
\[z\] at distance \[1/1000\] from \[z_0\], then the distance from
\[z_0\] will roughly double with each iteration. Hence, after only ten
iterations the image of \[z\] will move substantially away from \[z_0\].
The result will be a computer picture which is quite sharp and accurate
near \[z_0\]. (Figures 3, 4.)

Now suppose that we try to construct a picture for \[z\mapsto z+z^4\] by
the same method. Again start with a point \[z\in J\] at a distance of
\[\epsilon\approx 1/1000\] from the fixed point at zero. Examining
the proof of 7.2, we see that the associated coordinate \[w_0=- 1/
3z_0^{\,3}\], which increases by one under each iteration, is roughly
equal to \[-1/3\epsilon^3\approx 300,000,000\]. {\it Thus
we would have to follow such an orbit for some \[300,000,000\]
iterations in order
to escape from a neighborhood of \[z=0\].} The result in practice will be a
false picture which shows everything near the origin to be in the Fatou set.
(This difficulty was eliminated in Figure 8 by using a special computer
program, which extrapolated iterates of \[f\] in order to distinguish
different Fatou components. Similarly, Figures
10, 12 were plotted with special purpose programs.)

For the fixed points of Cremer type, the situation is much worse. As far
as I know, no useful
computer picture of such a point has ever been produced!\smallskip

Many Julia sets are made up of very fine filaments. For such sets, it
is essential to make some kind of distance estimate in order to obtain
a sharp picture. For all of the center points of our pixels will quite likely
lie outside of these filaments, and hence correspond to escaping orbits.
But a good distance estimate can tell us that our pixel intersects the
\break

\noindent
set \[J(f)\], even though its center point is outside. In the case of
the Mandelbrot set \[M\], this procedure is even more important, since \[M\]
contains both large regions and also very fine filaments. Indeed,
it was precisely the difficulty of seeing such filaments which led
to Mandelbrot's incorrect belief that \[M\] has many components.

Here is an example of how first derivatives can be used to make such
distance estimates, following Fisher.
Consider a rational map \[f:\hat\C\to\hat\C\] with a superattractive
fixed point at the origin. Let \[\Omega\] be the basin of attraction
of this fixed point. Let us assume, to simplify the discussion, that this
basin is connected, simply connected,\break and
contains no other critical point of \[f\]. Then it is not difficult to
show that the B\"otk$\!$her coordinate of \S6.7 can be defined throughout
\[\Omega\], and yields a conformal isomorphism
\[\;\phi:\Omega\to D\;\] with \[\;\phi(f(z))=\phi(z)^n\,\].
Define the {\bit canonical potential function}\break
\[G:\Omega-\{0\}\to \R\;\] by
$$	G(z)=\log|\phi(z)|<0\,.	$$
(Compare \S17.) Denoting the gradient vector of \[G\] by \[G'\], we will show
that:\smallskip

(1) {\it the function \[G\] and the norm \[\| G' \|= |\phi'(z)/
\phi(z)|\] are easy to compute, and}\smallskip

(2) {\it the distance of \[z\] from the boundary of \[\Omega\]
can be computed, up to a factor\break
\indent of two, from a knowledge of \[G\] and \[\| G' \|\].}\smallskip

\noindent
In fact, for any orbit \[z_0\mapsto z_1\mapsto\cdots\] in \[\Omega\], it is easy
to check that
$$	G(z_0)=\lim_{k\to\infty}\log|z_k|/n^k\,.	$$
Since the convergence is locally uniform, we can also write
$$	\|G'(z_0)\|=\lim_{k\to\infty}{|dz_k/dz_0|\over n^k|z_k|}\,.	$$
In both cases, the successive terms can easily be computed
inductively, and
we obtain good approximations by iterating until \[|z_k|\]
is small. (If many iterations do not yield any small \[z_k\], then we must
assume that \[z_0\ne \Omega\], and set \[G=0\].)

Setting \[\phi(z)=w\], a brief computation shows that
the \Poincare metric on \[\Omega\] can be written as
$$	{2|dw|\over 1-|w|^2}\;=\; {2|\phi'(z)dz|\over 1-|\phi(z)|^2}\;
	=\;{\|G'(z)dz\|\over |\sinh G(z)|}\;.	$$
As an immediate consequence of the Quarter Theorem, \S A.7, A.8, we
obtain the following.

{\QP {\bf Corollary.} {\it The distance between \[z\] and the
boundary of \[\Omega\] is equal to\break
\[|\sinh( G)|/\|G'\|\], up to a factor of two.}\medskip}

\noindent If \[z\] is very close to \[\partial \Omega\], then \[G\] is small and
this distance estimate is very close to the ratio \[|G|/\|G'\|\]. It is
interesting to note that this is just the step size which would be
prescribed if we tried to solve the equation \[G(z)=0\] by Newton's method.
\smallskip

There are similar estimates in the case of a superattractive fixed point
at infinity, or in other words for the basin \[\C\ssm K(f)\] of a polynomial
map. For further information, the reader is referred to
Fisher, to Milnor [M5], and to Peitgen.

\vfill
\rightline{jwm, version of 9-4-91}
\end

%% file: cdref.tex
\input header.tex
\footline={\hss\tenrm R-\folio\hss}\pageno=1

\def\3i{\indent\indent\indent}
\font\tencyr=wncyr10
\input cyracc.def
\def\cyr{\tencyr\cyracc}

\centerline{\bf REFERENCES$\qquad$}
\bigskip

\indent\indent\indent{\bf For general background, see for example:}\smallskip

\ref [A1] L. Ahlfors, Complex Analysis, McGraw-Hill 1966.

\ref [Mu] J. Munkres, Topology: A First Course, Prentice-Hall 1975.

\ref [Wi] T. Willmore, An Introduction to Differential Geometry, Clarendon, 1959.
\medskip

\3i{\bf Surveys of conformal dynamics.}
\smallskip
\ref [Bl1] P. Blanchard, Complex analytic dynamics on the Riemann sphere,
Bull. Amer. Math. Soc. {\bf 11} (1984), 85-141.

\ref [BlC] P. Blanchard and A. Chiu, Conformal dynamics: an informal
discussion, Lecture Notes, Boston University 1990.

\ref [Br] H. Brolin, Invariant sets under iteration of rational functions,
Arkiv f\"or Mat. {\bf 6} (1965) 103-144. % (See p. 118.)

\ref [C] L. Carleson, Complex Dynamics, Lecture Notes UCLA 1990.

\ref [Dv1] R. Devaney, An Introduction to Chaotic Dynamical Systems,
$2^{nd}$ ed., Addison-Wesley, 1989.

\ref [D1] A. Douady, Syst\`emes dynamiques holomorphes, S\'eminar Bourbaki, $35^e$
ann\'ee 1982\--83, ${\rm n}^o$ 599.

\ref [Fa] K. Falconer, Fractal Geometry: Mathematical Foundations and
Applications, Wiley 1990. (Ch. 14.)

\ref [K1] L. Keen, Julia sets, pp. 57-74 of ``Chaos and Fractals, the Mathematics
behind the Computer Graphics'', edit. Devaney and Keen, Proc. Symp. Appl.
Math. {\bf 39}, Amer. Math. Soc. 1989.

\ref [L1] M. Lyubich, The dynamics of rational transforms: the topological picture,
Russian Math. Surveys {\bf 41:4} (1986), 43-117.

\medskip
\3i{\bf \S1. Simply connected surfaces; uniformization.}\smallskip

\ref [A2] L. Ahlfors, Conformal Invariants, McGraw-Hill 1973.

\ref [AS] L. Ahlfors and L. Sario, Riemann Surfaces, Princeton U. Press 1960.

\ref [Be] A. Beardon, A  Primer on Riemann Surfaces, Cambridge U. Press 1984.

\ref [FK] H. Farkas and I. Kra, Riemann Surfaces, Springer 1980.

\ref [Sp] G. Springer, Introduction to Riemann Surfaces, Addison-Wesley 1957.

\medskip
\3i{\bf \S2. Montel's theorem.}\smallskip

\ref [Mo] P. Montel, Le\c cons sur les Familles Normales, Gauthier-Villars 1927.
\medskip

\3i{\bf \S3. Fatou and Julia.}\smallskip

\ref [Ca] A. Cayley, Application of the Newton-Fourier method to an
immaginary root of an equation, Quart. J. Pure Appl. Math. {\bf 16} (1879)
179-185.

\ref [F1] P. Fatou, Sur les solutions uniformes de certaines \'equations
fonctionnelle, C. R. Acad. Sci. Paris {\bf 143} (1906) 546-548.

\ref [F2] P. Fatou, Sur les \'equations fonctionnelles, Bull. Soc. math. France
{\bf 47} (1919) 161-271, and {\bf 48} (1920) 33-94, 208-314.

\ref [J] G. Julia, Memoire sur l'iteration des fonctions rationelles,
J. Math. Pure Appl. {\bf 8} (1918) 47-245.

\ref [Ri] J. F. Ritt, On the iteration of rational functions, Trans. Amer.
Math. Soc. {\bf 21} (1920) 348-356.
\medskip

\3i{\bf \S4. Iterated maps (Hyperbolic and Euclidean cases).}
\smallskip

\ref [Ba1]
I. N. Baker, Repulsive fixedpoints of entire functions, Math. Zeit. {\bf 104}
(1968) 252-256.

\ref [Ba2], ---------, An entire function which has wandering domains,
J. Austral. Math. Soc. {\bf 22} (1976) 173-176.

\ref [Dn1] A. Denjoy, Sur l'it\'eration des fonctions analytiques, C.R.Acad.Sci.
Paris {\bf 182} (1926) 255-257.

\ref [Dv2] R. Devaney, Exploding Julia sets, pp. 141-154 of ``Chaotic
Dynamics and Fractals'', edit. Barnsley and Demko, Academic Press 1986.

\ref [EL] A. Eremenko and M. Lyubich, Dynamical properties of some classes of
entire functions, Stony Brook Institute for Mathematical Sciences
Preprint 1990\#4. (See also Sov. Math. Dokl. {\bf 30} (1984) 592-594;
Func. Anal. Appl. {\bf 19} (1985) 323-324; and
J. Lond. Math. Soc. {\bf 36} (1987) 458-468.)

\ref [F3] P. Fatou, Sur l'it\'eration des fonctions transcendantes enti\`eres,
Acta Math. {\bf 47} (1926) 337-370.

\ref [GK] L. Goldberg and L. Keen, A finiteness theorem for a dynamical class
of entire functions, Erg. Th. \& Dy. Sy. {\bf 6} (1986) 183-192.

\ref [K2] L. Keen, The dynamics of holomorphic self-maps of \[\C^*\], in
``Holomorphic Functions and Moduli'', edit. Drasin et al., Springer 1988.

\ref [L2] M. Lyubich, The measurable dynamics of the exponential map,
Siber. J. Math. {\bf 28} (1987) 111-127. (See also Sov. Math. Dokl.
{\bf 35} (1987) 223-226.)

\ref [R1] M. Rees, The exponential map is not recurrent, Math. Zeit. {\bf
191} (1986) 593-598.
\medskip

\3i{\bf \S5. Smooth Julia sets.}\smallskip

\ref [La] S. Latt\`es, Sur l'iteration des substitutions rationelles et les
fonctions de \Poincare, C. R. Acad. Sci. Paris {\bf 16} (1918) 26-28.

\ref [He1] M. Herman, Exemples de fractions rationnelles ayant une orbite dense
sur la sph\`ere de Riemann, Bull. Soc. Math. France {\bf 112} (1984) 93-142.

\ref [R2] M. Rees, Ergodic rational maps with dense critical point forward
orbit, Erg. Th. \& Dy. Sy. {\bf 4} (1984) 311-322.

\ref [R3] M. Rees, Positive measure sets of ergodic rational maps, Ann. Sci.
\'Ecole Norm. Sup. (4) {\bf 19} (1986) 383-407.

\ref [UvN] S. Ulam and J. von Neumann, On combinations of stochastic and
deterministic processes, Bull. Amer. Math. Soc. {\bf 53} (1947) 1120.
\medskip
\eject
\3i{\bf \S6. Attracting and repelling fixed points.}\smallskip

\ref [Sch] E. Schr\"oder, Ueber iterirte Functionen, Math. Ann. {\bf 3} (1871);
see p. 303.

\ref [K\oe ] G. K\oe nigs, Recherches sur les integrals de certains equations
fonctionelles, Ann. Sci. \'Ec. Norm. Sup. ($3^e$ ser.)
{\bf 1} (1884) suppl\'em. 1-41.

\ref [B\"o] L. E. B\"ottcher (~{\cyr B\"etkher\cdprime}),
The principal laws of convergence of iterates and
their application to analysis (Russian), Izv. Kazan. Fiz.-Mat. Obshch.
{\bf 14} (1904) 155-234.

\medskip

\3i{\bf \S7. Parabolic fixed points.}\smallskip

\ref [Le] L. Leau, \'Etude sur les equations fonctionelles \`a une ou plusi\`ers
variables, Ann. Fac. Sci. Toulouse {\bf 11} (1897).

\ref [\'Ec] J. \'Ecalle, Th\'eorie it\'erative: introduction a la th\'eorie des
invariants holomorphes, J. Math. Pure Appl. {\bf 54} (1975) 183-258.

\ref [Ca] C. Camacho, On the local structure of conformal mappings and holomorphic
vector fields, Ast\'erisque {\bf 59-60} (1978) 83-94.

\ref [MT] J. Milnor and W. Thurston, Iterated maps of the interval,
pp. 465-563 of ``Dynamical Systems (Maryland 1986-87)'', edit. J.C.Alexander,
Lect. Notes Math. {\bf 1342}, Springer 1988.

\ref({\bit See also~} [B\"o, pp.201-234] {\bit and~} [F3].)
\medskip

\3i{\bf \S8. Cremer points and Siegel disks.}\smallskip

\ref [Pf] G. A. Pfeifer, On the conformal mapping of curvilinear angles; the
functional equation \[\phi[f(x)]=a_1\phi(x)\], Trans. A. M. S. {\bf 18} (1917)
185-198.

\ref [Cr1] H. Cremer, Zum Zentrumproblem, Math. Ann. {\bf 98} (1927) 151-163.

\ref [Cr2] ------------,  \"Uber die H\"aufigkeit der Nichtzentren, Math. Ann. {\bf
115} (1938) 573-580.

\ref [Si] C. L. Siegel, Iteration of analytic functions, Ann. of Math. {\bf 43}
(1942) 607-612.

\ref [HW] G. H. Hardy and E. M. Wright, An Introduction to the Theory of Numbers,
Clarendon Press 1938, 1945, 1954.

\ref [Br] A. D. Bryuno, Convergence of transformations of differential equations
to normal forms, Dokl. Akad. Nauk USSR {\bf 165} (1965) 987-989.

\ref [SiM] C. L. Siegel and J. Moser, Lectures on Celestial Mechanics, Springer
1971.

\ref [He2] M. Herman, Recent results and some open questions on Siegel's
linearization  theorem of germs of complex analytic diffeomorphisms of
\[C^n\] near a fixed point, pp. 138-198 of Proc \[8^{th}\] Int. Cong. Math.
Phys., World Sci. 1986.

\ref [D2] A. Douady, Disques de Siegel et anneaux de Herman, S\'em. Bourbaki
39$^e$ ann\'ee (1986-87) n$^o$ 677.

\ref [Y1] J.-C. Yoccoz, Lin\'earisation des germes de diff\'eomorphismes
holomorphes de \[(\C,0)\], C. R. Acad. Sci. Paris {\bf 306} (1988)
55-58.
\medskip\eject

\3i{\bf \S9. Holomorphic fixed point formula.}\smallskip

\ref [AB] M. Atiyah and R. Bott, A Lefschetz fixed point formula for elliptic
differential operators, Bull.A.M.S. {\bf 72} (1966) 245-250.
\medskip

\3i{\bf \S10. Most periodic orbits repel.}\smallskip

\ref [Shi] M. Shishikura, On the quasiconformal surgery of rational functions,
Ann. Sci. \'Ec. Norm. Sup. {\bf 20} (1987) 1-29.
\medskip

% {\bf \S11. Repelling cycles are dense in the Julia set.}\smallskip

\3i{\bf \S12. Herman rings}\smallskip

\ref [Ar] V. Arnold, Small denominators I, on the mappings of the circumference
into itself, Amer. Math. Soc. Transl. (2) {\bf 46} (1965) 213-284.

\ref [CL] E. Coddington and N. Levinson, Theory of Ordinary Differential
Equations, McGraw-Hill 1955.

\ref [Dn2] A. Denjoy, Sur les courbes d\'efinies par les \'equations
diff\'erentielles \`a la surface du tore, Journ. de Math. {\bf 11} (1932)
333-375.

\ref [He3] M. Herman, Sur la conjugation diff\'erentiables des diff\'eomorphismes
du cercle \`a les rotations, Pub. I.H.E.S. {\bf 49} (1979) 5-233.

\ref [dM] W. de Melo, Lectures on One-Dimensional Dynamics, 17$^o$ Col.
Brasil. Mat., IMPA 1990.

\ref [Y2] J.-C. Yoccoz, Conjugation diff\'erentiable des diff\'eomorphismes du
cercle dont le nombre de rotation v\'erifie une condition diophantienne,
Ann. Sci. E.N.S. Paris (4) {\bf 17} (1984) 333-359.
 
\ref ({\bit See also~} [He1], [D2], [Shi].)

\medskip

\3i{\bf \S13. Classification of Fatou components.}\smallskip

\ref [S1] D. Sullivan, Conformal dynamical systems, pp. 725-752 of ``Geometric
Dynamics'', edit. Palis, Lecture Notes Math. {\bf 1007} Springer 1983.

\ref [S2] D. Sullivan, Quasiconformal homeomorphisms and dynamics I, solution
of the Fatou-Julia problem on wandering domains, Ann. Math. {\bf 122}
(1985) 401-418.

\noindent{(\bit See also~} [DH1] ?? {\bit and~} [DH2] {\bit below.})
\medskip

\3i{\bf \S14. Sub-hyperbolic and hyperbolic maps.}\smallskip

\ref [Sm] S. Smale, Differentiable dynamical systems, Bull. Amer. Math. Soc. {\bf
73} (1967) 747-817.

\ref [MSS] R. Ma\~n\'e, P. Sad and D. Sullivan, On the dynamics of rational maps,
Ann. Sci. \'Ec. Norm. Sup. Paris (4) {\bf 16} (1983) 193-217.

\ref [L3] M. Lyubich, Some typical properties of the dynamics of rational
maps, Russian Math. Surveys {\bf 38} (1983) 154-155. (See also Sov. Math.
Dokl. {\bf 27} (1983) 22-25/

\ref [L4] M. Lyubich, An analysis of the stability of the dynamics of rational
functions, Selecta Math. Sovietica {\bf 9} (1990) 69-90. (Russian original
published in 1984.)

\ref [DH1] A. Douady and J. H. Hubbard, A proof of Thurston's topological
characterization of rational functions, preprint, Mittag-Leffler 1984.

\ref [Shu] M. Shub, Global Stability of Dynamical Systems, Springer 1987.
\smallskip
{\bf Basic reference for \S\S14 and 16-19:}

\ref [DH2] A. Douady and J. H. Hubbard, \'Etude dynamique des polyn\^omes
complexes I \& II, Publ. Math. Orsay (1984-85).

\medskip

\3i{\bf \S15. Prime ends.}\smallskip

\ref [Ca] C. Carath\'eodory, \"Uber die Begrenzung einfach zusammenh\"angender
Gebiete, Math. Ann. {\bf 73} (1913) 323-370. (Gesam. Math. Schr., v. 4.)

\ref [Ep] D.B.A.Epstein, Prime Ends, preprint, Univ. Warwick 1978.

\ref [Ho] K. Hoffman, Banach Spaces of Analytic Functions, Prentice-Hall 1962.

\ref [Ma] J. Mather, Topological proofs of some purely topological consequences
of Carath\'eo\-dory's theory of prime ends, pp. 225-255 of ``Selected Studies'',
edit. T and G. Rassias, North-Holland 1982.

\ref [Oh] M. Ohtsuka, Dirichlet Problem, Extremal Length and Prime Ends, van
Nostrand 1970.

\noindent({\bit See also~} [A2].)
\medskip

\3i{\bf \S16. Local connectivity.}\smallskip

% \ref (The next few sections are based on the Douady-Hubbard Orsay Notes.
% See \S13 above.)\smallskip

\ref [Ku] K. Kuratowski, Topologie, Warsaw 1958-61.

\ref [HY] Hocking and Young, Topology, Addison-Wesley 1961.

\noindent({\bit See also~} [Mu].)
\medskip

\3i{\bf \S17. The filled Julia set.}\smallskip

\ref [Bro] M. Brown, A proof of the generalized Schoenflies theorem, Bulletin
A.M.S. {\bf 66} (1960) 74-76. (See also ``The monotone union of open
n-cells is an open n-cell'', Proc.A.M.S. {\bf 12} (1961) 812-814.)

\ref [Bl2] P. Blanchard, Disconnected Julia sets, pp. 181-201 of ``Chaotic
Dynamics and\break Fractals'', edit. Barnsley and Demko, Academic Press 1986.

\ref [BH1] B. Branner and J. H. Hubbard, The iteration of cubic polynomials,
Part I: the global topology of parameter space, Acta Math. {\bf 160} (1988)
143-206.

\ref [BH2] B. Branner and J. H. Hubbard, The iteration of cubic polynomials,
Part II: patterns and parapatterns, Acta Math., to appear.\medskip

\3i{\bf \S18. External rays and periodic points.}\smallskip

\ref [Pe] C. Petersen, On the Pommerenke-Levin-Yoccoz inequality,
preprint, IHES 1991.

\ref [Po] C. Pommerenke, On conformal mapping and iteration of rational
functions, Complex Var. Th. \& Appl. {\bf 5}, 1986.

\ref [G] L. Goldberg, Rotation cycles on the unit circle, Stony Brook
IMS preprint 1990/14.

\ref [GM] L. Goldberg and J. Milnor, Fixed point portraits of polynomial maps,
Stony Brook IMS preprint 1990/14.\medskip
% \eject

\3i{\bf Appendix A. Theorems from classical analysis.}\smallskip

\ref [Je] J. L. W. V. Jensen, Sur un nouvel et important th\'eoreme de la th\'eorie
des fonctions, Acta Math. {\bf 22} (1899) 219-251.

\ref [RR] F. and M. Riesz, \"Uber Randwerte einer analytischen Funktionen,
Quatr. Congr. Math. Scand. Stockholm 1916, 27-44.

\ref [Ko] P. Koebe, \"Uber die Uniformizierung beliebiger analytischer Kurven,
Nachr. Akad. Wiss. G\"ottingen, Math.-Phys. Kl. (1907) 191-210.

\ref [Gr] T. H. Gronwall, Some remarks on conformal representation, Ann. of
Math. {\bf 16} (1914-15) 72-76.

\ref [Bi] L. Bieberbach, \"Uber die Koeffizienten derjenigen Potenzreihen
welche eine schlichte Abbildung des Einheitskreises vermitteln, S.-B.
Preuss. Akad. Wiss. (1916) 940-955.

\ref [dB] L. de Brange, A proof of the Bieberbach conjecture, Acta Math. {\bf
154} (1985) 137-152.
\medskip

\3i{\bf Appendix B. Length-area-modulus inequalities.}\smallskip

\ref [CR] R. Courant and H. Robbins, What is Mathematics?, Oxford U. Press
 1941.\smallskip

\ref ({\bit See also~} [A2], [BH2].)\medskip

\3i{\bf Appendix C. Continued fractions.}\smallskip

\ref ({\bit See~} Hardy and Wright [HW].)\medskip

\3i{\bf Appendix D. Two complex variables.}\smallskip

\ref [Ze] E. Zehnder, A simple proof of a generalization of a theorem by C. L.
Siegel, in ``Geometry and Topology III'', edit. do Carmo and Palis, Lecture
Notes Math. {\bf 597}, Springer 1977.

\ref [Ne] S. Newhouse, Diffeomorphisms with infinitely many sinks, Topology
{\bf 16} (1974) 9-18.

\ref [H] J. Hubbard, The H\'enon mapping in the complex domain, pp. 101-111 of
Chaotic Dynamics and Fractals, M. Barnsley and S. Demko, (ed.),
Academic Press 1986.

\ref [FM] S. Friedland and J. Milnor, Dynamical properties of plane polynomial
automorphisms, Erg. Th. \& Dy. Sy. {\bf 9} (1989) 67-99.

\ref [Be] E. Bedford, Iteration of polynomial automorphisms of \[\C^2\],
Preprint, Purdue 1990.

\medskip

\3i{\bf Appendix E. Branched coverings and orbifolds}\smallskip

\ref [Mas] B. Maskit, Kleinian Groups, Grundl. math. Wiss.  {\bf 287} Springer
1987.

\ref [M1] J. Milnor, On the 3-dimensional Brieskorn manifolds \[M(p,q,r)\],
pp. 175-225 of ``Knots, Groups, and 3-Manifolds'', edit. Neuwirth, Ann.
Math. Studies {\bf 84}, Princeton U. Press 1975. (See \S2.)

\ref [T] W. Thurston, Three-dimensional Geometry
and Topology, to appear.
\medskip\eject

\medskip\3i{\bf Appendix F. Parameter space.}\smallskip

\ref [BM] R. Brooks and P. Matelski, The dynamics of 2-generator subgroups
of \[PSL(2,\C)\], pp. 65-71 of ``Riemann Surfaces and Related Topics'',
Proceedings 1978 Stony Brook Conference, edit. Kra and Maskit,
Ann. Math. Stud. {\bf 97} Princeton U. Press 1981.

\ref [Man] B. Mandelbrot, Fractal aspects of the iteration of \[z\mapsto\lambda
z(1-z)\] for complex \[\lambda,\,z\], Annals NY Acad. Sci. {\bf 357} (1980)
249-259.

\ref [DH3] A. Douady and J. H. Hubbard, It\'eration des polyn\^omes quadratiques
complexes, CRAS Paris {\bf 294} (1982) 123-126.

\ref [DH4] A. Douady and J. H. Hubbard, On the dynamics of polynomial-like
mappings, Ann. Sci. Ec. Norm. Sup. (Paris) {\bf 18} (1985), 287-343.

\ref [D3] A. Douady, Julia sets and the Mandelbrot  set, pp. 161-173
of ``The Beauty of Fractals'', edit. Peitgen and Richter, Springer 1986.

\ref [B1] B. Branner, The Mandelbrot set, pp. 75-105
of ``Chaos and Fractals'', edit.
Devaney and Keen, Proc. Symp. Applied Math. {\bf 39}, Amer. Math. Soc. 1989.

\ref [B2] B. Branner, The parameter space for complex cubic polynomials,
pp. 169-179 of ``Chaotic Dynamics and fractals'', edit. Barnsley and Demko,
Academic Press 1986.

%\ref [BH2] B. Branner and J. H. Hubbard, The iteration of cubic polynomials,
%Part I: the global topology of parameter space, Acta Math. {\bf 160} (1988)
%143-206.

\ref [M2] J. Milnor, { Remarks on iterated cubic maps}, Stony Brook IMS
preprint 1990/6.

\ref [M3] J. Milnor, { Hyperbolic components in spaces of polynomial maps},
in preparation.

\ref [M4] J. Milnor, { On cubic polynomials with periodic critical point},
in preparation.

\ref [R4] M. Rees, Components of degree two hyperbolic rational maps,
Invent. Math. {\bf 100} 1990, 357-382.

\ref [R5] M. Rees, { A partial description of parameter space of rational
maps of degree two}, Part I: preprint, Univ. Liverpool 1990;
and Part II: Stony Brook IMS preprints 1991/4.

\ref({\bit See also~} [BH1].)\medskip

\3i{\bf Appendix G. Computer graphics.}\smallskip

\ref [M5] J. Milnor, Self-similarity and hairiness in the Mandelbrot set,
pp. 211-257 of ``Computers in Geometry and Topology'', edit. Tangora,
Lect. Notes Pure Appl. Math. {\bf 114}, Dekker 1989 (cf. p. 218 and \S5).

\ref [Fi] Y. Fisher, Exploring the Mandelbrot set, pp. 287-296 of
The Science of Fractal Images,  edit. Peitgen and Saupe,
Springer 1989.

\ref [Pei] H.-O. Peitgen, Fantastic deterministic fractals, ibid. pp. 169-218.
\medskip

\vfill\bigskip

%\rightline{jwm, version of 8-30-91}
\eject

\footline={\hss\tenrm I-\folio\hss}\pageno=1

\centerline{\bf INDEX}\smallskip
$$\vbox{\settabs 2 \columns
\+Abel & 7-4\cr
\+annulus & 2-2, 2-5, 4-2, 12-1\cr
\+Arnold & 12-1, 12-2 \cr
\+attracting & 4-1, 3-2, 6-1, 7-1\cr
\+Baire's Theorem & 3-7, 8-1 \cr
\+ Baker & 4-1, 11-5\cr
\+basin & 4-2, 3-2, 7-1, 10-1\cr
\+Bieberbach & A-3ff, D-2\cr
\+Blaschke product & 5-3, 12-2, 12-4\cr
\+B\"otk$\!$her (= B\"ottcher) & 6-6, 17-1\cr
\+Branner & 17-3, F-2\cr
\+Brooks \& Matelski & F-1-F-2\cr
\+Brown & 17-2\cr
\+Bryuno & 8-5ff, 12-4\cr
\+Carath\'eodory & 1-1, \S15, \S16, 17-3, 18-1\cr
\+cellular set & 17-2, F-2\cr
\+complete surface & 1-4\cr
\+conformal automorphism, metric & 1-1,1-3\cr
\+continued fraction & 8-5, 12-4, \S C\cr
\+covering surface & 2-1, 2-5, E-1\cr
\+Cremer & \S8, 11-4, 18-3\cr
\+curvature & 1-3, 1-5, 1-9, 2-5\cr
\+cycle & 3-2, 10-1\cr
\+deck transformation & 2-1, E-1\cr
\+Denjoy & 4-1, 12-1\cr
\+Diophantine number & 8-3, 12-2 \cr
\+Douady & 14-1, 17-3, \S18, F-2, G-1\cr
\+\'Ecalle & 7-4\cr
\+Euclidean surface & 2-1\cr
\+Euler characteristic & 5-2, 14-4, E-2ff\cr
\+expanding map & 14-1\cr
\+external ray & 17-2, \S18\cr
\+Fatou & \S3, 4-1, 6-4, \S7, \S10, 11-2, 13-1, 15-3, D-2\cr
\+Flower Theorem & 7-1\cr
\+fractal & 5-1, 6-2, F-2\cr
\+fractional linear transformation & 1-2, 1-5, 1-9\cr
\+fully invariant & 3-2\cr
\+generic & 3-7, 8-1\cr
\+grand orbit & 3-2, 3-6\cr
\+Green's function & 17-1, G-2\cr
}$$ 
\vfil\eject

$$\vbox{\settabs 2 \columns
\+half-plane & 1-2\cr
\+H\'enon map & D-1\cr
\+Herman & 5-3, 8-8, \S12, 13-1\cr
\+homoclinic point & 11-1  \cr
\+Hubbard & 14-1, 17-3, \S18, \S F-1, G-1\cr
\+hyperbolic, sub-hyperbolic map & 14-1, 14-3, 17-3, F-1\cr
\+Hyperbolic surface & 2-2\cr
\+Jensen & A-1\cr
\+Jordan curve & 3-3, 6-4, 7-7, 15-4, 16-2\cr
\+Julia & \S3, 6-4, 11-1, 13-1\cr
\+Koebe & 1-1, A-4\cr
\+K\oe nigs & 6-1, D-1\cr
\+Latt\`es & 5-1ff\cr
\+Leau & \S7\cr
\+Lebesgue measure & 8-2, 8-8, 15-2, 15-3, A-2\cr
\+Liouville & 1-2, 8-3, 8-10\cr
\+local connectivity & \S16, 17-3, 18-1ff, F-2\cr
\+locally uniform convergence & 2-3\cr
\+Lyubich & 4-1, 8-2, F-2\cr
\+Mandelbrot & \S F, \S G \cr
\+maximum modulus principle & 1-1, 17-1\cr
\+Montel's Theorem& 2-4, 3-6, 11-2\cr
\+multiplier & 3-2\cr
\+neutral & 3-2\cr
\+Newton's method & 5-3, F-1\cr
\+normal family & 2-4\cr
\+$o(),\;{\cal O}()$ & 6-6, 7-3\cr
\+orbifold, orbifold metric & 14-3ff, \S E\cr
\+parabolic & 1-7, 1-9, 2-2, \S7\cr
\+periodic orbit & 3-2, \S10, \S11\cr
\+petal & 7-1, 7-8, 18-6\cr
\+Pfeiffer & 8-1\cr
\+Picard's Theorem & 2-3\cr
\+\Poincare  & 1-3ff, 2-2, 12-1\cr
\+prime end & 15-5\cr
\+proper map & 5-2, 12-4, 13-1, E-1\cr
\+PSL & 1-2, 1-6, 1-9\cr
\+punctured plane, disk & 2-1, 2-2, 4-2\cr
\+quasi-conformal deformation, surgery & 10-1, 12-1, 12-4, 13-4\cr
}$$
\vfil\eject

$$\vbox{\settabs 2 \columns
\+ramification & 14-4, E-1\cr
\+rational function & 3-1, 5-3\cr
\+Rees & 4-1, 5-3, 8-2, F-2\cr
\+repelling & 3-2, 6-1, 7-1\cr
\+residue & 9-1 \cr
\+Riemann surface, metric, sphere, map& 1-1, 1-3, 1-5, 15-1\cr
\+Riemann-Hurwitz &5-2, 17-2, E-2 \cr
\+Riesz & 8-9, 15-3, 16-4, A-2\cr
\+root of unity & 7-1, 7-4\cr
\+rotation, rotation number & 1-2, 4-1ff, 12-1ff, 13-1\cr
\+Schr\"oder & 6-1\cr
\+Schwarz Lemma, & 1-1ff, A-4\cr
\+Schwarz inequality&  15-2, B-1\cr
\+self-similarity & 3-2, 3-9\cr
\+Shishikura & 10-1, 12-1, 12-4\cr
\+Siegel & 8-2ff, 11-3, 12-1\cr
\+size (of disk) & 8-8\cr
\+Smale & 14-3\cr
\+Snail Lemma & 13-2, 18-5 \cr
\+spherical metric, surface & 1-5, 2-1\cr
\+spiral & 6-2\cr
\+Sullivan & \S13-1, 14-3, 18-2, F-2\cr
\+superattracting & 3-2, 6-6\cr
\+Tchebycheff polynomial & 5-3\cr
\+thrice punctured sphere & 2-3\cr
\+Thurston & 11-4, 13-4, 14-3, \S E \cr
\+torus & 2-1, 4-1, 4-6, 5-1ff\cr
\+transcendental function & 4-1\cr
\+transitive action & 1-2, 1-8\cr
\+Ulam \& Von Neumann & 5-1\cr
\+uniformization & 1-1, 2-1\cr
\+univalent & 7-4, 15-3, A-1ff\cr
\+wandering component, point & 13-4, 14-3\cr
\+Weierstrass Theorem, \[\wp\] & 1-2, 5-2\cr
\+Yoccoz & 8-6ff, 12-2, 18-2, 18-6, F-2\cr
}$$

\end